%% file: Paper.tex
\documentclass[a4paper,twoside,leqno]{article}

\newcommand\CustomTitle{%
    Scaling-robust built-in a~posteriori error estimation for discontinuous least-squares finite element methods
}
\newcommand\CustomShortTitle{%
    Scaling-robust discontinuous LSFEM%
}
\newcommand\CustomAbstract{%
    A convincing feature of least-squares finite element methods is the built-in a posteriori error estimator for any conforming discretization.
    In order to generalize this property to discontinuous finite element ansatz functions, this paper introduces a least-squares principle
    on piecewise Sobolev functions by the example of the Poisson model problem with mixed boundary conditions.
    It allows for fairly general discretizations including standard piecewise polynomial ansatz spaces on triangular and polygonal meshes.
    The presented scheme enforces the interelement continuity of the piecewise polynomials by additional least-squares residuals.
    A side condition on the normal jumps of the flux variable requires a vanishing integral mean and
    enables the penalization of the jump with the natural power of the mesh size in the least-squares functional.
    This avoids over-penalization with additional regularity assumptions on the exact solution
    as usually present in the literature on discontinuous LSFEM.
    The proof of the built-in a~posteriori error estimation for the over-penalized scheme is presented as well.
    All results in this paper are robust with respect to the size of the domain guaranteed by
    a suitable weighting of the residuals in the least-squares functional.
    Numerical experiments illustrate the importance of the proposed weighting and
    exhibit optimal convergence rates of the adaptive mesh-refining algorithm for various polynomial degrees.
}
\newcommand\CustomKeywords{%
    least-squares finite element method,
    discontinuous Galerkin method,
    a posteriori error analysis,
    mixed boundary conditions,
    weighting of least-squares functional,
    scaling-robust estimates
}
\newcommand\CustomClassification{%
    65N15, 
    65N30, 
    65N50  
}
\newcommand\CustomFunding{%
    This research has been supported by the Austrian Science Fund (FWF) through the project \emph{Computational nonlinear PDEs} (grant DOI \href{https://doi.org/10.55776/P33216}{10.55776/P33216}).
    The author gratefully acknowledges the many fruitful discussions on this topic
    in the working group of Prof.\ Carsten Carstensen at Humboldt-Universit\"at zu Berlin and
    in the working group of Prof.\ Dirk Praetorius at the Institute of Analysis and Scientific Computing of TU Wien.
}
\newcommand\AuthorA{Philipp Bringmann}
\newcommand\ShortAuthorA{P.~Bringmann}
\newcommand\AffiliationA{%
    TU Wien,
    Institute of Analysis and Scientific Computing,
    Wiedner Hauptstr.\ 8--10,
    1040 Wien,
    Austria
}

\input{minarthd}

\numberwithin{equation}{section}

\usepackage{hyperref}

\newcommand\vvvert{|\mkern-1.5mu|\mkern-1.5mu|}

\newcommand\E{\mathcal{E}}

\newcommand\M{\mathcal{M}}
\newcommand\N{\mathbb N}
\newcommand\R{\mathbb{R}}
\newcommand\T{\mathcal{T}}
\newcommand\V{\mathcal{V}}

\DeclareMathOperator\conv{conv}
\DeclareMathOperator\curl{curl}
\DeclareMathOperator\Curl{Curl}

\DeclareMathOperator\diam{{diam}}
\DeclareMathOperator\dist{{dist}}
\DeclareMathOperator\Div{div}

\DeclareMathOperator\Int{int}
\DeclareMathOperator\Mid{mid}
\DeclareMathOperator\osc{osc}
\DeclareMathOperator\supp{supp}

\newcommand\dx{\;\textup{d}x}
\newcommand\ds{\;\textup{d}s}
\newcommand\CO{{\textup{C}}}
\newcommand\CR{{\textup{CR}}}

\newcommand\DIR{{\textup{D}}}

\newcommand\FR{{\textup{F}}}
\newcommand\NEU{{\textup{N}}}
\newcommand\LS{{\textup{LS}}}

\newcommand\PW{{\textup{pw}}}
\newcommand\RT{{\textup{RT}}}

\begin{document}
    \maketitle
    %
    %
    \CustomAbstractFormat

\section{Introduction}
\label{sec:introduction}
Minimal residual methods with discontinuous ansatz functions have gained much attention in the past decade.
Besides discontinuous Petrov--Galerkin methods (DPG) for ultraweak formulations \cite{DemkowiczGopalakrishnan2017}, 
this includes least-squares finite element methods (LSFEMs) with discontinuous polynomial function spaces.
A key feature of minimal residual methods is the built-in a~posteriori error estimation;
see \cite{MR3215064} for the DPG method and \cite{MR4782072} for an overview for LSFEMs.
This makes these classes of methods a powerful choice for adaptive mesh-refining algorithms; see also \cite{MR3671598,MR4138307,MR4216839}
for convergence results of adaptive LSFEMs driven by the built-in error estimator.

This paper considers the model example of the Poisson problem to establish a novel approach in the a~posteriori error analysis
for LSFEMs with discontinuous ansatz spaces.
Given a right-hand side \(f \in L^2(\Omega)\) on a bounded polygonal Lipschitz domain \(\Omega \subseteq \R^2\)
and an \(\Omega\)-dependent weight factor \(c_\Omega > 0\), 
the conforming least-squares formulation minimizes the squared \(L^2\) norms of the residuals
\begin{equation}
    \label{eq:least_squares_functional}
    c_\Omega^2\,
    \Vert f + \Div \sigma \Vert_{L^2(\Omega)}^2
    +
    \Vert \sigma - \nabla u \Vert_{L^2(\Omega)}^2
\end{equation}
over Sobolev functions \(\sigma \in H_\NEU(\Div, \Omega)\) and \(u \in H^1_\DIR(\Omega)\).
The well-posedness of this formulation follows from the fundamental equivalence~\cite[Lem.~4.3]{MR0461948}
between the homogeneous least-squares functional with \(f \equiv 0\) and the norms of the underlying function spaces
which reads, for all \(\tau \in H_\NEU(\Div, \Omega)\) and \(v \in H^1_\DIR(\Omega)\),
\begin{equation}
    \label{eq:least_squares_principle}
    c_\Omega^2\,
    \Vert \Div \tau \Vert_{L^2(\Omega)}^2
    +
    \Vert \tau - \nabla v \Vert_{L^2(\Omega)}^2
    \approx
    c_\Omega^2\,
    \Vert \Div \tau \Vert_{L^2(\Omega)}^2
    +
    \Vert \tau \Vert_{L^2(\Omega)}^2
    +
    \Vert \nabla v \Vert_{L^2(\Omega)}^2.
\end{equation}
The LSFEM minimizes~\eqref{eq:least_squares_functional} over discrete conforming subspaces of \(H_\NEU(\Div, \Omega)\) and \(H^1_\DIR(\Omega)\).

In case of nonconforming discrete spaces in the LSFEM, the interelement continuity conditions need to be enforced
for the normal jump \([\sigma \cdot n_E]_E = 0\) and the full jump \([u]_E = 0\) along every interior edge \(E \in \E(\Omega)\).
One possible realization introduces additional weighted residual terms of discontinuous Galerkin type.
For some discrete subspaces \(\Sigma(\T) \subseteq H(\Div, \T)\) and \(U(\T) \subseteq H^1(\T)\) of
piecewise Sobolev spaces and an exponent \(\alpha \in \{-1, 1\}\),
the discontinuous LSFEM seeks minimizers \(\sigma_h \in \Sigma(\T)\) and \(u_h \in U(\T)\) of the functional
\begin{equation}
    \label{eq:intro_discontinuous_LSFEM}
    \begin{split}
        (\sigma_h, u_h)
        &\mapsto
        c_\Omega^2\,
        \Vert f + \Div_\PW \sigma_h \Vert_{L^2(\Omega)}^2
        +
        \Vert
        \sigma_h - \nabla_\PW u_h
        \Vert_{L^2(\Omega)}^2
        \\
        &\phantom{{}\mapsto{}}
        +
        \sum_{E \in \E(\Omega) \cup \E(\Gamma_\NEU)}
        c_\Omega^{(1 - \alpha)}\,
        h_E^\alpha\,
        \Vert [\sigma_h \cdot n_E]_E \Vert_{L^2(E)}^2
        +
        \sum_{E \in \E(\Omega) \cup \E(\Gamma_\DIR)}
        h_E^{-1}\,
        \Vert [u_h]_E \Vert_{L^2(E)}^2.
    \end{split}
\end{equation}
Since the normal jump of an \(H(\Div)\) function belongs to \(H^{-1/2}\),
the \(L^2\) norm of the jump is not well-defined in general and the natural choice for the exponent is \(\alpha = 1\).
However, the presence of the divergence term leads to the exponent \(\alpha = -1\) in the common discontinuous least-squares
discretizations in the literature; see, e.g., \cite{MR2149925,MR3895838}.
The latter will be called \emph{over-penalized} discontinuous LSFEM throughout this paper.
Indeed, in the case of polygonal discretization of piecewise polynomial degree \(k \in \N\),
the over-penalization necessitates additional regularity assumptions
on the exact solution \(\sigma \in H^{k+1}(\Omega; \R^2)\) and, thus, \(u \in H^{k+2}(\Omega)\)
in order to guarantee a~priori convergence of rate \(h^k\) in \cite[Lem.~4.1]{MR3895838}.
In contrast to that, the novel discontinuous least-squares scheme proposed in this paper
follows a more general approach by establishing a fundamental equivalence for piecewise Sobolev functions
including suitable measures for the interelement jumps in \(H^{-1/2}\).
As a consequence, \emph{any} discontinuous discrete subspace of the piecewise Sobolev spaces leads
to a well-posed discontinuous LSFEM providing a consistent extension of a conforming LSFEM
and satisfying a quasi-bestapproximation property without any additional regularity assumption.
A side condition on the lowest moments of the normal jumps of the flux \(\sigma\) enables the natural power of \(h\) 
with \(\alpha = 1\) in the discontinuous least-squares functional~\eqref{eq:intro_discontinuous_LSFEM} and,
thereby, optimal convergence rates without additional regularity assumptions.
On triangular meshes, the resulting saddle-point problem can be replaced by a symmetric and positive definite system matrix
using suitable basis functions for the piecewise Raviart--Thomas finite element space.
A detailed construction can be found in Appendix~\ref{app:Sigma_basis} below.
The results carry over to polygonal meshes as well and allow for variable polynomial degrees on different triangles,
although the analysis in this paper is not \(p\)-robust.

This paper validates the novel approach in the analysis of discontinuous LSFEMs for a simple model problem.
However, the application of the underlying methodology is by no means restricted to the Poisson problem.
Once a fundamental equivalence on the piecewise Sobolev spaces is established,
the key results of well-posedness of the discrete formulation as well as the a~priori and a~posteriori error analysis
can be directly deduced for all discontinuous subspaces.
In this sense, the novel approach mimics the analysis for conforming LSFEMs.
Moreover, the presentation is restricted to the two-dimensional case for simplicity,
although the arguments can be generalized to three dimensions as well.
The necessary adaptations of the dimension-dependent arguments are addressed
in Remarks~\ref{rem:3D_Helmholtz} and~\ref{rem:3D} below.

A common criticism towards least-squares formulations concerns the dependence of the fundamental equivalence
constants in~\eqref{eq:least_squares_principle} on the diameter \(\diam(\Omega)\) of the domain \(\Omega\)
caused by the different scaling of the divergence and the \(L^2\) contribution of the flux variable \(\sigma\).
This results in a~priori and a~posteriori error estimates that might degenerate for large domains.
As a remedy, all results of this paper include explicit weights \(c_\Omega\) in the least-squares functionals
in terms of the domain size to guarantee scaling-robust estimates.
Note that the over-penalization requires an additional weight factor \(c_\Omega^2\) in front of the normal jump term
in~\eqref{eq:intro_discontinuous_LSFEM} as shown in Section~\ref{sec:overpenalised_DLSFEM} below.
The reader is referred to \cite{MR4253887} for related results in the context of DPG methods.

Typically, the stability of discontinuous Galerkin schemes critically hinges on the choice of sufficiently large penalty parameters~\cite{MR2882148}.
Although, these parameters can be chosen automatically in terms of local geometric quantities (see, e.g., \cite{MR4796047}),
they are still limited to moderately large values.
In contrast to that, the discontinuous LSFEM may be considered as a discontinuous Galerkin method with stability independent of any penalty parameter.

Contributions on discontinuous LSFEMs date back to the late 1990s.
Cao and Gunzburger employed penalty terms of discontinuous Galerkin type for an LSFEM for the interface problem \cite{MR1618488}
and Day and Bochev for a least-squares method for consistent mesh tying \cite{MR2431595}.
Gerritsma and Proot introduced the first least-squares formulation with discontinuous ansatz functions in \cite{MR1910570}
for a one-dimensional model problem discretized using Lagrangian polynomials for Gauss--Lobatto points.
This approach has been applied to 2D viscoelastic fluid flow \cite{MR2285779} and 
the Stokes equations \cite{MR4058304,MR4414910,MohapatraKumarJoshi2023} as well.
Bensow and Larson employed discontinuous finite element functions for a least-squares formulation of the Poisson model problem \cite{MR2149925}
and the div-curl problem \cite{MR2195400}.
Their so-called discontinuous/continuous LSFEM follows a domain decomposition ansatz to employ the computationally more demanding
discontinuous functions solely close to expected singularities of the exact solution.
Bochev, Lai, and Olson showed that discontinuous ansatz functions for the velocity allow for exact local conservation of mass
in the context of the Stokes equations \cite{MR2878612,MR3049438}.
The a~priori error analysis in \cite{MR3895838} establishes optimal convergence rates for uniform refinement of polygonal meshes
under increased regularity assumptions.
This has been applied to div-curl systems \cite{MR4018037}, the time-harmonic Maxwell equations \cite{MR4367673},
and the Helmholtz equation \cite{MR4563176}.
Alternative discretizations enforce the interelement continuity by Lagrange multipliers similar
to hybridized discontinuous Galerkin schemes for linear elasticity and hyperelasticity in \cite{MR4433564,IgelbuscherDissertation}.

While the idea of discontinuous ansatz spaces for LSFEMs has been applied to several problems, the a~posteriori analysis is still widely open.
To the best of the author's knowledge, none of the mentioned publications provide a rigorous a~posteriori error analysis
except for \cite[Thm.~4.1]{MR4018037} and the recent publication \cite{MR4593282}.
The result in the former reference is restricted to the norm induced by the least-squares functional.
The latter considers the discretization with lowest-order nonconforming Crouzeix--Raviart functions for the primal variable \(u\)
and \(H(\Div)\)-conforming Raviart--Thomas functions for the dual variable \(\sigma\).
The analysis in~\cite{MR4593282} employs suitable averaging operators and guides the proofs of a built-in a~posteriori estimator
for the (over-penalized) fully discontinuous LSFEM in Section~\ref{sec:overpenalised_DLSFEM} below.

The outline of this paper reads as follows.
Section~\ref{sec:piecewise_Sobolev} introduces the piecewise Sobolev spaces with and without boundary conditions.
The main result is a scaling-robust least-squares principle for these piecewise Sobolev spaces in Section~\ref{sec:discontinuous_least_squares}.
Section~\ref{sec:piecewise_polynomials} defines piecewise polynomial functions and characterizes their jumps 
to allow for the formulation of discontinuous LSFEMs in Section~\ref{sec:discontinuous_LSFEM}.
The following Section~\ref{sec:error_analysis} presents a~priori and a~posteriori error estimates.
Section~\ref{sec:overpenalised_DLSFEM} establishes an a~posteriori result for the standard over-penalized discontinuous LSFEM.
The paper concludes with a numerical investigation of the considered schemes in Section~\ref{sec:experiments},
a required quasi-interpolation result in Appendix~\ref{app:quasi_interpolation},
and the explicit construction of a Raviart--Thomas basis with vanishing mean value of the normal jumps in Appendix~\ref{app:Sigma_basis}.

\section{Piecewise Sobolev spaces and jumps}
\label{sec:piecewise_Sobolev}

Let \(\Omega \subseteq \R^2\) denote a bounded and simply connected polygonal Lipschitz domain 
with outward unit normal vector \(n\colon \partial\Omega \to \R^2\)
and \(t\colon \partial\Omega \to \R^2\) its counter-clockwise rotation by \(\pi / 2\).
The boundary \(\partial\Omega\) is decomposed into the closed Dirichlet part \(\Gamma_\DIR \subseteq \partial\Omega\)
with positive surface measure \(\vert \Gamma_\DIR \vert > 0\) and the relatively open Neumann part
\(\Gamma_\NEU \coloneqq \partial\Omega \setminus \Gamma_\DIR\).
Without loss of generality, assume that \(\vert \Gamma_\NEU \vert > 0\).
The Neumann boundary solely affects the definition of the jump term \(s_t\) in~\eqref{eq:jump_terms} below.
All arguments in this paper apply analogously in the case of pure Dirichlet boundary data.
The definition of piecewise Sobolev spaces employs the notion of a regular triangulation \(\T\) of \(\Omega\) into closed triangles.
Let \(A \lesssim B\) abbreviate the relation \(A \leq C\,B\) with a positive generic constant \(C > 0\) which solely depends on
the polynomial degree of the discretization and the interior angles of the triangles in \(\T\), 
but is independent of the underlying piecewise constant mesh-size function \(h_\T \in P_0(\T)\) 
defined by the diameter \(h_\T \vert_T \coloneqq h_T \coloneqq \diam(T)\) of the triangle \(T \in \T\) and the size of the domain \(\Omega\).
The equivalence \(A \approx B\) means \(A \lesssim B\) and \(B \lesssim A\).
For any triangle \(T \in \T\), let \(n_T\colon \partial T \to \R^2\) denote the outward unit normal vector
and \(t_T\colon \partial T \to \R^2\) its counter-clockwise rotation by \(\pi/2\).
Moreover, define the enlarged patch \(\Omega_T \subseteq \Omega\) of a triangle \(T \in \T\) by
\[
    \Omega_T
    \coloneqq
    \Int\Big(
        \bigcup
        \big\{
            K \in \T
            \;:\;
            \dist(T, K)
            \coloneqq
            \inf_{x \in T, y \in K}
            \vert x - y \vert
            = 0
        \big\}
    \Big).
\]

Let \(\E(T)\) denote the set of the three edges of \(T\).
This leads to the set \(\E \coloneqq \bigcup_{T \in \T} \E(T)\) of all edges in \(\T\).
Throughout the paper, assume that every triangulation \(\T\) reflects the dissection of the boundary of \(\Omega\) 
in that the Dirichlet edges \(\E(\Gamma_\DIR) \coloneqq \{E \in \E(\partial\Omega) : E \subseteq \Gamma_\DIR\}\)
and the Neumann edges \(\E(\Gamma_\NEU) \coloneqq \{E \in \E(\partial\Omega) : E \subseteq \overline\Gamma_\NEU\}\)
partition the set \(\E(\partial\Omega)\) of boundary edges.
It remains the set of interior edges \(\E(\Omega) \coloneqq \E \setminus \E(\partial\Omega)\).
The set of adjacent triangles \(\T(E) \subseteq \T\) of an edge \(E \in \E\) consists either of exactly two triangles \(T_+\) and \(T_-\)
if \(E \in \E(\Omega)\) belongs to the interior or of exactly one triangle \(T_+\) if \(E \in \E(\partial\Omega)\) lies on the boundary.
The index \(\pm\) is determined by the fixed orientation of the unit normal vector \(n_E \in \R^2\) of \(E\) such that
\(n_E \cdot n_{T_\pm} = \pm 1\) as illustrated in Figure~\ref{fig:edge_patch}.
The fixed orientation of \(n_E\) induces an orientation of the tangential vector \(t_E \in \R^2\) of \(E\).
Let \(\omega_E \coloneqq \Int(\bigcup_{T \in \T(E)} T) \subseteq \Omega\)
and \(\Omega_E \coloneqq \bigcup_{T \in \T(E)} \Omega_T\) denote the standard and the enlarged edge patch,
\(h_E \coloneqq \vert E \vert\) the length of \(E\), and \(\Mid(E) \in \overline\Omega\) its midpoint.
\begin{figure}
    \centering
    \begin{tikzpicture}[line width=1.5,scale=5,line join=round,line cap=round,>=stealth]
        \path[draw,TUblue,fill=TUblue!25!white] (0.2,0) coordinate (A1) -- (0.8,0) coordinate (A2) --
        (1.1,.5) coordinate (A4) -- (.5,.5) coordinate (A3) --cycle;
        \path[draw,TUmagenta] (A2)--(A3) coordinate[midway,label=below left:\(E\)] (A23);
        \path[TUblue] ($(A23)!.5!(A1)$) node {\(T_+\)};
        \path[TUblue] ($(A23)!.5!(A4)$) node {\(T_-\)};
        \draw ($(A2)!.25!(A3)$) node[shape=coordinate] (nuBase) {};
        \draw[TUmagenta,->] (nuBase) -- ($(nuBase)!1!90:(A2)$) node[shape=coordinate,label=above:{$n_E$}] (nuTop) {nuTop};
        \draw ($(A2)!.5!(A4)$) node[shape=coordinate] (nuBaseMinus) {};
        \draw[TUblue,->] (nuBaseMinus) -- ($(nuBaseMinus)!0.5!90:(A2)$) node[shape=coordinate,label=below:{$n_{T_-}$}] (nuTop) {nuTop};
        \path ($(A1)!.5!(A3)$) node[coordinate] (nuBasePlus) {};
        \draw[TUblue,->] (nuBasePlus) -- ($(nuBasePlus)!0.5!-90:(A1)$) node[shape=coordinate,label=above:{$n_{T_+}$}] (nuTop) {nuTop};
    \end{tikzpicture}
    \caption{Edge patch \(\omega_E\)}
    \label{fig:edge_patch}
\end{figure}

Let \(\V\) denote the set of vertices of \(\T\) analogously partitioned into the set \(\V(\Omega)\) of interior vertices,
the set \(\V(\Gamma_\DIR)\) of vertices on the Dirichlet boundary, and
the set \(\V(\Gamma_\NEU)\) of vertices on the Neumann boundary.
The set \(\V(E)\) consists of the two vertices of an edge \(E \in \E\).

The context-sensitive measure \(\vert {}\cdot{} \vert\) denotes not only the Lebesgue measure of Lebesgue sets in \(\R^2\)
or the arc length of one-dimensional sets but also the modulus of real numbers, the cardinality of finite sets,
and the Euclidean norm of vectors in \(\R^2\).

Given any open domain \(\omega \subseteq \R^2\), this paper employs standard notation for Sobolev and Lebesgue spaces
\(H^1(\omega)\), \(H(\Div, \omega)\), and \(L^2(\omega)\).
Appropriate subscripts designate their usual norms \(\Vert {}\cdot{} \Vert_{H^1(\omega)}\), \(\Vert {}\cdot{} \Vert_{H(\Div, \omega)}\),
and \(\Vert {}\cdot{} \Vert_{L^2(\omega)}\).
For \(T \in \T\), write \(H(\Div, T) \coloneqq H(\Div, \operatorname{int}(T))\) and \(H^1(T) \coloneqq H^1(\operatorname{int}(T))\).
Define the piecewise Sobolev spaces by \cite{MR3521055}
\begin{align*}
    H(\Div, \T)
    &\coloneqq
    \{
        \tau \in L^2(\Omega; \R^2)
        \;:\;
        \forall T \in \T,\: \tau\vert_T \in H(\Div, T)
    \},
    \\
    H^1(\T)
    &\coloneqq
    \{
        v \in L^2(\Omega)
        \;:\;
        \forall T \in \T,\: v\vert_T \in H^1(T)
    \}.
\end{align*}
The corresponding spaces with partial homogeneous boundary conditions read
\begin{align*}
    H_\NEU(\Div, \T)
    &\coloneqq
    \left\{
        \tau \in H(\Div, \T)
        \;:\;
        \begin{aligned}
            &\forall E \in \mathcal{E}(\Gamma_\NEU)\:
            \forall w \in H^1(\omega_E) \text{ with }
            \\
            &w\vert_{\partial \omega_E \setminus \Gamma_\NEU} = 0,
            \langle \tau \cdot n_E,\: w \rangle_{\partial \omega_E} = 0
        \end{aligned}
    \right\},
    \\
    H^1_\DIR(\T)
    &\coloneqq
    \{
        v \in L^2(\Omega)
        \;:\;
        \forall E \in \mathcal{E}(\Gamma_\DIR),\:
        v\vert_E = 0
    \}.
\end{align*}
The differential operators \(\Div_\PW\) and \(\nabla_\PW\) apply piecewise, i.e., for \(\tau \in H(\Div, \T)\),
the function \(\Div_\PW \tau \in L^2(\Omega)\) is defined by \((\Div_\PW\tau)\vert_T \coloneqq \Div(\tau \vert_T)\) for all \(T \in \T\)
and analogously for \(\nabla_\PW v\) with \(v \in H^1(\T)\).
For two dimensional domains, the curl operators are defined, for \(v \in H^1(\Omega)\) and \(\beta \in H^1(\Omega; \R^2)\), by
\[
    \Curl v
    \coloneqq
    \begin{pmatrix}
        -\partial v / \partial x_2
        \\
        \partial v / \partial x_1
    \end{pmatrix}
    \quad\text{and}\quad
    \curl \beta
    \coloneqq
    \partial \beta_2 / \partial x_1
    - \partial \beta_1 / \partial x_2.
\]

In contrast to the global counterparts
\begin{align*}
    H_\NEU(\Div, \Omega)
    &\coloneqq
    H(\Div, \Omega) \cap H_\NEU(\Div, \T)
    =
    \{
        \tau \in H(\Div, \Omega)
        \;:\;
        \forall w \in H^1_\DIR(\Omega),\:
        \langle \tau \cdot n, w \rangle_{\partial\Omega}
        =
        0
    \},
    \\
    H^1_\DIR(\Omega)
    &\coloneqq
    H^1(\Omega)
    \cap
    H^1_\DIR(\T)
    =
    \{
        v \in H^1(\Omega)
        \;:\;
        v\vert_{\Gamma_\DIR} \equiv 0
        \text{ in the sense of traces}
    \},
\end{align*}
the piecewise Sobolev spaces allow for non-vanishing interelement jumps.
For all \(v \in L^2(\Omega)\) with Lebesgue square-integrable traces \((v\vert_T)\vert_E \in L^2(E)\) for all \(E \in \mathcal{E}(T)\) and \(T \in \T\),
define the jump \([v]_E \in L^2(E)\) and the average \(\langle v \rangle_E \in L^2(E)\) by
\begin{align*}
    [v]_E
    &\coloneqq
    (v\vert_{T_+})\vert_E
    -
    (v\vert_{T_-})\vert_E,
    &
    \langle v \rangle_E
    &\coloneqq
    \frac12
    \big(
        (v\vert_{T_+})\vert_E
        +
        (v\vert_{T_-})\vert_E
    \big),
    &
    &\text{if } E \in \mathcal{E}(\Omega),
    \\
    [v]_E
    &\coloneqq
    \langle v \rangle_E
    \coloneqq
    (v\vert_{T_+})\vert_E,
    &&&
    &\text{if } E \in \mathcal{E}(\partial\Omega).
\end{align*}

The Friedrichs constant \(C_\FR > 0\) equals the inverse square root of the smallest eigenvalue
of the Laplace operator on \(H^1_\DIR(\Omega)\), i.e.,
\begin{equation}
    \label{eq:Friedrichs}
    C_\FR^{-2}
    \coloneqq
    \min_{v \in H^1_\DIR(\Omega) \setminus \{0\}}
    \frac{\Vert \nabla v \Vert_{L^2(\Omega)}^2}%
    {\Vert v \Vert_{L^2(\Omega)}^2}.
\end{equation}
In the case of an isotropic domain \(\Omega\) with full Dirichlet boundary \(\Gamma_\DIR = \partial\Omega\),
it is well-known that \(C_\FR \approx \diam(\Omega)\).
However, the scaling of \(C_\FR\) is much more involved for general domains and depends on the geometry of \(\Omega\) and \(\Gamma_\DIR\).
Nevertheless, the analysis in this paper employs the (unknown) \(\Omega\)-dependent weighting factor \(c_\Omega \coloneqq C_\FR \)
in the least-squares functional~\eqref{eq:intro_discontinuous_LSFEM} for the sake of simplified equivalence constants.
In practice, upper bounds of the Laplace eigenvalue can be computed by any conforming discretization
and lower bounds by (a postprocessing) of suitable nonconforming methods; see, e.g., \cite{MR3246802,MR4700405}.
Alternatively, the weight \(c_\Omega\) can be chosen by theoretical upper bounds for the Friedrichs constant \(C_\FR\) in practice.
In the case of full Dirichlet boundary \(\Gamma_\DIR = \partial\Omega\), for example, it is well-known that
\(C_\FR \leq \operatorname{width}(\Omega) / \pi\) for the width of the domain \(\Omega\) defined as the smallest possible distance 
of two parallel hyperplanes (lines in 2D) enclosing \(\Omega\)
\begin{equation}
    \label{eq:width}
    \operatorname{width}(\Omega)
    \coloneqq
    \inf
    \left\{
        \ell > 0
        :
        \begin{gathered}
            \exists H_1, H_2 \subseteq \R^2 \text{ lines with }
            \Omega \subseteq \conv(H_1 \cup H_2) \text{ and }
            \\
            \operatorname{dist}(H_1, H_2) \coloneqq \inf\{\vert x_1 - x_2 \vert : x_1 \in H_1, x_2 \in H_2\} = \ell,\,
        \end{gathered}
    \right\}.
\end{equation}
In the general case, the Friedrichs constant \(C_\FR \leq C \diam(\Omega)\) is bounded by the diameter of the domain \(\Omega\)
and a constant \(C > 0\) independent of the scaling of the domain \cite[Lem.~3.30]{MR4242224}.
The piecewise differential operators and the jumps provide the (weighted) norms on \(H(\Div, \T)\) and \(H^1(\T)\),
for \(\tau \in H(\Div, \T)\) and \(v \in H^1(\T)\),
\begin{align*}
    \Vert \tau \Vert_{H(\Div, \T)}^2
    &\coloneqq
    c_\Omega^2\,
    \Vert \Div_\PW \tau \Vert_{L^2(\Omega)}^2
    +
    \Vert \tau \Vert_{L^2(\Omega)}^2,
    \\
    \Vert v \Vert_{H^1(\T)}^2
    &\coloneqq
    \Vert \nabla_\PW v \Vert_{L^2(\Omega)}^2
    +
    j^2(v)
    \quad\text{with}\quad
    j^2(v)
    \coloneqq
    \sum_{E \in \E(\Omega) \cup \E(\Gamma_\DIR)}
    h_E^{-1}\, \Vert [v]_E \Vert_{L^2(E)}^2.
\end{align*}

The traces \((v \vert_T)\vert_{\partial T} \in H^{1/2}(\partial T) \subseteq L^2(\partial T)\) of functions \(v \in H^1(\T)\)
are Lebesgue square-integrable on edges \(E \in \E\) and, thus, the jump term \(j^2(v)\) is well-defined.
The particular weight~\(h_E^{-1}\) therein ensures that the norm \(\Vert \cdot \Vert_{H^1(\T)}\) bounds the \(L^2\) norm
with a generic multiplicative constant independent of the mesh size \cite[Thm.~10.6.12]{MR2373954};
see also Lemma~\ref{lem:U_jump_equivalence} below for the mesh-size-robust relation between \(j^2(v)\) 
and the jumps of the tangential derivatives of \(v \in H^1(\T)\). 
On the contrary, the piecewise normal traces of vector fields \(\tau \in H(\Div, \T)\) are not square-integrable in general,
but \((\tau \vert_T \cdot n_T)\vert_{\partial T} \in H^{-1/2}(\partial T)\).
The same holds for the piecewise tangential derivatives 
\((\nabla v\vert_T \cdot t_T)\vert_{\partial T} \in H^{-1/2}(\partial T)\) of \(v \in H^1(\T)\) as well.
This prevents the measurement of the jumps by \(L^2\) norms on the edges.
Using the duality pairing \(\langle \cdot, \cdot \rangle_{\partial T}\) in \(H^{-1/2}(\partial T)\) and \(H^{1/2}(\partial T)\),
the appropriate terms measuring interelement jumps of piecewise Sobolev functions \(\tau \in H(\Div, \T)\) and \(v \in H^1(\T)\)
read \cite{MR3521055}
\begin{equation}
    \label{eq:jump_terms}
    \begin{split}
        s_n(\tau)
        &\coloneqq
        \sup_{w \in H^1_\DIR(\Omega) \setminus \{0\}}
        \sum_{T \in \T}
        \frac{
            \langle \tau \vert_T \cdot n_T, w \rangle_{\partial T}
        }{
            \Vert \nabla w \Vert_{L^2(\Omega)}
        },
        \\
        s_t(\nabla_\PW v)
        &\coloneqq
        \sup_{w \in H^1_\NEU(\Omega) \setminus \{0\}}
        \sum_{T \in \T}
        \frac{
            \langle
            \nabla_\PW v \vert_T \cdot t_T, w
            \rangle_{\partial T}
        }{
            \Vert \nabla w \Vert_{L^2(\Omega)}
        }.
    \end{split}
\end{equation}
In the case of pure Dirichlet boundary conditions,
the space \(H^1_\NEU(\Omega)\) is replaced by \(\{w \in H^1(\Omega) \,:\, \int_\Omega w \dx = 0\}\).
The jump terms in~\eqref{eq:jump_terms} are consistent in the sense that they vanish for global Sobolev functions.
\begin{lemma}[consistency of jump terms]
    \label{lem:consistency_jump}
    Let \(\sigma, \tau \in H(\Div, \T)\) and \(u, v \in H^1(\T)\).
    If \(\sigma \in H_\NEU(\Div, \Omega)\) and \(u \in H^1_\DIR(\Omega)\),
    then \(s_n(\sigma + \tau) = s_n(\tau)\), \(s_t(\nabla u + \nabla_\PW v) = s_t(\nabla_\PW v)\), and \(j(u) = 0\).
\end{lemma}

\begin{proof}
    It is well-known that any \(u \in H^1(\T)\) belongs to \(H^1(\Omega)\) if and only if 
    the jump \([u]_E \in L^2(E)\) vanishes on all interior edges \(E \in \E(\Omega)\).
    Hence, the Dirichlet boundary conditions ensure \(j(u) = 0\) for \(u \in H^1_\DIR(\Omega)\).
    The duality pairings in the definition of the jump terms~\eqref{eq:jump_terms} 
    cannot be split into edge-oriented jump terms as in \(j(v)\).
    However, the definition of the normal trace operator and a global integration by parts show,
    for \(\sigma \in H_\NEU(\Div, \Omega)\) and \(w \in H^1_\DIR(\Omega)\),
    \begin{align*}
        \sum_{T \in \T}
        \langle \sigma \vert_T \cdot n_T, w \rangle_{\partial T}
        &=
        \sum_{T \in \T}
        \big(
            (\Div \sigma, w)_{L^2(T)}
            +
            (\sigma, \nabla w)_{L^2(T)}
        \big)
        =
        (\Div \sigma, w)_{L^2(\Omega)}
        +
        (\sigma, \nabla w)_{L^2(\Omega)}
        \\
        &=
        \langle \sigma \cdot n, w \rangle_{\partial\Omega}
        =
        0
    \end{align*}
    Since this holds for all \(w \in H^1_\DIR(\Omega)\), it follows that \(s_n(\sigma + \tau) = s_n(\tau)\).
    Analogously, for \(u \in H^1_\DIR(\Omega)\) and \(w \in H^1_\NEU(\Omega)\),
    \begin{align*}
        \sum_{T \in \T}
        \langle \nabla u\vert_T \cdot t_T, w \rangle_{\partial T}
        &=
        \sum_{T \in \T}
        \big(
            (\curl \nabla u, w)_{L^2(T)}
            +
            (\nabla u, \Curl w)_{L^2(T)}
        \big)
        \\
        &=
        (\curl \nabla u, w)_{L^2(\Omega)}
        +
        (\nabla u, \Curl w)_{L^2(\Omega)}
        =
        \langle \nabla u \cdot t, w \rangle_{\partial\Omega}
        =
        0
    \end{align*}
    and, thus, \(s_t(\nabla u + \nabla_\PW v) = s_t(\nabla_\PW v)\).
\end{proof}

Abbreviate the weighted norm in the product space \(H(\Div, \T) \times H^1(\T)\) by
\begin{equation}
    \label{eq:weighted_norm}
    \Vert (\tau, v) \Vert_\T^2
    \coloneqq
    \Vert \tau \Vert_{H(\Div, \Omega)}^2
    +
    s_n^2(\tau)
    +
    \Vert v \Vert_{H^1(\T)}^2
    +
    s_t^2(\nabla_\PW v).
\end{equation}
The following lemma justifies that the generalized jump terms \(s_n(\tau)\) and \(s_t(\nabla_\PW v)\)
are suitable for the norm on \(H(\Div, \T) \times H^1(\T)\).
\begin{lemma}[equivalence of norms]
    \label{lem:norm_equivalence}
    For all \(\tau \in H(\Div, \T)\) and \(v \in H^1(\T)\),
    \[
        \Vert \tau \Vert_{H(\Div, \T)}^2
        +
        \Vert v \Vert_{H^1(\T)}^2
        \leq
        \Vert (\tau, v) \Vert_\T^2
        \leq
        3\,
        \big(
            \Vert \tau \Vert_{H(\Div, \T)}^2
            +
            \Vert v \Vert_{H^1(\T)}^2
        \big).
    \]
\end{lemma}

\begin{proof}
    The definition of the trace operator, a Cauchy--Schwarz inequality, and the Friedrichs inequality show,
    for \(\tau \in H(\Div, \T)\) and \(w \in H^1_\DIR(\Omega)\) with \(\Vert \nabla w \Vert_{L^2(\Omega)} = 1\), that
    \begin{align*}
        \sum_{T \in \T}
        \langle
            \tau \vert_T \cdot n_T, w
        \rangle_{\partial T}
        &=
        \sum_{T \in \T}
        \big(
             (\Div \tau, w)_{L^2(T)}
             +
             (\tau, \nabla w)_{L^2(T)}
        \big)
        =
        (\Div_\PW \tau, w)_{L^2(\Omega)}
        +
        (\tau, \nabla w)_{L^2(\Omega)}
        \\
        &\leq
        \Vert \Div_\PW \tau \Vert_{L^2(\Omega)}
        \Vert w \Vert_{L^2(\Omega)}
        +
        \Vert \tau \Vert_{L^2(\Omega)}
        \\
        &\leq
        C_\FR\,
        \Vert \Div_\PW \tau \Vert_{L^2(\Omega)}
        +
        \Vert \tau \Vert_{L^2(\Omega)}
        \leq
        \sqrt{2}\,
        \Vert \tau \Vert_{H(\Div, \T)}.
    \end{align*}
    Since this holds for all such \(w \in H^1_\DIR(\Omega)\), it follows that \(s_n(\tau) \leq \sqrt{2}\, \Vert \tau \Vert_{H(\Div, \T)}\).
    Analogously, for \(v \in H^1(\T)\) and \(w \in H^1_\NEU(\T)\) with \(\Vert \nabla w \Vert_{L^2(\Omega)} = 1\), the estimate
    \begin{align*}
        \sum_{T \in \T}
        \langle
            \nabla_\PW v \vert_T \cdot t_T, w
        \rangle_{\partial T}
        &=
        \sum_{T \in \T}
        (\nabla v, \Curl w)_{L^2(T)}
        =
        (\nabla_\PW v, \Curl w)_{L^2(\Omega)}
        \\
        &\leq
        \Vert \nabla_\PW v \Vert_{L^2(\Omega)}
        \Vert \Curl w \Vert_{L^2(\Omega)}
        =
        \Vert \nabla_\PW v \Vert_{L^2(\Omega)}
        \Vert \nabla w \Vert_{L^2(\Omega)}
        =
        \Vert \nabla_\PW v \Vert_{L^2(\Omega)}
    \end{align*}
    implies \(s_t(\nabla_\PW v) \leq \Vert \nabla_\PW v \Vert_{L^2(\Omega)}\).
\end{proof}

\section{Least-squares principle for piecewise Sobolev spaces}
\label{sec:discontinuous_least_squares}

This section introduces a consistent generalization of the least-squares principle~\eqref{eq:least_squares_principle}
to piecewise Sobolev functions for the weighted functional \(LS(f)\colon H(\Div, \T) \times H^1(\T) \to \R\) defined,
for \(\sigma \in H(\Div, \T)\) and \(u \in H^1(\T)\), by
\[
    LS(f; \sigma, u)
    \coloneqq
    c_\Omega^2\,
    \Vert f + \Div_\PW \sigma \Vert_{L^2(\Omega)}^2
    +
    \Vert \sigma - \nabla_\PW u \Vert_{L^2(\Omega)}^2
    +
    s_n^2(\sigma)
    +
    j^2(u)
    +
    s_t^2(\nabla_\PW u).
\]
The weights in the norm and the least-squares functional ensure 
that the equivalence constants in the following main result do not depend on the size of the domain.
\begin{theorem}[fundamental equivalence]
    \label{thm:fundamental_equivalence}
    For all \(\tau \in H(\Div, \T)\) and \(v \in H^1(\T)\),
    \begin{equation}
        \label{eq:fundamental_equivalence}
        \frac18\,
        \Vert (\tau, v) \Vert_\T^2
        \leq
        LS(0; \tau, v)
        \leq
        2\,
        \Vert (\tau, v) \Vert_\T^2.
    \end{equation}
\end{theorem}
Its proof relies on stability estimates for the standard Helmholtz decomposition applied to piecewise \(H(\Div)\) functions.
It characterizes \(L^2\) vector fields as a gradient and a Curl part.
The additional piecewise \(H(\Div)\) regularity allows for the desired stability estimates.
\begin{lemma}[Helmholtz decomposition for piecewise Sobolev functions]
    \label{lem:Helmholtz}
    Assume that the bounded polygonal Lipschitz domain \(\Omega \subseteq \R^2\) is simply connected and recall \(c_\Omega = C_\FR\).
    Given \(\tau \in L^2(\Omega; \R^2)\), there exist \(a \in H^1_\DIR(\Omega)\) and \(b \in H^1_\NEU(\Omega)\)
    such that the decomposition \(\tau = \nabla a + \Curl b\) is \(L^2\) orthogonal.
    For all \(v \in H^1(\T)\), the following quasi-orthogonality holds
    \begin{equation}
        \label{eq:Helmholtz_quasi_orthogonality}
        (\nabla_\PW v, \Curl b)_{L^2(\Omega)}
        \leq
        s_t(\nabla_\PW v)
        \Vert \tau \Vert_{L^2(\Omega)}.
    \end{equation}
    If \(\tau \in H(\Div, \T)\), then
    \begin{equation}
        \label{eq:Helmholtz_stability}
        \Vert \nabla a \Vert_{L^2(\Omega)}
        \leq
        C_\FR\,
        \Vert \Div_\PW \tau \Vert_{L^2(\Omega)}
        +
        s_n(\tau).
    \end{equation}
\end{lemma}

\begin{proof}
    Given a function \(\tau \in L^2(\Omega; \R^2)\), recall that \(a \in H^1_\DIR(\Omega)\) is determined by
    \[
        (\nabla a, \nabla w)_{L^2(\Omega)}
        =
        (\tau, \nabla w)_{L^2(\Omega)}
        \quad\text{for all }
        w \in H^1_\DIR(\Omega).
    \]
    Hence, the weak divergence of the difference \(\tau - \nabla a\) vanishes
    and there exists \(b \in H^1_\NEU(\Omega)\) such that \(\tau - \nabla a = \Curl b\) \cite[Thm.~3.2]{MR851383}.
    The orthogonality follows from the integration by parts
    \begin{align*}
        (\nabla a, \Curl b)_{L^2(\Omega)}
        &=
        -(\curl \nabla a, b)_{L^2(\Omega)}
        +
        \langle \nabla a \cdot t, b \rangle_{\partial\Omega}
        =
        0.
    \end{align*}
    Another integration by parts proves, for all \(v \in H^1(\T)\), that
    \begin{align*}
        (\nabla_\PW v, \Curl b)_{L^2(\Omega)}
        &=
        \sum_{T \in \T}
        (\nabla v, \Curl b)_{L^2(T)}
        =
        \sum_{T \in \T}
        \Big(
            -(\curl \nabla v, b)_{L^2(T)}
            + \langle \nabla v \cdot t_T, b \rangle_{\partial T}
        \Big)
        \\
        &=
        \sum_{T \in \T}
        \langle \nabla v \cdot t_T, b \rangle_{\partial T}
        \leq
        s_t(\nabla_\PW v) \,
        \Vert \nabla b \Vert_{L^2(\Omega)}
        =
        s_t(\nabla_\PW v) \,
        \Vert \Curl b \Vert_{L^2(\Omega)}.
    \end{align*}
    This and the estimate
    \(\Vert \Curl b \Vert_{L^2(\Omega)}^2 = (\tau, \Curl b)_{L^2(\Omega)} \leq \Vert \tau \Vert_{L^2(\Omega)} \Vert \Curl b \Vert_{L^2(\Omega)}\)
    verifies the quasi-orthogonality estimate~\eqref{eq:Helmholtz_quasi_orthogonality}.
    Assume that \(\tau \in H(\Div, \T)\).
    A piecewise integration by parts, the Cauchy--Schwarz inequality, and the Friedrichs inequality~\eqref{eq:Friedrichs}
    prove the stability estimate~\eqref{eq:Helmholtz_stability}
    \begin{align*}
        \Vert \nabla a \Vert_{L^2(\Omega)}^2
        &=
        (\tau, \nabla a)_{L^2(\Omega)}
        =
        \sum_{T \in \T}
        (\tau, \nabla a)_{L^2(T)}
        =
        - (\Div_\PW \tau, a)_{L^2(\Omega)}
        +
        \sum_{T \in \T}
        \langle \tau \cdot n_T, a \rangle_{\partial T}
        \\
        &\leq
        \Vert \Div_\PW \tau \Vert_{L^2(\Omega)}
        \Vert a \Vert_{L^2(\Omega)}
        +
        s_n(\tau)\, \Vert \nabla a \Vert_{L^2(\Omega)}
        \\
        &\leq
        \big(
            C_\FR\,
            \Vert \Div_\PW \tau \Vert_{L^2(\Omega)}
            +
            s_n(\tau)
        \big)
        \Vert \nabla a \Vert_{L^2(\Omega)}
    \end{align*}
    and conclude the proof.
\end{proof}

\begin{proof}[Proof of Theorem~\ref{thm:fundamental_equivalence}]
    \emph{Step 1}.\enskip
    The proof of the stability \(LS(0; \tau, v) \leq 2\, \Vert (\tau, v) \Vert_\T^2\)
    employs the triangle and Young inequality in
    \(
        \Vert \tau - \nabla_\PW v \Vert_{L^2(\Omega)}^2
        \leq
        2\,
        (
            \Vert \tau \Vert_{L^2(\Omega)}^2
            +
            \Vert \nabla_\PW v \Vert_{L^2(\Omega)}^2
        )
    \).
    The remaining terms coincide on both sides of the stability estimate.

    \emph{Step 2}.\enskip
    The proof of ellipticity \(\Vert (\tau, v) \Vert_\T^2 \leq 8 \, LS(0; \tau, v)\) departs with
    the Helmholtz decomposition \(\tau = \nabla a + \Curl b\) for \(a \in H^1_\DIR(\Omega)\) and \(b \in H^1_\NEU(\Omega)\)
    from Lemma~\ref{lem:Helmholtz}.
    This plus an algebraic identity result in
    \begin{equation}
        \label{eq:Helmholtz_split}
        \begin{split}
            \Vert \tau \Vert_{L^2(\Omega)}^2
            +
            \Vert \nabla_\PW v \Vert_{L^2(\Omega)}^2
            &=
            \Vert \tau - \nabla_\PW v \Vert_{L^2(\Omega)}^2
            + 2\, (\nabla_\PW v, \tau)_{L^2(\Omega)}
            \\
            &=
            \Vert \tau - \nabla_\PW v \Vert_{L^2(\Omega)}^2
            + 2\, (\nabla_\PW v, \nabla a)_{L^2(\Omega)}
            + 2\, (\nabla_\PW v, \Curl b)_{L^2(\Omega)}.
        \end{split}
    \end{equation}
    The weighted Young inequality and the stability estimate~\eqref{eq:Helmholtz_stability} from Lemma~\ref{lem:Helmholtz} prove
    \begin{align*}
        (\nabla_\PW v, \nabla a)_{L^2(\Omega)}
        &\leq
        \Vert \nabla_\PW v \Vert_{L^2(\Omega)}
        \Vert \nabla a \Vert_{L^2(\Omega)}
        \leq
        \frac14
        \Vert \nabla_\PW v \Vert_{L^2(\Omega)}^2
        +
        \Vert \nabla a \Vert_{L^2(\Omega)}^2
        \\
        &\leq
        \frac14
        \Vert \nabla_\PW v \Vert_{L^2(\Omega)}^2
        +
        2 \,
        \big(
            c_\Omega^2\,
            \Vert \Div_\PW \tau \Vert_{L^2(\Omega)}^2
            +
            s_n^2(\tau)
        \big).
    \end{align*}
    The quasi-orthogonality estimate~\eqref{eq:Helmholtz_quasi_orthogonality} from Lemma~\ref{lem:Helmholtz} and another Young inequality show
    \[
        (\nabla_\PW v, \Curl b)_{L^2(\Omega)}
        \leq
        s_t(\nabla_\PW v)
        \Vert \tau \Vert_{L^2(\Omega)}
        \leq
        s_t^2(\nabla_\PW v)
        +
        \frac{1}{4}
        \Vert \tau \Vert_{L^2(\Omega)}^2.
    \]
    The combination of the two previously displayed formulas with identity~\eqref{eq:Helmholtz_split} reads
    \begin{align*}
        \Vert \tau \Vert_{L^2(\Omega)}^2
        +
        \Vert \nabla_\PW v \Vert_{L^2(\Omega)}^2
        &=
        \Vert \tau - \nabla_\PW v \Vert_{L^2(\Omega)}^2
        + 2\, (\nabla_\PW v, \tau)_{L^2(\Omega)}
        \\
        &\leq
        \Vert \tau - \nabla_\PW v \Vert_{L^2(\Omega)}^2
        +
        4 \,
        \big(
            c_\Omega^2\,
            \Vert \Div_\PW \tau \Vert_{L^2(\Omega)}^2
            +
            s_n^2(\tau)
        \big)
        \\
        &\phantom{{}\leq{}}
        +
        2 \,
        s_t^2(\nabla_\PW v)
        +
        \frac12
        \Vert \nabla_\PW v \Vert_{L^2(\Omega)}^2
        +
        \frac12
        \Vert \tau \Vert_{L^2(\Omega)}^2.
    \end{align*}
    The absorption of \(\Vert \nabla_\PW v \Vert_{L^2(\Omega)}^2\) and \(\Vert \tau \Vert_{L^2(\Omega)}^2\) on the left-hand side
    concludes the proof of the ellipticity.
\end{proof}

\begin{remark}[scaling-robust conforming LSFEM]
    For conforming functions \(\sigma \in H_\NEU(\Div, \Omega) \subseteq H(\Div, \T)\) and \(u \in H^1_\DIR(\Omega) \subseteq H^1(\T)\),
    the consistency of \(s_n\) and \(s_t\) from Lemma~\ref{lem:consistency_jump} ensures 
    that all jump contributions vanish in~\eqref{eq:fundamental_equivalence}.
    This results in the fundamental equivalence~\eqref{eq:least_squares_principle} with \(\Omega\)-independent generic constants
    and corresponding a~priori and a~posteriori error estimates.
    To the best of the author's knowledge, this scaling-robust conforming LSFEM for the Poisson model problem
    has not been stated explicitly in the literature.
    However, the reader is referred to \cite[Thm.~3.1]{MR3715170} for a similar approach applied
    to the conforming LSFEM for the Stokes equations.
\end{remark}

\begin{remark}[three space dimensions]
    \label{rem:3D_Helmholtz}
    A Helmholtz decomposition in three dimensions allows to generalize Theorem~\ref{thm:fundamental_equivalence}
    to contractible domains \(\Omega \subseteq \R^3\).
    However, this leads to the Sobolev space \(H_\NEU(\curl, \Omega)\) in the definition of \(s_t(\nabla_\PW v)\).
    The remaining arguments in this section apply analogously.
\end{remark}

\section{Computable jump terms for piecewise polynomials}
\label{sec:piecewise_polynomials}

Given any subset \(\omega \subseteq \R^2\) and \(k \in \N_0\), let \(P^k(\omega)\) denote the space of polynomials on
\(\omega\) of total degree at most \(k\).
For a regular triangulation \(\T\) of \(\Omega\), define the set of piecewise polynomials as
\[
    P^k(\T)
    \coloneqq
    \{
        v_h \in L^\infty(\Omega)
        \;:\;
        \forall T \in \T,\:
        v_h\vert_{T} \in P^k(T)
    \}
\]
and analogously for the vector-valued polynomials \(P^k(\T; \R^2)\).
Appropriate indices specify partial homogeneous boundary conditions \(P^k_\DIR(\T) \coloneqq P^k(\T) \cap H^1_\DIR(\T)\)
and analogously for \(P^k_\NEU(\T)\).
For \(k \geq 1\), the conforming counterparts read \(S^k(\T) \coloneqq P^k(\T) \cap H^1(\Omega)\),
\(S^k_\DIR(\T) \coloneqq P^k(\T) \cap H^1_\DIR(\Omega)\), and analogously for \(S^k_\NEU(\T)\).

The jump terms \(s_n\) and \(s_t\) from~\eqref{eq:jump_terms} cannot be computed exactly and
need to be replaced by equivalent computable terms for discrete functions in order to derive a discontinuous LSFEM.
This bases on the following lemma rewriting the \(H^{-1/2}\) norm for piecewise polynomial functions
on the skeleton as weighted least-squares residuals.
For an alternative approach, the reader is referred to \cite{MR4708036,MR4803195}
for residual minimization in fractional or negative Sobolev norms.
\begin{lemma}[equivalence of edge-oriented measures]
    \label{lem:sup_control}
    For any edgewise polynomial \(q \in P^k(\E) \coloneqq \prod_{E \in \E} P^k(E)\) 
    with \(\int_E q \ds = 0\) for all \(E \in \E(\Omega) \cup \E(\Gamma_\NEU)\), it holds that
    \begin{equation}
        \label{eq:sup_control}
        \sup_{w \in H^1_\DIR(\Omega) \setminus \{0\}}
        \sum_{E \in \E(\Omega) \cup \E(\Gamma_\NEU)}
        \frac{(q, w)_{L^2(E)}}{\Vert \nabla w \Vert_{L^2(\Omega)}}
        \approx
        \bigg(
            \sum_{E \in \E(\Omega) \cup \E(\Gamma_\NEU)}
            h_E\:
            \Vert q \Vert_{L^2(E)}^2
        \bigg)^{1/2}.
    \end{equation}
    The equivalence constants depend on the polynomial degree \(k\) and the interior angles in \(\T\),
    but not on the size of the domain \(\Omega\).
\end{lemma}

\begin{remark}[generalization of upper bound]
    The proof of equivalence~\eqref{eq:sup_control} below reveals that the upper bound \(\lesssim\)
    even holds for arbitrary \(q \in L^2(\E) \coloneqq \prod_{E \in \E} L^2(E)\)
    with \(\int_E q \ds = 0\) for all \(E \in \E(\Omega) \cup \E(\Gamma_\NEU)\).
\end{remark}

\begin{remark}[necessity of a side condition]
    \label{rem:counter_example}
    The following counterexample justifies that the side condition enforcing integral mean zero
    on the interior and Neumann edges cannot be neglected.
    To illustrate this, let \(\T\) be a quasi-uniform mesh, i.e.,
    there exists \(h > 0\) such that \(h_E \approx h\) for all \(E \in \E\).
    For \(q \equiv 1 \in P^0(\E)\), the right-hand side reads
    \[
        \bigg(
            \sum_{E \in \E(\Omega) \cup \E(\Gamma_\NEU)}
            h_E \Vert q \Vert_{L^2(E)}^2
        \bigg)^{1/2}
        =
        \bigg(
            \sum_{E \in \E(\Omega) \cup \E(\Gamma_\NEU)}
            h_E^2
        \bigg)^{1/2}
        \approx
        \vert \E(\Omega) \cup \E(\Gamma_\NEU) \vert^{1/2}
        \,h
        \approx
        1.
    \]
    Assume that there exists an edge 
    \(E \in \E' \coloneqq \{ E \in \E \;:\; \dist(E, \Gamma_\DIR) > 0 \} \subsetneq \E(\Omega) \cup \E(\Gamma_\NEU)\)
    with positive distance \(\dist(E, \Gamma_\DIR) > 0\) to the Dirichlet boundary.
    Define \(w_h \in S^1_\DIR(\T) \subseteq H^1_\DIR(\Omega)\) by
    \(w_h(z) \coloneqq 1\) for all \(z \in \V(\Omega) \cup \V(\Gamma_\NEU)\)
    and \(w_h(z) \coloneqq 0\) for all \(z \in \V(\Gamma_\DIR)\).
    Then \(w_h \vert_E \equiv 1\) on exactly all edges \(E \in \E'\) and
    \[
        \Vert \nabla w_h \Vert_{L^2(\Omega)}^2
        =
        \sum_{
            \substack{T \in \T \\
            \V(T) \cap \V(\Gamma_\DIR) \neq \emptyset}
        }
        \Vert \nabla w_h \Vert_{L^2(T)}^2
        \approx
        \vert
            \{
                T \in \T \;:\;
                \V(T) \cap \V(\Gamma_\DIR) \neq \emptyset
            \}
        \vert
        \approx
        \vert \V(\Gamma_\DIR) \vert
        \approx
        h^{-1}.
    \]
    The insertion of \(w_h\) into the left-hand side of~\eqref{eq:sup_control} provides a lower bound of the supremum, namely
    \[
        \sum_{E \in \E(\Omega) \cup \E(\Gamma_\NEU)}
        \frac{(q, w_h)_{L^2(E)}}{\Vert \nabla w_h \Vert_{L^2(\Omega)}}
        \geq
        \!
        \sum_{E \in \E'}
        \frac{(q, w_h)_{L^2(E)}}{\Vert \nabla w_h \Vert_{L^2(\Omega)}}
        =
        \!
        \sum_{E \in \E'}
        \frac{h_E}{\Vert \nabla w_h \Vert_{L^2(\Omega)}}
        \approx
        \!
        \sum_{E \in \E'}
        h_E^{3/2}
        =
        \vert \E' \vert \,
        h^{3/2}
        \approx
        h^{-1/2}.
    \]
    Asymptotically, this leads to a contradiction and illustrates that the equivalence~\eqref{eq:sup_control}
    cannot hold for arbitrary \(q \in P^k(\E)\).
\end{remark}

\begin{proof}[Proof of Lemma~\ref{lem:sup_control}]
    \emph{Step 1}.\enskip
    For the proof of the upper bound ``\(\lesssim\)'' in~\eqref{eq:sup_control},
    let \(w \in H^1_\DIR(\Omega)\) be arbitrary with \(\Vert \nabla w \Vert_{L^2(\Omega)} = 1\).
    For an edge \(E \in \E(\Omega) \cup \E(\Gamma_\NEU)\) with adjacent triangle \(T \in \T(E)\),
    let \(w_E \coloneqq \int_E w \ds / h_E \in \R\) abbreviate the integral mean of \(w\) on \(E\).
    Using the characteristic function \(\chi_T \in P^0(\T)\) of the triangle \(T\) with \(\chi_T\vert_T \equiv 1\) and 
    \(\chi_T\vert_{T'} \equiv 0\) for all \(T' \in \T \setminus \{T\}\),
    this integral mean can be considered as a function \(w_E \chi_T \in P^0(\T)\).
    The trace inequality and a Poincar\'e inequality show
    \begin{equation}
        \label{eq:trace_Poincare}
        h_E^{-1/2}\:
        \Vert w - w_E \Vert_{L^2(E)}
        \lesssim
        \Vert \nabla (w - w_E \chi_T) \Vert_{L^2(T)}
        +
        h_E^{-1}\:
        \Vert w - w_E \chi_T \Vert_{L^2(T)}
        \lesssim
        \Vert \nabla w \Vert_{L^2(T)}.
    \end{equation}
    This, the integral mean property \(\int_E q \ds = 0\) for all \(E \in \E(\Omega) \cup \E(\Gamma_\NEU)\),
    and a Cauchy--Schwarz inequality in \(\R^N\) for \(N \coloneqq \vert \E(\Omega) \cup \E(\Gamma_\NEU) \vert\) prove
    \begin{align*}
        \sum_{E \in \E(\Omega) \cup \E(\Gamma_\NEU)}
        (q, w)_{L^2(E)}
        &=
        \sum_{E \in \E(\Omega) \cup \E(\Gamma_\NEU)}
        (q, w - w_E)_{L^2(E)}
        \\
        &\leq
        \sum_{E \in \E(\Omega) \cup \E(\Gamma_\NEU)}
        h_E^{1/2}\:
        \Vert q \Vert_{L^2(E)}
        \:
        h_E^{-1/2}\:
        \Vert w - w_E \Vert_{L^2(E)}
        \\
        &\lesssim
        \sum_{E \in \E(\Omega) \cup \E(\Gamma_\NEU)}
        h_E^{1/2}\:
        \Vert q \Vert_{L^2(E)}
        \:
        \Vert \nabla w \Vert_{L^2(\omega_E)}
        \\
        &\leq
        \bigg(
            \sum_{E \in \E(\Omega) \cup \E(\Gamma_\NEU)}
            h_E\:
            \Vert q \Vert_{L^2(E)}^2
        \bigg)^{1/2}
        \bigg(
            \sum_{E \in \E(\Omega) \cup \E(\Gamma_\NEU)}
            \Vert \nabla w \Vert_{L^2(\omega_E)}^2
        \bigg)^{1/2}.
    \end{align*}
    The finite overlap of the edge patches \(\omega_E\) solely depends on the spatial dimension and shows
    \[
        \sum_{E \in \E(\Omega) \cup \E(\Gamma_\NEU)}
        \Vert \nabla w \Vert_{L^2(\omega_E)}^2
        \lesssim
        \Vert \nabla w \Vert_{L^2(\Omega)}^2
        =
        1.
    \]
    This concludes the proof of the upper bound with a generic constant involving the interior angles in \(\T\) from the
    estimate~\eqref{eq:trace_Poincare} only.

    \emph{Step 2}.\enskip
    The proof of the lower bound ``\(\gtrsim\)'' in~\eqref{eq:sup_control} employs the bubble-function technique by Verf\"urth.
    For \(E \in \E(\Omega) \cup \E(\Gamma_\NEU)\), the edge-bubble function
    \(0 \leq b_E \coloneqq 4\prod_{z \in \V(E)} \varphi_z \in S^2_\DIR(\T)\)
    satisfies \(\supp(b_E) = \overline\omega_E\),
    \begin{equation}
        \label{eq:bubble_estimates}
        \Vert q \Vert_{L^2(E)}^2
        \lesssim
        (q, q b_E)_{L^2(E)},
        \quad\text{and}\quad
        \Vert \nabla (q\, b_E) \Vert_{L^2(\omega_E)}^2
        \lesssim
        h_E^{-1}\:
        \Vert q \Vert_{L^2(E)}^2.
    \end{equation}
    Hence, for \(v \coloneqq \sum_{E \in \E(\Omega) \cup \E(\Gamma_\NEU)} h_E\,q\,b_E \in H^1_\DIR(\Omega)\),
    \begin{align*}
        \sum_{E \in \E(\Omega) \cup \E(\Gamma_\NEU)}
        h_E\:
        \Vert q \Vert_{L^2(E)}^2
        &\lesssim
        \sum_{E \in \E(\Omega) \cup \E(\Gamma_\NEU)}
        (q, v)_{L^2(E)}
        \\
        &\leq
        \Vert \nabla v \Vert_{L^2(\Omega)}
        \sup_{w \in H^1_\DIR(\Omega) \setminus \{0\}}
        \sum_{E \in \E(\Omega) \cup \E(\Gamma_\NEU)}
        \frac{(q, w)_{L^2(E)}}{\Vert \nabla w \Vert_{L^2(\Omega)}}.
    \end{align*}
    The local support \(\omega_E\) of the bubble functions \(b_E\) enables the localization in the estimate
    \[
        \Vert \nabla v \Vert_{L^2(\Omega)}
        \lesssim
        \bigg(
            \sum_{E \in \E(\Omega) \cup \E(\Gamma_\NEU)}
            \hspace{-3pt}
            h_E^2\:
            \Vert \nabla (q\,b_E) \Vert_{L^2(\omega_E)}^2
        \bigg)^{1/2}
        \lesssim
        \bigg(
            \sum_{E \in \E(\Omega) \cup \E(\Gamma_\NEU)}
            \hspace{-3pt}
            h_E\:
            \Vert q \Vert_{L^2(E)}^2
        \bigg)^{1/2}.
    \]
    The constants in the local estimates~\eqref{eq:bubble_estimates} include the polynomial degree \(k\) and the interior angles in \(\T\).
    The global overlapping argument of the patches \(\omega_E\) depend on the spatial dimension only.
    Hence, the resulting lower bound is independent of the size of the domain.
\end{proof}

The side condition in Lemma~\ref{lem:sup_control} motivates the inclusion of the integral mean property
for the jumps into the discrete subspace \(\Sigma(\T)\) of \(H(\Div, \T)\) such that
\begin{equation}
    \label{eq:Sigma_condition}
    \Sigma(\T)
    \subseteq
    \left\{
        \tau_h \in H(\Div, \T)
        \;:\;
        \begin{aligned}
            &\forall T \in \T\: \forall E \in \E(T),\:
            (\tau_h \vert_T \cdot n_T)\vert_E \in P^{k+1}(E)
            \text{ and }
            \\
            &\forall E \in \E(\Omega) \cup \E(\Gamma_\NEU),\:
            \textstyle
            \int_E [\tau_h \cdot n_E]_E \ds = 0
        \end{aligned}
    \right\}.
\end{equation}
For \(\tau_h \in \Sigma(\T)\), the reformulation
\[
    s_n^2(\tau_h)
    =
    \sup_{w \in H^1_\DIR(\Omega) \setminus \{0\}}
    \sum_{T \in \T}
    \frac{(\tau_h \vert_T \cdot n_T, w)_{L^2(\partial T)}}{\Vert \nabla w \Vert_{L^2(\Omega)}}
    =
    \sup_{w \in H^1_\DIR(\Omega) \setminus \{0\}}
    \sum_{E \in \E(\Omega) \cup \E(\Gamma_\NEU)}
    \frac{([\tau_h \cdot n_E]_E, w)_{L^2(E)}}{\Vert \nabla w \Vert_{L^2(\Omega)}}
\]
and Lemma~\ref{lem:sup_control} immediately result in
the following equivalence.
\begin{corollary}
    \label{cor:Sigma_jump_equivalence}
    For every
    \(\tau_h \in \Sigma(\T)\),
    \begin{align*}
        s_n^2(\tau_h)
        &\approx
        \sum_{E \in \E(\Omega) \cup \E(\Gamma_\NEU)}
        h_E\:
        \Vert [\tau_h \cdot n_E]_E \Vert_{L^2(E)}^2.
    \end{align*}
\end{corollary}

In addition to~\eqref{eq:Sigma_condition},
consider the discrete subspace \(U(\T)\) of \(H^1(\T)\) with polynomial traces such that
\begin{equation}
    \label{eq:U_condition}
    U(\T)
    \subseteq
    \big\{
        v_h \in H^1(\T)
        \;:\;
        \forall T \in \T\: \forall E \in \E(T),\:
        (v_h\vert_T)\vert_E \in P^{k+1}(E)
    \big\}.
\end{equation}
The control of the jump term \(s_t^2(\nabla_\PW v_h)\) for \(v_h \in U(\T)\) employs the following quasi-interpolation.
\begin{lemma}
    \label{lem:quasi_interpolation}
    Given \(w \in H^1(\Omega)\), there exists \(w_h \in S^{k+2}(\T)\) satisfying, for all \(E \in \E\) and \(v_h \in P^k(E)\),
    \begin{equation}
        \label{eq:quasi_interpolation}
        (w - w_h, v_h)_{L^2(E)}
        =
        0
        \quad\text{and}\quad
        \Vert \nabla w_h \Vert_{L^2(\Omega)}
        \lesssim
        \Vert \nabla w \Vert_{L^2(\Omega)}.
    \end{equation}
    The generic constant in the stability estimate solely depends on the polynomial degree \(k\) and the interior angles in \(\T\),
    but not on the size of the domain \(\Omega\).
    Moreover, if \(w \vert_E \equiv 0\) vanishes along a Neumann boundary edge \(E \in \E(\Gamma_\NEU)\), then \(w_h \vert_E \equiv 0\) also vanishes.
\end{lemma}

The author expects the result from Lemma~\ref{lem:quasi_interpolation} to be well-known;
see, e.g., the interpolation operator in \cite[Sect.~3.2]{MR3696081} for continuous arguments \(w \in C^0(\overline\Omega)\).
Nevertheless, a direct proof is given in Appendix~\ref{app:quasi_interpolation} for the sake of the explicit dependencies
of the generic constants.

\begin{lemma}
    \label{lem:U_jump_equivalence}
    For all \(v_h \in U(\T)\), it holds that \(j^2(v_h) + s_t^2(\nabla_\PW v_h) \approx j^2(v_h)\).
    The equivalence constants depend on the polynomial degree \(k\) and the interior angles in \(\T\),
    but not on the size of the domain \(\Omega\).
\end{lemma}

\begin{proof}
    Given \(w \in H^1_\NEU(\Omega)\), let \(w_h \in S^{k + 2}_\NEU(\T)\) denote
    the quasi-interpolation of \(w\) from Lemma~\ref{lem:quasi_interpolation}.
    For each edge \(E \in \E\), let \(z_j(E) \in \V\) for \(j = 1,2\) denote
    the two vertices of \(E = \conv\{z_1(E), z_2(E)\}\).
    For all \(T \in \T\) and \(w \in H^1_\NEU(\Omega)\) with \(\Vert \nabla w \Vert_{L^2(\Omega)} = 1\),
    the integration by parts for the arc-length derivative shows
    \begin{equation}
        \label{eq:ibp_arc_length}
        \begin{aligned}
            \langle \nabla v_h \vert_T \cdot t_T, w \rangle_{\partial T}
            &=
            (\nabla v_h \vert_T \cdot t_T, w_h)_{L^2(\partial T)}
            +
            (\nabla v_h \vert_T \cdot t_T, w - w_h)_{L^2(\partial T)}
            \\
            &=
            (\nabla v_h \vert_T \cdot t_T, w_h)_{L^2(\partial T)}
            =
            -(v_h \vert_T, \nabla w_h \cdot t_T)_{L^2(\partial T)}.
        \end{aligned}
    \end{equation}
    The sum over every \(T \in \T\) reads
    \begin{align*}
        \sum_{T \in \T}
        \langle \nabla v_h\vert_T \cdot t_T, w \rangle_{\partial T}
        &=
        - \sum_{T \in \T}
        (v_h\vert_T, \nabla w_h \cdot t_T)_{L^2(\partial T)}
        =
        - \sum_{E \in \E(\Omega) \cup \E(\Gamma_\DIR)}
        ([v_h]_E, \nabla w_h \cdot t_E)_{L^2(E)}
        \\
        &\leq
        \sum_{E \in \E(\Omega) \cup \E(\Gamma_\DIR)}
        h_E^{-1/2}
        \Vert [v_h]_E \Vert_{L^2(E)}
        \,
        h_E^{1/2}
        \Vert \nabla w_h \cdot t_E \Vert_{L^2(E)}
        \\
        &\leq
        \bigg(
            \sum_{E \in \E(\Omega) \cup \E(\Gamma_\DIR)}
            \hspace{-6pt}
            h_E^{-1}
            \Vert [v_h]_E \Vert_{L^2(E)}^2
        \bigg)^{1/2}
        \bigg(
            \sum_{E \in \E(\Omega) \cup \E(\Gamma_\DIR)}
            \hspace{-6pt}
            h_E\,
            \Vert \nabla w_h \cdot t_E \Vert_{L^2(E)}^2
        \bigg)^{1/2}.
    \end{align*}
    The discrete trace inequality \(h_E^{1/2} \Vert \nabla w_h \cdot t_E \Vert_{L^2(E)} \lesssim \Vert \nabla w_h \Vert_{L^2(\omega_E)}\)
    and the finite overlap of the patches \(\omega_E\) for all \(E \in \E(\Omega) \cup \E(\Gamma_\DIR)\), lead to
    \[
        \sum_{T \in \T}
        \langle \nabla v_h\vert_T \cdot t_T, w \rangle_{\partial T}
        \lesssim
        j(v_h)
        \Vert \nabla w_h \Vert_{L^2(\Omega)}
        \lesssim
        j(v_h)
        \Vert \nabla w \Vert_{L^2(\Omega)}
        =
        j(v_h).
    \]
    The supremum over all \(w \in H^1_\NEU(\Omega)\) with \(\Vert \nabla w \Vert_{L^2(\Omega)} = 1\) concludes
    \(s_t(\nabla_\PW v_h) \lesssim j(v_h)\).
    The generic constant solely depends on the stability estimate~\eqref{eq:quasi_interpolation} and the finite overlap of \(\omega_E\)
    and, thus, is independent of the size of the domain \(\Omega\).
\end{proof}

\begin{remark}[polygonal meshes]
    \label{rem:polygonal_mesh}
    Given a polygonal mesh \(\M\), assume that there exists a regular subtriangulation \(\T\) into triangles.
    Since the contributions from any interior edge \(E \in \E(\T) \setminus \E(\M)\) inside the polygons vanish on both sides of the equivalences,
    the results from Corollary~\ref{cor:Sigma_jump_equivalence} and Lemma~\ref{lem:U_jump_equivalence} generalize 
    to discrete spaces \(\Sigma(\M)\) and \(U(\M)\).
\end{remark}

\begin{remark}[three space dimensions]
    \label{rem:3D}
    All arguments in the proof of Lemma~\ref{lem:sup_control}, in particular the bubble-function technique, 
    generalize to domains \(\Omega \subseteq \R^d\) of any dimension \(d \geq 2\); see, e.g., \cite[Sect.~3.6]{MR3059294}.
    In the proof of Lemma~\ref{lem:U_jump_equivalence}, the analogous argumentation with an integration by parts
    of a surface derivative to obtain~\eqref{eq:ibp_arc_length} in higher dimensions is not straight-forward.
    However, the control of the jump term \(s_t^2(\nabla_\PW v_h)\) (modified as mentioned in Remark~\ref{rem:3D_Helmholtz})
    by \(j^2(v_h)\) is in fact not necessary.
    Instead, Lemma~\ref{lem:sup_control} can be employed to transform \(s_t^2(\nabla_\PW v_h)\)
    into a computable jump term of discontinuous Galerkin type as shown in Corollary~\ref{cor:Sigma_jump_equivalence} for \(s_n^2(\tau_h)\).
    This procedure leads to a modified discontinuous LSFEM in three dimensions with an additional side constraint
    on the discrete space \(U(\T)\).
\end{remark}

\section{Discontinuous least-squares FEM}
\label{sec:discontinuous_LSFEM}

The following section presents the discretization of the least-squares minimization in~\eqref{eq:fundamental_equivalence}.
Let \(\Sigma(\T) \subseteq H(\Div, \T)\) and \(U(\T) \subseteq H^1(\T)\) be subspaces satisfying the conditions~\eqref{eq:Sigma_condition}
and~\eqref{eq:U_condition}.
Since the jump terms \(s_n\) and \(s_t\) are not computable in general, they are replaced by the least-squares jump terms introduced
in Section~\ref{sec:piecewise_polynomials}.
This leads to the discontinuous least-squares functional \(LS_h(f)\colon \Sigma(\T) \times U(\T) \to \R\) defined,
for \(\sigma_h \in \Sigma(\T)\) and \(u_h \in U(\T)\), by
\begin{equation}
    \label{eq:LS_functional}
    \begin{split}
        LS_h(f; \sigma_h, u_h)
        &\coloneqq
        c_\Omega^2\,
        \Vert f + \Div_\PW \sigma_h \Vert_{L^2(\Omega)}^2
        +
        \Vert \sigma_h - \nabla_\PW u_h \Vert_{L^2(\Omega)}^2
        \\
        &\hphantom{{}\coloneqq{}}
        +
        \sum_{E \in \E(\Omega) \cup \E(\Gamma_\NEU)}
        h_E\,
        \Vert [\sigma_h \cdot n_E]_E \Vert_{L^2(E)}^2
        +
        \sum_{E \in \E(\Omega) \cup \E(\Gamma_\DIR)}
        h_E^{-1}\,
        \Vert [u_h]_E \Vert_{L^2(E)}^2.
    \end{split}
\end{equation}
The discontinuous LSFEM seeks solutions \(\sigma_h \in \Sigma(\T)\) and \(u_h \in U(\T)\) to the minimization problem
\begin{equation}
    \label{eq:minimisation}
    \text{Minimize }
    LS_h(f; \sigma_h, u_h)
    \text{ over all }
    \sigma_h \in \Sigma(\T)
    \text{ and }
    u_h \in U(\T).
\end{equation}
The well-posedness of this variational formulation follows from a discrete version of the fundamental equivalence~\eqref{eq:fundamental_equivalence}.
\begin{theorem}[discrete fundamental equivalence]
    \label{thm:discrete_fundamental_equivalence}
    For all \(\tau_h \in \Sigma(\T)\) and \(v_h \in U(\T)\),
    \begin{equation}
        \label{eq:discrete_fundamental_equivalence}
        \Vert (\tau_h, v_h) \Vert_\T^2
        \approx
        LS_h(0; \tau_h, v_h).
    \end{equation}
    The equivalence constants depend on the polynomial degree \(k\) and the interior angles in \(\T\),
    but not on the size of the domain \(\Omega\).
\end{theorem}

\begin{proof}
    The discontinuous fundamental equivalence~\eqref{eq:fundamental_equivalence} and the equivalences
    from Corollary~\ref{cor:Sigma_jump_equivalence} and Lemma~\ref{lem:U_jump_equivalence} show,
    for \(\tau_h \in \Sigma(\T)\) and \(v_h \in U(\T)\), that
    \[
        \Vert (\tau_h, v_h) \Vert_\T^2
        \approx
        LS(0; \tau_h, v_h)
        \approx
        LS_h(0; \tau_h, v_h).
        \qedhere
    \]
    Note that all involved estimates are independent of the size of the domain \(\Omega\).
\end{proof}

The proof of Theorem~\ref{thm:discrete_fundamental_equivalence} applies to any admissible discrete subspaces \(\Sigma(\T)\) and \(U(\T)\)
satisfying the conditions in~\eqref{eq:Sigma_condition} and \eqref{eq:U_condition}.
This includes piecewise polynomial spaces \(P^k(\T)\) or piecewise versions of standard \(H(\Div)\)-based
finite element spaces such as Raviart--Thomas elements \cite[Sect.~2.3.1]{MR3097958}.
The remaining part of this section investigates the particular choice of \(U(\T) \coloneqq P^{k+1}(\T)\)
with polynomial degree \(k \in \N_0\) and
\begin{equation}
    \label{eq:Sigma_RT}
    \Sigma(\T)
    \coloneqq
    \bigg\{
        \tau_h \in RT^{k, \PW}(\T)
        \;:\;
        \forall E \in \E(\Omega) \cup \E(\Gamma_\NEU),\:
        \int_E [\tau_h \cdot n_E]_E \ds = 0
    \bigg\}
\end{equation}
employing the space of piecewise Raviart--Thomas functions
\[
    RT^{k,\PW}(\T)
    \coloneqq
    \{
        a +
        b\,\operatorname{id}
        \;:\;
        a \in P^k(\T; \R^2),\,
        b \in P^k(\T)
    \}.
\]
The first variation of the least-squares functional \(LS_h\) from~\eqref{eq:LS_functional} leads to the following forms
\begin{align}
    \label{eq:bilinear_form}
    B(\sigma_h, u_h; \tau_h, v_h)
    &\coloneqq
    c_\Omega^2\,
    (\Div_\PW \sigma_h, \Div_\PW \tau_h)_{L^2(\Omega)}
    +
    (\sigma_h - \nabla_\PW u_h, \tau_h - \nabla_\PW v_h)_{L^2(\Omega)}
    \\
    \notag
    &\phantom{\coloneqq{}}
    +
    \sum_{E \in \E(\Omega) \cup \E(\Gamma_\NEU)}
    h_E\,
    ([\sigma_h \cdot n_E]_E, [\tau_h \cdot n_E]_E)_{L^2(E)}
    \\
    \notag
    &\phantom{\coloneqq{}}
    +
    \sum_{E \in \E(\Omega) \cup \E(\Gamma_\DIR)}
    h_E^{-1}\,
    ([u_h]_E, [v_h]_E)_{L^2(E)},
    \\
    \notag
    F(\tau_h, v_h)
    &\coloneqq
    -
    c_\Omega^2\,
    (f, \Div_\PW \tau_h)_{L^2(\Omega)}.
\end{align}
A straight-forward realization of the constrained minimization problem~\eqref{eq:minimisation} employs a Lagrange multiplier
\[
    \lambda_h
    \in
    \Lambda(\T)
    \coloneqq
    P^0(\E(\Omega) \cup \E(\Gamma_\NEU))
    \coloneqq
    \prod_{E \in \E(\Omega) \cup \E(\Gamma_\NEU)}
    P^0(E)
\]
and the bilinear form \(C\colon \Lambda(\T) \times (\Sigma(\T) \times U(\T)) \to \R\),
\[
    C(\lambda_h; \tau_h, v_h)
    \coloneqq
    \sum_{E \in \E(\Omega) \cup \E(\Gamma_\NEU)}
    (\lambda_h, [\tau_h \cdot n_E]_E)_{L^2(E)}.
\]
This leads to the following saddle-point problem characterizing the solution \(\sigma_h \in \Sigma(\T)\)
with vanishing integral mean of the normal jumps and \(u_h \in U(\T)\) to problem~\eqref{eq:minimisation}:
Seek \(\sigma_h \in RT^{k,\PW}(\T)\), \(u_h \in U(\T)\), and \(\lambda_h \in \Lambda(\T)\) satisfying,
for all \(\tau_h \in RT^{k,\PW}(\T)\), \(v_h \in U(\T)\), and \(\mu_h \in \Lambda(\T)\),
\begin{equation}
    \label{eq:saddle_point}
    \begin{aligned}
        B(\sigma_h, u_h; \tau_h, v_h)
        + C(\lambda_h; \tau_h, v_h)
        &=
        F(\tau_h, v_h),
        \\
        C(\mu_h; \sigma_h, u_h)
        \:\hphantom{+ B(\sigma_h, u_h; \tau_h, v_h)}
        &= 0.
    \end{aligned}
\end{equation}
The proof of the well-posedness of this saddle-point problem employs the injectivity of \(C\):
If \(\mu_h \in \Lambda(\T)\) satisfies \(C(\mu_h; \tau_h, v_h) = 0\) for all \(\tau_h \in RT^{k,\PW}(\T)\) and \(v_h \in U(\T)\),
then \(\mu_h = 0\) \cite[Lem.~7.2.2]{MR3097958}.
Standard arguments for hybridization in the context of mixed FEMs \cite[Subsect.~7.2.2]{MR3097958} prove well-posedness
of problem~\eqref{eq:saddle_point}:
There exists a unique \(\lambda_h \in \Lambda(\T)\) such that the unique discrete solution \(\sigma_h \in \Sigma(\T)\) and \(u_h \in U(\T)\)
to~\eqref{eq:minimisation} and \(\lambda_h\) solve the saddle-point problem~\eqref{eq:saddle_point}.

\begin{remark}[elliptic discretization]
    Although hybridization is a standard technique in the context of mixed FEM \cite{MR3097958}, the saddle-point problem~\eqref{eq:saddle_point}
    breaks the Rayleigh--Ritz setting which is typical for LSFEMs.
    As a remedy, the edge-oriented basis functions of \(RT^k(\T)\) enable the definition of an explicit basis for
    \(\Sigma(\T)\) from~\eqref{eq:Sigma_RT} leading to a symmetric and positive definite system matrix.
    Further details of the construction are postponed to the Appendix~\ref{app:Sigma_basis}.
\end{remark}

\begin{remark}[semi-discontinuous Crouzeix--Raviart discretization]
    The Raviart--Thomas space \(\Sigma(\T) \coloneqq RT^0_\NEU(\T)\) and the Crouzeix--Raviart space
    \[
        U(\T)
        \coloneqq
        CR^1_\DIR(\T)
        \coloneqq
        \bigg\{
            v_\CR \in P^1(\T)
            \;:\;
            \forall E \in \E(\Omega) \cup \E(\Gamma_\DIR),\:
            \int_E [v_\CR]_E \ds = 0
        \bigg\}
        \subseteq
        H^1(\T)
    \]
    are an admissible choice for the discretization in~\eqref{eq:minimisation}.
    However, the (discrete) least-squares functional must explicitly include the full jump terms,
    here \(j^2(v_\CR)\) for \(v_\CR \in U(\T)\), in order to guarantee well-posedness.
    In a counterexample from~\cite[Lem.~7.1]{MR4593282}, the Helmholtz decomposition leads to \(\tau \in H_\NEU(\Div, \Omega)\),
    \(v \in H^1_\DIR(\Omega)\), and \(w_\CR \in CR^1_\DIR(\T)\) such that
    \[
        \Vert \Div \tau \Vert_{L^2(\Omega)}
        =
        \Vert \tau - \nabla_\PW (v + w_\CR) \Vert_{L^2(\Omega)}
        =
        0
        \quad\text{and}\quad
        0
        <
        \Vert \tau \Vert_{L^2(\Omega)}
        =
        \Vert \nabla_\PW (v + w_\CR) \Vert_{L^2(\Omega)}.
    \]
    This contradicts the fundamental equivalence~\eqref{eq:least_squares_principle} without an underlying least-squares principle
    for piecewise Sobolev spaces.
\end{remark}

\begin{remark}[polygonal discretization]
    Let \(\M\) denote a polygonal mesh with a regular subtriangulation \(\T\).
    Remark~\ref{rem:polygonal_mesh} justifies that the well-posedness of the discontinuous LSFEM also holds for the discrete subspaces
    \begin{align*}
        \Sigma(\M)
        &\coloneqq
        \bigg\{
            \tau_h \in P^{k+1}(\M; \R^2)
            \;:\;
            \forall E \in \E(\Omega) \cup \E(\Gamma_\NEU),\:
            \int_E [\tau_h \cdot n_E]_E \ds = 0
        \bigg\},
        \;
        U(\M)
        \coloneqq
        P^{k+1}(\M).
    \end{align*}
    Note that the number of side conditions in \(\Sigma(\M)\) implicitly poses a condition on the polynomial degree 
    \(k \in \N_0\) in order to guarantee that the space is nontrivial.

    Virtual element spaces \cite{MR4586821} may provide another choice for the discrete spaces \(\Sigma(\M)\) and \(U(\M)\)
    on polygonal meshes.
    The virtual element methodology covers approximations of \(H(\Div)\) functions as well as of discontinuous functions.
    However, the explicit construction and analysis of virtual element spaces for the use in a discontinuous least-squares scheme
    goes beyond the scope of this paper.
\end{remark}

\section{Error analysis}
\label{sec:error_analysis}
This section derives an a~posteriori and an a~priori error estimate 
for the discontinuous least-squares solutions~\(\sigma_h \in \Sigma(\T)\) and \(u_h \in U(\T)\) to~\eqref{eq:minimisation}
and the exact solutions \(\sigma \in H_\NEU(\Div, \Omega)\) and \(u \in H^1_\DIR(\Omega)\) characterized,
for all \(\tau \in H_\NEU(\Div, \Omega)\) and \(v \in H^1_\DIR(\Omega)\), by
\[
    (\Div \sigma, \Div \tau)_{L^2(\Omega)}
    +
    (\sigma - \nabla u, \tau - \nabla v)_{L^2(\Omega)}
    =
    - (f, \Div \tau)_{L^2(\Omega)}.
\]
The analysis employs the fundamental equivalence~\eqref{eq:fundamental_equivalence} and 
the equivalences of the edge terms from Section~\ref{sec:discontinuous_LSFEM} 
and, thus, applies to \emph{any} subspaces \(\Sigma(\T) \subseteq H(\Div, \T)\) and \(U(\T) \subseteq H^1(\T)\)
satisfying the conditions~\eqref{eq:Sigma_condition} and~\eqref{eq:U_condition}.
The independence of the generic constants from the size of the domain~\(\Omega\) transfers to all subsequent statements.
\begin{theorem}[a~posteriori error estimate]
    \label{thm:a_posteriori}
    For any function \(\tau_h \in \Sigma(\T)\) and \(v_h \in U(\T)\),
    \begin{equation}
        \label{eq:a_posteriori}
        \Vert (\sigma - \tau_h, u - v_h) \Vert_\T^2
        \approx
        LS_h(f; \tau_h, v_h).
    \end{equation}
\end{theorem}

\begin{proof}
    The fundamental equivalence~\eqref{eq:fundamental_equivalence} shows, for \(\tau = \sigma - \tau_h \in H(\Div, \T)\)
    and \(v = u - v_h \in H^1(\T)\), that
    \begin{align*}
        \Vert (\sigma - \tau_h, u - v_h) \Vert_\T^2
        &\approx
        c_\Omega^2\,
        \Vert f + \Div_\PW\tau_h \Vert_{L^2(\Omega)}^2
        +
        \Vert \tau_h - \nabla_\PW v_h \Vert_{L^2(\Omega)}^2
        \\
        &\phantom{{}\approx{}}
        +
        s_n^2(\sigma - \tau_h)
        + j^2(u - v_h)
        + s_t^2(\nabla_\PW(u - v_h)).
    \end{align*}
    The consistency of the jump terms from Lemma~\ref{lem:consistency_jump} and 
    the equivalences from Corollary~\ref{cor:Sigma_jump_equivalence} and Lemma~\ref{lem:U_jump_equivalence} prove
    \begin{equation}
        \label{eq:jump_equivalence}
        s_n^2(\sigma - \tau_h)
        + j^2(u - v_h)
        + s_t^2(\nabla_\PW(u - v_h))
        \approx
        \sum_{E \in \E(\Omega) \cup \E(\Gamma_\NEU)}
        h_E\,
        \Vert [\tau_h \cdot n_E]_E \Vert_{L^2(E)}^2
        + j^2(v_h).
    \end{equation}
    The combination of the two previously displayed formulas concludes the proof of~\eqref{eq:a_posteriori}.
\end{proof}

\begin{remark}[inhomogeneous boundary conditions]
    Suppose that inhomogeneous Neumann and Dirichlet boundary conditions satisfy
    the mild additional regularity assumptions \(g_\NEU \in L^2(\Gamma_\NEU)\) and \(g_\DIR \in H^1(\Gamma_\DIR)\).
    Since the fundamental equivalence~\eqref{eq:fundamental_equivalence} holds independently of the boundary conditions,
    the first equivalence in the proof of Theorem~\ref{thm:a_posteriori} remains valid for inhomogeneous boundary conditions.
    The equivalence of the jump terms in~\eqref{eq:jump_equivalence}, however, holds only for functions with polynomial boundary conditions
    (with vanishing piecewise integral mean of the Neumann boundary data).
    In order to approximate the boundary data, consider the \(L^2(\Gamma_\NEU)\)-orthogonal projection
    \(\Pi_\NEU^k\colon L^2(\Gamma_\NEU) \to P^k(\E(\Gamma_\NEU))\) onto the piecewise polynomial functions on \(\Gamma_\NEU\)
    and the nodal interpolation operator \(I_\DIR^{k+1}\colon H^1(\Gamma_\DIR) \to S^{k+1}(\E(\Gamma_\DIR))\)
    in the Lagrange nodes on \(\Gamma_\DIR\). 
    Recall from \cite[Lem.~2.2]{MR3757107} and \cite[Lem.~2.3]{MR3715170} that there exist conforming extensions 
    \(\widehat\sigma_{\NEU,h} \in H(\Div, \Omega)\) and \(\widehat u_{\DIR,h} \in H^1(\Omega)\) with 
    \(\widehat\sigma_{\NEU,h}\vert_{\Gamma_\NEU} \cdot n = \Pi_\NEU^k g_\NEU\)
    and \(\widehat u_{\DIR,h}\vert_{\Gamma_\DIR} = I_\DIR^{k+1} g_\DIR\) satisfying
    \begin{align*}
        \Vert \sigma - \widehat\sigma_{\NEU,h} \Vert_{L^2(\Omega)}^2
        &\lesssim
        \osc_k^2(g_\NEU, \E(\Gamma_\NEU))
        \coloneqq
        \sum_{E \in \E(\Gamma_\NEU)}
        h_E\,
        \Vert (1 - \Pi_\NEU^k) g_\NEU \Vert_{L^2(E)}^2,
        \\
        \Vert \nabla(u - \widehat u_{\DIR,h}) \Vert_{L^2(\Omega)}^2
        &\lesssim
        \osc_k^2(\partial_s g_\DIR, \E(\Gamma_\DIR))
        \coloneqq
        \sum_{E \in \E(\Gamma_\DIR)}
        h_E\,
        \Vert (1 - \Pi_\DIR^k) \partial_s g_\DIR \Vert_{L^2(E)}^2.
    \end{align*}
    This and Lemma~\ref{lem:norm_equivalence} show that \(s_n(\sigma - \widehat\sigma_{\NEU,h})\) and \(s_t(\nabla(u - \widehat u_{\DIR,h}))\)
    are controlled by the oscillation terms \(\osc_k^2(g_\NEU, \E(\Gamma_\NEU))\) and \(\osc_k^2(\partial_s g_\DIR, \E(\Gamma_\DIR))\) as well.
    If \(\tau_h\) belongs to the modified space
    \[
        \tau_h
        \in
        \Sigma(\T; g_\NEU)
        \coloneqq
        \left\{
            \rho_h \in H(\Div, \T)
            \;:\;
            \begin{aligned}
                &\forall T \in \T\: \forall E \in \E(T),\:
                (\tau_h \vert_T \cdot n_T)\vert_E \in
                P^{k+1}(E),
                \\
                &\forall E \in \E(\Omega),\:
                \textstyle
                \int_E [\tau_h \cdot n_E]_E \ds
                = 0,
                \\
                &\forall E \in \E(\Gamma_\NEU),\:
                \textstyle
                \int_E \tau_h \cdot n_E \ds
                = 
                \int_E \Pi_\NEU^k g_\NEU \ds
                =
                \int_E g_\NEU \ds
            \end{aligned}
        \right\},
    \]
    the remaining terms \(s_n(\widehat\sigma_{\NEU,h} - \tau_h)\) and \(s_t(\nabla(\widehat u_{\DIR,h} - v_h))\)
    can be estimated analogously to~\eqref{eq:jump_equivalence}.
    Accordingly, the jump terms in the discontinuous least-squares functional~\eqref{eq:LS_functional}
    are evaluated with respect to the approximated boundary data \(\Pi_\NEU^k g_N\) and \(I_\DIR^{k+1} g_\DIR\) on boundary edges.
    Altogether, this validates the a~posteriori error estimate~\eqref{eq:a_posteriori} (up to oscillations) 
    in the case of inhomogeneous boundary conditions.

    For \(\Omega \subseteq \R^3\), the nodal interpolation is replaced by the Scott--Zhang quasi-interpolation \cite{MR3082295}.
    The reader is referred to \cite[Sect.~3.3--3.4]{BringmannDissertation} for the analogous extension results in three dimensions.
\end{remark}

The a~posteriori estimate from Theorem~\ref{thm:a_posteriori} enables the derivation of the following quasi-best\-approxi\-mation result.
The proof mimics the \emph{medius analysis} in \cite{MR2684360} in that
the control of the a~priori error \(\Vert (\sigma - \sigma_h, u - u_h) \Vert_\T\)
employs the a~posteriori estimate~\eqref{eq:a_posteriori} from Theorem~\ref{thm:a_posteriori}.
In the classical theory of nonconforming methods in \cite[Chap.~10]{MR2373954} or \cite{MR2684360},
the proof of a quasi-bestapproximation estimate requires additional regularity of the exact solution \((\sigma, u)\).
However, such regularity assumptions can be avoided for the discontinuous least-squares FEM at hand.
\begin{theorem}[quasi-bestapproximation]
    \label{thm:quasi_bestapproximation}
    The discrete solutions \(\sigma_h \in \Sigma(\T)\) and \(u_h \in U(\T)\) to the discontinuous LSFEM~\eqref{eq:minimisation} satisfy
    \[
        \Vert (\sigma - \sigma_h, u - u_h) \Vert_\T
        \lesssim
        \inf_{\substack{\tau_h \in \Sigma(\T) \\ v_h \in U(\T)}}
        \Vert (\sigma - \tau_h, u - v_h) \Vert_\T
    \]
\end{theorem}

\begin{proof}
    For the discrete solutions \(\sigma_h \in \Sigma(\T)\) and \(u_h \in U(\T)\) as well as
    arbitrary \(\tau_h \in \Sigma(\T)\) and \(v_h \in U(\T)\),
    set \(\rho_h \coloneqq \tau_h - \sigma_h \in \Sigma(\T)\) and \(w_h \coloneqq v_h - u_h \in \Sigma(\T)\).
    The discrete fundamental equivalence~\eqref{eq:discrete_fundamental_equivalence}
    from Theorem~\ref{thm:discrete_fundamental_equivalence} and the discrete equation~\eqref{eq:saddle_point} yield
    \begin{align*}
        \Vert (\rho_h, w_h) \Vert_\T^2
        &\lesssim
        B(\tau_h - \sigma_h, v_h - u_h; \rho_h, w_h)
        =
        B(\tau_h, v_h; \rho_h, w_h) - F(\rho_h, w_h)
        \\
        &=
        c_\Omega^2\,
        (f + \Div_\PW \tau_h, \Div_\PW \rho_h)_{L^2(\Omega)}
        +
        (\tau_h - \nabla_\PW v_h, \rho_h - \nabla_\PW w_h)_{L^2(\Omega)}
        \\
        &\phantom{={}}
        +
        \sum_{E \in \E(\Omega) \cup \E(\Gamma_\NEU)}
        h_E\, ([\tau_h \cdot n_E]_E, [\rho_h \cdot n_E]_E)_{L^2(E)}
        \\
        &\phantom{={}}
        +
        \sum_{E \in \E(\Omega) \cup \E(\Gamma_\DIR)}
        h_E^{-1}\, ([v_h]_E, [w_h]_E)_{L^2(E)}.
    \end{align*}
    Cauchy--Schwarz inequalities in \(L^2\) and in \(\R^N\) and the efficiency in~\eqref{eq:a_posteriori} result in
    \begin{align*}
        \Vert (\rho_h, w_h) \Vert_\T^2
        &\lesssim
        LS_h(f; \tau_h, v_h)^{1/2}\;
        LS_h(0; \rho_h, w_h)^{1/2}
        \lesssim
        \Vert (\sigma - \tau_h, u - v_h) \Vert_\T
        \,
        \Vert (\rho_h, w_h) \Vert_\T.
    \end{align*}
    This, the division by \(\Vert (\rho_h, w_h) \Vert_\T\), and the triangle inequality
    \[
        \Vert (\sigma - \sigma_h, u - u_h) \Vert_\T
        \leq
        \Vert (\sigma - \tau_h, u - v_h) \Vert_\T
        +
        \Vert (\rho_h, w_h) \Vert_\T
    \]
    conclude the proof.
\end{proof}

\begin{remark}[quasi-bestapproximation results in the literature]
    Using standard approximation results of the Raviart--Thomas and Lagrange finite element spaces \cite[Chap.~2]{MR3097958},
    the quasi-bestapproximation result from Theorem~\ref{thm:quasi_bestapproximation} leads to the optimal a~priori convergence rates,
    for exact solutions \(\sigma \in H^s(\Omega; \R^2)\) with \(\Div \sigma \in H^s(\Omega)\) and \(u \in H^{1+s}(\Omega)\) for \(s > 0\),
    \begin{equation}
        \label{eq:apriori_rate}
        \Vert (\sigma - \sigma_h, u - u_h) \Vert_\T
        \lesssim
        h^{\min\{k+1, s\}} \,
        \big( \vert \sigma \vert_{H^s(\Omega)} + \vert \Div \sigma \vert_{H^s(\Omega)} + \vert u \vert_{H^{1+s}(\Omega)}\big).
    \end{equation}%
    The argumentation for the proof of Theorem~\ref{thm:quasi_bestapproximation} is essentially included in~\cite[Thm.~3.1]{MR2149925}.
    But therein, quasi-bestapproximation is shown in a conforming discrete subspace.
    As a consequence, the normal jump term vanishes and the negative power of \(h\) does not decrease the order of convergence.
    However, this approach is limited to situations where a conforming subspace is available, e.g., for triangular meshes.
    If one aims at quasi-bestapproximation in the discontinuous space on general polygonal meshes,
    the over-penalization of the normal jumps in \cite[Lem.~4.1 and Thm.~4.1]{MR3895838}
    requires increased regularity assumptions on \(\nabla u = \sigma \in H^{1+s}(\Omega, \R^2)\) and, thus, 
    \(u \in H^{2+s}(\Omega)\) to achieve the same a~priori convergence rate as in \eqref{eq:apriori_rate}.
    In this sense, Theorem~\ref{thm:quasi_bestapproximation} generalizes both references.
\end{remark}

\begin{remark}[quasi-bestapproximation theory of nonconforming methods]
    The framework by Veeser and Zanotti \cite{MR3816182} characterizes quasi-best\-approxi\-mation of nonconforming methods
    without any regularity assumptions for formulations including smoothing operators
    \(J_1\colon \Sigma(\T) \to H_\NEU(\Div, \Omega)\) and \(J_2\colon U(\T) \to H^1_\DIR(\Omega)\)
    \[
        B(\sigma_h, u_h; \tau_h, u_h)
        =
        F(J_1 \tau_h, J_2 v_h)
        \quad\text{for all }
        \tau_h \in \Sigma(\T)
        \text{ and }
        v_h \in U(\T).
    \]
    However, this variational formulation does \emph{not} comply with a discrete least-squares minimization.
\end{remark}

\section{A~posteriori error estimation for over-penalized discontinuous LSFEM}
\label{sec:overpenalised_DLSFEM}

Established discontinuous LSFEMs from the literature such as \cite{MR2149925} include an a~posteriori error estimator as well,
albeit it is subject to the existence of suitable enrichment operators,
as they are used in the proof of Theorem~\ref{thm:a_posteriori_classical} below.
In contrast to that, the analysis in Sections~\ref{sec:discontinuous_LSFEM}--\ref{sec:error_analysis} applies to any discrete subspace.
Recall the weight factor \(c_\Omega = C_\FR\) in terms of the Friedrichs constant from~\eqref{eq:Friedrichs}.
Under the additional assumption of a quasi-uniform \emph{initial} triangulation,
there exists an initial mesh size \(h_0 > 0\) such that \(h_E \leq h_0\) for all \(E \in \E\) and, thus, \(h_E \leq h_0^2\, h_E^{-1}\).
As the initial triangulation resolves the domain and the Dirichlet boundary, it follows \(h_0 \lesssim C_\FR\) and hence, for all \(E \in \E\),
\begin{equation}
    \label{eq:mesh_size_estimate}
    h_E \lesssim C_\FR^2\, h_E^{-1}.
\end{equation}
All generic constants in this section are independent of the size of the domain \(\Omega\), but may depend
on the polynomial degree \(k\) and the interior angles in \(\T\).

Discontinuous least-squares formulations from the literature such as in \cite{MR2149925} usually employ the
(unweighted version of the) least-squares functional \(LS_h^\star(f)\colon RT^{k,\PW}(\T) \times P^{k+1}(\T) \to \R\)
\begin{equation}
    \label{eq:op_LS_functional}
    \begin{aligned}
        LS_h^\star(f; \tau_h, v_h)
        &\coloneqq
        c_\Omega^2 \,
        \Vert f + \Div_\PW \tau_h \Vert_{L^2(\Omega)}^2
        +
        \Vert \tau_h - \nabla_\PW v_h \Vert_{L^2(\Omega)}^2
        \\
        &\phantom{{}\coloneqq{}}
        +
        \sum_{E \in \E(\Omega) \cup \E(\Gamma_\NEU)}
        c_\Omega^2 \,
        h_E^{-1}
        \Vert [\tau_h \cdot n_E]_E \Vert_{L^2(E)}^2
        +
        \sum_{E \in \E(\Omega) \cup \E(\Gamma_\DIR)}
        h_E^{-1}
        \Vert [v_h]_E \Vert_{L^2(E)}^2.
    \end{aligned}
\end{equation}
For \(\tau + \tau_h \in H_\NEU(\Div, \Omega) + RT^{k, \PW}(\T)\) and \(v + v_h \in H^1_\DIR(\Omega) + P^{k+1}(\T)\),
abbreviate the jump terms
\[
    j_{n,\alpha}^2(\tau_h)
    \coloneqq
    \sum_{E \in \E(\Omega) \cup \E(\Gamma_\NEU)}
    h_E^{\alpha}
    \Vert [\tau_h \cdot n_E]_E \Vert_{L^2(E)}^2
\]
and the norm
\begin{align*}
    \Vert (\tau + \tau_h, v + v_h) \Vert_\star^2
    &\coloneqq
    c_\Omega^2 \,
    \Vert \Div_\PW(\tau + \tau_h) \Vert_{L^2(\Omega)}^2
    +
    \Vert \tau + \tau_h \Vert_{L^2(\Omega)}^2
    +
    c_\Omega^2 \,
    j_{n,-1}^2(\tau_h)
    \\
    &\phantom{{}\coloneqq \Big({}}
    +
    \Vert \nabla_\PW (v + v_h) \Vert_{L^2(\Omega)}^2
    +
    j^2(v_h).
\end{align*}
The cause for the negative power of \(h_E\) in \(j_{n,-1}\) lies in the imbalance of the scaling of
the natural Sobolev jumps \(s_n\) and the computable jumps \(j_{n,-1}\) in absence of additional side conditions
as addressed in Remark~\ref{rem:counter_example}.
Nevertheless, the over-penalization of the normal jumps is a common choice for discontinuous LSFEMs which justifies the investigation
of its built-in error estimation property in the following result.
\begin{theorem}[built-in a~posteriori estimation]
    \label{thm:a_posteriori_classical}
    For every \(\tau_h \in RT^{k, \PW}(\T)\) and \(v_h \in P^{k+1}(\T)\),
    \begin{equation}
        \label{eq:a_posteriori_overpenalised}
        \Vert (\sigma - \tau_h, u - v_h) \Vert_\star^2
        \approx
        LS_h^\star(f; \tau_h, v_h).
    \end{equation}
    The equivalence constants solely depend on the polynomial degree \(k\) and the interior angles in \(\T\),
    but not on the size of the domain \(\Omega\).
\end{theorem}
Before the remaining part of this section addresses the proof of Theorem~\ref{thm:a_posteriori_classical},
the following lemma recalls suitable averaging operators.
They link the discontinuous spaces with their conforming counterparts and enable the application of the fundamental
equivalence~\eqref{eq:least_squares_principle} for the conforming LSFEM.
This argumentation solely relies on the availability of conforming subspaces of the discontinuous ansatz spaces with
corresponding averaging operators and the conforming fundamental equivalence.
Thus, the analogous arguments enable the generalization to higher spatial dimensions and various applications
with a well-posed conforming LSFEM.
This includes, but is not limited to, the linear elasticity problem \cite{MR2084237,MR4065377},
the Helmholtz equation \cite{MR1302685,MR4563176}, and the time-harmonic Maxwell equations \cite{MR2189548,MR4367673}.
\begin{lemma}[averaging operators]
    \label{lem:averaging}
    There exist operators \(J_\RT\colon RT^{k,\PW}(\T) \to RT^k_\NEU(\T)\) and \(J_\CO\colon P^{k+1}(\T) \to S^{k+1}_\DIR(\T)\) satisfying
    \begin{align}
        \label{eq:averaging_estimates_RT}
        c_\Omega^2\,
        \Vert \Div_\PW (1 - J_\RT) \tau_h \Vert_{L^2(\Omega)}^2
        +
        \Vert (1 - J_\RT) \tau_h \Vert_{L^2(\Omega)}^2
        &\lesssim
        c_\Omega^2\,
        j_{n,-1}^2(\tau_h)
        \\
        \label{eq:averaging_estimates_Sk}
        \Vert \nabla_\PW (1 - J_\CO) v_h
        \Vert_{L^2(\Omega)}^2
        &\lesssim
        j^2(v_h).
    \end{align}
    The generic constants solely depend on the polynomial degree \(k\) and the interior angles in \(\T\),
    but not on the size of the domain \(\Omega\).
\end{lemma}

\begin{proof}
    The operators in \cite[Sect.~22.2]{MR4242224} are defined by 
    averaging the coefficients at the degrees of freedom associated to edges and vertices.
    In order to meet the boundary conditions for the averaging \(J_\RT \tau_h\) of \(\tau_h \in RT^{k,\PW}(\T)\)
    and \(J_\CO v_h\) of \(v_h \in P^{k+1}(\T)\), set the degrees of freedom on the boundary explicitly
    to \((J_\RT\tau_h) \cdot n_E \vert_E \equiv 0\) for all \(E \in \E(\Gamma_\NEU)\)
    and \((J_\CO v_h) \vert_E \equiv 0\) for all \(E \in \E(\Gamma_\DIR)\).
    This leads to some operators \(J_\RT\colon RT^{k,\PW}(\T) \to RT^k_\NEU(\T)\) and \(J_\CO\colon P^{k+1}(\T) \to S^{k+1}_\DIR(\T)\)
    with \cite[Lem.~22.3]{MR4242224}
    \begin{align*}
        \Vert (1 - J_\RT) \tau_h \Vert_{L^2(T)}^2
        &\lesssim
        \sum_{E \in \E(T) \setminus \E(\Gamma_\DIR)}
        h_E
        \Vert [\tau_h \cdot n_E]_E \Vert_{L^2(E)}^2,
        \\
        \Vert \Div (1 - J_\RT) \tau_h \Vert_{L^2(T)}^2
        &\lesssim
        \sum_{E \in \E(T) \setminus \E(\Gamma_\DIR)}
        h_E^{-1}
        \Vert [\tau_h \cdot n_E]_E \Vert_{L^2(E)}^2,
        \\
        \Vert \nabla (1 - J_\CO) v_h \Vert_{L^2(T)}^2
        &\lesssim
        \sum_{E \in \E(T) \setminus \E(\Gamma_\NEU)}
        h_E^{-1}
        \Vert [v_h]_E \Vert_{L^2(E)}^2.
    \end{align*}
    Summing over \(T \in \T\), a finite overlapping argument of the edge patches, and the estimate~\eqref{eq:mesh_size_estimate}
    prove the claimed estimates with generic constants independent of the domain size.
\end{proof}

\begin{proof}[Proof of Theorem~\ref{thm:a_posteriori_classical}]
    The proof of the equivalence~\eqref{eq:a_posteriori_overpenalised} follows arguments from \cite[Thm.~7.3]{MR4593282}.
    It bases on the fundamental equivalence~\eqref{eq:least_squares_principle} of the conforming LSFEM
    \begin{equation}
        \label{eq:averaging_cfe}
        \Vert (\sigma - J_\RT \tau_h, u - J_\CO v_h) \Vert_\star^2 \approx LS_h^\star(f; J_\RT \tau_h, J_\CO v_h)
    \end{equation}
    and on the estimates~\eqref{eq:averaging_estimates_RT}--\eqref{eq:averaging_estimates_Sk} 
    for the averaging operators \(J_\RT\) and \(J_\CO\) from Lemma~\ref{lem:averaging}
    \begin{equation}
        \label{eq:averaging_error_estimate}
        \begin{split}
            &\hspace{-2em}
            \Vert ( (1 - J_\RT) \tau_h, (1 - J_\CO) v_h) \Vert_\star^2
            \lesssim
            c_\Omega^2\,
            j_{n,-1}^2(\tau_h)
            +
            j^2(v_h).
        \end{split}
    \end{equation}
    This and several triangle inequalities show
    \begin{align}
        \label{eq:triangle_ls_averaging_error}
        \begin{split}
            LS_h^\star(0; (1 - J_\RT) \tau_h, (1 - J_\CO) v_h)
            &\lesssim
            \Vert ( (1 - J_\RT) \tau_h, (1 - J_\CO) v_h) \Vert_\star^2
            \\
            &\lesssim
            c_\Omega^2 \,
            j_{n,-1}^2(\tau_h)
            +
            j^2(v_h),
        \end{split}
        \\[\smallskipamount]
        \label{eq:triangle_ls_functional}
        \begin{split}
            LS_h^\star(f; J_\RT \tau_h, J_\CO v_h)
            &\lesssim
            LS_h^\star(f; \tau_h, v_h)
            +
            LS_h^\star(0; (1 - J_\RT) \tau_h, (1 - J_\CO)
            v_h)
            \\
            &\lesssim
            LS_h^\star(f; \tau_h, v_h),
        \end{split}
        \\[\smallskipamount]
        \label{eq:triangle_error}
        \begin{split}
            \Vert (\sigma - J_\RT \tau_h, u - J_\CO v_h) \Vert_\star^2
            &\lesssim
            \Vert (\sigma - \tau_h, u - v_h) \Vert_\star^2
            +
            \Vert ( (1 - J_\RT) \tau_h, (1 - J_\CO) v_h) \Vert_\star^2
            \\
            &\lesssim
            \Vert (\sigma - \tau_h, u - v_h) \Vert_\star^2.
        \end{split}
    \end{align}
    The triangle inequality, reliability from \eqref{eq:averaging_cfe}, and the estimates
    \eqref{eq:averaging_error_estimate} and \eqref{eq:triangle_ls_functional} prove the reliability
    \begin{align*}
        \Vert (\sigma - \tau_h, u - v_h) \Vert_\star^2
        &\lesssim
        \Vert (\sigma - J_\RT \tau_h, u - J_\CO v_h) \Vert_\star^2
        +
        \Vert ( (1 - J_\RT) \tau_h, (1 - J_\CO) v_h) \Vert_\star^2
        \\
        &\lesssim
        LS_h^\star(f; J_\RT \tau_h, J_\CO v_h)
        +
        c_\Omega^2\,
        j_{n,-1}^2(\tau_h)
        +
        j^2(v_h)
        \lesssim
        LS_h^\star(f; \tau_h, v_h).
    \end{align*}
    The triangle inequality, efficiency from \eqref{eq:averaging_cfe}, and the estimates
    \eqref{eq:triangle_ls_averaging_error} and \eqref{eq:triangle_error} verify
    \begin{align*}
        LS_h^\star(f; \tau_h, v_h)
        &\lesssim
        LS_h^\star(f; J_\RT \tau_h, J_\CO v_h)
        +
        LS_h^\star(0; (1 - J_\RT) \tau_h, (1 - J_\CO) v_h)
        \\
        &\lesssim
        \Vert (\sigma - J_\RT \tau_h, u - J_\CO v_h) \Vert_\star^2
        +
        c_\Omega^2\,
        j_{n,-1}^2(\tau_h)
        +
        j^2(v_h)
        \\
        &\lesssim
        \Vert (\sigma - \tau_h, u - v_h) \Vert_\star^2
    \end{align*}
    and conclude the proof of the efficiency.
\end{proof}

\begin{remark}[scaling-robust discrete fundamental equivalence]
    Replacing \(\sigma\) and \(u\) by \(0\) in the proof of Theorem~\ref{thm:a_posteriori_classical},
    an analogous argumentation establishes the discrete fundamental equivalence, 
    for all \(\tau_h \in RT^{k, \PW}(\T)\) and \(v_h \in P^{k+1}(\T)\),
    \(\Vert (\tau_h, v_h) \Vert_\star^2 \approx LS_h^\star(0; \tau_h, v_h)\) with equivalence constants 
    independent of the size of the domain \(\Omega\).
\end{remark}

\section{Numerical experiments}
\label{sec:experiments}

This section presents numerical experiments with the naturally penalized discontinuous LSFEM from
Section~\ref{sec:discontinuous_LSFEM} and the over-penalized method from Section~\ref{sec:overpenalised_DLSFEM}.
The implementation utilizes the NGSolve software package (version 6.2.2405) \cite{ngsolve}
and is provided as a code capsule \cite{DLSFEM} on the Code Ocean platform for the sake of reproducibility.
The Subsections~\ref{sec:scaling} and~\ref{sec:scaling_anisotropic} investigate the influence of the weighting factor \(c_\Omega\)
for large domains using benchmark problems from~\cite{MR4253887}
while Subsection~\ref{sec:lshape} focusses on the adaptive mesh refinement for the L-shaped domain.

The exponent \(\alpha \in \{-1, 1\}\) depends on the chosen method, i.e.,
\(\alpha = 1\) for the naturally penalized formulation and \(\alpha = -1\) for the over-penalized scheme.
The homogeneous boundary conditions are weakly enforced in the minimization
of the discontinuous least-squares functional~\eqref{eq:intro_discontinuous_LSFEM}.
Recall that, for \(\alpha = 1\), the discrete flux space incorporates the integral mean side condition in~\eqref{eq:Sigma_RT}.
\begin{figure}
    \centering
    \hfil
    \subfloat[Square (Subsection~\ref{sec:scaling})]{%
        \label{fig:example:scaling_initial_triangulations:square}
        \begin{tikzpicture}[line width=2, scale=0.8, line join=round, line cap=round, >=stealth]
            \draw[TUblue, fill=TUblue!25!white] (0,0) node[coordinate] (A1) {} -- (4,0) node[coordinate] (A2) {}
                -- (4,4) node[coordinate] (A3) {} -- (0,4) node[coordinate] (A4) {} -- cycle;
            \draw[TUblue] (A1) -- (A3) (A2) -- (A4) ($(A1)!0.5!(A2)$) -- ($(A3)!0.5!(A4)$) ($(A1)!0.5!(A4)$) -- ($(A2)!0.5!(A3)$);
            \draw[TUgray, line width=1.3, decorate, decoration={brace, amplitude=4mm, mirror, raise=1mm}]
                (A1) -- (A2) node[midway,yshift=-10mm] {edge length \(\ell\)};
        \end{tikzpicture}
    }
    \hfil
    \subfloat[Rectangle (Subsection~\ref{sec:scaling_anisotropic})]{%
        \label{fig:example:scaling_initial_triangulations:rectangle}
        \centering
        \hspace{6mm}
        \begin{tikzpicture}[line width=1.3, scale=0.8, line join=round, line cap=round, >=stealth]
            \draw[TUblue, fill=TUblue!25!white] (0,0) node[coordinate] (A1) {} -- (1,0) node[coordinate] (A2) {}
                -- (2,0) node[coordinate] (A3) {} -- (2,1) node[coordinate] (A4) {} -- (1,1) node[coordinate] (A5) {}
                -- (0,1) node[coordinate] (A6) {} -- cycle;
            \draw[TUblue] (A1) -- (A5) (A2) -- (A6) (A2) -- (A5) (A2) -- (A4) (A3) -- (A5) 
                ($(A1)!0.5!(A2)$) -- ($(A5)!0.5!(A6)$) ($(A2)!0.5!(A3)$) -- ($(A4)!0.5!(A5)$) ($(A1)!0.5!(A6)$) -- ($(A3)!0.5!(A4)$);
            \draw[TUblue, fill=TUblue!25!white] (3,0) node[coordinate] (B1) {} -- (4,0) node[coordinate] (B2) {} 
                -- (4,1) node[coordinate] (B3) {} -- (3,1) node[coordinate] (B4) {} -- cycle;
            \draw[TUblue] (B1) -- (B3) (B2) -- (B4) ($(B1)!0.5!(B2)$) -- ($(B3)!0.5!(B4)$) ($(B1)!0.5!(B4)$) -- ($(B2)!0.5!(B3)$);
            \draw[TUblue] (2.5,0.5) node {\(\cdots\)};
            \draw[TUgray, line width=1.3, decorate, decoration={brace, amplitude=4mm, mirror, raise=1mm}]
                (A1) -- (B2) node[midway,yshift=-10mm] {\(\ell\) unit squares};
        \end{tikzpicture}
        \hspace{6mm}
    }
    \hfil
    \subfloat[L-shape (Subsection~\ref{sec:lshape})]{%
        \label{fig:example:scaling_initial_triangulations:lshape}
        \centering
        \begin{tikzpicture}[line width=1.6, scale=0.5, line join=round, line cap=round, >=stealth]
            \draw[TUblue, fill=TUblue!25!white] (0,0) node[coordinate] (A1) {} -- (0,4) node[coordinate] (A2) {}
                -- (-4,4) node[coordinate] (A3) {} -- (-4,0) node[coordinate] (A4) {}
                -- (-4,-4) node[coordinate] (A5) {} -- (0,-4) node[coordinate] (A6) {}
                -- (4,-4) node[coordinate] (A7) {} -- (4,0) node[coordinate] (A8) {} -- cycle;
            \draw[TUblue] (A7) -- (A1) -- (A6) -- (A4) -- (A1) -- (A3) (A2) -- (A4) (A1) -- (A5) (A6) -- (A8)
                ($(A1)!0.5!(A2)$) -- ($(A3)!0.5!(A4)$) ($(A4)!0.5!(A5)$) -- ($(A7)!0.5!(A8)$)
                ($(A2)!0.5!(A3)$) -- ($(A5)!0.5!(A6)$) ($(A6)!0.5!(A7)$) -- ($(A8)!0.5!(A1)$);
            \draw[TUgray, line width=1.3, decorate, decoration={brace, amplitude=4mm, mirror, raise=1mm}]
                (A5) -- (A6) node[midway,yshift=-10mm] {edge length \(\ell\)};
        \end{tikzpicture}
    }
    \caption{%
        Illustration of initial triangulations \(\T_0\) of the scaled square domain, the anisotropic rectangle domain,
        and the scaled L-shaped domain.
    }
    \label{fig:example:scaling_initial_triangulations}
\end{figure}

\subsection{Scaled square example}
\label{sec:scaling}

Consider the square domain \(\Omega = \Omega_\ell \coloneqq (0, \ell)^2 \subseteq \R^2\) with edge length \(\ell > 0\).
The smooth exact solution \(u(x) \coloneqq \sin(\pi x_1/ \ell) \, \sin(\pi x_2/\ell)\)
solves the Poisson model problem with right-hand side \(f(x) = 2 \pi^2 / \ell^2 \, u(x)\)
\[
    -\Delta u = f \quad \text{in } \Omega_\ell
    \quad \text{with} \quad 
    u = 0 \quad \text{on } \partial \Omega_\ell.
\]
The initial triangulation \(\T_0\) of \(\Omega_\ell\) into 8 triangles is displayed in Figure~\ref{fig:example:scaling_initial_triangulations:square}.
The convergence history plot in Figure~\ref{fig:example:scaling:one} for the relative energy error
\(\Vert \nabla(u - u_\LS) \Vert_{L^2(\Omega_\ell)} / \Vert \nabla u \Vert_{L^2(\Omega_\ell)}\)
of the unweighted lowest-order discontinuous LSFEM exhibits a drastic preasymptotic effect with almost no error reduction at all.
In contrast to that, the proper scaling with \(c_\Omega = C_\FR = \ell/\pi\) completely avoids this undesired behavior
and ensures the identical convergence of the relative energy error independently of the scaling of the domain
as displayed in Figure~\ref{fig:example:scaling:friedrichs}.

Although the exact Friedrichs constant \(C_\FR = \lambda(1)^{-1/2}\) with the principal Dirichlet eigenvalue \(\lambda(1) > 0\)
is usually unknown for nontrivial domains,
an upper bound of \(C_\FR\) can be computed by guaranteed lower eigenvalue bounds on \(\lambda(1)\) from nonconforming discretizations
such as the hybrid high-order method in \cite{MR4332791,MR4700405}.
Alternatively, the Figures~\ref{fig:example:scaling:diam} and~\ref{fig:example:scaling:width} justify the choice of every equivalent constant
\(c_\Omega\) with the same scaling in \(\ell\) as the Friedrichs constant \(C_\FR\), e.g., the width of \(\Omega_\ell\)
from~\eqref{eq:width} or the diameter of \(\Omega_\ell\).
Figure~\ref{fig:example:weighting:k0} investigates the fixed domain \(\Omega_\ell\) for \(\ell = \pi\) with varying weightings.
It illustrates that overly weighting the equilibrium residual does not deteriorate the convergence behavior.
However, very large weights might result in disadvantageously large condition numbers of the system matrix.
The discontinuous LSFEM for higher-order discretization turns out to be more robust with respect to the choice of the weighting
as shown by Figure~\ref{fig:example:weighting:k1}.
Nevertheless, an under-weighted equilibrium residual still exhibits a preasymptotic increase of the relative energy error by one order of magnitude.

Similarly dramatic preasymptotic effects are observed for the over-penalized discontinuous LSFEM with the unweighted functional
as shown in Figure~\ref{fig:example:scaling:overpenalized}.
Undisplayed numerical experiments exhibit comparable results for the relative error contributions of the dual variable
\(\Vert \tau - \tau_\LS \Vert_{L^2(\Omega_\ell)} / \Vert \tau \Vert_{L^2(\Omega_\ell)}\)
and \(\Vert \Div (\tau - \tau_\LS) \Vert_{L^2(\Omega_\ell)} / \Vert \Div \tau \Vert_{L^2(\Omega)}\) and are therefore omitted.
\begin{figure}[p]
    \centering
    \subfloat[\(c_\Omega = 1\)]{\label{fig:example:scaling:one}\input{figures/plot_dls_square_error_scalingOne.tex}}
    \hfill
    \subfloat[\(c_\Omega = \operatorname{diam}(\Omega_\ell) = \sqrt{2}\,\ell\)]{\label{fig:example:scaling:diam}\input{figures/plot_dls_square_error_scalingDiam.tex}}

    \subfloat[\(c_\Omega = \operatorname{width}(\Omega_\ell) = \ell\)]{\label{fig:example:scaling:width}\input{figures/plot_dls_square_error_scalingWidth.tex}}
    \hfill
    \subfloat[\(c_\Omega = C_\FR = \ell/\pi\)]{\label{fig:example:scaling:friedrichs}\input{figures/plot_dls_square_error_scalingFriedrichs.tex}}
    \caption{%
        Convergence history plot of the relative energy error
        \(\Vert \nabla(u - u_\LS) \Vert_{L^2(\Omega)} / \Vert \nabla u \Vert_{L^2(\Omega)}\)
        under uniform refinement of the scaled square \(\Omega_\ell\) from Subsection~\ref{sec:scaling}.
        The experiments are carried out with the naturally penalized discontinuous LSFEM (\(\alpha = 1\)) 
        for lowest-order ansatz spaces (\(k = 0\)) and different weighting factors \(c_\Omega\).
    }
    \label{fig:example:scaling}
\end{figure}

\begin{figure}[p]
    \centering
    \subfloat[\(k = 0\)]{\label{fig:example:weighting:k0}\input{figures/plot_dls_square_error_weighting_k0.tex}}
    \hfill
    \subfloat[\(k = 1\)]{\label{fig:example:weighting:k1}\input{figures/plot_dls_square_error_weighting_k1.tex}}
    \caption{%
        Convergence history plot of the relative energy error
        \(\Vert \nabla(u - u_\LS) \Vert_{L^2(\Omega)} / \Vert \nabla u \Vert_{L^2(\Omega)}\)
        under uniform refinement of the square \(\Omega_\pi\) with edge length \(\ell = \pi\) from Subsection~\ref{sec:scaling}.
        The experiments are carried out with the naturally penalized discontinuous LSFEM (\(\alpha = 1\)) 
        for different polynomial degrees \(k\) and weighting factors \(c_\Omega\).
    }
    \label{fig:example:weighting}
\end{figure}

\begin{figure}
    \centering
    \subfloat[\(c_\Omega = 1\)]{\input{figures/plot_opdls_square_error_scalingOne.tex}}
    \hfill
    \subfloat[\(c_\Omega = C_\FR = \ell/\pi\)]{\input{figures/plot_opdls_square_error_scalingFriedrichs.tex}}
    \caption{%
        Convergence history plot of the relative energy error
        \(\Vert \nabla(u - u_\LS) \Vert_{L^2(\Omega)} / \Vert \nabla u \Vert_{L^2(\Omega)}\)
        under uniform refinement of the scaled square \(\Omega_\ell\) from Subsection~\ref{sec:scaling}.
        The experiments are carried out with the over-penalized discontinuous LSFEM (\(\alpha = -1\)) 
        for lowest-order ansatz spaces (\(k = 0\)) and different weighting factors \(c_\Omega\).
    }
    \label{fig:example:scaling:overpenalized}
\end{figure}

\subsection{Anisotropic rectangle example}
\label{sec:scaling_anisotropic}

For the anisotropic rectangle \(\Omega_\ell \coloneqq (0, \ell) \times (0, 1) \subseteq \R^2\) with \(\ell > 0\),
set the Dirichlet boundary \(\Gamma_\DIR \coloneqq \{0\} \times [0, 1] \cup \{\ell\} \times [0,1]\) as the left and right boundary
and the Neumann boundary \(\Gamma_\NEU \coloneqq \partial\Omega_\ell \setminus \Gamma_\DIR\).
For this domain, the Friedrichs constant reads \(c_\Omega = C_\FR = \ell / \pi\).
The smooth function \(u(x) \coloneqq \sin(\pi x_1 / \ell)\)
satisfies the Poisson model problem with right-hand side \(f(x) = \pi^2/\ell^2 \, u(x)\)
and mixed homogeneous boundary conditions \(u = 0\) on \(\Gamma_\DIR\) and \(\partial u / \partial n = 0\) on \(\Gamma_\NEU\).
In order to avoid degenerate triangles with small interior angles in the case of large \(\ell \in \N\),
the initial triangulation \(\T_0\) consists of \(8 \ell\) congruent triangles as displayed 
in Figure~\ref{fig:example:scaling_initial_triangulations:rectangle}.
The convergence graphs in Figure~\ref{fig:example:scaling_anisotropic} show the similar preasymptotic effect
as in Subsection~\ref{sec:scaling} for the unweighted discontinuous LSFEM.
Although the initial triangulations for large \(\ell > 0\) have significantly more degrees of freedom,
the relative error is almost equal.
Again, the weighting with \(c_\Omega = C_\FR(\Omega_\ell) = \ell/\pi\) or similarly scaling constants improves the convergence
and the high number of initial degrees of freedom yield the expected improvement of the accuracy.
\begin{figure}
    \centering
    \subfloat[\(c_\Omega = 1 = \operatorname{width}(\Omega)\)]{\input{figures/plot_dls_rectangle_error_scalingOne.tex}}
    \hfill
    \subfloat[\(c_\Omega = \operatorname{diam}(\Omega) = \sqrt{2}\,\ell\)]{\input{figures/plot_dls_rectangle_error_scalingDiam.tex}}

    \subfloat[\(c_\Omega = C_\FR = \ell/\pi\)]{\input{figures/plot_dls_rectangle_error_scalingFriedrichs.tex}}
    \caption{%
        Convergence history plot of the relative energy error
        \(\Vert \nabla(u - u_\LS) \Vert_{L^2(\Omega)} / \Vert \nabla u \Vert_{L^2(\Omega)}\)
        under uniform refinement of the anisotropic rectangle \(\Omega_\ell\) from Subsection~\ref{sec:scaling_anisotropic}.
        The experiments are carried out with the naturally penalized discontinuous LSFEM (\(\alpha = 1\)) 
        for lowest-order ansatz spaces (\(k = 0\)) and different weighting factors \(c_\Omega\).
    }
    \label{fig:example:scaling_anisotropic}
\end{figure}

\subsection{L-shaped domain example}
\label{sec:lshape}

The benchmark problem on the scaled L-shaped domain \(\Omega_\ell \coloneqq (-\ell, \ell)^2 \setminus [0, \ell)^2\)
with full Dirichlet boundary \(\Gamma_\DIR \coloneqq \partial\Omega_{\ell}\)
employs the exact solution \(u \in H^1_0(\Omega_\ell)\) given in polar coordinates by
\[
    u(r, \varphi)
    \coloneqq
    (r/\ell)^{2/3} \sin(2\varphi/3) \,
    \big(1 - (r/\ell)^2 \cos^2(\varphi)\big) \,
    \big(1 - (r/\ell)^2 \sin^2(\varphi)\big).
\]
The exact solution determines the right-hand side \(f \coloneqq - \Delta u \in L^2(\Omega_{\ell})\) and 
the homogeneous Dirichlet boundary conditions \(u\vert_{\Gamma_\DIR} = 0\).
The weight factor \(c_\Omega \coloneqq \ell \lambda(1)^{-1/2}\) is chosen with the approximation
\(\lambda(1) = 9.6397238389738806\) of the principal eigenvalue on standard L-shaped domain \(\Omega_1\) from~\cite{MR4332791}.
The adaptive mesh-refining algorithm employs the triangle-wise contributions to the least-squares functional~\eqref{eq:intro_discontinuous_LSFEM}
as the built-in error estimator and the standard D\"orfler marking strategy \cite{MR1393904} with some bulk parameter \(0 < \theta \leq 1\).

In accordance with the investigation in Subsections~\ref{sec:scaling} and~\ref{sec:scaling_anisotropic},
Figure~\ref{fig:example:scaling_lshape} shows that the weighting with correctly scaled \(c_\Omega\) is required
to avoid a large preasymptotic effect for the discontinuous LSFEM.
In particular, the adaptive mesh refinement seems not to be able to compensate 
for the lack of a proper weighting in Figure~\ref{fig:example:scaling_lshape:one}.

For the unscaled L-shape \(\Omega_1\) with weighting \(c_\Omega = \lambda(1)^{-1/2}\), the
adaptive algorithm results in an increased refinement towards the singularity of the exact solution at the
reentrant corner as displayed in Figure~\ref{fig:mesh}.
Indeed, the convergence history plots in Figure~\ref{fig:convergence} exhibit the optimal convergence rate \(-(k+1)/2\)
for the polynomial degrees \(k = 0, \dots, 3\) in both schemes.
As expected for high-order discontinuous Galerkin methods, the condition number of the linear system increases
drastically during adaptive mesh refinement.
Therefore, the Figures~\ref{fig:convergence}--\ref{fig:efficiency} omit the values corresponding to 
unreliable solutions in case of very large condition numbers for \(k = 2,3\).
Undisplayed graphs confirm the same rate of convergence for the exact error in the weighted norm
\begin{equation}
    \label{eq:computation_norm}
    \vvvert (\tau, v) \vvvert^2 \coloneqq
    c_\Omega^2\,
    \Vert \Div_\PW \tau \Vert_{L^2(\Omega)}^2
    +
    \Vert \tau \Vert_{L^2(\Omega)}^2
    +
    \Vert \nabla_\PW v \Vert_{L^2(\Omega)}^2
    +
    j^2(v).
\end{equation}
According to Figure~\ref{fig:efficiency}, both schemes do not provide an asymptotically exact error estimator in the sense of~\cite{MR3820383}
as the efficiency indices appear to be uniformly bounded away from one.
This may be due to the facts that the minimization of the discrete least-squares functional~\eqref{eq:intro_discontinuous_LSFEM}
does not immediately follow from a continuous least-squares principle and that the equivalence constants from
Corollary~\ref{cor:Sigma_jump_equivalence} and Lemma~\ref{lem:U_jump_equivalence} introduce additional (\(k\)-dependent) factors.

\begin{figure}
    \centering
    \subfloat[\(c_\Omega = 1\)]{\label{fig:example:scaling_lshape:one}\input{figures/plot_dls_lshape_error_scalingOne.tex}}
    \hfill
    \subfloat[\(c_\Omega = \operatorname{diam}(\Omega) = 2\sqrt{2}\,\ell\)]{\input{figures/plot_dls_lshape_error_scalingDiam.tex}}

    \subfloat[\(c_\Omega = \operatorname{width}(\Omega) = 2\ell\)]{\input{figures/plot_dls_lshape_error_scalingWidth.tex}}
    \hfill
    \subfloat[\(c_\Omega = C_\FR = \ell/\pi\)]{\input{figures/plot_dls_lshape_error_scalingFriedrichs.tex}}

    \subfloat[Legend for Subfigures (a)--(d)]{\hspace*{2cm}\input{figures/legend_scal.tex}\hspace*{2cm}}
    \caption{%
        Convergence history plot for the relative energy error
        \(\Vert \nabla(u - u_\LS) \Vert_{L^2(\Omega)} / \Vert \nabla u \Vert_{L^2(\Omega)}\)
        under uniform refinement of the scaled L-shaped domain \(\Omega_\ell\) from Subsection~\ref{sec:lshape}.
        The experiments are carried out with the naturally penalized discontinuous LSFEM (\(\alpha = 1\)) 
        for lowest-order ansatz spaces (\(k = 0\)) and different weighting factors \(c_\Omega\).
    }
    \label{fig:example:scaling_lshape}
\end{figure}

\begin{figure}[p]
    \begin{minipage}[b]{0.49\textwidth}
        \centering
        \input{figures/plot_dls_lshape_mesh.tex}
        \caption{Adaptively refined mesh with \(854\)~triangles from discontinuous LSFEM with \(\alpha = 1\),
            polynomial degree \(k = 1\), and bulk parameter \(\theta = 0.5\).}
        \label{fig:mesh}
    \end{minipage}
    \hfill
    \begin{minipage}[b]{0.49\textwidth}
        \centering
        \input{figures/legend.tex}
        \caption{Legend for Figures~\ref{fig:convergence}--\ref{fig:efficiency}.
        \vspace{2\baselineskip}}
        \label{fig:legend}
    \end{minipage}
\end{figure}

\begin{figure}[p]
    \centering
    \subfloat[Naturally penalized discontinuous LSFEM \(\alpha = 1\).]{%
        \label{fig:dls:estimator}
        \input{figures/plot_dls_lshape_estimator_degree.tex}
    }
    \hfil
    \subfloat[Over-penalized discontinuous LSFEM \(\alpha = -1\).]{%
        \label{fig:opdls:estimator}
        \input{figures/plot_opdls_lshape_estimator_degree.tex}
    }
    \caption{%
        Convergence history plot of the built-in error estimator \eqref{eq:intro_discontinuous_LSFEM}
        for various choices of the polynomial degree \(k \in \N_0\) under uniform (\(\theta = 1\)) 
        and adaptive (\(\theta = 0.5\)) mesh refinement.
        Both graphs employ the line styles according to the legend in Figure~\ref{fig:legend}.
    }
    \label{fig:convergence}
\end{figure}

\begin{figure}[p]
    \centering
    \subfloat[Naturally penalized discontinuous LSFEM \(\alpha = 1\).]{%
        \label{fig:dls:efficiency}
        \input{figures/plot_dls_lshape_efficiency_degree.tex}
    }
    \hfil
    \subfloat[Over-penalized discontinuous LSFEM \(\alpha = -1\).]{%
        \label{fig:opdls:efficiency}
        \input{figures/plot_opdls_lshape_efficiency_degree.tex}
    }
    \caption{Plot of efficiency indices of the built-in error estimator \eqref{eq:intro_discontinuous_LSFEM}
        divided by the error in the weighted norm \eqref{eq:computation_norm} for various choices of the
        polynomial degree \(k \in \N_0\) under uniform (\(\theta = 1\)) and adaptive (\(\theta = 0.5\)) mesh refinement.
        Both graphs employ the line styles according to the legend in Figure~\ref{fig:legend}.}
    \label{fig:efficiency}
\end{figure}

{
    \bibliography{Paper}
}

\appendix

\section{Proof of Lemma~{\ref{lem:quasi_interpolation}}}
\label{app:quasi_interpolation}

\emph{Step~1}.\enskip
Let \(J\colon H^1(\Omega) \to S^1(\T)\) denote the lowest-order Scott--Zhang quasi-interpolation operator from~\cite{MR1011446}.
This operator can be constructed in such a way that it preserves homogeneous boundary conditions on \(\Gamma_\NEU\) \cite{MR3082295}.
Hence, \(w \in H^1_\NEU(\Omega)\) implies \(Jw \in S^1_\NEU(\T)\).

\emph{Step~2}.\enskip
Let \(Q_j \in P^k(-1, 1)\) for \(j = 0, \dots, k\) denote the Jacobi polynomials on \([-1,1]\) with respect to the weight \(r \mapsto (1+r)(1-r)\).
For each edge \(E = \conv\{z_1(E), z_2(E)\} \in \E\) with the two vertices \(z_1(E), z_2(E) \in \V\) and for all \(j,\ell = 0, \dots, k\),
the transformation formula verifies
\begin{equation}
    \label{eq:Jacobi_orthonormality}
    \delta_{j\ell}
    =
    \int_E
    \frac{8}{h_E^3}\,
    \vert x - z_1(E) \vert\,
    \vert z_2(E) - x \vert
    \,
    Q_j\Big(\frac{2 \vert x - z_1(E) \vert}{h_E}\Big)
    \,
    Q_\ell\Big(\frac{2 \vert x - z_1(E) \vert}{h_E}\Big)
    \ds.
\end{equation}
The functions \(q_{E, j} \in P^k(E)\) defined by \(q_{E, j}(x) \coloneqq Q_j(2 h_E^{-1} \, \vert x - z_1(E) \vert )\)
for \(j = 0, \dots, k\) form a basis of \(P^k(E)\).
Let \(b_{E, j} \in S^{k+2}(\T)\) denote an edge-bubble functions such that, for all \(x \in E\),
\[
    b_{E, j}(x)
    =
    \frac{2}{h_E}\,
    \frac{2 \vert x - z_1(E) \vert}{h_E}\,
    \frac{2 \vert z_2(E) - x \vert}{h_E}
    \,
    q_{E, j}(x)
\]
and \(b_{E, j}(z) = 0\) in all Lagrange interpolation nodes \(z\) which belong to \(\overline\Omega \setminus E\).
Thus, for all \(E \in \E\) and \(j, \ell = 0, \dots, k\), the orthonormality~\eqref{eq:Jacobi_orthonormality} reads
\(\delta_{jk} = (b_{E, j}, q_{E, \ell})_{L^2(E)}\).
Moreover,
\begin{equation}
    \label{eq:scaling_bubbles}
    \Vert q_{E, j} \Vert_{L^2(E)}
    \lesssim
    h_E^{1/2}
    \quad\text{and}\quad
    \Vert \nabla b_{E, j} \Vert_{L^2(\omega_E)}
    \lesssim
    h_E^{-1}.
\end{equation}

\emph{Step~3}.\enskip
Define the quasi-interpolation \(w_h \in S^{k+2}(\T)\) by
\[
    w_h
    \coloneqq
    Jw
    +
    \sum_{E \in \E} \sum_{j = 0}^k
    ((1 - J)w, q_{E, j})_{L^2(E)} \,
    b_{E, j}.
\]
It satisfies, for all \(E \in \E\) and \(\ell = 0, \dots, k\),
\[
    (w - w_h, q_{E, \ell})_{L^2(E)}
    =
    ((1 - J)w, q_{E, \ell})_{L^2(E)}
    -
    \sum_{E \in \E} \sum_{j = 0}^k
    ((1 - J)w, q_{E, j})_{L^2(E)}\,
    (b_{E, j},q_{E, \ell})_{L^2(E)}
    =
    0.
\]
Since the function \(q_{E, \ell}\) for \(\ell = 0, \dots, k\) form a basis of \(P^k(E)\),
this implies the orthogonality in~\eqref{eq:quasi_interpolation}.
If \(w \in H^1_\NEU(\Omega)\), then the inclusions \(Jw, (1 - J)w \in H^1_\NEU(\Omega)\) lead to \(w_h \in S^{k+2}_\NEU(\T)\).

\emph{Step~4}.\enskip
For \(E \in \E\) and \(j = 0, \dots, k\), abbreviate the coefficients \(\alpha_{E, j} \coloneqq ((1 - J)w, q_{E, j})_{L^2(E)}\).
The finite overlap of the edge patches \(\omega_E\) and the second estimate in~\eqref{eq:scaling_bubbles} verify the localization
\[
    \Big\Vert
        \nabla
        \Big(
            \sum_{E \in \E}
            \sum_{j = 0}^k
            \alpha_{E, j}\, b_{E, j}
        \Big)
    \Big\Vert_{L^2(\Omega)}^2
    \lesssim
    \sum_{E \in \E}
    \sum_{j = 0}^k
    \vert \alpha_{E, j} \vert^2\,
    \Vert \nabla b_{E, j} \Vert_{L^2(\omega_E)}^2
    \lesssim
    \sum_{E \in \E}
    \sum_{j = 0}^k
    h_E^{-2}
    \vert \alpha_{E, j} \vert^2.
\]
The first estimate in~\eqref{eq:scaling_bubbles}, the trace inequality,
and the local first-order approximation and stability properties of the Scott--Zhang operator \cite[Thm.~3.1 and Eqn.~(4.3)]{MR1011446}
prove, for all \(E \in \E\) and \(j = 0, \dots, k\),
\begin{align*}
    h_E^{-1}\, \vert \alpha_{E, j} \vert
    \leq
    h_E^{-1} \,
    \Vert (1 - J) w \Vert_{L^2(E)} \,
    \Vert q_{E, j} \Vert_{L^2(E)}
    &\lesssim
    h_E^{-1} \,
    \Vert (1 - J) w \Vert_{L^2(\omega_E)}
    +
    \Vert \nabla (1 - J) w \Vert_{L^2(\omega_E)}
    \\
    &\lesssim
    \Vert \nabla w \Vert_{L^2(\Omega_E)}.
\end{align*}
The latter estimate and a finite overlap argument also prove the stability estimate
\(\Vert \nabla J w \Vert_{L^2(\Omega)} \lesssim \Vert \nabla w \Vert_{L^2(\Omega)}\).
The combination of the two previously displayed formulas with the finite overlap of the enlarged patches \(\Omega_E\) results in
\begin{align*}
    \Vert \nabla w_h \Vert_{L^2(\Omega)}
    &\leq
    \Vert \nabla J w \Vert_{L^2(\Omega)}
    +
    \Big\Vert
        \nabla
        \Big(
            \sum_{E \in \E}
            \sum_{j = 0}^k
            \alpha_{E, j}\, b_{E, j}
        \Big)
    \Big\Vert_{L^2(\Omega)}
    \\
    &\lesssim
    \Vert \nabla w \Vert_{L^2(\Omega)}
    +
    \bigg(
        \sum_{E \in \E}
        \sum_{j = 0}^k
        h_E^{-2}
        \vert \alpha_{E, j} \vert^2
    \bigg)^{1/2}
    \lesssim
    \Vert \nabla w \Vert_{L^2(\Omega)}.
\end{align*}
The constants in the estimates~\eqref{eq:scaling_bubbles} depend exclusively on the maximum norm of \(Q_j\) and \(Q_j'\)
for \(j = 0, \dots, k\) on \([-1, 1]\).
The overlap of the patches \(\omega_E\) solely relates to the spatial dimension, whereas the overlap of the patches \(\Omega_E\)
depends on the interior angles in \(\T\) as well.
The remaining estimates in this proof are local with generic constants depending on the polynomial degree \(k\) and the interior angles.
This concludes the proof of Lemma~\ref{lem:quasi_interpolation}.
\hfill \qed

\section{Piecewise Raviart--Thomas basis with vanishing mean normal jump}
\label{app:Sigma_basis}

This section constructs an explicit basis for the space \(\Sigma(\T)\) from~\eqref{eq:Sigma_RT} employing a correction in terms of
the lowest-order edge-oriented basis functions \(\psi_E\) of \(RT^0(\T)\).
To this end, recall from \cite[Def.~4.1]{MR2194203} that, for \(E \in \E\),
\begin{equation}
    \label{eq:edge_basis_function}
    \psi_E(x)
    \coloneqq
    \begin{cases}
        \pm\frac{\vert E \vert}{2\vert T_\pm \vert}
        (x - P_{\pm})
        &\text{for }
        x \in T_\pm,
        \\
        0 &\text{otherwise}.
    \end{cases}
\end{equation}
Abbreviate its restriction to one of the two triangles \(T_+, T_- \in \T\) adjacent to an interior edge \(E \in \E(\Omega)\) by
\(\psi_{E,\pm} \in RT^{k,\PW}(\T)\) with \(\psi_{E,\pm} = \psi_E\) on \(T_\pm\) and \(\psi_{E,\pm} = 0\) otherwise.
For any boundary edge \(E \in \E(\partial\Omega)\), identify \(\psi_{E,+} = \psi_E\) and set \(\psi_{E,-} \coloneqq 0\).
They satisfy \([\psi_{E,\pm} \cdot n_E]_E \equiv \pm 1\) on \(E \in \E(\Omega)\) and
\([\psi_{E,+} \cdot n_E]_E = \psi_{E,+} \cdot n_E \equiv 1\) on \(E \in \E(\partial\Omega)\)
(using the convention that \(n_E = n\) on \(E \in \E(\partial\Omega)\)).

The inclusion \(\psi_E \in RT^0_\NEU(\T) \subseteq \Sigma(\T)\) motivates the direct sum
\(\Sigma(\T) = RT^0_\NEU(\T) \oplus Y(\T)\) with \(Y(\T) \cong \Sigma(\T) / RT^0_\NEU(\T)\).
Since \(\{\psi_E \;:\; E \in \E(\Omega) \cup \E(\Gamma_\DIR)\}\) provides a basis for \(RT^0_\NEU(\T)\),
it suffices to characterize the remaining basis functions in \(Y(\T)\).
For \(k = 0\), the space \(Y(\T) = \operatorname{span}(\{ \psi_E \;:\; E \in \E(\Gamma_\NEU) \})\)
solely consists of the basis functions belonging to the Neumann boundary.
For \(k \geq 1\), let \(\psi_{E, j} \in RT^{k, \PW}(\T)\) for \(j = 1,\dots,2(k+1)\)
denote the \(2(k + 1)\) basis functions associated to the edge \(E \in \E(\Omega)\) \cite{MR2970849}
with the enumeration such that \(\supp(\psi_{E, j}) = T_+\), \(\supp(\psi_{E, k + 1 + j}) = T_-\), and
\(\psi_{E, j} + \psi_{E, k + 1 + j} \in RT^k(\T)\) for \(j = 1, \dots, k + 1\).
The definition of the basis functions as the product of the lowest-order basis functions and Lagrange polynomials
with respect to Gaussian quadrature nodes \cite{MR2970849} yields
\begin{equation}
    \label{eq:sum_local_basis}
    \sum_{j = 1}^{k + 1} \psi_{E, j} = \psi_{E,+}
    \quad\text{and}\quad
    \sum_{j = k + 2}^{2(k + 1)} \psi_{E, j} = \psi_{E,-}.
\end{equation}
For \(j = 1, \dots, k + 1\), abbreviate
\[
    \alpha_{E,j}
    \coloneqq
    \frac1{\vert E \vert}
    \int_{E} [\psi_{E, j} \cdot n_E]_E \ds
    \quad\text{and}\quad
    \alpha_{E,k + 1 + j}
    \coloneqq
    - \frac1{\vert E \vert}
    \int_{E} [\psi_{E, k + 1 + j} \cdot n_E]_E \ds.
\]
Since \(\psi_{E, j} + \psi_{E, k + 1 + j} \in RT^k(\T)\), it holds that \(\alpha_{E,j} = \alpha_{E,j + k + 1}\).
The relation~\eqref{eq:sum_local_basis} ensures that
\[
    \sum_{j = 1}^{k + 1}
    \alpha_{E,j}
    =
    \sum_{j = 1}^{k + 1}
    \frac1{\vert E \vert}
    \int_{E} [\psi_{E, j} \cdot n_E]_E \ds
    =
    \frac1{\vert E \vert}
    \int_{E} [\psi_{E,+} \cdot n_E]_E \ds
    =
    1
    =
    \sum_{j = k + 1}^{2(k + 1)}
    \alpha_{E,j}.
\]
Define the corrected basis function \(\widetilde\psi_{E, j} \in RT^{k, \PW}(\T)\) by
\begin{equation}
    \label{eq:edge_basis_correction}
    \widetilde\psi_{E, j}
    \coloneqq
    \begin{cases}
        \psi_{E, j}
        -
        \alpha_{E, j}
        \psi_{E,+}
        &\text{for }
        j = 1, \dots, k+1,
        \\
        \psi_{E, j}
        -
        \alpha_{E, j}
        \psi_{E,-}
        &\text{for }
        j = k+2, \dots, 2(k+1).
    \end{cases}
\end{equation}
Since the functions \(\widetilde\psi_{E,j} \in RT^{k, \PW}(\T)\) for all \(j = 1,\dots, 2(k+1)\) satisfy the side condition
\[
    \int_E [\widetilde\psi_{E, j} \cdot n_E]_E \ds
    =
    \int_E [\psi_{E, j} \cdot n_E]_E \ds
    -
    \alpha_{E, j}\,
    \vert E \vert
    =
    0,
\]
they belong to the space \(\widetilde\psi_{E,j} \in Y(\T) \subseteq \Sigma(T)\).
However, the two equalities
\[
    \sum_{j = 1}^{k+1} \widetilde\psi_{E, j}
    =
    \bigg(
        1
        -
        \sum_{j = 1}^{k + 1} \alpha_{E, j}
    \bigg)
    \psi_{E,+}
    =
    0
    \quad\text{and}\quad
    \sum_{j = k + 1}^{2(k + 1)}
    \widetilde\psi_{E, j}
    =
    \bigg(
        1
        -
        \sum_{j = k + 1}^{2(k + 1)} \alpha_{E, j}
    \bigg)
    \psi_{E,-}
    =
    0
\]
show that there are two independent nontrivial linear combinations of zero.
Consequently, exactly one of the functions \(\widetilde\psi_{E,1}, \dots, \widetilde\psi_{E,k + 1}\)
and one of the functions \(\widetilde\psi_{E,k + 2}, \dots, \widetilde\psi_{E,2(k + 1)}\)
need to be neglected for the basis of \(Y(\T)\).
An analogous procedure for the basis functions belonging to a Neumann edge \(E \in \E(\Gamma_\NEU)\) leads to the \(k+1\) functions
\(\widetilde\psi_{E,j}\) of the form for \(j = 1,\dots, k+1\) from which one functions needs to be neglected.
Additionally, \(Y(\T)\) contains the basis functions \(\psi_{E,j} \in RT^{k, \PW}(\T)\) for \(j = 1,\dots, k+1\)
belonging to any Dirichlet edge \(E \in \E(\Gamma_\DIR)\) where again one basis function is redundant
due to \(\psi_E \in RT^0_\NEU(\T)\) already included in \(RT^0_\NEU(\T)\).
The \(k(k+1)\) basis functions \(\psi_{T, j} \in RT^{k,\PW}(\T)\) for \(j = 1,\dots, k(k+1)\)
associated with the interior of each triangle \(T \in \T\) complete the basis of \(Y(\T)\), namely
\begin{align*}
    Y(\T)
    =
    \operatorname{span}
    \big(
        &\{
            \widetilde\psi_{E,j}
            \;:\;
            E \in \E(\Omega),
            j = 1,\dots, k, k+2, \dots, 2k + 1
        \}
        \\
        &\cup
        \{
            \widetilde\psi_{E,j}
            \;:\;
            E \in \E(\Gamma_\NEU),
            j = 1,\dots, k
        \}
        \cup
        \{
            \psi_{E,j}
            \;:\;
            E \in \E(\Gamma_\DIR),
            j = 1,\dots, k
        \}
        \\
        &\cup
        \{
            \psi_{T,j}
            \;:\;
            T \in \T,
            j = 1, \dots, k (k + 1)
        \}
    \big).
\end{align*}
The combination with the basis for \(RT^0_\NEU(\T)\) provides a basis for \(\Sigma(\T)\) from~\eqref{eq:Sigma_RT}.

\end{document}

%% file: minarthd.tex
\usepackage[%
    includehead,%
    inner  = 3cm,%
    outer  = 3cm,%
    top    = 2.5cm,%
    bottom = 3.5cm%
]{geometry}
\usepackage{lmodern}
\usepackage{microtype}
\usepackage{titlesec}
\titleformat{\section}[runin]{\bfseries}{\thesection.}{0.5em}{}[.]
\titleformat{\subsection}[runin]{\bfseries}{\thesubsection.}{0.5em}{}[.]
\titleformat{\subsubsection}[runin]{\bfseries}{\thesubsection.}{0.5em}{}[.]
\renewcommand\maketitle{
    \begin{center}
        {\bfseries\MakeUppercase{\CustomTitle}} \par\medskip
        {\footnotesize\MakeUppercase{\AuthorA}} \par
        {\footnotesize\AffiliationA} \par
    \end{center}
    \thispagestyle{plain}
}
\newcommand\CustomAbstractFormat{%
    {
        \footnotesize
        \hspace\parindent
        \textbf{Abstract.}
        \CustomAbstract
        \medskip

        \textbf{Keywords.}
        \CustomKeywords
        \medskip

        \textbf{AMS subject classification.}
        \CustomClassification
        \medskip

        \textbf{Acknowledgement.}
        \CustomFunding
    }
}
\usepackage[british]{babel}
\renewenvironment{thebibliography}[1]{%
    \titleformat{\section}{\centering}{}{}{\MakeUppercase}%
    \begin{oldthebibliography}{#1}%
    \footnotesize
    \setlength{\parskip}{0ex}%
    \setlength{\itemsep}{0ex plus 0.25ex}%
}{%
    \end{oldthebibliography}%
}

\usepackage{amsmath}
\usepackage{amssymb}
\usepackage{mathtools}
\usepackage{amsthm}
\theoremstyle{plain}
\newtheorem{theorem}{Theorem}[section]
\newtheorem{lemma}[theorem]{Lemma}
\newtheorem{corollary}[theorem]{Corollary}
\theoremstyle{definition}

\theoremstyle{remark}
\newtheorem{remark}[theorem]{Remark}
\makeatletter
\thm@headfont{\scshape}
\makeatother
\makeatletter
\let\thm@indent\indent
\renewenvironment{proof}[1][\proofname]{\par
    \normalfont
    \topsep6\p@\@plus6\p@ \trivlist
    \item[\hskip4\labelsep\itshape
        #1\@addpunct{.}]\ignorespaces
}{%
    \qed\endtrivlist
}
\makeatother
\usepackage{tikz} 
\usetikzlibrary{%
    calc,%
    matrix,%
}
\definecolor{TUblue}{rgb}{0,0.4,0.6}
\definecolor{TUgray}{rgb}{0.3922,0.3882,0.3882}
\definecolor{TUgreen}{rgb}{0,0.4941,0.4431}
\definecolor{TUmagenta}{rgb}{0.7294,0.2745,0.5098}
\definecolor{TUyellow}{rgb}{0.8824,0.5373,0.1333}

\usepackage[font=footnotesize,labelfont=bf]{caption}
\usepackage[caption=false]{subfig}
\usepackage{algorithm}
\usepackage{algpseudocode}
\algrenewcommand\algorithmicrequire{\textbf{Input:}}
\algrenewcommand\algorithmicensure{\textbf{Output:}}
\usepackage{pgfplotstable}
\usepackage{pgfplots}
\pgfplotsset{%
    compat            = newest,%
    every axis/.style = {scale only axis},%
    grid style        = {densely dotted, semithick},%
}
\newcommand\drawslopetriangle[4][ST]{
    \pgfplotsextra
    {
        \pgfkeys{/pgf/fpu=true}
        \pgfmathsetmacro\leftcoord{#3}
        \pgfmathsetmacro\rightcoord{10*#3}
        \pgfmathsetmacro\bottomcoord{#4}
        \pgfmathsetmacro\topcoord{10^(#2)*#4}
        \pgfkeys{/pgf/fpu=false}
        \coordinate (#1-BL) at (axis cs:\leftcoord,\bottomcoord);
        \coordinate (#1-BR) at (axis cs:\rightcoord,\bottomcoord);
        \coordinate (#1-TL) at (axis cs:\leftcoord,\topcoord);
        \shadedraw[%
            bottom color = black!20,%
            middle color = black!5,%
            top color    = white,%
            draw         = black,%
            font         = \footnotesize%
        ]
        (#1-TL) -- (#1-BL) node[midway, left] {\(#2\)}
        -- (#1-BR) node[midway, below] {\(1\)} -- (#1-TL);
    }
}
\newcommand\drawswappedslopetriangle[4][SST]{
    \pgfplotsextra
    {
        \pgfkeys{/pgf/fpu=true}
        \pgfmathsetmacro\leftcoord{#3/10}
        \pgfmathsetmacro\rightcoord{#3}
        \pgfmathsetmacro\topcoord{#4}
        \pgfmathsetmacro\bottomcoord{10^(-#2)*#4}
        \pgfkeys{/pgf/fpu=false}
        \coordinate (#1-TR) at (axis cs:\rightcoord,\topcoord);
        \coordinate (#1-BR) at (axis cs:\rightcoord,\bottomcoord);
        \coordinate (#1-TL) at (axis cs:\leftcoord,\topcoord);
        \shadedraw[%
            bottom color = black!20,%
            middle color = black!5,%
            top color    = white,%
            draw         = black,%
            font         = \footnotesize%
        ]
        (#1-BR) -- (#1-TR) node[midway, right] {\(#2\)}
        -- (#1-TL) node[midway, above] {\(1\)} -- (#1-BR);
    }
}

\pgfplotsset{
    colormap={parula}{
        rgb=(0.2081,0.1663,0.5292)
        rgb=(0.2116,0.1898,0.5777)
        rgb=(0.2123,0.2138,0.627)
        rgb=(0.2081,0.2386,0.6771)
        rgb=(0.1959,0.2645,0.7279)
        rgb=(0.1707,0.2919,0.7792)
        rgb=(0.1253,0.3242,0.8303)
        rgb=(0.0591,0.3598,0.8683)
        rgb=(0.0117,0.3875,0.882)
        rgb=(0.006,0.4086,0.8828)
        rgb=(0.0165,0.4266,0.8786)
        rgb=(0.0329,0.443,0.872)
        rgb=(0.0498,0.4586,0.8641)
        rgb=(0.0629,0.4737,0.8554)
        rgb=(0.0723,0.4887,0.8467)
        rgb=(0.0779,0.504,0.8384)
        rgb=(0.0793,0.52,0.8312)
        rgb=(0.0749,0.5375,0.8263)
        rgb=(0.0641,0.557,0.824)
        rgb=(0.0488,0.5772,0.8228)
        rgb=(0.0343,0.5966,0.8199)
        rgb=(0.0265,0.6137,0.8135)
        rgb=(0.0239,0.6287,0.8038)
        rgb=(0.0231,0.6418,0.7913)
        rgb=(0.0228,0.6535,0.7768)
        rgb=(0.0267,0.6642,0.7607)
        rgb=(0.0384,0.6743,0.7436)
        rgb=(0.059,0.6838,0.7254)
        rgb=(0.0843,0.6928,0.7062)
        rgb=(0.1133,0.7015,0.6859)
        rgb=(0.1453,0.7098,0.6646)
        rgb=(0.1801,0.7177,0.6424)
        rgb=(0.2178,0.725,0.6193)
        rgb=(0.2586,0.7317,0.5954)
        rgb=(0.3022,0.7376,0.5712)
        rgb=(0.3482,0.7424,0.5473)
        rgb=(0.3953,0.7459,0.5244)
        rgb=(0.442,0.7481,0.5033)
        rgb=(0.4871,0.7491,0.484)
        rgb=(0.53,0.7491,0.4661)
        rgb=(0.5709,0.7485,0.4494)
        rgb=(0.6099,0.7473,0.4337)
        rgb=(0.6473,0.7456,0.4188)
        rgb=(0.6834,0.7435,0.4044)
        rgb=(0.7184,0.7411,0.3905)
        rgb=(0.7525,0.7384,0.3768)
        rgb=(0.7858,0.7356,0.3633)
        rgb=(0.8185,0.7327,0.3498)
        rgb=(0.8507,0.7299,0.336)
        rgb=(0.8824,0.7274,0.3217)
        rgb=(0.9139,0.7258,0.3063)
        rgb=(0.945,0.7261,0.2886)
        rgb=(0.9739,0.7314,0.2666)
        rgb=(0.9938,0.7455,0.2403)
        rgb=(0.999,0.7653,0.2164)
        rgb=(0.9955,0.7861,0.1967)
        rgb=(0.988,0.8066,0.1794)
        rgb=(0.9789,0.8271,0.1633)
        rgb=(0.9697,0.8481,0.1475)
        rgb=(0.9626,0.8705,0.1309)
        rgb=(0.9589,0.8949,0.1132)
        rgb=(0.9598,0.9218,0.0948)
        rgb=(0.9661,0.9514,0.0755)
        rgb=(0.9763,0.9831,0.0538)
    }
}

%% file: figures/plot_dls_square_error_scalingOne.tex
\begin{tikzpicture}[>=stealth]
    %
    %
    \colorlet{col1}{TUblue}
    \colorlet{col2}{TUgreen}
    \colorlet{col3}{TUmagenta}
    \colorlet{col4}{TUyellow}
    \colorlet{col5}{purple}
    \colorlet{col6}{green}
    \pgfplotsset{%
        linedefault/.style = {%
            mark = *,%
            mark size = 2pt,%
            every mark/.append style = {solid},%
            gray,%
            every mark/.append style = {fill = gray!60!white}%
        },%
        line1/.style = {%
            linedefault,%
            col1,%
            every mark/.append style = {fill = col1!60!white}%
        },%
        line2/.style = {%
            linedefault,%
            mark = triangle*,%
            mark size = 2.75pt,%
            col2,%
            every mark/.append style = {fill = col2!60!white}%
        },%
        line3/.style = {%
            linedefault,%
            mark = square*,%
            mark size = 1.66pt,%
            col3,%
            every mark/.append style = {fill = col3!60!white}%
        },%
        line4/.style = {%
            linedefault,%
            mark = pentagon*,%
            mark size = 2.2pt,%
            col4,%
            every mark/.append style = {fill = col4!60!white}%
        },%
        line5/.style = {%
            linedefault,%
            mark = diamond*,%
            mark size = 2.75pt,%
            col5,%
            every mark/.append style = {fill = col5!60!white}%
        },%
        line6/.style = {%
            linedefault,%
            mark = halfsquare*,%
            mark size = 1.66pt,%
            col6,%
            every mark/.append style = {fill = col6!60!white}%
        },%
        minorline/.style = {%
            dashed,%
            every mark/.append style = {fill = black!20!white}%
        },%
        majorline/.style = {%
            solid%
        }%
    }

    %
    %
    \begin{loglogaxis}[%
            width            = 0.36\textwidth,%
            xlabel           = ndof,%
            ylabel           = {relative energy error},%
            ymajorgrids      = true,%
            font             = \footnotesize,%
            grid style       = {densely dotted, semithick},%
            legend style     = {legend pos  = south west}%
        ]

        %
        %
        \pgfplotstableread[col sep=comma]{
ndof,nelem,res,err,res1,res2,res3,res4,err1,err2,err3,relerr1,relerr2,relerr3,cond
6.400000000000000000e+01,8.000000000000000000e+00,5.824570832028575396e+00,5.977244475212094343e+00,5.782112461670300796e+00,6.633255354617408894e-01,1.512260111074443762e-16,2.297827058755128882e-01,1.065049253138570906e+00,5.782112461670299908e+00,1.077170871050412426e+00,2.158447365525732464e-01,5.935338977916125580e-02,2.183013247498441223e-01,nan
2.480000000000000000e+02,3.200000000000000000e+01,2.602371996235329110e+00,2.685125338769466907e+00,2.535070454618536306e+00,4.269606730335537281e-01,5.123007166607173057e-16,4.043048114885007260e-01,5.095882467238075364e-01,2.535070454618536751e+00,7.236267641698759245e-01,1.032641686565600458e-01,2.602498830148856565e-02,1.466374405219338373e-01,nan
9.760000000000000000e+02,1.280000000000000000e+02,1.327125105375399627e+00,1.364790322148153212e+00,1.285757090183150764e+00,2.251127377995592316e-01,8.983284403256719854e-16,2.396121984411689054e-01,2.532835614963079296e-01,1.285757090183150986e+00,3.812201017112842871e-01,5.132598049590578193e-02,1.319955947165485539e-02,7.725134386727980362e-02,nan
3.872000000000000000e+03,5.120000000000000000e+02,6.675928494082248310e-01,6.856738648501695677e-01,6.451905433792359812e-01,1.147501483277262452e-01,2.932740124430430338e-14,1.274432374579519966e-01,1.262431509775446825e-01,6.451905433792360922e-01,1.947831576493090078e-01,2.558221096756643428e-02,6.623514669221424421e-03,3.947132017324025965e-02,nan
1.542400000000000000e+04,2.048000000000000000e+03,3.344676124072866741e-01,3.433546575684364721e-01,3.228852943694279776e-01,5.781886798415108303e-02,5.485731744465843133e-15,6.535020225240387171e-02,6.303238677978414972e-02,3.228852943694281441e-01,9.830781134146816813e-02,1.277303207265883568e-02,3.314734702295095078e-03,1.992132761281052364e-02,nan
6.156800000000000000e+04,8.192000000000000000e+03,1.673581016659725917e-01,1.717684039048081035e-01,1.614789286726169282e-01,2.900610726176727772e-02,1.044091967505671859e-13,3.304444778628871926e-02,3.149634503287283144e-02,1.614789286726171502e-01,4.936517420529724226e-02,6.382493918276365905e-03,1.657739816258531340e-03,1.000347576238196437e-02,nan
2.460160000000000000e+05,3.276800000000000000e+04,8.370461695986376738e-02,8.590215287799810029e-02,8.074400021075507639e-02,1.452536059883421070e-02,2.692954524005542241e-12,1.660973298511842985e-02,1.574358018206599777e-02,8.074400021075506251e-02,2.473309487491211905e-02,3.190316357629961134e-03,8.289164733358427061e-04,5.011972895728674347e-03,nan
        }\tableUniformOne

        \pgfplotstableread[col sep=comma]{
ndof,nelem,res,err,res1,res2,res3,res4,err1,err2,err3,relerr1,relerr2,relerr3,cond
6.400000000000000000e+01,8.000000000000000000e+00,7.857406276729371752e-01,2.218510490995226014e+00,6.795362753173194159e-01,3.727537301046064955e-01,1.021328880308320395e-16,1.291256798503988001e-01,1.470394155499268773e+00,6.795362753173194159e-01,1.515902462778306647e+00,2.979926404219448433e-01,6.975440496073854746e-01,3.072154332333111859e-01,nan
2.480000000000000000e+02,3.200000000000000000e+01,5.610885180281114382e-01,1.203364742039569224e+00,3.423348877958233838e-01,3.256404013458841229e-01,3.801621560167472951e-16,3.026309281486193603e-01,7.401190586389244652e-01,3.423348877958234393e-01,8.849391538525188672e-01,1.499794781150977918e-01,3.514403883271292606e-01,1.793261650395727713e-01,nan
9.760000000000000000e+02,1.280000000000000000e+02,3.375144326323609723e-01,5.438613833934129449e-01,1.506719687649808181e-01,2.072868266538815119e-01,6.654142625707036514e-16,2.196499978588501434e-01,3.121468927721309239e-01,1.506719687649808181e-01,4.191031847979508895e-01,6.325418529190070249e-02,1.546795757619652822e-01,8.492806150400403209e-02,nan
3.872000000000000000e+03,5.120000000000000000e+02,1.809498642876982699e-01,2.518188696736184395e-01,6.811287206303717756e-02,1.122392604669124605e-01,1.300634829062521559e-15,1.245224415713943877e-01,1.360179626673898856e-01,6.811287206303719144e-02,2.006800777726470419e-01,2.756300194815859347e-02,6.992455359147266292e-02,4.066628602672212045e-02,nan
1.542400000000000000e+04,2.048000000000000000e+03,9.273646775303205636e-02,1.226349340923030135e-01,3.277970162358899431e-02,5.749170756708794905e-02,2.058309173436487412e-15,6.496342951074915439e-02,6.437977379306429604e-02,3.277970162358898043e-02,9.909623074388780595e-02,1.304606976667855439e-02,3.365158351816125576e-02,2.008109478693000799e-02,nan
6.156800000000000000e+04,8.192000000000000000e+03,4.680189085613303823e-02,6.093178511744803960e-02,1.621119189984216585e-02,2.896456799108897487e-02,4.182860168604945501e-15,3.299497001037377669e-02,3.167042072448816165e-02,1.621119189984217279e-02,4.946578763063684547e-02,6.417769028512013314e-03,1.664238083710630270e-02,1.002386430507340269e-02,nan
2.460160000000000000e+05,3.276800000000000000e+04,2.349117431346321380e-02,3.043406005337259690e-02,8.082402781182289678e-03,1.452013444999486133e-02,4.491945273503703007e-14,1.660348552148985593e-02,1.576560218255703555e-02,8.082402781182261922e-03,2.474575851406417629e-02,3.194778947941809364e-03,8.297380352682756297e-03,5.014539085544886428e-03,nan
        }\tableUniformTwo

        \pgfplotstableread[col sep=comma]{
ndof,nelem,res,err,res1,res2,res3,res4,err1,err2,err3,relerr1,relerr2,relerr3,cond
6.400000000000000000e+01,8.000000000000000000e+00,9.829653030123151891e-02,3.118536848558413510e+00,9.790142890413658372e-02,8.319409737669449662e-03,3.429770089157427023e-18,2.881928070925343036e-03,2.203322324403350851e+00,9.790142890413658372e-02,2.204780742119493819e+00,4.465291396146996439e-01,1.004958258457235232e+01,4.468247050890631900e-01,nan
2.480000000000000000e+02,3.200000000000000000e+01,9.745797204608264153e-02,3.050956295641871208e+00,9.573464064389757000e-02,1.398966905882345797e-02,1.675629933906204639e-17,1.171426718763997864e-02,2.154700421632864238e+00,9.573464064389757000e-02,2.157877496410004703e+00,4.366335942289736316e-01,9.828101220088379719e+00,4.372774041828538238e-01,nan
9.760000000000000000e+02,1.280000000000000000e+02,9.462194481254627976e-02,2.810917514110204962e+00,8.823545353024880855e-02,2.465315136075517388e-02,6.171186099853503710e-17,2.366514965943703960e-02,1.985049750362527066e+00,8.823545353024878080e-02,1.988227669171082290e+00,4.022551805913183731e-01,9.058236001755613387e+00,4.028991615817225092e-01,nan
3.872000000000000000e+03,5.120000000000000000e+02,8.579094740043224920e-02,2.157456887177632954e+00,6.790366052029141175e-02,3.627243075850075660e-02,2.618492854665533285e-16,3.786146207792982082e-02,1.523302710789681491e+00,6.790366052029137012e-02,1.526289017263426251e+00,3.086856674055855065e-01,6.970977739294243136e+00,3.092908196188626779e-01,nan
1.542400000000000000e+04,2.048000000000000000e+03,6.571909110909364415e-02,1.134163783712789808e+00,3.589319434524300412e-02,3.694577154365034227e-02,8.975543460860748686e-16,4.081283500253767016e-02,8.002443739200585382e-01,3.589319434524297636e-02,8.028998124888930832e-01,1.621634143374048043e-01,3.684788962122008460e+00,1.627015186951724646e-01,nan
6.156800000000000000e+04,8.192000000000000000e+03,4.030423109595750869e-02,3.963381202398773007e-01,1.257472837928853317e-02,2.535209375935540410e-02,2.935917186000901098e-15,2.869805903611455725e-02,2.789140632280114329e-01,1.257472837928850715e-02,2.813053982686919974e-01,5.651980604151264631e-02,1.290919383992520597e+00,5.700439183512695518e-02,nan
2.460160000000000000e+05,3.276800000000000000e+04,2.156819303766166485e-02,1.127085917635530926e-01,3.559974838846860501e-03,1.401648656882630017e-02,8.383968939389572952e-15,1.600161349122413909e-02,7.852861744622204043e-02,3.559974838846863103e-03,8.077010262748807179e-02,1.591322493889512496e-02,3.654663852272512781e-01,1.636744480226039536e-02,nan
        }\tableUniformThree

        \pgfplotstableread[col sep=comma]{
ndof,nelem,res,err,res1,res2,res3,res4,err1,err2,err3,relerr1,relerr2,relerr3,cond
6.400000000000000000e+01,8.000000000000000000e+00,9.869667710478499320e-03,3.141210289127430944e+00,9.869265136414114709e-03,8.423172980568547567e-05,2.575242544887425464e-20,2.917872712657214977e-05,2.221152719966857259e+00,9.869265136414114709e-03,2.221167547194388803e+00,4.501426786333805330e-01,1.013080157742658542e+02,4.501456835451355110e-01,nan
2.480000000000000000e+02,3.200000000000000000e+01,9.868337131701748061e-03,3.140652958340188761e+00,9.866531926725615950e-03,1.451255902374350355e-04,1.965288740626901341e-19,1.206815929867755194e-04,2.220749909762239849e+00,9.866531926725617685e-03,2.220782180802458150e+00,4.500180188614515409e-01,1.012896416750420059e+02,4.500245583414361206e-01,nan
9.760000000000000000e+02,1.280000000000000000e+02,9.865219739475631494e-03,3.137849592594695647e+00,9.857723911485670454e-03,2.809190901820016488e-04,7.413604557887979381e-19,2.625342983966443390e-04,2.218767267603718452e+00,9.857723911485668719e-03,2.218800284666537781e+00,4.496162515609709032e-01,1.011992187468897555e+02,4.496229422167414280e-01,nan
3.872000000000000000e+03,5.120000000000000000e+02,9.853377012979254793e-03,3.126886337180413200e+00,9.823292263769744723e-03,5.512728031890763410e-04,3.501765449461282533e-18,5.367177628265374289e-04,2.211015115364882710e+00,9.823292263769748192e-03,2.211048128875229946e+00,4.480453370797411461e-01,1.008457440624385413e+02,4.480520270156302942e-01,nan
1.542400000000000000e+04,2.048000000000000000e+03,9.807530709808281663e-03,3.084210173597199400e+00,9.689378464029481733e-03,1.076967669073022969e-03,2.165303877831281515e-17,1.069459790925144423e-03,2.180838750335053700e+00,9.689378464029483468e-03,2.180871581708214713e+00,4.419303270340517731e-01,9.947098737064762020e+01,4.419369800612302779e-01,nan
6.156800000000000000e+04,8.192000000000000000e+03,9.633218291688285320e-03,2.926130027568388137e+00,9.194707356821357064e-03,2.016700497103481819e-03,1.603255604266407968e-16,2.046746291731230214e-03,2.069059997602616274e+00,9.194707356821344921e-03,2.069092342529608164e+00,4.192792159682247055e-01,9.439270256214339838e+01,4.192857704208000080e-01,nan
2.460160000000000000e+05,3.276800000000000000e+04,9.026969530177987588e-03,2.436085291149299170e+00,7.668239454398222479e-03,3.273109803154469321e-03,1.331849533357026990e-15,3.459918320606191963e-03,1.722548424832619718e+00,7.668239454398223347e-03,1.722579364779543765e+00,3.490612905705714475e-01,7.872201016352275360e+01,3.490675603146597972e-01,nan
        }\tableUniformFour

        %
        %
        \addlegendimage{line1}
        \addlegendentry{\(\ell = \hphantom{0}1\hphantom{{}^2}\)}
        \addlegendimage{line2}
        \addlegendentry{\(\ell = 10\hphantom{{}^2}\)}
        \addlegendimage{line3}
        \addlegendentry{\(\ell = 10^2\)}
        \addlegendimage{line4}
        \addlegendentry{\(\ell = 10^3\)}

        %
        %
        \addplot+ [line1, forget plot] table [x=ndof, y=relerr1] {\tableUniformOne};
        \addplot+ [line2, forget plot] table [x=ndof, y=relerr1] {\tableUniformTwo};
        \addplot+ [line3, forget plot] table [x=ndof, y=relerr1] {\tableUniformThree};
        \addplot+ [line4, forget plot] table [x=ndof, y=relerr1] {\tableUniformFour};

        %
        %
        \drawslopetriangle[ST1]{0.5}{7e3}{4e-3}
    \end{loglogaxis}
\end{tikzpicture}

%% file: figures/plot_dls_square_error_scalingDiam.tex
\begin{tikzpicture}[>=stealth]
    %
    %
    \colorlet{col1}{TUblue}
    \colorlet{col2}{TUgreen}
    \colorlet{col3}{TUmagenta}
    \colorlet{col4}{TUyellow}
    \colorlet{col5}{purple}
    \colorlet{col6}{green}
    \pgfplotsset{%
        linedefault/.style = {%
            mark = *,%
            mark size = 2pt,%
            every mark/.append style = {solid},%
            gray,%
            every mark/.append style = {fill = gray!60!white}%
        },%
        line1/.style = {%
            linedefault,%
            col1,%
            every mark/.append style = {fill = col1!60!white}%
        },%
        line2/.style = {%
            linedefault,%
            mark = triangle*,%
            mark size = 2.75pt,%
            col2,%
            every mark/.append style = {fill = col2!60!white}%
        },%
        line3/.style = {%
            linedefault,%
            mark = square*,%
            mark size = 1.66pt,%
            col3,%
            every mark/.append style = {fill = col3!60!white}%
        },%
        line4/.style = {%
            linedefault,%
            mark = pentagon*,%
            mark size = 2.2pt,%
            col4,%
            every mark/.append style = {fill = col4!60!white}%
        },%
        line5/.style = {%
            linedefault,%
            mark = diamond*,%
            mark size = 2.75pt,%
            col5,%
            every mark/.append style = {fill = col5!60!white}%
        },%
        line6/.style = {%
            linedefault,%
            mark = halfsquare*,%
            mark size = 1.66pt,%
            col6,%
            every mark/.append style = {fill = col6!60!white}%
        },%
        minorline/.style = {%
            dashed,%
            every mark/.append style = {fill = black!20!white}%
        },%
        majorline/.style = {%
            solid%
        }%
    }

    %
    %
    \begin{loglogaxis}[%
            width            = 0.36\textwidth,%
            xlabel           = ndof,%
            ylabel           = {relative energy error},%
            ymajorgrids      = true,%
            font             = \footnotesize,%
            grid style       = {densely dotted, semithick},%
            legend style     = {legend pos  = south west}%
        ]

        %
        %
        \pgfplotstableread[col sep=comma]{
ndof,nelem,res,err,res1,res2,res3,res4,err1,err2,err3,relerr1,relerr2,relerr3,cond
6.400000000000000000e+01,8.000000000000000000e+00,8.206663850902652158e+00,8.314799755258867719e+00,8.176345598487792188e+00,6.659473755228546565e-01,1.990907266437116664e-16,2.306909379145465622e-01,1.062815445331813935e+00,8.176345598487792188e+00,1.074100066709976486e+00,2.153920291720106128e-01,4.196509757404362495e-02,2.176789901940358873e-01,nan
2.480000000000000000e+02,3.200000000000000000e+01,3.632982795565297796e+00,3.692511692315125771e+00,3.584922160348372078e+00,4.276412934161094581e-01,6.134453514609902047e-16,4.049938533879294811e-01,5.093667500000186310e-01,3.584922160348371634e+00,7.235476598825351102e-01,1.032192840361056546e-01,1.840137364128058978e-02,1.466214106418847196e-01,nan
9.760000000000000000e+02,1.280000000000000000e+02,1.847814312591796515e+00,1.875016863124936206e+00,1.818304432680139859e+00,2.252108941164258638e-01,9.432699448614904730e-16,2.397222915920180775e-01,2.532454942913133777e-01,1.818304432680139193e+00,3.812137808628918711e-01,5.131826646736976916e-02,9.333340519754131093e-03,7.725006299559832035e-02,nan
3.872000000000000000e+03,5.120000000000000000e+02,9.284128198659448161e-01,9.414937481932699193e-01,9.124330996374908631e-01,1.147631256314201736e-01,1.856383008823780383e-13,1.274583373995707147e-01,1.262378856443633901e-01,9.124330996374910852e-01,1.947827322967140840e-01,2.558114398798603736e-02,4.683511004732584561e-03,3.947123397878349316e-02,nan
1.542400000000000000e+04,2.048000000000000000e+03,4.648908050067093733e-01,4.713244714081266107e-01,4.566282307318228484e-01,5.782053010681648120e-02,6.665070859672262899e-15,6.535216734491797186e-02,6.303170362414495242e-02,4.566282307318223488e-01,9.830778496697788194e-02,1.277289363638288433e-02,2.343868656840408706e-03,1.992132226822137392e-02,nan
6.156800000000000000e+04,8.192000000000000000e+03,2.325600492378970818e-01,2.357536692180298432e-01,2.283656234745806424e-01,2.900631736943806716e-02,1.835927920033552518e-12,3.304469805071025973e-02,3.149625833305461087e-02,2.283656234745809477e-01,4.936517270148697845e-02,6.382476349219958281e-03,1.172198719084990501e-03,1.000347545764629364e-02,nan
2.460160000000000000e+05,3.276800000000000000e+04,1.163015749887662126e-01,1.178930496290545443e-01,1.141892516783707123e-01,1.452538700335275802e-02,1.144326036017877051e-11,1.660976454993120383e-02,1.574356927250954510e-02,1.141892516783708927e-01,2.473309480299789795e-02,3.190314146891615708e-03,5.861324157029215807e-04,5.011972881155806261e-03,nan
        }\tableUniformOne

        \pgfplotstableread[col sep=comma]{
ndof,nelem,res,err,res1,res2,res3,res4,err1,err2,err3,relerr1,relerr2,relerr3,cond
6.400000000000000000e+01,8.000000000000000000e+00,8.206663850902652158e+00,8.314799755258869496e+00,8.176345598487792188e+00,6.659473755228544345e-01,3.798502056565034185e-16,2.306909379145464511e-01,1.062815445331813935e+00,8.176345598487792188e+00,1.074100066709976709e+00,2.153920291720106128e-01,4.196509757404362495e-02,2.176789901940359151e-01,nan
2.480000000000000000e+02,3.200000000000000000e+01,3.632982795565298240e+00,3.692511692315126215e+00,3.584922160348372522e+00,4.276412934161095136e-01,8.588096331930980362e-16,4.049938533879292590e-01,5.093667500000187420e-01,3.584922160348372522e+00,7.235476598825349992e-01,1.032192840361056824e-01,1.840137364128059672e-02,1.466214106418846919e-01,nan
9.760000000000000000e+02,1.280000000000000000e+02,1.847814312591796959e+00,1.875016863124935318e+00,1.818304432680140081e+00,2.252108941164265299e-01,8.469862223466421594e-16,2.397222915920200481e-01,2.532454942913091589e-01,1.818304432680139193e+00,3.812137808628914271e-01,5.131826646736892955e-02,9.333340519754129358e-03,7.725006299559823708e-02,nan
3.872000000000000000e+03,5.120000000000000000e+02,9.284128198659450382e-01,9.414937481932723617e-01,9.124330996374918623e-01,1.147631256314188275e-01,1.707044695473087845e-13,1.274583373995668012e-01,1.262378856443726605e-01,9.124330996374918623e-01,1.947827322967158326e-01,2.558114398798791433e-02,4.683511004732588030e-03,3.947123397878384704e-02,nan
1.542400000000000000e+04,2.048000000000000000e+03,4.648908050067094289e-01,4.713244714081269438e-01,4.566282307318230149e-01,5.782053010681771632e-02,4.012496867575791185e-15,6.535216734491620938e-02,6.303170362414425854e-02,4.566282307318227374e-01,9.830778496697777091e-02,1.277289363638274382e-02,2.343868656840410875e-03,1.992132226822135310e-02,nan
6.156800000000000000e+04,8.192000000000000000e+03,2.325600492378959439e-01,2.357536692180162707e-01,2.283656234745809199e-01,2.900631736944774691e-02,1.161702849110912529e-13,3.304469805069178145e-02,3.149625833300409572e-02,2.283656234745813085e-01,4.936517270145271419e-02,6.382476349209721678e-03,1.172198719084992452e-03,1.000347545763934955e-02,nan
2.460160000000000000e+05,3.276800000000000000e+04,1.163015749887029576e-01,1.178930496289714441e-01,1.141892516783160338e-01,1.452538700339341300e-02,7.911180815443812496e-12,1.660976454982873371e-02,1.574356927248731289e-02,1.141892516783162004e-01,2.473309480286844594e-02,3.190314146887110198e-03,5.861324157026409892e-04,5.011972881129573773e-03,nan
        }\tableUniformTwo

        \pgfplotstableread[col sep=comma]{
ndof,nelem,res,err,res1,res2,res3,res4,err1,err2,err3,relerr1,relerr2,relerr3,cond
6.400000000000000000e+01,8.000000000000000000e+00,8.206663850902653934e+00,8.314799755258869496e+00,8.176345598487792188e+00,6.659473755228548786e-01,3.293684730871717475e-16,2.306909379145464234e-01,1.062815445331814157e+00,8.176345598487793964e+00,1.074100066709976709e+00,2.153920291720106128e-01,4.196509757404362495e-02,2.176789901940358873e-01,nan
2.480000000000000000e+02,3.200000000000000000e+01,3.632982795565298240e+00,3.692511692315126215e+00,3.584922160348372966e+00,4.276412934161097357e-01,7.363165360142812954e-16,4.049938533879294256e-01,5.093667500000185200e-01,3.584922160348372522e+00,7.235476598825353323e-01,1.032192840361056130e-01,1.840137364128058978e-02,1.466214106418847474e-01,nan
9.760000000000000000e+02,1.280000000000000000e+02,1.847814312591797403e+00,1.875016863124937094e+00,1.818304432680140525e+00,2.252108941164264466e-01,9.249603372038634746e-16,2.397222915920174391e-01,2.532454942913131002e-01,1.818304432680140081e+00,3.812137808628920377e-01,5.131826646736968589e-02,9.333340519754132827e-03,7.725006299559830647e-02,nan
3.872000000000000000e+03,5.120000000000000000e+02,9.284128198659458153e-01,9.414937481932709185e-01,9.124330996374919733e-01,1.147631256314203124e-01,1.099615659356772676e-13,1.274583373995695490e-01,1.262378856443634456e-01,9.124330996374920844e-01,1.947827322967143060e-01,2.558114398798604777e-02,4.683511004732588898e-03,3.947123397878353479e-02,nan
1.542400000000000000e+04,2.048000000000000000e+03,4.648908050067095954e-01,4.713244714081269993e-01,4.566282307318231815e-01,5.782053010681845184e-02,4.432767373739244201e-15,6.535216734491512691e-02,6.303170362414418915e-02,4.566282307318227929e-01,9.830778496697790969e-02,1.277289363638272474e-02,2.343868656840411308e-03,1.992132226822137392e-02,nan
6.156800000000000000e+04,8.192000000000000000e+03,2.325600492378950834e-01,2.357536692180263738e-01,2.283656234745828628e-01,2.900631736944081843e-02,6.233605678988960828e-13,3.304469805067847266e-02,3.149625833303889427e-02,2.283656234745830849e-01,4.936517270147063041e-02,6.382476349216773329e-03,1.172198719085000909e-03,1.000347545764298032e-02,nan
2.460160000000000000e+05,3.276800000000000000e+04,1.163015749889107220e-01,1.178930496294187807e-01,1.141892516783493128e-01,1.452538700392427481e-02,8.406426856172798522e-12,1.660976455059038834e-02,1.574356927427026862e-02,1.141892516783490630e-01,2.473309480371403690e-02,3.190314147248412332e-03,5.861324157028093658e-04,5.011972881300925421e-03,nan
        }\tableUniformThree

        \pgfplotstableread[col sep=comma]{
ndof,nelem,res,err,res1,res2,res3,res4,err1,err2,err3,relerr1,relerr2,relerr3,cond
6.400000000000000000e+01,8.000000000000000000e+00,8.206663850902652158e+00,8.314799755258869496e+00,8.176345598487792188e+00,6.659473755228546565e-01,1.023882225369320495e-16,2.306909379145465622e-01,1.062815445331814157e+00,8.176345598487792188e+00,1.074100066709976486e+00,2.153920291720106128e-01,4.196509757404361801e-02,2.176789901940358318e-01,nan
2.480000000000000000e+02,3.200000000000000000e+01,3.632982795565298240e+00,3.692511692315125771e+00,3.584922160348372522e+00,4.276412934161096802e-01,6.576999104044124367e-16,4.049938533879294256e-01,5.093667500000184090e-01,3.584922160348372078e+00,7.235476598825352212e-01,1.032192840361056130e-01,1.840137364128058978e-02,1.466214106418847474e-01,nan
9.760000000000000000e+02,1.280000000000000000e+02,1.847814312591796959e+00,1.875016863124936872e+00,1.818304432680140303e+00,2.252108941164263634e-01,1.031525423232391930e-15,2.397222915920173836e-01,2.532454942913132112e-01,1.818304432680139859e+00,3.812137808628921487e-01,5.131826646736974834e-02,9.333340519754131093e-03,7.725006299559838974e-02,nan
3.872000000000000000e+03,5.120000000000000000e+02,9.284128198659452602e-01,9.414937481932705854e-01,9.124330996374915292e-01,1.147631256314202153e-01,1.088282897719513504e-13,1.274583373995694657e-01,1.262378856443632791e-01,9.124330996374916403e-01,1.947827322967143615e-01,2.558114398798600961e-02,4.683511004732586296e-03,3.947123397878354173e-02,nan
1.542400000000000000e+04,2.048000000000000000e+03,4.648908050067090403e-01,4.713244714081262776e-01,4.566282307318225708e-01,5.782053010681841715e-02,8.652047693031471238e-15,6.535216734491595958e-02,6.303170362414360628e-02,4.566282307318221823e-01,9.830778496697788194e-02,1.277289363638260851e-02,2.343868656840407839e-03,1.992132226822137045e-02,nan
6.156800000000000000e+04,8.192000000000000000e+03,2.325600492378985806e-01,2.357536692180295101e-01,2.283656234745796987e-01,2.900631736944914857e-02,6.167821486416770957e-13,3.304469805071744842e-02,3.149625833305428474e-02,2.283656234745798652e-01,4.936517270149062137e-02,6.382476349219892361e-03,1.172198719084984429e-03,1.000347545764703090e-02,nan
2.460160000000000000e+05,3.276800000000000000e+04,1.163015749887800904e-01,1.178930496289479768e-01,1.141892516782774120e-01,1.452538700370485311e-02,1.198320819925216164e-11,1.660976455036190444e-02,1.574356927239316944e-02,1.141892516782772593e-01,2.473309480299630894e-02,3.190314146868032576e-03,5.861324157024408454e-04,5.011972881155483603e-03,nan
        }\tableUniformFour

        %
        %
        \addlegendimage{line1}
        \addlegendentry{\(\ell = \hphantom{0}1\hphantom{{}^2}\)}
        \addlegendimage{line2}
        \addlegendentry{\(\ell = 10\hphantom{{}^2}\)}
        \addlegendimage{line3}
        \addlegendentry{\(\ell = 10^2\)}
        \addlegendimage{line4}
        \addlegendentry{\(\ell = 10^3\)}

        %
        %
        \addplot+ [line1, forget plot] table [x=ndof, y=relerr1] {\tableUniformOne};
        \addplot+ [line2, forget plot] table [x=ndof, y=relerr1] {\tableUniformTwo};
        \addplot+ [line3, forget plot] table [x=ndof, y=relerr1] {\tableUniformThree};
        \addplot+ [line4, forget plot] table [x=ndof, y=relerr1] {\tableUniformFour};

        %
        %
        \drawslopetriangle[ST1]{0.5}{7e3}{4e-3}
    \end{loglogaxis}
\end{tikzpicture}

%% file: figures/plot_dls_square_error_scalingWidth.tex
\begin{tikzpicture}[>=stealth]
    %
    %
    \colorlet{col1}{TUblue}
    \colorlet{col2}{TUgreen}
    \colorlet{col3}{TUmagenta}
    \colorlet{col4}{TUyellow}
    \colorlet{col5}{purple}
    \colorlet{col6}{green}
    \pgfplotsset{%
        linedefault/.style = {%
            mark = *,%
            mark size = 2pt,%
            every mark/.append style = {solid},%
            gray,%
            every mark/.append style = {fill = gray!60!white}%
        },%
        line1/.style = {%
            linedefault,%
            col1,%
            every mark/.append style = {fill = col1!60!white}%
        },%
        line2/.style = {%
            linedefault,%
            mark = triangle*,%
            mark size = 2.75pt,%
            col2,%
            every mark/.append style = {fill = col2!60!white}%
        },%
        line3/.style = {%
            linedefault,%
            mark = square*,%
            mark size = 1.66pt,%
            col3,%
            every mark/.append style = {fill = col3!60!white}%
        },%
        line4/.style = {%
            linedefault,%
            mark = pentagon*,%
            mark size = 2.2pt,%
            col4,%
            every mark/.append style = {fill = col4!60!white}%
        },%
        line5/.style = {%
            linedefault,%
            mark = diamond*,%
            mark size = 2.75pt,%
            col5,%
            every mark/.append style = {fill = col5!60!white}%
        },%
        line6/.style = {%
            linedefault,%
            mark = halfsquare*,%
            mark size = 1.66pt,%
            col6,%
            every mark/.append style = {fill = col6!60!white}%
        },%
        minorline/.style = {%
            dashed,%
            every mark/.append style = {fill = black!20!white}%
        },%
        majorline/.style = {%
            solid%
        }%
    }

    %
    %
    \begin{loglogaxis}[%
            width            = 0.36\textwidth,%
            xlabel           = ndof,%
            ylabel           = {relative energy error},%
            ymajorgrids      = true,%
            font             = \footnotesize,%
            grid style       = {densely dotted, semithick},%
            legend style     = {legend pos  = south west}%
        ]

        %
        %
        \pgfplotstableread[col sep=comma]{
ndof,nelem,res,err,res1,res2,res3,res4,err1,err2,err3,relerr1,relerr2,relerr3,cond
6.400000000000000000e+01,8.000000000000000000e+00,5.824570832028575396e+00,5.977244475212094343e+00,5.782112461670300796e+00,6.633255354617408894e-01,1.512260111074443762e-16,2.297827058755128882e-01,1.065049253138570906e+00,5.782112461670299908e+00,1.077170871050412426e+00,2.158447365525732464e-01,5.935338977916125580e-02,2.183013247498441223e-01,nan
2.480000000000000000e+02,3.200000000000000000e+01,2.602371996235329110e+00,2.685125338769466907e+00,2.535070454618536306e+00,4.269606730335537281e-01,5.123007166607173057e-16,4.043048114885007260e-01,5.095882467238075364e-01,2.535070454618536751e+00,7.236267641698759245e-01,1.032641686565600458e-01,2.602498830148856565e-02,1.466374405219338373e-01,nan
9.760000000000000000e+02,1.280000000000000000e+02,1.327125105375399627e+00,1.364790322148153212e+00,1.285757090183150764e+00,2.251127377995592316e-01,8.983284403256719854e-16,2.396121984411689054e-01,2.532835614963079296e-01,1.285757090183150986e+00,3.812201017112842871e-01,5.132598049590578193e-02,1.319955947165485539e-02,7.725134386727980362e-02,nan
3.872000000000000000e+03,5.120000000000000000e+02,6.675928494082248310e-01,6.856738648501695677e-01,6.451905433792359812e-01,1.147501483277262452e-01,2.932740124430430338e-14,1.274432374579519966e-01,1.262431509775446825e-01,6.451905433792360922e-01,1.947831576493090078e-01,2.558221096756643428e-02,6.623514669221424421e-03,3.947132017324025965e-02,nan
1.542400000000000000e+04,2.048000000000000000e+03,3.344676124072866741e-01,3.433546575684364721e-01,3.228852943694279776e-01,5.781886798415108303e-02,5.485731744465843133e-15,6.535020225240387171e-02,6.303238677978414972e-02,3.228852943694281441e-01,9.830781134146816813e-02,1.277303207265883568e-02,3.314734702295095078e-03,1.992132761281052364e-02,nan
6.156800000000000000e+04,8.192000000000000000e+03,1.673581016659725917e-01,1.717684039048081035e-01,1.614789286726169282e-01,2.900610726176727772e-02,1.044091967505671859e-13,3.304444778628871926e-02,3.149634503287283144e-02,1.614789286726171502e-01,4.936517420529724226e-02,6.382493918276365905e-03,1.657739816258531340e-03,1.000347576238196437e-02,nan
2.460160000000000000e+05,3.276800000000000000e+04,8.370461695986376738e-02,8.590215287799810029e-02,8.074400021075507639e-02,1.452536059883421070e-02,2.692954524005542241e-12,1.660973298511842985e-02,1.574358018206599777e-02,8.074400021075506251e-02,2.473309487491211905e-02,3.190316357629961134e-03,8.289164733358427061e-04,5.011972895728674347e-03,nan
        }\tableUniformOne

        \pgfplotstableread[col sep=comma]{
ndof,nelem,res,err,res1,res2,res3,res4,err1,err2,err3,relerr1,relerr2,relerr3,cond
6.400000000000000000e+01,8.000000000000000000e+00,5.824570832028575396e+00,5.977244475212094343e+00,5.782112461670299908e+00,6.633255354617407784e-01,1.936287918639931947e-16,2.297827058755127216e-01,1.065049253138570906e+00,5.782112461670299908e+00,1.077170871050412648e+00,2.158447365525732464e-01,5.935338977916125580e-02,2.183013247498441778e-01,nan
2.480000000000000000e+02,3.200000000000000000e+01,2.602371996235329110e+00,2.685125338769466907e+00,2.535070454618537195e+00,4.269606730335535061e-01,9.226501774606658982e-16,4.043048114884999489e-01,5.095882467238084246e-01,2.535070454618536751e+00,7.236267641698758135e-01,1.032641686565602263e-01,2.602498830148856912e-02,1.466374405219338095e-01,nan
9.760000000000000000e+02,1.280000000000000000e+02,1.327125105375400071e+00,1.364790322148153434e+00,1.285757090183150764e+00,2.251127377995592593e-01,8.954220383283036462e-16,2.396121984411694605e-01,2.532835614963078186e-01,1.285757090183151208e+00,3.812201017112841761e-01,5.132598049590576805e-02,1.319955947165485539e-02,7.725134386727978975e-02,nan
3.872000000000000000e+03,5.120000000000000000e+02,6.675928494082256082e-01,6.856738648501705669e-01,6.451905433792364253e-01,1.147501483277259815e-01,1.026132585397630841e-14,1.274432374579533567e-01,1.262431509775466532e-01,6.451905433792367583e-01,1.947831576493091188e-01,2.558221096756683674e-02,6.623514669221431360e-03,3.947132017324028047e-02,nan
1.542400000000000000e+04,2.048000000000000000e+03,3.344676124072866741e-01,3.433546575684368607e-01,3.228852943694280886e-01,5.781886798415141610e-02,9.552309269503031965e-15,6.535020225240328884e-02,6.303238677978459381e-02,3.228852943694282551e-01,9.830781134146873712e-02,1.277303207265892415e-02,3.314734702295096379e-03,1.992132761281063813e-02,nan
6.156800000000000000e+04,8.192000000000000000e+03,1.673581016659741461e-01,1.717684039047924216e-01,1.614789286726119877e-01,2.900610726178922544e-02,6.429880705941083651e-13,3.304444778630142437e-02,3.149634503283511161e-02,1.614789286726120432e-01,4.936517420528349631e-02,6.382493918268721846e-03,1.657739816258478865e-03,1.000347576237918014e-02,nan
2.460160000000000000e+05,3.276800000000000000e+04,8.370461695982803207e-02,8.590215287807560773e-02,8.074400021081869216e-02,1.452536059861641270e-02,3.089692813276274212e-12,1.660973298481951965e-02,1.574358018217631230e-02,8.074400021081869216e-02,2.473309487490336564e-02,3.190316357652315648e-03,8.289164733364960463e-04,5.011972895726900593e-03,nan
        }\tableUniformTwo

        \pgfplotstableread[col sep=comma]{
ndof,nelem,res,err,res1,res2,res3,res4,err1,err2,err3,relerr1,relerr2,relerr3,cond
6.400000000000000000e+01,8.000000000000000000e+00,5.824570832028577172e+00,5.977244475212095232e+00,5.782112461670300796e+00,6.633255354617411115e-01,2.337755210490751476e-16,2.297827058755128327e-01,1.065049253138571128e+00,5.782112461670301684e+00,1.077170871050412648e+00,2.158447365525732464e-01,5.935338977916126274e-02,2.183013247498441223e-01,nan
2.480000000000000000e+02,3.200000000000000000e+01,2.602371996235329554e+00,2.685125338769467351e+00,2.535070454618537195e+00,4.269606730335538947e-01,7.219855308028877787e-16,4.043048114885005040e-01,5.095882467238077584e-01,2.535070454618537195e+00,7.236267641698760356e-01,1.032641686565600736e-01,2.602498830148856565e-02,1.466374405219338095e-01,nan
9.760000000000000000e+02,1.280000000000000000e+02,1.327125105375400294e+00,1.364790322148153656e+00,1.285757090183151208e+00,2.251127377995596479e-01,8.743565343235763321e-16,2.396121984411683226e-01,2.532835614963079296e-01,1.285757090183151652e+00,3.812201017112845092e-01,5.132598049590575418e-02,1.319955947165485886e-02,7.725134386727980362e-02,nan
3.872000000000000000e+03,5.120000000000000000e+02,6.675928494082254971e-01,6.856738648501706779e-01,6.451905433792367583e-01,1.147501483277273276e-01,1.632616760209468172e-14,1.274432374579506921e-01,1.262431509775450433e-01,6.451905433792370914e-01,1.947831576493093131e-01,2.558221096756650714e-02,6.623514669221433962e-03,3.947132017324032210e-02,nan
1.542400000000000000e+04,2.048000000000000000e+03,3.344676124072866186e-01,3.433546575684370272e-01,3.228852943694284772e-01,5.781886798415200590e-02,1.049917809332941811e-14,6.535020225240065206e-02,6.303238677978352522e-02,3.228852943694287547e-01,9.830781134146861222e-02,1.277303207265870384e-02,3.314734702295101584e-03,1.992132761281060344e-02,nan
6.156800000000000000e+04,8.192000000000000000e+03,1.673581016659693166e-01,1.717684039048222033e-01,1.614789286726250050e-01,2.900610726173910581e-02,1.169398264023716218e-12,3.304444778625745260e-02,3.149634503289364812e-02,1.614789286726250606e-01,4.936517420530718570e-02,6.382493918280583885e-03,1.657739816258611571e-03,1.000347576238398012e-02,nan
2.460160000000000000e+05,3.276800000000000000e+04,8.370461695983663630e-02,8.590215287802992206e-02,8.074400021081394596e-02,1.452536059866131601e-02,2.112086940767491040e-12,1.660973298484670624e-02,1.574358018201784185e-02,8.074400021081394596e-02,2.473309487486104880e-02,3.190316357620202013e-03,8.289164733364469319e-04,5.011972895718324987e-03,nan
        }\tableUniformThree

        \pgfplotstableread[col sep=comma]{
ndof,nelem,res,err,res1,res2,res3,res4,err1,err2,err3,relerr1,relerr2,relerr3,cond
6.400000000000000000e+01,8.000000000000000000e+00,5.824570832028576284e+00,5.977244475212094343e+00,5.782112461670300796e+00,6.633255354617408894e-01,2.380401256515800723e-16,2.297827058755128882e-01,1.065049253138570906e+00,5.782112461670299908e+00,1.077170871050412426e+00,2.158447365525731909e-01,5.935338977916124192e-02,2.183013247498440945e-01,nan
2.480000000000000000e+02,3.200000000000000000e+01,2.602371996235329110e+00,2.685125338769466907e+00,2.535070454618536751e+00,4.269606730335538392e-01,6.062680048551588790e-16,4.043048114885005040e-01,5.095882467238075364e-01,2.535070454618536751e+00,7.236267641698759245e-01,1.032641686565600458e-01,2.602498830148856565e-02,1.466374405219338373e-01,nan
9.760000000000000000e+02,1.280000000000000000e+02,1.327125105375400071e+00,1.364790322148153656e+00,1.285757090183151208e+00,2.251127377995595369e-01,8.819553947937163033e-16,2.396121984411684058e-01,2.532835614963079851e-01,1.285757090183151652e+00,3.812201017112844537e-01,5.132598049590580275e-02,1.319955947165485886e-02,7.725134386727984526e-02,nan
3.872000000000000000e+03,5.120000000000000000e+02,6.675928494082250531e-01,6.856738648501704558e-01,6.451905433792363143e-01,1.147501483277272721e-01,6.576230226168422522e-14,1.274432374579499982e-01,1.262431509775448768e-01,6.451905433792368694e-01,1.947831576493092298e-01,2.558221096756646898e-02,6.623514669221431360e-03,3.947132017324029435e-02,nan
1.542400000000000000e+04,2.048000000000000000e+03,3.344676124072861745e-01,3.433546575684364721e-01,3.228852943694279776e-01,5.781886798415242917e-02,3.789910400700843764e-15,6.535020225240033287e-02,6.303238677978402482e-02,3.228852943694279776e-01,9.830781134146887590e-02,1.277303207265880793e-02,3.314734702295093344e-03,1.992132761281066242e-02,nan
6.156800000000000000e+04,8.192000000000000000e+03,1.673581016659745346e-01,1.717684039048100464e-01,1.614789286726127926e-01,2.900610726178009385e-02,2.085346776280678598e-13,3.304444778630749591e-02,3.149634503288635534e-02,1.614789286726127648e-01,4.936517420530970451e-02,6.382493918279105900e-03,1.657739816258485370e-03,1.000347576238449013e-02,nan
2.460160000000000000e+05,3.276800000000000000e+04,8.370461695983397177e-02,8.590215287802405175e-02,8.074400021081548640e-02,1.452536059863524312e-02,2.619566980731968840e-12,1.660973298484865607e-02,1.574358018197992079e-02,8.074400021081536150e-02,2.473309487486015368e-02,3.190316357612517622e-03,8.289164733364614603e-04,5.011972895718143708e-03,nan
        }\tableUniformFour

        %
        %
        \addlegendimage{line1}
        \addlegendentry{\(\ell = \hphantom{0}1\hphantom{{}^2}\)}
        \addlegendimage{line2}
        \addlegendentry{\(\ell = 10\hphantom{{}^2}\)}
        \addlegendimage{line3}
        \addlegendentry{\(\ell = 10^2\)}
        \addlegendimage{line4}
        \addlegendentry{\(\ell = 10^3\)}

        %
        %
        \addplot+ [line1, forget plot] table [x=ndof, y=relerr1] {\tableUniformOne};
        \addplot+ [line2, forget plot] table [x=ndof, y=relerr1] {\tableUniformTwo};
        \addplot+ [line3, forget plot] table [x=ndof, y=relerr1] {\tableUniformThree};
        \addplot+ [line4, forget plot] table [x=ndof, y=relerr1] {\tableUniformFour};

        %
        %
        \drawslopetriangle[ST1]{0.5}{7e3}{4e-3}
    \end{loglogaxis}
\end{tikzpicture}

%% file: figures/plot_dls_square_error_scalingFriedrichs.tex
\begin{tikzpicture}[>=stealth]
    %
    %
    \colorlet{col1}{TUblue}
    \colorlet{col2}{TUgreen}
    \colorlet{col3}{TUmagenta}
    \colorlet{col4}{TUyellow}
    \colorlet{col5}{purple}
    \colorlet{col6}{green}
    \pgfplotsset{%
        linedefault/.style = {%
            mark = *,%
            mark size = 2pt,%
            every mark/.append style = {solid},%
            gray,%
            every mark/.append style = {fill = gray!60!white}%
        },%
        line1/.style = {%
            linedefault,%
            col1,%
            every mark/.append style = {fill = col1!60!white}%
        },%
        line2/.style = {%
            linedefault,%
            mark = triangle*,%
            mark size = 2.75pt,%
            col2,%
            every mark/.append style = {fill = col2!60!white}%
        },%
        line3/.style = {%
            linedefault,%
            mark = square*,%
            mark size = 1.66pt,%
            col3,%
            every mark/.append style = {fill = col3!60!white}%
        },%
        line4/.style = {%
            linedefault,%
            mark = pentagon*,%
            mark size = 2.2pt,%
            col4,%
            every mark/.append style = {fill = col4!60!white}%
        },%
        line5/.style = {%
            linedefault,%
            mark = diamond*,%
            mark size = 2.75pt,%
            col5,%
            every mark/.append style = {fill = col5!60!white}%
        },%
        line6/.style = {%
            linedefault,%
            mark = halfsquare*,%
            mark size = 1.66pt,%
            col6,%
            every mark/.append style = {fill = col6!60!white}%
        },%
        minorline/.style = {%
            dashed,%
            every mark/.append style = {fill = black!20!white}%
        },%
        majorline/.style = {%
            solid%
        }%
    }

    %
    %
    \begin{loglogaxis}[%
            width            = 0.36\textwidth,%
            xlabel           = ndof,%
            ylabel           = {relative energy error},%
            ymajorgrids      = true,%
            font             = \footnotesize,%
            grid style       = {densely dotted, semithick},%
            legend style     = {legend pos  = south west}%
        ]

        %
        %
        \pgfplotstableread[col sep=comma]{
ndof,nelem,res,err,res1,res2,res3,res4,err1,err2,err3,relerr1,relerr2,relerr3,cond
6.400000000000000000e+01,8.000000000000000000e+00,1.964012006227397134e+00,2.435195810553486773e+00,1.851157089722178073e+00,6.200234656969803027e-01,2.557749207902330937e-16,2.147824288944214510e-01,1.106583467230553897e+00,1.851157089722178073e+00,1.130870945316297282e+00,2.242621327172892309e-01,1.875434983974247627e-01,2.291842753257133936e-01,nan
2.480000000000000000e+02,3.200000000000000000e+01,9.931180565166085605e-01,1.206606461999659707e+00,8.122934449900176546e-01,4.152531943598093966e-01,6.693400313719406793e-16,3.924635243937316109e-01,5.168207882201721892e-01,8.122934449900176546e-01,7.273065283833670680e-01,1.047297879868680365e-01,8.230253331130109318e-02,1.473831166532051096e-01,nan
9.760000000000000000e+02,1.280000000000000000e+02,5.240209889521627584e-01,6.153812373155924664e-01,4.101301086247215455e-01,2.233862199479774713e-01,7.575587710238683397e-16,2.376760139333177757e-01,2.545278906518143591e-01,4.101301086247216010e-01,3.817105108847301809e-01,5.157813430165894553e-02,4.155486805321752936e-02,7.735072154312500059e-02,nan
3.872000000000000000e+03,5.120000000000000000e+02,2.674215403715322070e-01,3.101135001170129968e-01,2.054887345699082102e-01,1.145204463430959940e-01,4.416572415349288141e-15,1.271759709572485020e-01,1.264165169110532116e-01,2.054887345699081547e-01,1.948425702246812385e-01,2.561734225073912172e-02,2.082036181178931536e-02,3.948335967816016256e-02,nan
1.542400000000000000e+04,2.048000000000000000e+03,1.348039265099734152e-01,1.555896879551097456e-01,1.027929339957385085e-01,5.778939969788898096e-02,2.032818477661997758e-15,6.531536270749545292e-02,6.305495843224699637e-02,1.027929339957384530e-01,9.831498327123786019e-02,1.277760604574743400e-02,1.041510174251695370e-02,1.992278094963705060e-02,nan
6.156800000000000000e+04,8.192000000000000000e+03,6.763856749378296951e-02,7.791920791811607350e-02,5.140229119474090669e-02,2.900238064162486556e-02,1.108088471549749617e-14,3.304000892343096629e-02,3.149921425633222388e-02,5.140229119474098302e-02,4.936605031548839395e-02,6.383075344510445552e-03,5.208140985779264633e-03,1.000365329942510227e-02,nan
2.460160000000000000e+05,3.276800000000000000e+04,3.387362840959455557e-02,3.898959606869041583e-02,2.570185941855715472e-02,1.452489221878952227e-02,5.332871477539866470e-13,1.660917307070294141e-02,1.574394152504849514e-02,2.570185941855711656e-02,2.473320297260788933e-02,3.190389581027336125e-03,2.604142818097149748e-03,5.011994800901643267e-03,nan
        }\tableUniformOne

        \pgfplotstableread[col sep=comma]{
ndof,nelem,res,err,res1,res2,res3,res4,err1,err2,err3,relerr1,relerr2,relerr3,cond
6.400000000000000000e+01,8.000000000000000000e+00,1.964012006227397134e+00,2.435195810553487217e+00,1.851157089722178073e+00,6.200234656969798586e-01,3.052454996149280033e-16,2.147824288944212290e-01,1.106583467230554341e+00,1.851157089722178073e+00,1.130870945316297727e+00,2.242621327172893142e-01,1.875434983974247627e-01,2.291842753257134768e-01,nan
2.480000000000000000e+02,3.200000000000000000e+01,9.931180565166085605e-01,1.206606461999659707e+00,8.122934449900177656e-01,4.152531943598094522e-01,7.223394904411071674e-16,3.924635243937314444e-01,5.168207882201723002e-01,8.122934449900178766e-01,7.273065283833668460e-01,1.047297879868680504e-01,8.230253331130110706e-02,1.473831166532050541e-01,nan
9.760000000000000000e+02,1.280000000000000000e+02,5.240209889521628694e-01,6.153812373155922444e-01,4.101301086247212679e-01,2.233862199479778043e-01,9.336874404051448201e-16,2.376760139333180810e-01,2.545278906518141371e-01,4.101301086247213235e-01,3.817105108847300698e-01,5.157813430165891083e-02,4.155486805321749466e-02,7.735072154312498671e-02,nan
3.872000000000000000e+03,5.120000000000000000e+02,2.674215403715322070e-01,3.101135001170141625e-01,2.054887345699090151e-01,1.145204463430948422e-01,2.120354816709280293e-15,1.271759709572482244e-01,1.264165169110544329e-01,2.054887345699090428e-01,1.948425702246813773e-01,2.561734225073936805e-02,2.082036181178940557e-02,3.948335967816019726e-02,nan
1.542400000000000000e+04,2.048000000000000000e+03,1.348039265099734707e-01,1.555896879551078582e-01,1.027929339957368987e-01,5.778939969789177733e-02,2.051378797953550485e-15,6.531536270749564721e-02,6.305495843224520613e-02,1.027929339957369959e-01,9.831498327123754100e-02,1.277760604574706971e-02,1.041510174251680625e-02,1.992278094963698815e-02,nan
6.156800000000000000e+04,8.192000000000000000e+03,6.763856749378296951e-02,7.791920791811458857e-02,5.140229119473969932e-02,2.900238064162669743e-02,1.649608099377380905e-14,3.304000892343126466e-02,3.149921425633116223e-02,5.140229119473980340e-02,4.936605031548795680e-02,6.383075344510230446e-03,5.208140985779146671e-03,1.000365329942501380e-02,nan
2.460160000000000000e+05,3.276800000000000000e+04,3.387362840959397270e-02,3.898959606869305261e-02,2.570185941856031886e-02,1.452489221878787082e-02,5.591217233826742225e-13,1.660917307069830276e-02,1.574394152504827657e-02,2.570185941856036743e-02,2.473320297260880180e-02,3.190389581027291890e-03,2.604142818097478478e-03,5.011994800901828015e-03,nan
        }\tableUniformTwo

        \pgfplotstableread[col sep=comma]{
ndof,nelem,res,err,res1,res2,res3,res4,err1,err2,err3,relerr1,relerr2,relerr3,cond
6.400000000000000000e+01,8.000000000000000000e+00,1.964012006227397356e+00,2.435195810553487661e+00,1.851157089722178295e+00,6.200234656969801916e-01,1.968823092040905334e-16,2.147824288944212290e-01,1.106583467230554563e+00,1.851157089722178295e+00,1.130870945316297949e+00,2.242621327172893420e-01,1.875434983974247627e-01,2.291842753257134768e-01,nan
2.480000000000000000e+02,3.200000000000000000e+01,9.931180565166086716e-01,1.206606461999660151e+00,8.122934449900178766e-01,4.152531943598096187e-01,8.108935025185472134e-16,3.924635243937313334e-01,5.168207882201724113e-01,8.122934449900178766e-01,7.273065283833672900e-01,1.047297879868680642e-01,8.230253331130109318e-02,1.473831166532051373e-01,nan
9.760000000000000000e+02,1.280000000000000000e+02,5.240209889521628694e-01,6.153812373155927995e-01,4.101301086247214900e-01,2.233862199479781652e-01,7.996369159770887145e-16,2.376760139333174982e-01,2.545278906518143036e-01,4.101301086247215455e-01,3.817105108847306250e-01,5.157813430165890389e-02,4.155486805321751548e-02,7.735072154312505610e-02,nan
3.872000000000000000e+03,5.120000000000000000e+02,2.674215403715317629e-01,3.101135001170156613e-01,2.054887345699096812e-01,1.145204463430955083e-01,1.456377764846994925e-15,1.271759709572456432e-01,1.264165169110547937e-01,2.054887345699095980e-01,1.948425702246829594e-01,2.561734225073944091e-02,2.082036181178945414e-02,3.948335967816051645e-02,nan
1.542400000000000000e+04,2.048000000000000000e+03,1.348039265099723605e-01,1.555896879551120493e-01,1.027929339957392302e-01,5.778939969789066017e-02,2.099117773099914603e-15,6.531536270749067896e-02,6.305495843224787067e-02,1.027929339957391608e-01,9.831498327124020553e-02,1.277760604574760574e-02,1.041510174251702309e-02,1.992278094963751897e-02,nan
6.156800000000000000e+04,8.192000000000000000e+03,6.763856749377936128e-02,7.791920791812696756e-02,5.140229119474688801e-02,2.900238064162297472e-02,3.750692608499959903e-14,3.304000892341594359e-02,3.149921425633932237e-02,5.140229119474697822e-02,4.936605031549482631e-02,6.383075344511884505e-03,5.208140985779872653e-03,1.000365329942640505e-02,nan
2.460160000000000000e+05,3.276800000000000000e+04,3.387362840959225879e-02,3.898959606868622474e-02,2.570185941855310241e-02,1.452489221880214586e-02,4.523125551333918796e-13,1.660917307069349064e-02,1.574394152504163952e-02,2.570185941855314751e-02,2.473320297260977324e-02,3.190389581025946178e-03,2.604142818096747292e-03,5.011994800902024039e-03,nan
        }\tableUniformThree

        \pgfplotstableread[col sep=comma]{
ndof,nelem,res,err,res1,res2,res3,res4,err1,err2,err3,relerr1,relerr2,relerr3,cond
6.400000000000000000e+01,8.000000000000000000e+00,1.964012006227397134e+00,2.435195810553487217e+00,1.851157089722178073e+00,6.200234656969803027e-01,3.070507784663998085e-16,2.147824288944213955e-01,1.106583467230553897e+00,1.851157089722178295e+00,1.130870945316297504e+00,2.242621327172892032e-01,1.875434983974247627e-01,2.291842753257133936e-01,nan
2.480000000000000000e+02,3.200000000000000000e+01,9.931180565166086716e-01,1.206606461999659929e+00,8.122934449900178766e-01,4.152531943598095632e-01,6.935911916584082573e-16,3.924635243937312779e-01,5.168207882201721892e-01,8.122934449900177656e-01,7.273065283833671790e-01,1.047297879868680365e-01,8.230253331130107930e-02,1.473831166532051373e-01,nan
9.760000000000000000e+02,1.280000000000000000e+02,5.240209889521627584e-01,6.153812373155934656e-01,4.101301086247220451e-01,2.233862199479775268e-01,8.616459835404618220e-16,2.376760139333169153e-01,2.545278906518148587e-01,4.101301086247222116e-01,3.817105108847307915e-01,5.157813430165905655e-02,4.155486805321758487e-02,7.735072154312513937e-02,nan
3.872000000000000000e+03,5.120000000000000000e+02,2.674215403715318184e-01,3.101135001170146066e-01,2.054887345699089318e-01,1.145204463430961467e-01,2.455061853408083980e-15,1.271759709572463926e-01,1.264165169110540443e-01,2.054887345699089318e-01,1.948425702246824598e-01,2.561734225073928478e-02,2.082036181178939169e-02,3.948335967816040543e-02,nan
1.542400000000000000e+04,2.048000000000000000e+03,1.348039265099720274e-01,1.555896879551131318e-01,1.027929339957401045e-01,5.778939969788945280e-02,2.103240203826982961e-15,6.531536270748966588e-02,6.305495843224900865e-02,1.027929339957400767e-01,9.831498327124024716e-02,1.277760604574783819e-02,1.041510174251711850e-02,1.992278094963753285e-02,nan
6.156800000000000000e+04,8.192000000000000000e+03,6.763856749377977762e-02,7.791920791812606550e-02,5.140229119474651331e-02,2.900238064162268328e-02,7.598830985201369646e-14,3.304000892341764362e-02,3.149921425633885746e-02,5.140229119474653413e-02,4.936605031549418099e-02,6.383075344511789963e-03,5.208140985779827550e-03,1.000365329942627494e-02,nan
2.460160000000000000e+05,3.276800000000000000e+04,3.387362840957731241e-02,3.898959606876716694e-02,2.570185941860862397e-02,1.452489221874287903e-02,4.657866869660552010e-13,1.660917307062891729e-02,1.574394152510336445e-02,2.570185941860855111e-02,2.473320297264049866e-02,3.190389581038454402e-03,2.604142818102360857e-03,5.011994800908250829e-03,nan
        }\tableUniformFour

        %
        %
        \addlegendimage{line1}
        \addlegendentry{\(\ell = \hphantom{0}1\hphantom{{}^2}\)}
        \addlegendimage{line2}
        \addlegendentry{\(\ell = 10\hphantom{{}^2}\)}
        \addlegendimage{line3}
        \addlegendentry{\(\ell = 10^2\)}
        \addlegendimage{line4}
        \addlegendentry{\(\ell = 10^3\)}

        %
        %
        \addplot+ [line1, forget plot] table [x=ndof, y=relerr1] {\tableUniformOne};
        \addplot+ [line2, forget plot] table [x=ndof, y=relerr1] {\tableUniformTwo};
        \addplot+ [line3, forget plot] table [x=ndof, y=relerr1] {\tableUniformThree};
        \addplot+ [line4, forget plot] table [x=ndof, y=relerr1] {\tableUniformFour};

        %
        %
        \drawslopetriangle[ST1]{0.5}{7e3}{4e-3}
    \end{loglogaxis}
\end{tikzpicture}

%% file: figures/plot_dls_square_error_weighting_k0.tex
\begin{tikzpicture}[>=stealth]
    %
    %
    \colorlet{col0}{TUblue}
    \colorlet{col1}{TUgreen}
    \colorlet{col2}{TUmagenta}
    \colorlet{col3}{TUyellow}
    \colorlet{col4}{purple}
    \colorlet{col5}{green}
    \pgfplotsset{%
        linedefault/.style = {%
            mark = *,%
            mark size = 2pt,%
            every mark/.append style = {solid},%
            gray,%
            every mark/.append style = {fill = gray!60!white}%
        },%
        line1/.style = {%
            linedefault,%
            col0,%
            every mark/.append style = {fill = col0!60!white}%
        },%
        line2/.style = {%
            linedefault,%
            mark = triangle*,%
            mark size = 2.75pt,%
            col1,%
            every mark/.append style = {fill = col1!60!white}%
        },%
        line3/.style = {%
            linedefault,%
            mark = square*,%
            mark size = 1.66pt,%
            col2,%
            every mark/.append style = {fill = col2!60!white}%
        },%
        line4/.style = {%
            linedefault,%
            mark = pentagon*,%
            mark size = 2.2pt,%
            col3,%
            every mark/.append style = {fill = col3!60!white}%
        },%
        line5/.style = {%
            linedefault,%
            mark = diamond*,%
            mark size = 2.75pt,%
            col4,%
            every mark/.append style = {fill = col4!60!white}%
        },%
        line6/.style = {%
            linedefault,%
            mark = halfsquare*,%
            mark size = 1.66pt,%
            col5,%
            every mark/.append style = {fill = col5!60!white}%
        },%
        minorline/.style = {%
            dashed,%
            every mark/.append style = {fill = black!20!white}%
        },%
        majorline/.style = {%
            solid%
        }%
    }

    %
    %
    \begin{loglogaxis}[%
            width            = 0.36\textwidth,%
            xlabel           = ndof,%
            ylabel           = {relative energy error},%
            ymin             = 1e-3,%
            ymajorgrids      = true,%
            font             = \footnotesize,%
            grid style       = {densely dotted, semithick},%
            legend style     = {legend pos  = south west}%
        ]

        %
        %
        \pgfplotstableread[col sep=comma]{
ndof,nelem,res,err,res1,res2,res3,res4,err1,err2,err3,relerr1,relerr2,relerr3,cond
6.400000000000000000e+01,8.000000000000000000e+00,3.140426790102585769e-02,3.139092834487416273e+00,3.139130679113536632e-02,8.524650936622455528e-04,2.422140795617594229e-19,2.953025707803935435e-04,2.219487838960685799e+00,3.139130679113536632e-02,2.219637840300487053e+00,4.498052709490791834e-01,3.180300325435872111e+01,4.498356704818486929e-01,nan
2.480000000000000000e+02,3.200000000000000000e+01,3.137514448725216981e-02,3.132103790075521754e+00,3.131717726402560159e-02,1.465354987099717595e-03,1.957077733888844758e-18,1.219358262091326103e-03,2.214457864438242662e+00,3.131717726402560159e-02,2.214784358090258198e+00,4.487429839019342337e-01,3.173093468727996935e+01,4.488091453484820970e-01,nan
9.760000000000000000e+02,1.280000000000000000e+02,3.127553832771842041e-02,3.104128557182525228e+00,3.103758744638523623e-02,2.809369846531066962e-03,6.980954078857991383e-18,2.633420639012266310e-03,2.194673994611542334e+00,3.103758744638523970e-02,2.195007249431626395e+00,4.447339336862004000e-01,3.144765097470330772e+01,4.448014652318512674e-01,nan
3.872000000000000000e+03,5.120000000000000000e+02,3.090608750763604654e-02,2.999296236118383252e+00,2.999207723969524522e-02,5.314510440327749677e-03,3.631709143958952127e-17,5.236174811583208313e-03,2.120551433361919091e+00,2.999207723969525216e-02,2.120881893402708940e+00,4.297135624053726177e-01,3.038832765818356307e+01,4.297805276103312533e-01,nan
1.542400000000000000e+04,2.048000000000000000e+03,2.958868685119422567e-02,2.648921288223409665e+00,2.651313544678659748e-02,9.148117283962975410e-03,2.127641654965740308e-16,9.426345448647647715e-03,1.872816492607403793e+00,2.651313544678660442e-02,1.873136255019813179e+00,3.795119675517499847e-01,2.686342265538060659e+01,3.795767649639692487e-01,nan
6.156800000000000000e+04,8.192000000000000000e+03,2.566328475354108865e-02,1.820598600865565242e+00,1.829906374482356129e-02,1.217831989476989875e-02,1.096813678284408472e-15,1.324526235847633786e-02,1.287145394244711305e+00,1.829906374482354742e-02,1.287439763178978014e+00,2.608301897293187066e-01,1.854082798172110813e+01,2.608898413470135402e-01,nan
2.460160000000000000e+05,3.276800000000000000e+04,1.813960372440420049e-02,8.167980413513724613e-01,8.256979154483560721e-03,1.072237516816562999e-02,5.006343819822642096e-15,1.207883228105878677e-02,5.774035142958292566e-01,8.256979154483541639e-03,5.776643007227526549e-01,1.170064150153977550e-01,8.366069012423810847e+00,1.170592613943052251e-01,nan
        }\tableUniformOne

        \pgfplotstableread[col sep=comma]{
ndof,nelem,res,err,res1,res2,res3,res4,err1,err2,err3,relerr1,relerr2,relerr3,cond
6.400000000000000000e+01,8.000000000000000000e+00,3.022657632911757708e-01,2.935315339791428002e+00,2.914580206370028237e-01,7.569212702240828150e-02,3.214715071135365497e-17,2.622052194715360565e-02,2.058925511127845276e+00,2.914580206370028237e-01,2.071703189611911800e+00,4.172654299518494536e-01,2.952804877000247963e+00,4.198549765272959755e-01,nan
2.480000000000000000e+02,3.200000000000000000e+01,2.809043367425435034e-01,2.448142612184258216e+00,2.421490815761406512e-01,1.072926491818104872e-01,1.393929808323488978e-16,9.359143189495920667e-02,1.708287645836542801e+00,2.421490815761406235e-01,1.736812998271960895e+00,3.461714525554837052e-01,2.453483156319967495e+00,3.519518975006265538e-01,nan
9.760000000000000000e+02,1.280000000000000000e+02,2.300210422155115619e-01,1.497320493975264677e+00,1.484982218281948596e-01,1.224777952597883302e-01,3.302319274925034468e-16,1.259251668415640735e-01,1.040369661384129563e+00,1.484982218281948874e-01,1.066558909513987352e+00,2.108229710340373542e-01,1.504601560441521713e+00,2.161300222724775177e-01,nan
3.872000000000000000e+03,5.120000000000000000e+02,1.522481228981040491e-01,6.134365203659553067e-01,6.063770737366683261e-02,9.395971958401669322e-02,8.075759499924403848e-16,1.033156835439552057e-01,4.200717986384981395e-01,6.063770737366683261e-02,4.429075715674475888e-01,8.512434370563709551e-02,6.143884284457660216e-01,8.975183879073472815e-02,nan
1.542400000000000000e+04,2.048000000000000000e+03,8.480008957331372499e-02,2.044971771948399242e-01,1.976081141657253654e-02,5.472316743899156150e-02,2.045971218409882197e-15,6.169230478142943841e-02,1.339727321342580046e-01,1.976081141657254348e-02,1.532315661255825301e-01,2.714855159128195042e-02,2.002188802460157790e-01,3.105120730242655683e-02,nan
6.156800000000000000e+04,8.192000000000000000e+03,4.386326441306044599e-02,7.344249793261388437e-02,6.790616114479536281e-03,2.859783727783100221e-02,7.093663075258849341e-15,3.255821249372461285e-02,4.423388075155389870e-02,6.790616114479532811e-03,5.823273846407254345e-02,8.963658309683042447e-03,6.880332623798023073e-02,1.180041997582919472e-02,nan
2.460160000000000000e+05,3.276800000000000000e+04,2.216263584411782769e-02,3.147906424089075977e-02,2.805980988673616519e-03,1.447348564612275373e-02,2.663701619764608460e-14,1.654772224512053927e-02,1.760486813441195994e-02,2.805980988673623024e-03,2.594468296478675412e-02,3.567492154491786784e-03,2.843053150503099466e-02,5.257491974434778947e-03,nan
        }\tableUniformTwo

        \pgfplotstableread[col sep=comma]{
ndof,nelem,res,err,res1,res2,res3,res4,err1,err2,err3,relerr1,relerr2,relerr3,cond
6.400000000000000000e+01,8.000000000000000000e+00,1.964012006227397356e+00,2.435195810553487217e+00,1.851157089722178295e+00,6.200234656969800806e-01,2.926329773543867200e-16,2.147824288944212012e-01,1.106583467230554341e+00,1.851157089722178295e+00,1.130870945316297727e+00,2.242621327172893142e-01,1.875434983974247904e-01,2.291842753257134768e-01,nan
2.480000000000000000e+02,3.200000000000000000e+01,9.931180565166084495e-01,1.206606461999660151e+00,8.122934449900178766e-01,4.152531943598094522e-01,4.948046936343726307e-16,3.924635243937309448e-01,5.168207882201724113e-01,8.122934449900178766e-01,7.273065283833672900e-01,1.047297879868680781e-01,8.230253331130110706e-02,1.473831166532051651e-01,nan
9.760000000000000000e+02,1.280000000000000000e+02,5.240209889521626474e-01,6.153812373155927995e-01,4.101301086247215455e-01,2.233862199479778876e-01,9.141287406841202494e-16,2.376760139333170818e-01,2.545278906518143591e-01,4.101301086247216010e-01,3.817105108847305694e-01,5.157813430165895247e-02,4.155486805321752936e-02,7.735072154312509773e-02,nan
3.872000000000000000e+03,5.120000000000000000e+02,2.674215403715315964e-01,3.101135001170151062e-01,2.054887345699092926e-01,1.145204463430955916e-01,5.731112490467163133e-15,1.271759709572457819e-01,1.264165169110547104e-01,2.054887345699091539e-01,1.948425702246826263e-01,2.561734225073943050e-02,2.082036181178940903e-02,3.948335967816045400e-02,nan
1.542400000000000000e+04,2.048000000000000000e+03,1.348039265099724715e-01,1.555896879551110779e-01,1.027929339957386196e-01,5.778939969789162467e-02,7.586730460492429892e-15,6.531536270749102591e-02,6.305495843224703800e-02,1.027929339957385502e-01,9.831498327123984471e-02,1.277760604574744094e-02,1.041510174251696584e-02,1.992278094963745305e-02,nan
6.156800000000000000e+04,8.192000000000000000e+03,6.763856749378012456e-02,7.791920791812274871e-02,5.140229119474341163e-02,2.900238064162711377e-02,2.230844509753342579e-13,3.304000892341928813e-02,3.149921425633601252e-02,5.140229119474348796e-02,4.936605031549391731e-02,6.383075344511213167e-03,5.208140985779516167e-03,1.000365329942622117e-02,nan
2.460160000000000000e+05,3.276800000000000000e+04,3.387362840956339993e-02,3.898959606882903411e-02,2.570185941864954610e-02,1.452489221870304978e-02,3.449105056253131216e-13,1.660917307057204265e-02,1.574394152515188119e-02,2.570185941864952875e-02,2.473320297266456622e-02,3.190389581048280309e-03,2.604142818106514219e-03,5.011994800913118463e-03,nan
        }\tableUniformThree

        \pgfplotstableread[col sep=comma]{
ndof,nelem,res,err,res1,res2,res3,res4,err1,err2,err3,relerr1,relerr2,relerr3,cond
6.400000000000000000e+01,8.000000000000000000e+00,1.841576170976315652e+01,1.846387285679245949e+01,1.840218488719548873e+01,6.680667258782667650e-01,3.379987689399953195e-16,2.314251024134690304e-01,1.061034164300986848e+00,1.840218488719548873e+01,1.071633683427948291e+00,2.150310316559878032e-01,1.864352923402529463e-02,2.171791486625966050e-01,nan
2.480000000000000000e+02,3.200000000000000000e+01,8.090267747992456293e+00,8.117102343278704168e+00,8.068743183709022659e+00,4.281892118718177431e-01,7.185075617526670592e-16,4.055486017639908236e-01,5.092044873873511346e-01,8.068743183709020883e+00,7.234947661416836207e-01,1.031864027561529285e-01,8.175346098794463931e-03,1.466106921290240284e-01,nan
9.760000000000000000e+02,1.280000000000000000e+02,4.105805766962500769e+00,4.118108399431804401e+00,4.092600193576081402e+00,2.252897397005093061e-01,8.358242299853964777e-16,2.398107269269356134e-01,2.532174716681399040e-01,4.092600193576081402e+00,3.812103842418186339e-01,5.131258789667238196e-02,4.146670957879906340e-03,7.724937469626327879e-02,nan
3.872000000000000000e+03,5.120000000000000000e+02,2.060843698474333241e+00,2.066768095020991680e+00,2.053693030919703677e+00,1.147735435215623961e-01,6.433828775484165171e-13,1.274704593257322549e-01,1.262340096362830344e-01,2.053693030919705009e+00,1.947826177041495355e-01,2.558035854453302380e-02,2.080826087308056647e-03,3.947121075747494351e-02,nan
1.542400000000000000e+04,2.048000000000000000e+03,1.031471935056957889e+00,1.034387368128164786e+00,1.027774232451958314e+00,5.782186421127862713e-02,2.567613100099886908e-15,6.535374462925140626e-02,6.303120087821882656e-02,1.027774232451958758e+00,9.830779301127844094e-02,1.277279175875827184e-02,1.041353017491282148e-03,1.992132389833748199e-02,nan
6.156800000000000000e+04,8.192000000000000000e+03,5.158804362250402065e-01,5.173279556503191934e-01,5.140031934772832933e-01,2.900648600588436737e-02,6.187913215746760111e-12,3.304489891772533633e-02,3.149619454312154593e-02,5.140031934772828492e-01,4.936517519259252945e-02,6.382463422676913113e-03,5.207941195906017322e-04,1.000347596244981901e-02,nan
2.460160000000000000e+05,3.276800000000000000e+04,2.579615346509456231e-01,2.586829458616488053e-01,2.570161099100783386e-01,1.452540823190759216e-02,3.128061437136879751e-11,1.660978994021222038e-02,1.574356137372161235e-02,2.570161099100771729e-01,2.473309526977430714e-02,3.190312546262523766e-03,2.604117647123820110e-04,5.011972975744473022e-03,nan
        }\tableUniformFour

        \pgfplotstableread[col sep=comma]{
ndof,nelem,res,err,res1,res2,res3,res4,err1,err2,err3,relerr1,relerr2,relerr3,cond
6.400000000000000000e+01,8.000000000000000000e+00,1.840211525124394711e+02,1.840259653257881212e+02,1.840197922119121188e+02,6.685847868008387707e-01,2.138043153680955691e-16,2.316045639813310575e-01,1.060602074636052050e+00,1.840197922119121188e+02,1.071032965006385895e+00,2.149434636119559594e-01,1.864332087071481547e-03,2.170574060210593015e-01,nan
2.480000000000000000e+02,3.200000000000000000e+01,8.068946885226203847e+01,8.069216209161049846e+01,8.068731216104573889e+01,4.283228413096957787e-01,7.395946695250566293e-16,4.056839037225830613e-01,5.091670866609260271e-01,8.068731216104573889e+01,7.234833256978820959e-01,1.031788237844119688e-01,8.175333973075981243e-04,1.466083738104093093e-01,nan
9.760000000000000000e+02,1.280000000000000000e+02,4.092731384803398953e+01,4.092854959416605709e+01,4.092599094923934899e+01,2.253089456534433377e-01,8.803399420097704106e-16,2.398322690182585482e-01,2.532109905630927638e-01,4.092599094923933478e+01,3.812097836477758017e-01,5.131127455020301209e-02,4.146669844712552256e-04,7.724925299046465910e-02,nan
3.872000000000000000e+03,5.120000000000000000e+02,2.053764537210790309e+01,2.053824066173526575e+01,2.053692904042138423e+01,1.147760803562660886e-01,1.160050659280821413e-15,1.274734111099221634e-01,1.262331131140574936e-01,2.053692904042138423e+01,1.947826203903541276e-01,2.558017687114691463e-02,2.080825958754200285e-04,3.947121130181380222e-02,nan
1.542400000000000000e+04,2.048000000000000000e+03,1.027811260008056493e+01,1.027840558614067490e+01,1.027774216704494670e+01,5.782218904697869477e-02,2.862193283056301823e-15,6.535412867562880035e-02,6.303108461604632617e-02,1.027774216704493782e+01,9.830779890921921560e-02,1.277276819911640944e-02,1.041353001535762473e-04,1.992132509351016503e-02,nan
6.156800000000000000e+04,8.192000000000000000e+03,5.140219985294321781e+00,5.140365469665326259e+00,5.140031914979858563e+00,2.900693720274111556e-02,3.617370185260173932e-10,3.304554132957532309e-02,3.149751604577658709e-02,5.140031914979858563e+00,4.936578976432819121e-02,6.382731215103303220e-03,5.207941175851512255e-05,1.000360050072120026e-02,nan
2.460160000000000000e+05,3.276800000000000000e+04,2.570277446607707184e+00,2.570380373940063468e+00,2.570161096618854835e+00,1.582205802303749895e-02,3.048741740329522965e-09,1.864812181197303756e-02,2.003357718172786919e-02,2.570161096618861496e+00,2.694178814586630985e-02,4.059651505285591204e-03,2.604117644609134331e-05,5.459547728760554865e-03,nan
        }\tableUniformFive

        %
        %
        \addlegendimage{line1}
        \addlegendentry{\(c_\Omega = 10^{-4}\)}
        \addlegendimage{line2}
        \addlegendentry{\(c_\Omega = 10^{-2}\)}
        \addlegendimage{line3}
        \addlegendentry{\(c_\Omega = \hphantom{0}1\hphantom{{}^{-2}}\)}
        \addlegendimage{line4}
        \addlegendentry{\(c_\Omega = 10^2\hphantom{{}^{-}}\)}
        \addlegendimage{line5}
        \addlegendentry{\(c_\Omega = 10^4\hphantom{{}^{-}}\)}

        %
        %
        \addplot+ [line1, forget plot] table [x=ndof, y=relerr1] {\tableUniformOne};
        \addplot+ [line2, forget plot] table [x=ndof, y=relerr1] {\tableUniformTwo};
        \addplot+ [line3, forget plot] table [x=ndof, y=relerr1] {\tableUniformThree};
        \addplot+ [line4, forget plot] table [x=ndof, y=relerr1] {\tableUniformFour};
        \addplot+ [line5, forget plot] table [x=ndof, y=relerr1] {\tableUniformFive};

        %
        %
        \drawslopetriangle[ST1]{0.5}{9e3}{4e-3}
    \end{loglogaxis}
\end{tikzpicture}

%% file: figures/plot_dls_square_error_weighting_k1.tex
\begin{tikzpicture}[>=stealth]
    %
    %
    \colorlet{col0}{TUblue}
    \colorlet{col1}{TUgreen}
    \colorlet{col2}{TUmagenta}
    \colorlet{col3}{TUyellow}
    \colorlet{col4}{purple}
    \colorlet{col5}{green}
    \pgfplotsset{%
        linedefault/.style = {%
            mark = *,%
            mark size = 2pt,%
            every mark/.append style = {solid},%
            gray,%
            every mark/.append style = {fill = gray!60!white}%
        },%
        line1/.style = {%
            linedefault,%
            col0,%
            every mark/.append style = {fill = col0!60!white}%
        },%
        line2/.style = {%
            linedefault,%
            mark = triangle*,%
            mark size = 2.75pt,%
            col1,%
            every mark/.append style = {fill = col1!60!white}%
        },%
        line3/.style = {%
            linedefault,%
            mark = square*,%
            mark size = 1.66pt,%
            col2,%
            every mark/.append style = {fill = col2!60!white}%
        },%
        line4/.style = {%
            linedefault,%
            mark = pentagon*,%
            mark size = 2.2pt,%
            col3,%
            every mark/.append style = {fill = col3!60!white}%
        },%
        line5/.style = {%
            linedefault,%
            mark = diamond*,%
            mark size = 2.75pt,%
            col4,%
            every mark/.append style = {fill = col4!60!white}%
        },%
        line6/.style = {%
            linedefault,%
            mark = halfsquare*,%
            mark size = 1.66pt,%
            col5,%
            every mark/.append style = {fill = col5!60!white}%
        },%
        minorline/.style = {%
            dashed,%
            every mark/.append style = {fill = black!20!white}%
        },%
        majorline/.style = {%
            solid%
        }%
    }

    %
    %
    \begin{loglogaxis}[%
            width            = 0.36\textwidth,%
            xlabel           = ndof,%
            ylabel           = {relative energy error},%
            ymin             = 3e-6,%
            ymajorgrids      = true,%
            font             = \footnotesize,%
            grid style       = {densely dotted, semithick},%
            legend style     = {legend pos  = south west}%
        ]

        %
        %
        \pgfplotstableread[col sep=comma]{
ndof,nelem,res,err,res1,res2,res3,res4,err1,err2,err3,relerr1,relerr2,relerr3,cond
1.280000000000000000e+02,8.000000000000000000e+00,3.127444551713214060e-02,3.112228605801071524e+00,3.112118408478225379e-02,2.448126492259920014e-03,8.491025771854680566e-04,1.687763329520671710e-03,2.200409233894064620e+00,3.112118408478225379e-02,2.200726600480254369e+00,4.459387676185898597e-01,3.152933789323236624e+01,4.460030856836792457e-01,nan
5.040000000000000000e+02,3.200000000000000000e+01,1.349059871075149013e-02,6.495947810757677354e-01,1.324367321342360265e-02,2.032324401958702717e-03,6.342109862485657400e-04,1.438334919516565851e-03,4.590957985493445137e-01,1.324367321342360265e-02,4.593789639012982873e-01,9.303225942848802887e-02,1.341864645756408336e+01,9.308964072574094084e-02,nan
2.000000000000000000e+03,1.280000000000000000e+02,6.731320080567550733e-03,1.563949894119044115e-01,6.294055615508064444e-03,1.883279121861911260e-03,5.720312203382983598e-04,1.349656981629276204e-03,1.103536343511310053e-01,6.294055615508065311e-03,1.106429074209253766e-01,2.236232170338078862e-02,6.377211648739828931e+00,2.242094068303562490e-02,nan
7.968000000000000000e+03,5.120000000000000000e+02,3.038274983597794889e-03,3.095475417421657272e-02,2.467170485869195004e-03,1.403226813346256053e-03,4.292449653044802960e-04,9.954335435773751638e-04,2.168432439219721009e-02,2.467170485869193703e-03,2.195221973819038233e-02,4.394162827803690247e-03,2.499766338757094974e+00,4.448449774900294258e-03,nan
3.180800000000000000e+04,2.048000000000000000e+03,1.158356575513208630e-03,4.609723587879824280e-03,7.139825307461958488e-04,7.290250898613020005e-04,2.282272654716716454e-04,4.984512368045536292e-04,3.105524588433945365e-03,7.139825307461957404e-04,3.330990473279242439e-03,6.293108542610267213e-04,7.234155511515690806e-01,6.749997949080121552e-04,nan
1.271040000000000000e+05,8.192000000000000000e+03,3.642235833837260918e-04,6.351681268992295966e-04,1.463593103098212542e-04,2.695403364251538732e-04,8.492692200078895545e-05,1.771248032823249279e-04,3.588483320686073567e-04,1.463593103098213084e-04,5.032349116326169164e-04,7.271787550654089697e-05,1.482929855766729754e-01,1.019767137935293297e-04,nan
        }\tableUniformOne

        \pgfplotstableread[col sep=comma]{
ndof,nelem,res,err,res1,res2,res3,res4,err1,err2,err3,relerr1,relerr2,relerr3,cond
1.280000000000000000e+02,8.000000000000000000e+00,2.375582897223976753e-01,1.826813220988875663e+00,1.858848670459484209e-01,1.189490308530601614e-01,4.040145419231600671e-02,7.809998577083511473e-02,1.270970971082695522e+00,1.858848670459484487e-01,1.298971189537385662e+00,2.575771905485951252e-01,1.883227439664203029e+00,2.632517635863744210e-01,nan
5.040000000000000000e+02,3.200000000000000000e+01,7.639579290037777637e-02,2.251652437385122385e-01,4.368131067025483871e-02,5.123326677053109679e-02,1.571414530015809510e-02,3.250351187966010746e-02,1.463200465622337409e-01,4.368131067025483177e-02,1.654743909491074583e-01,2.965064061657510588e-02,4.425842100158915993e-01,3.353212230691714896e-02,nan
2.000000000000000000e+03,1.280000000000000000e+02,2.372367075564015340e-02,3.589604235027741425e-02,8.560160406122438448e-03,1.814518408301150909e-02,5.613800684067880853e-03,1.134785197368445853e-02,1.886111933452701284e-02,8.560160406122440183e-03,2.931736154027327762e-02,3.822061871587336887e-03,8.673255845166046030e-02,5.940939545061707690e-03,nan
7.968000000000000000e+03,5.120000000000000000e+02,6.414960750285819517e-03,8.011998063363964459e-03,1.399488411679003985e-03,5.120038029111360391e-03,1.593141614551505996e-03,3.231139735406294514e-03,3.731519764844570472e-03,1.399488411679004419e-03,6.950489579698571063e-03,7.561639987176620821e-04,1.417978223650525650e-02,1.408463662217592600e-03,nan
3.180800000000000000e+04,2.048000000000000000e+03,1.644725050349019101e-03,1.983400249728644194e-03,2.475619877032677902e-04,1.326355007806635745e-03,4.141879872839005178e-04,8.444313222198483303e-04,9.117908856187558440e-04,2.475619877032677902e-04,1.743911406514281988e-03,1.847674635303767760e-04,2.508327361894504247e-03,3.533903357508046844e-04,nan
1.271040000000000000e+05,8.192000000000000000e+03,4.145715838735529934e-04,4.955688913316846174e-04,5.251842167587991744e-05,3.349250520004459645e-04,1.048076090555402423e-04,2.143641485105670535e-04,2.272639325712182557e-04,5.251842167587995132e-05,4.372430061106229410e-04,4.605330129465714579e-05,5.321228647227591079e-04,8.860395783692563828e-05,nan
        }\tableUniformTwo

        \pgfplotstableread[col sep=comma]{
ndof,nelem,res,err,res1,res2,res3,res4,err1,err2,err3,relerr1,relerr2,relerr3,cond
1.280000000000000000e+02,8.000000000000000000e+00,4.871087757425438536e-01,6.730697811503696393e-01,2.536673945166825206e-01,3.571091853888883061e-01,1.030899235930625990e-01,1.864760299045930547e-01,3.292539624805464471e-01,2.536673945166825761e-01,5.294030708711304767e-01,6.672718147172666814e-02,2.569942381505733661e-02,1.072897483619334319e-01,nan
5.040000000000000000e+02,3.200000000000000000e+01,1.555119313674683534e-01,1.750823018626533889e-01,1.226468963129830753e-01,8.004968256135694404e-02,2.386884883206959532e-02,4.651909285111180403e-02,6.274289711044697682e-02,1.226468963129830891e-01,1.080503591502847416e-01,1.271436920075985112e-02,1.242672870449056631e-02,2.189558056417318516e-02,nan
2.000000000000000000e+03,1.280000000000000000e+02,4.021441193566580030e-02,4.426432484022040786e-02,3.111951858324811526e-02,2.101442338065660081e-02,6.451318380281016035e-03,1.286658430248901856e-02,1.505644111020733485e-02,3.111951858324811873e-02,2.764434115334071101e-02,3.051072869454737966e-03,3.153066457234439054e-03,5.601914733338167474e-03,nan
7.968000000000000000e+03,5.120000000000000000e+02,1.016722582488613705e-02,1.109839424116029760e-02,7.808590452094870023e-03,5.333061169316734296e-03,1.655897181385254397e-03,3.348859390104458322e-03,3.693007060643157050e-03,7.808590452094867421e-03,6.968641814135141492e-03,7.483596931677531042e-04,7.911756271845110048e-04,1.412142073975325776e-03,nan
3.180800000000000000e+04,2.048000000000000000e+03,2.552647844936707742e-03,2.776470844509669293e-03,1.953940654926306694e-03,1.340391767924414834e-03,4.183470029092952662e-04,8.523277720005651998e-04,9.144670492768303084e-04,1.953940654926305826e-03,1.747757501264212386e-03,1.853097676692902073e-04,1.979755799240186112e-04,3.541697174957317631e-04,nan
1.271040000000000000e+05,8.192000000000000000e+03,6.393035049182512267e-04,6.942067855335728142e-04,4.885971204718436668e-04,3.358164401747438881e-04,1.050726547364468484e-04,2.148694510682737076e-04,2.275170273457742423e-04,4.885971204718438836e-04,4.375293329769109414e-04,4.610458901892055671e-05,4.950523857044502338e-05,8.866197979092743816e-05,nan
        }\tableUniformThree

        \pgfplotstableread[col sep=comma]{
ndof,nelem,res,err,res1,res2,res3,res4,err1,err2,err3,relerr1,relerr2,relerr3,cond
1.280000000000000000e+02,8.000000000000000000e+00,2.477436157319054644e+00,2.519131816916768241e+00,2.440338381255072520e+00,3.674387645898135601e-01,1.055354746603028859e-01,1.904980746099869393e-01,3.300319907637220829e-01,2.440338381255072076e+00,5.308979007759280222e-01,6.688485803862498358e-02,2.472343378285480307e-03,1.075926931938659681e-01,nan
5.040000000000000000e+02,3.200000000000000000e+01,1.228066820938175763e+00,1.230671057125228396e+00,1.224291586092678852e+00,8.058572754255033388e-02,2.401586588941609715e-02,4.676972967518463337e-02,6.288643010098576847e-02,1.224291586092678630e+00,1.081973196412657212e-01,1.274345506574608121e-02,1.240466726262658050e-03,2.192536098596288910e-02,nan
2.000000000000000000e+03,1.280000000000000000e+02,3.120916731609919559e-01,3.126375577659990523e-01,3.110471704071128274e-01,2.104972179247748335e-02,6.461476366043795820e-03,1.288521531781008575e-02,1.506813781487985685e-02,3.110471704071127719e-01,2.765594927694120866e-02,3.053443117379007424e-03,3.151566747424871195e-04,5.604267030984304981e-03,nan
7.968000000000000000e+03,5.120000000000000000e+02,7.834760217636968505e-02,7.847375489982999497e-02,7.807633290028756501e-02,5.335317286586622500e-03,1.656558961260621915e-03,3.350105892278353124e-03,3.693756247018875269e-03,7.807633290028756501e-02,6.969405549232493384e-03,7.485115100684660780e-04,7.910786463910372303e-05,1.412296839063431533e-03,nan
3.180800000000000000e+04,2.048000000000000000e+03,1.960773726969222697e-02,1.963812209302915820e-02,1.953879935109807328e-02,1.340534227829549493e-03,4.183891651276177106e-04,8.524079101746796696e-04,9.145131835401144665e-04,1.953879935109804900e-02,1.747805737312255799e-03,1.853191164256150583e-04,1.979694277203405991e-05,3.541794921626957278e-04,nan
1.271040000000000000e+05,8.192000000000000000e+03,4.903298091552881106e-03,4.910757823413826675e-03,4.885932993185726318e-03,3.358255123969451769e-04,1.050753749068850161e-04,2.148751999110262419e-04,2.275214299212678897e-04,4.885932993185729788e-03,4.375330676026845588e-04,4.610548116724073534e-05,4.950485140667314751e-06,8.866273658433513748e-05,nan
        }\tableUniformFour

        \pgfplotstableread[col sep=comma]{
ndof,nelem,res,err,res1,res2,res3,res4,err1,err2,err3,relerr1,relerr2,relerr3,cond
1.280000000000000000e+02,8.000000000000000000e+00,2.440661546819113781e+01,2.441088215056095123e+01,2.440287565732292663e+01,3.675457002110618654e-01,1.055607616814322330e-01,1.905396264106859194e-01,3.300651015893625861e-01,2.440287565732292663e+01,5.309275379222327951e-01,6.689156833621577591e-02,2.472291896317969874e-04,1.075986995095527532e-01,nan
5.040000000000000000e+02,3.200000000000000000e+01,1.224328983719081698e+01,1.224355130627389698e+01,1.224291168917549300e+01,8.059113671423764891e-02,2.401734878276611321e-02,4.677225834412194333e-02,6.288826805237120177e-02,1.224291168917549122e+01,1.081990117920886146e-01,1.274382751256577899e-02,1.240466303575878484e-04,2.192570388741144072e-02,nan
2.000000000000000000e+03,1.280000000000000000e+02,3.110576150728308065e+00,3.110630966365126149e+00,3.110471523453970377e+00,2.105007560396272501e-02,6.461578167360341671e-03,1.288540206314039809e-02,1.506826220073978137e-02,3.110471523453970821e+00,2.765606924026212804e-02,3.053468323224119065e-03,3.151566564421434165e-05,5.604291340635643784e-03,nan
7.968000000000000000e+03,5.120000000000000000e+02,7.807904955709963346e-01,7.808031710993930119e-01,7.807633189084690617e-01,5.335432969824040715e-03,1.656621019509744346e-03,3.350609326170610105e-03,3.694862191269649571e-03,7.807633189084689507e-01,6.969923908319593416e-03,7.487356212295279362e-04,7.910786361632631651e-06,1.412401880575945038e-03,nan
3.180800000000000000e+04,2.048000000000000000e+03,1.953949220576949863e-01,1.953980304779120436e-01,1.953879928956225587e-01,1.341297582086592227e-03,4.188381005185239959e-04,8.563256358959044122e-04,9.236067726756199857e-04,1.953879928956223921e-01,1.751998115502595172e-03,1.871618628551480960e-04,1.979694270968521658e-06,3.550290456037387540e-04,nan
1.271040000000000000e+05,8.192000000000000000e+03,4.889251856176605321e-02,4.897046995392578123e-02,4.885932992579716111e-02,8.858970223196182322e-04,3.581182673234899389e-04,1.526822939045137818e-03,2.389450075158271094e-03,4.885932992579721662e-02,2.272297306421662728e-03,4.842038197386242557e-04,4.950485140053268168e-07,4.604637053478834933e-04,nan
        }\tableUniformFive

        %
        %
        \addlegendimage{line1}
        \addlegendentry{\(c_\Omega = 10^{-4}\)}
        \addlegendimage{line2}
        \addlegendentry{\(c_\Omega = 10^{-2}\)}
        \addlegendimage{line3}
        \addlegendentry{\(c_\Omega = \hphantom{0}1\hphantom{{}^{-2}}\)}
        \addlegendimage{line4}
        \addlegendentry{\(c_\Omega = 10^2\hphantom{{}^{-}}\)}
        \addlegendimage{line5}
        \addlegendentry{\(c_\Omega = 10^4\hphantom{{}^{-}}\)}

        %
        %
        \addplot+ [line1, forget plot] table [x=ndof, y=relerr1] {\tableUniformOne};
        \addplot+ [line2, forget plot] table [x=ndof, y=relerr1] {\tableUniformTwo};
        \addplot+ [line3, forget plot] table [x=ndof, y=relerr1] {\tableUniformThree};
        \addplot+ [line4, forget plot] table [x=ndof, y=relerr1] {\tableUniformFour};
        \addplot+ [line5, forget plot] table [x=ndof, y=relerr1] {\tableUniformFive};

        %
        %
        \drawslopetriangle[ST1]{1}{1e4}{3e-5}
    \end{loglogaxis}
\end{tikzpicture}

%% file: figures/plot_opdls_square_error_scalingOne.tex
\begin{tikzpicture}[>=stealth]
    %
    %
    \colorlet{col1}{TUblue}
    \colorlet{col2}{TUgreen}
    \colorlet{col3}{TUmagenta}
    \colorlet{col4}{TUyellow}
    \colorlet{col5}{purple}
    \colorlet{col6}{green}
    \pgfplotsset{%
        linedefault/.style = {%
            mark = *,%
            mark size = 2pt,%
            every mark/.append style = {solid},%
            gray,%
            every mark/.append style = {fill = gray!60!white}%
        },%
        line1/.style = {%
            linedefault,%
            col1,%
            every mark/.append style = {fill = col1!60!white}%
        },%
        line2/.style = {%
            linedefault,%
            mark = triangle*,%
            mark size = 2.75pt,%
            col2,%
            every mark/.append style = {fill = col2!60!white}%
        },%
        line3/.style = {%
            linedefault,%
            mark = square*,%
            mark size = 1.66pt,%
            col3,%
            every mark/.append style = {fill = col3!60!white}%
        },%
        line4/.style = {%
            linedefault,%
            mark = pentagon*,%
            mark size = 2.2pt,%
            col4,%
            every mark/.append style = {fill = col4!60!white}%
        },%
        line5/.style = {%
            linedefault,%
            mark = diamond*,%
            mark size = 2.75pt,%
            col5,%
            every mark/.append style = {fill = col5!60!white}%
        },%
        line6/.style = {%
            linedefault,%
            mark = halfsquare*,%
            mark size = 1.66pt,%
            col6,%
            every mark/.append style = {fill = col6!60!white}%
        },%
        minorline/.style = {%
            dashed,%
            every mark/.append style = {fill = black!20!white}%
        },%
        majorline/.style = {%
            solid%
        }%
    }

    %
    %
    \begin{loglogaxis}[%
            width            = 0.36\textwidth,%
            xlabel           = ndof,%
            ylabel           = {relative energy error},%
            ymajorgrids      = true,%
            font             = \footnotesize,%
            grid style       = {densely dotted, semithick},%
            legend style     = {legend pos  = south west}%
        ]

        %
        %
        \pgfplotstableread[col sep=comma]{
ndof,nelem,res,err,res1,res2,res3,res4,err1,err2,err3,relerr1,relerr2,relerr3,cond
4.800000000000000000e+01,8.000000000000000000e+00,5.824569740073050284e+00,5.978388014482380441e+00,5.782111036502226575e+00,6.632144186201187619e-01,3.287305622970130431e-02,2.277512143596326921e-01,1.068216734226136833e+00,5.782111036502225687e+00,1.080383368061995197e+00,2.164866637863290533e-01,5.935337514980965068e-02,2.189523749891618409e-01,nan
1.920000000000000000e+02,3.200000000000000000e+01,2.602226277743531657e+00,2.685085365994393847e+00,2.535069375812703640e+00,4.265970211133552348e-01,3.928216616372146724e-02,4.018417019103963650e-01,5.096194377425955491e-01,2.535069375812703640e+00,7.234602346838138187e-01,1.032704892784452738e-01,2.602497722648702480e-02,1.466036945926555146e-01,nan
7.680000000000000000e+02,1.280000000000000000e+02,1.327090836244198879e+00,1.364784603708738997e+00,1.285757038676186736e+00,2.250409508739079911e-01,1.192766113301835934e-02,2.391928645834421885e-01,2.533118841054591286e-01,1.285757038676186514e+00,3.811809824558484716e-01,5.133171985647185659e-02,1.319955894288527172e-02,7.724341664875188085e-02,nan
3.072000000000000000e+03,5.120000000000000000e+02,6.675874634407752417e-01,6.856729127348368813e-01,6.451905414811618389e-01,1.147394518641682293e-01,3.209566516471576945e-03,1.273842352695973290e-01,1.262480336710458395e-01,6.451905414811622830e-01,1.947766475625721250e-01,2.558320040813623700e-02,6.623514649735832462e-03,3.947000095385266860e-02,nan
1.228800000000000000e+04,2.048000000000000000e+03,3.344668691107521852e-01,3.433545167119145258e-01,3.228852943054964508e-01,5.781741999151951417e-02,8.281964356692758395e-04,6.534243100918565861e-02,6.303306929810956194e-02,3.228852943054965063e-01,9.830688196960070402e-02,1.277317037978801216e-02,3.314734701638773727e-03,1.992113928269521189e-02,nan
4.915200000000000000e+04,8.192000000000000000e+03,1.673580043991213240e-01,1.717683846873191145e-01,1.614789286705520244e-01,2.900591930008470515e-02,2.100834055041600229e-04,3.304345233759307249e-02,3.149643430427118834e-02,1.614789286705514970e-01,4.936505038606268242e-02,6.382512008443859235e-03,1.657739816237325430e-03,1.000345067135915804e-02,nan
1.966080000000000000e+05,3.276800000000000000e+04,8.370460453110137311e-02,8.590215036698857498e-02,8.074400021074099043e-02,1.452533666941699907e-02,5.288719279804675734e-05,1.660960707725690072e-02,1.574359156736353585e-02,8.074400021074086553e-02,2.473307890661073302e-02,3.190318664773606535e-03,8.289164733356968809e-04,5.011969659874280017e-03,nan
        }\tableUniformOne

        \pgfplotstableread[col sep=comma]{
ndof,nelem,res,err,res1,res2,res3,res4,err1,err2,err3,relerr1,relerr2,relerr3,cond
4.800000000000000000e+01,8.000000000000000000e+00,7.829055171052280437e-01,2.435514510113562459e+00,6.733766838797476950e-01,3.799700299976855877e-01,1.011009216965388924e-01,7.004477322817916074e-02,1.636483117141668941e+00,6.733766838797476950e-01,1.673385125394362793e+00,3.316525186523181512e-01,6.912212284257136519e-01,3.391311439140995976e-01,nan
1.920000000000000000e+02,3.200000000000000000e+01,5.343176916160587986e-01,1.582806095164434046e+00,3.222047208268953855e-01,3.246899547835546129e-01,1.922828919290280048e-01,1.981999903754766978e-01,1.045469225699386717e+00,3.222047208268953855e-01,1.143876457875618202e+00,2.118563588190012503e-01,3.307747946384443916e-01,2.317978333051255724e-01,nan
7.680000000000000000e+02,1.280000000000000000e+02,3.277225779072369116e-01,7.257148506209242234e-01,1.478780153591865432e-01,2.032736769476161109e-01,9.370907861689160845e-02,1.882354947627109809e-01,4.638604466907284563e-01,1.478780153591865432e-01,5.381706262711263289e-01,9.399777900713672529e-02,1.518113081535353848e-01,1.090561697106579714e-01,nan
3.072000000000000000e+03,5.120000000000000000e+02,1.791318552191530911e-01,2.942156127263246268e-01,6.795441470343259460e-02,1.113151935110564955e-01,2.999818894004793041e-02,1.190775185771132672e-01,1.739709820513511696e-01,6.795441470343260848e-02,2.273304240490678085e-01,3.525389164172570977e-02,6.976188154729018320e-02,4.606677528513238029e-02,nan
1.228800000000000000e+04,2.048000000000000000e+03,9.247453384659179432e-02,1.293412605826628203e-01,3.277360830246724144e-02,5.735215190456748047e-02,8.137382367087058962e-03,6.420236513640396880e-02,7.073277224582405487e-02,3.277360830246726919e-02,1.032081741838751648e-01,1.433345641250168175e-02,3.364532812551049301e-02,2.091434873974974337e-02,nan
4.915200000000000000e+04,8.192000000000000000e+03,4.676731410896414753e-02,6.183769593738255160e-02,1.621098738289421595e-02,2.894594511971810766e-02,2.091454876104925555e-03,3.289595164553909407e-02,3.255022455319542957e-02,1.621098738289421248e-02,5.001587156545786450e-02,6.596054559107378083e-03,1.664217088037044795e-02,1.013533461583080407e-02,nan
1.966080000000000000e+05,3.276800000000000000e+04,2.348675224692743993e-02,3.055020669572755532e-02,8.082396218013853870e-03,1.451774707487921839e-02,5.282766199032580068e-04,1.659091155320544311e-02,1.587933214921490727e-02,8.082396218013767134e-03,2.481605953856540897e-02,3.217825457616567582e-03,8.297373614945625919e-03,5.028785051572354052e-03,nan
        }\tableUniformTwo

        \pgfplotstableread[col sep=comma]{
ndof,nelem,res,err,res1,res2,res3,res4,err1,err2,err3,relerr1,relerr2,relerr3,cond
4.800000000000000000e+01,8.000000000000000000e+00,9.824089006612565789e-02,3.129928965513778305e+00,9.779055327021696220e-02,9.381015055927770588e-03,5.247920139868474327e-04,3.635865726524492889e-05,2.211449494609573296e+00,9.779055327021696220e-02,2.212777320392705249e+00,4.481762060831381289e-01,1.003820119972266767e+01,4.484453055689160972e-01,nan
1.920000000000000000e+02,3.200000000000000000e+01,9.621721840387577329e-02,3.100721695563716906e+00,9.374910410057528398e-02,2.136826909235328698e-02,3.481853913188010050e-03,3.638228176259539486e-04,2.189656131505477443e+00,9.374910410057528398e-02,2.193420197905622260e+00,4.437171020276763356e-01,9.624266390890641176e+00,4.444798613536183884e-01,nan
7.680000000000000000e+02,1.280000000000000000e+02,8.995785310646450816e-02,2.978932949346358772e+00,8.131900024461129195e-02,3.659767795230423754e-02,1.160820172620848599e-02,2.344174891988592014e-03,2.103979284384584147e+00,8.131900024461129195e-02,2.107296824792640244e+00,4.263553429056115007e-01,8.348194134798534449e+00,4.270276171474625682e-01,nan
3.072000000000000000e+03,5.120000000000000000e+02,7.649480397146880828e-02,2.607992871849256211e+00,5.582441132064394629e-02,4.504412323748387820e-02,2.469405238895495264e-02,9.814319829459057162e-03,1.842600111688690978e+00,5.582441132064393935e-02,1.844813075387871759e+00,3.733888485916050382e-01,5.730924159959331021e+00,3.738372887943220979e-01,nan
1.228800000000000000e+04,2.048000000000000000e+03,5.818799172245009421e-02,1.807470622681477401e+00,2.873477777925518076e-02,3.756353325342726146e-02,2.661124961349150761e-02,2.099945928341621221e-02,1.277227084630025100e+00,2.873477777925517729e-02,1.278598974929211707e+00,2.588203199895345863e-01,2.949907187741316328e+00,2.590983230874166487e-01,nan
4.915200000000000000e+04,8.192000000000000000e+03,3.791428782461976621e-02,8.332408701564942843e-01,1.121728468085511111e-02,2.459849173450751394e-02,1.426139677285637869e-02,2.243195292682496900e-02,5.886206512501700461e-01,1.121728468085513367e-02,5.896526512662909614e-01,1.192794822019834222e-01,1.151564454794018921e+00,1.194886091282879870e-01,nan
1.966080000000000000e+05,3.276800000000000000e+04,2.115063186876112167e-02,2.665907298629425215e-01,3.435872094424922713e-03,1.382665127379686422e-02,4.742875680991427007e-03,1.489539714898791348e-02,1.880321682864652955e-01,3.435872094424922280e-03,1.889516227008848359e-01,3.810328370723986557e-02,3.527260195073184113e-01,3.828960412638818256e-02,nan
        }\tableUniformThree

        \pgfplotstableread[col sep=comma]{
ndof,nelem,res,err,res1,res2,res3,res4,err1,err2,err3,relerr1,relerr2,relerr3,cond
4.800000000000000000e+01,8.000000000000000000e+00,9.869609383277208256e-03,3.141327018429738605e+00,9.869149311132413907e-03,9.529408433092013541e-05,5.390043773396054108e-07,3.734331868098888445e-09,2.221235876655103247e+00,9.869149311132412172e-03,2.221249472001913805e+00,4.501595313126478537e-01,1.013068268276417996e+02,4.501622865693104547e-01,nan
1.920000000000000000e+02,3.200000000000000000e+01,9.866951358272835551e-03,3.141177615258922806e+00,9.864299564703410753e-03,2.287111376316671980e-04,3.797138660043047574e-06,3.969225549899792963e-08,2.221116764416202649e+00,9.864299564703410753e-03,2.221157316586812591e+00,4.500923591569799265e-01,1.012667242861359966e+02,4.501005767448308936e-01,nan
7.680000000000000000e+02,1.280000000000000000e+02,9.858459630087607062e-03,3.139823111582714432e+00,9.847315530878770345e-03,4.683516840613853247e-04,1.580226134345673493e-05,3.193026239046140281e-07,2.220157856881814329e+00,9.847315530878768611e-03,2.220200732576166125e+00,4.498980438642049551e-01,1.010923664962793254e+02,4.499067322964073745e-01,nan
3.072000000000000000e+03,5.120000000000000000e+02,9.824834661931285129e-03,3.134423123080463114e+00,9.780067111422725140e-03,9.346985411242788233e-04,6.321191853589778922e-05,2.514145811479907740e-06,2.216339509413211672e+00,9.780067111422721671e-03,2.216382603207935453e+00,4.491242848940509846e-01,1.004019954154879315e+02,4.491330175226290744e-01,nan
1.228800000000000000e+04,2.048000000000000000e+03,9.695853039807432461e-03,3.113174434090232090e+00,9.522146708360899880e-03,1.810704937196767244e-03,2.434470770151383958e-04,1.921939381860080764e-05,2.201315603975961821e+00,9.522146708360901615e-03,2.201357308000463586e+00,4.460798051303872525e-01,9.775418913453424352e+01,4.460882561326375839e-01,nan
4.915200000000000000e+04,8.192000000000000000e+03,9.253708942882446845e-03,3.033298178127426858e+00,8.639766034073912207e-03,3.203838937660506182e-03,8.387947261007030206e-04,1.319528136925621545e-04,2.144838662157243725e+00,8.639766034073896595e-03,2.144875358893744366e+00,4.346351839432459618e-01,8.869568479042705178e+01,4.346426202568020591e-01,nan
1.966080000000000000e+05,3.276800000000000000e+04,8.148548419938150475e-03,2.771378800988218316e+00,6.483879339726128857e-03,4.419333315366084057e-03,2.096240808172882865e-03,6.583445767949370597e-04,1.959642365555702037e+00,6.483879339726120183e-03,1.959668394651211010e+00,3.971065679875494525e-01,6.656339024314299024e+01,3.971118425850816425e-01,nan
        }\tableUniformFour

        %
        %
        \addlegendimage{line1}
        \addlegendentry{\(\ell = \hphantom{0}1\hphantom{{}^2}\)}
        \addlegendimage{line2}
        \addlegendentry{\(\ell = 10\hphantom{{}^2}\)}
        \addlegendimage{line3}
        \addlegendentry{\(\ell = 10^2\)}
        \addlegendimage{line4}
        \addlegendentry{\(\ell = 10^3\)}

        %
        %
        \addplot+ [line1, forget plot] table [x=ndof, y=relerr1] {\tableUniformOne};
        \addplot+ [line2, forget plot] table [x=ndof, y=relerr1] {\tableUniformTwo};
        \addplot+ [line3, forget plot] table [x=ndof, y=relerr1] {\tableUniformThree};
        \addplot+ [line4, forget plot] table [x=ndof, y=relerr1] {\tableUniformFour};

        %
        %
        \drawslopetriangle[ST1]{0.5}{7e3}{4e-3}
    \end{loglogaxis}
\end{tikzpicture}

%% file: figures/plot_opdls_square_error_scalingFriedrichs.tex
\begin{tikzpicture}[>=stealth]
    %
    %
    \colorlet{col1}{TUblue}
    \colorlet{col2}{TUgreen}
    \colorlet{col3}{TUmagenta}
    \colorlet{col4}{TUyellow}
    \colorlet{col5}{purple}
    \colorlet{col6}{green}
    \pgfplotsset{%
        linedefault/.style = {%
            mark = *,%
            mark size = 2pt,%
            every mark/.append style = {solid},%
            gray,%
            every mark/.append style = {fill = gray!60!white}%
        },%
        line1/.style = {%
            linedefault,%
            col1,%
            every mark/.append style = {fill = col1!60!white}%
        },%
        line2/.style = {%
            linedefault,%
            mark = triangle*,%
            mark size = 2.75pt,%
            col2,%
            every mark/.append style = {fill = col2!60!white}%
        },%
        line3/.style = {%
            linedefault,%
            mark = square*,%
            mark size = 1.66pt,%
            col3,%
            every mark/.append style = {fill = col3!60!white}%
        },%
        line4/.style = {%
            linedefault,%
            mark = pentagon*,%
            mark size = 2.2pt,%
            col4,%
            every mark/.append style = {fill = col4!60!white}%
        },%
        line5/.style = {%
            linedefault,%
            mark = diamond*,%
            mark size = 2.75pt,%
            col5,%
            every mark/.append style = {fill = col5!60!white}%
        },%
        line6/.style = {%
            linedefault,%
            mark = halfsquare*,%
            mark size = 1.66pt,%
            col6,%
            every mark/.append style = {fill = col6!60!white}%
        },%
        minorline/.style = {%
            dashed,%
            every mark/.append style = {fill = black!20!white}%
        },%
        majorline/.style = {%
            solid%
        }%
    }

    %
    %
    \begin{loglogaxis}[%
            width            = 0.36\textwidth,%
            xlabel           = ndof,%
            ylabel           = {relative energy error},%
            ymajorgrids      = true,%
            font             = \footnotesize,%
            grid style       = {densely dotted, semithick},%
            legend style     = {legend pos  = south west}%
        ]

        %
        %
        \pgfplotstableread[col sep=comma]{
ndof,nelem,res,err,res1,res2,res3,res4,err1,err2,err3,relerr1,relerr2,relerr3,cond
4.800000000000000000e+01,8.000000000000000000e+00,1.963906731866401634e+00,2.465281827253487545e+00,1.850969649263351391e+00,6.195146866585253820e-01,8.957155191793357563e-02,1.975335391205203239e-01,1.139484790094576594e+00,1.850969649263351391e+00,1.163228377592000928e+00,2.309299721105272929e-01,1.875245085236940390e-01,2.357418889050717503e-01,nan
1.920000000000000000e+02,3.200000000000000000e+01,9.897182680353030504e-01,1.218506557525538270e+00,8.120196063671938358e-01,4.123382340308177629e-01,1.136234762786469121e-01,3.704501253242736825e-01,5.320997292967883929e-01,8.120196063671938358e-01,7.363778022797061729e-01,1.078259487762366214e-01,8.227478765790958992e-02,1.492213410698469411e-01,nan
7.680000000000000000e+02,1.280000000000000000e+02,5.231910176698021386e-01,6.177820082933358492e-01,4.101151425105433179e-01,2.227141830365041986e-01,3.659536074510880838e-02,2.336526976095885111e-01,2.578487622710572413e-01,4.101151425105432624e-01,3.833721370284854935e-01,5.225108358803082220e-02,4.155335166881445746e-02,7.768743739843768370e-02,nan
3.072000000000000000e+03,5.120000000000000000e+02,2.672899420862317466e-01,3.104661945547159130e-01,2.054881614391093625e-01,1.144162603766797626e-01,1.002075088152585766e-02,1.265979168944807720e-01,1.269220521433700177e-01,2.054881614391092792e-01,1.950760522190690405e-01,2.571978510696160819e-02,2.082030374149839227e-02,3.953067302222125490e-02,nan
1.228800000000000000e+04,2.048000000000000000e+03,1.347857638508218048e-01,1.556364298635636301e-01,1.027929144977703102e-01,5.777515582194953092e-02,2.597727238808860036e-03,6.523880784353577877e-02,6.312311871685116760e-02,1.027929144977704906e-01,9.834524347007496392e-02,1.279141820717442056e-02,1.041509976695976190e-02,1.992891294796390877e-02,nan
4.915200000000000000e+04,8.192000000000000000e+03,6.763619349765129241e-02,7.792518803863118793e-02,5.140229056168783839e-02,2.900052707395757434e-02,6.597315284061589869e-04,3.303018891266935181e-02,3.150800793034408409e-02,5.140229056168791472e-02,4.936987838884367152e-02,6.384857315429338757e-03,5.208140921637578746e-03,1.000442902927195510e-02,nan
1.966080000000000000e+05,3.276800000000000000e+04,3.387332533147739688e-02,3.899035114592436929e-02,2.570185939847360532e-02,1.452465609593290435e-02,1.661332436934150330e-04,1.660793056083211677e-02,1.574505643873186714e-02,2.570185939847349429e-02,2.473368358503709799e-02,3.190615509775452698e-03,2.604142816062253191e-03,5.012092193342044741e-03,nan
        }\tableUniformOne

        \pgfplotstableread[col sep=comma]{
ndof,nelem,res,err,res1,res2,res3,res4,err1,err2,err3,relerr1,relerr2,relerr3,cond
4.800000000000000000e+01,8.000000000000000000e+00,1.963906731866401412e+00,2.465281827253487990e+00,1.850969649263351391e+00,6.195146866585249379e-01,8.957155191793375604e-02,1.975335391205201574e-01,1.139484790094577260e+00,1.850969649263351391e+00,1.163228377592002039e+00,2.309299721105274317e-01,1.875245085236940390e-01,2.357418889050719724e-01,nan
1.920000000000000000e+02,3.200000000000000000e+01,9.897182680353031614e-01,1.218506557525539824e+00,8.120196063671940578e-01,4.123382340308172633e-01,1.136234762786508951e-01,3.704501253242724612e-01,5.320997292967902803e-01,8.120196063671939468e-01,7.363778022797071721e-01,1.078259487762369961e-01,8.227478765790960380e-02,1.492213410698471632e-01,nan
7.680000000000000000e+02,1.280000000000000000e+02,5.231910176698022497e-01,6.177820082933427326e-01,4.101151425105442616e-01,2.227141830365026998e-01,3.659536074513643211e-02,2.336526976095842922e-01,2.578487622710652349e-01,4.101151425105442061e-01,3.833721370284902674e-01,5.225108358803244590e-02,4.155335166881454073e-02,7.768743739843866902e-02,nan
3.072000000000000000e+03,5.120000000000000000e+02,2.672899420862318021e-01,3.104661945547119717e-01,2.054881614391090572e-01,1.144162603766802622e-01,1.002075088149233933e-02,1.265979168944836308e-01,1.269220521433656046e-01,2.054881614391090572e-01,1.950760522190658763e-01,2.571978510696071307e-02,2.082030374149837146e-02,3.953067302222061652e-02,nan
1.228800000000000000e+04,2.048000000000000000e+03,1.347857638508216938e-01,1.556364298635752874e-01,1.027929144977717119e-01,5.777515582194726190e-02,2.597727238977536306e-03,6.523880784352861784e-02,6.312311871686514253e-02,1.027929144977718229e-01,9.834524347008302692e-02,1.279141820717725163e-02,1.041509976695989548e-02,1.992891294796554288e-02,nan
4.915200000000000000e+04,8.192000000000000000e+03,6.763619349765159772e-02,7.792518803860612464e-02,5.140229056168487548e-02,2.900052707396155380e-02,6.597315276845963118e-04,3.303018891268549168e-02,3.150800793031539870e-02,5.140229056168487548e-02,4.936987838882559571e-02,6.384857315423525699e-03,5.208140921637271700e-03,1.000442902926829310e-02,nan
1.966080000000000000e+05,3.276800000000000000e+04,3.387332533147736219e-02,3.899035114592728363e-02,2.570185939847387940e-02,1.452465609593160852e-02,1.661332438898688074e-04,1.660793056083078104e-02,1.574505643873655089e-02,2.570185939847378226e-02,2.473368358503842332e-02,3.190615509776402025e-03,2.604142816062281814e-03,5.012092193342313623e-03,nan
        }\tableUniformTwo

        \pgfplotstableread[col sep=comma]{
ndof,nelem,res,err,res1,res2,res3,res4,err1,err2,err3,relerr1,relerr2,relerr3,cond
4.800000000000000000e+01,8.000000000000000000e+00,1.963906731866401856e+00,2.465281827253489766e+00,1.850969649263351613e+00,6.195146866585252710e-01,8.957155191793539362e-02,1.975335391205196578e-01,1.139484790094578592e+00,1.850969649263351613e+00,1.163228377592003371e+00,2.309299721105276537e-01,1.875245085236940112e-01,2.357418889050721944e-01,nan
1.920000000000000000e+02,3.200000000000000000e+01,9.897182680353032724e-01,1.218506557525541156e+00,8.120196063671942799e-01,4.123382340308173744e-01,1.136234762786532820e-01,3.704501253242716841e-01,5.320997292967912795e-01,8.120196063671941689e-01,7.363778022797086154e-01,1.078259487762371904e-01,8.227478765790961768e-02,1.492213410698474130e-01,nan
7.680000000000000000e+02,1.280000000000000000e+02,5.231910176698021386e-01,6.177820082933417334e-01,4.101151425105441506e-01,2.227141830365034492e-01,3.659536074512855647e-02,2.336526976095847641e-01,2.578487622710635696e-01,4.101151425105440396e-01,3.833721370284898788e-01,5.225108358803207120e-02,4.155335166881452685e-02,7.768743739843853024e-02,nan
3.072000000000000000e+03,5.120000000000000000e+02,2.672899420862316355e-01,3.104661945547165791e-01,2.054881614391095013e-01,1.144162603766804287e-01,1.002075088151905581e-02,1.265979168944802447e-01,1.269220521433697402e-01,2.054881614391094458e-01,1.950760522190700952e-01,2.571978510696155268e-02,2.082030374149840268e-02,3.953067302222147000e-02,nan
1.228800000000000000e+04,2.048000000000000000e+03,1.347857638508209444e-01,1.556364298635837806e-01,1.027929144977724613e-01,5.777515582194822641e-02,2.597727239065455128e-03,6.523880784352152629e-02,6.312311871687206755e-02,1.027929144977725584e-01,9.834524347009125644e-02,1.279141820717864982e-02,1.041509976695996834e-02,1.992891294796720128e-02,nan
4.915200000000000000e+04,8.192000000000000000e+03,6.763619349765090383e-02,7.792518803864569021e-02,5.140229056168943433e-02,2.900052707395550655e-02,6.597315288044361013e-04,3.303018891265992185e-02,3.150800793036070274e-02,5.140229056168949678e-02,4.936987838885430191e-02,6.384857315432706723e-03,5.208140921637739208e-03,1.000442902927410962e-02,nan
1.966080000000000000e+05,3.276800000000000000e+04,3.387332533147476010e-02,3.899035114611095615e-02,2.570185939849532059e-02,1.452465609590011981e-02,1.661332544598152826e-04,1.660793056071410354e-02,1.574505643895159068e-02,2.570185939849530671e-02,2.473368358516869411e-02,3.190615509819977411e-03,2.604142816064463229e-03,5.012092193368710910e-03,nan
        }\tableUniformThree

        \pgfplotstableread[col sep=comma]{
ndof,nelem,res,err,res1,res2,res3,res4,err1,err2,err3,relerr1,relerr2,relerr3,cond
4.800000000000000000e+01,8.000000000000000000e+00,1.963906731866401634e+00,2.465281827253487990e+00,1.850969649263351613e+00,6.195146866585253820e-01,8.957155191793422788e-02,1.975335391205201852e-01,1.139484790094577260e+00,1.850969649263351613e+00,1.163228377592001817e+00,2.309299721105273762e-01,1.875245085236940112e-01,2.357418889050718891e-01,nan
1.920000000000000000e+02,3.200000000000000000e+01,9.897182680353033835e-01,1.218506557525542044e+00,8.120196063671943909e-01,4.123382340308170413e-01,1.136234762786570568e-01,3.704501253242708514e-01,5.320997292967925008e-01,8.120196063671943909e-01,7.363778022797090594e-01,1.078259487762374541e-01,8.227478765790963156e-02,1.492213410698475240e-01,nan
7.680000000000000000e+02,1.280000000000000000e+02,5.231910176698022497e-01,6.177820082933428436e-01,4.101151425105442061e-01,2.227141830365036157e-01,3.659536074513263654e-02,2.336526976095840979e-01,2.578487622710640692e-01,4.101151425105443171e-01,3.833721370284910446e-01,5.225108358803220998e-02,4.155335166881455461e-02,7.768743739843882168e-02,nan
3.072000000000000000e+03,5.120000000000000000e+02,2.672899420862316355e-01,3.104661945547341761e-01,2.054881614391114997e-01,1.144162603766771674e-01,1.002075088165383342e-02,1.265979168944693090e-01,1.269220521433899740e-01,2.054881614391116940e-01,1.950760522190825297e-01,2.571978510696564663e-02,2.082030374149863167e-02,3.953067302222398188e-02,nan
1.228800000000000000e+04,2.048000000000000000e+03,1.347857638508210831e-01,1.556364298636016275e-01,1.027929144977745707e-01,5.777515582194431981e-02,2.597727239345779936e-03,6.523880784351079876e-02,6.312311871689377241e-02,1.027929144977746817e-01,9.834524347010337175e-02,1.279141820718305081e-02,1.041509976696018518e-02,1.992891294796966112e-02,nan
4.915200000000000000e+04,8.192000000000000000e+03,6.763619349764941890e-02,7.792518803866903265e-02,5.140229056169177968e-02,2.900052707395571125e-02,6.597315293615469908e-04,3.303018891264193624e-02,3.150800793038284475e-02,5.140229056169170335e-02,4.936987838887472307e-02,6.384857315437193585e-03,5.208140921637962120e-03,1.000442902927824694e-02,nan
1.966080000000000000e+05,3.276800000000000000e+04,3.387332533147476704e-02,3.899035114609177011e-02,2.570185939849308973e-02,1.452465609590515398e-02,1.661332532698127048e-04,1.660793056072506699e-02,1.574505643892660373e-02,2.570185939849313830e-02,2.473368358515661350e-02,3.190615509814914187e-03,2.604142816064243787e-03,5.012092193366262348e-03,nan
        }\tableUniformFour

        %
        %
        \addlegendimage{line1}
        \addlegendentry{\(\ell = \hphantom{0}1\hphantom{{}^2}\)}
        \addlegendimage{line2}
        \addlegendentry{\(\ell = 10\hphantom{{}^2}\)}
        \addlegendimage{line3}
        \addlegendentry{\(\ell = 10^2\)}
        \addlegendimage{line4}
        \addlegendentry{\(\ell = 10^3\)}

        %
        %
        \addplot+ [line1, forget plot] table [x=ndof, y=relerr1] {\tableUniformOne};
        \addplot+ [line2, forget plot] table [x=ndof, y=relerr1] {\tableUniformTwo};
        \addplot+ [line3, forget plot] table [x=ndof, y=relerr1] {\tableUniformThree};
        \addplot+ [line4, forget plot] table [x=ndof, y=relerr1] {\tableUniformFour};

        %
        %
        \drawslopetriangle[ST1]{0.5}{7e3}{4e-3}
    \end{loglogaxis}
\end{tikzpicture}

%% file: figures/plot_dls_rectangle_error_scalingOne.tex
\begin{tikzpicture}[>=stealth]
    %
    %
    \colorlet{col1}{TUblue}
    \colorlet{col2}{TUgreen}
    \colorlet{col3}{TUmagenta}
    \colorlet{col4}{TUyellow}
    \colorlet{col5}{purple}
    \colorlet{col6}{green}
    \pgfplotsset{%
        linedefault/.style = {%
            mark = *,%
            mark size = 2pt,%
            every mark/.append style = {solid},%
            gray,%
            every mark/.append style = {fill = gray!60!white}%
        },%
        line1/.style = {%
            linedefault,%
            col1,%
            every mark/.append style = {fill = col1!60!white}%
        },%
        line2/.style = {%
            linedefault,%
            mark = triangle*,%
            mark size = 2.75pt,%
            col2,%
            every mark/.append style = {fill = col2!60!white}%
        },%
        line3/.style = {%
            linedefault,%
            mark = square*,%
            mark size = 1.66pt,%
            col3,%
            every mark/.append style = {fill = col3!60!white}%
        },%
        line4/.style = {%
            linedefault,%
            mark = pentagon*,%
            mark size = 2.2pt,%
            col4,%
            every mark/.append style = {fill = col4!60!white}%
        },%
        line5/.style = {%
            linedefault,%
            mark = diamond*,%
            mark size = 2.75pt,%
            col5,%
            every mark/.append style = {fill = col5!60!white}%
        },%
        line6/.style = {%
            linedefault,%
            mark = halfsquare*,%
            mark size = 1.66pt,%
            col6,%
            every mark/.append style = {fill = col6!60!white}%
        },%
        minorline/.style = {%
            dashed,%
            every mark/.append style = {fill = black!20!white}%
        },%
        majorline/.style = {%
            solid%
        }%
    }

    %
    %
    \begin{loglogaxis}[%
            width            = 0.36\textwidth,%
            xlabel           = ndof,%
            ylabel           = {relative energy error},%
            ymajorgrids      = true,%
            font             = \footnotesize,%
            grid style       = {densely dotted, semithick},%
            legend style     = {legend pos  = south west}%
        ]

        %
        %
        \pgfplotstableread[col sep=comma]{
ndof,nelem,res,err,res1,res2,res3,res4,err1,err2,err3,relerr1,relerr2,relerr3,cond
6.400000000000000000e+01,8.000000000000000000e+00,2.747879563389379776e+00,2.933293319966482127e+00,2.510919688116741710e+00,8.151127593483262768e-01,4.602841695690459733e-16,7.627028282812968740e-01,9.745481625563479078e-01,2.510919688116741710e+00,1.161786512077739308e+00,1.974847466933499218e-01,5.155411392228821665e-02,2.354271690868403322e-01,nan
2.480000000000000000e+02,3.200000000000000000e+01,1.420527798690668364e+00,1.499536591330805635e+00,1.282717955549594002e+00,4.438391541559315590e-01,8.467313638123118331e-16,4.189757494196689191e-01,4.994026519837328482e-01,1.282717955549593780e+00,5.948458846140707301e-01,1.012001356262290097e-01,2.633671953887431125e-02,1.205409782277427622e-01,nan
9.760000000000000000e+02,1.280000000000000000e+02,7.163337311675977981e-01,7.538919706006652621e-01,6.448102943094657924e-01,2.268064475323365858e-01,8.530024755327154072e-16,2.142721029877520389e-01,2.512915846527825625e-01,6.448102943094656814e-01,2.990406781215357856e-01,5.092232159275750236e-02,1.323922207803757307e-02,6.059831092896293681e-02,nan
3.872000000000000000e+03,5.120000000000000000e+02,3.589320576877879088e-01,3.774599399383077314e-01,3.228377472450164598e-01,1.140276882470897674e-01,4.426441852002201100e-14,1.077297419600404815e-01,1.258468251756795220e-01,3.228377472450162378e-01,1.497143006035023816e-01,2.550189856886049197e-02,6.628493168719237424e-03,3.033846029066335562e-02,nan
1.542400000000000000e+04,2.048000000000000000e+03,1.795618219358548207e-01,1.887941006815728184e-01,1.614729843261882847e-01,5.709242371921573189e-02,4.582481997115731162e-15,5.393865389109854525e-02,6.294861166068546621e-02,1.614729843261885345e-01,7.488097271637528252e-02,1.275605568420510108e-02,3.315357583407147358e-03,1.517405757582556899e-02,nan
6.156800000000000000e+04,8.192000000000000000e+03,8.979288641544204419e-02,9.440505962037100407e-02,8.074325710821593505e-02,2.855605583517032212e-02,8.833569186228552088e-14,2.697852029877934044e-02,3.147745789958855678e-02,8.074325710821562974e-02,3.744344211824043173e-02,6.378666584876435328e-03,1.657817689316747577e-03,7.587627750127085013e-03,nan
2.460160000000000000e+05,3.276800000000000000e+04,4.489793998702961558e-02,4.720353053866453014e-02,4.037247420688632599e-02,1.427925940394449839e-02,4.227054015242156784e-12,1.349040742038853935e-02,1.573912302688278489e-02,4.037247420688609006e-02,1.872208930845060612e-02,3.189413149152439415e-03,8.289262075712468779e-04,3.793888497979739510e-03,nan
9.835520000000000000e+05,1.310720000000000000e+05,2.244915711882361151e-02,2.360189043644469326e-02,2.018634281207433434e-02,7.139783736805490125e-03,3.076245312464042103e-11,6.745347093418650047e-03,7.869610950441136357e-03,2.018634281207410536e-02,9.361090723631542571e-03,1.594716592606737699e-03,4.144652741914478500e-04,1.896953584603241301e-03,nan
        }\tableUniformOne

        \pgfplotstableread[col sep=comma]{
ndof,nelem,res,err,res1,res2,res3,res4,err1,err2,err3,relerr1,relerr2,relerr3,cond
6.220000000000000000e+02,8.000000000000000000e+01,3.996913869571407790e-02,5.566148027582590307e-02,9.756589010082951374e-03,2.817454284712745879e-02,3.122284048149736986e-16,2.661834244127413152e-02,3.586526985529987083e-02,9.756589010082949640e-03,4.143297932313103060e-02,7.267823186781709188e-02,2.003219392875203453e-01,8.396077013696320868e-02,nan
2.444000000000000000e+03,3.200000000000000000e+02,2.021568653856064754e-02,2.582763372176082040e-02,4.305428520072480773e-03,1.435763386662590631e-02,5.337387690122443796e-16,1.356449841445696476e-02,1.647190743011924322e-02,4.305428520072480773e-03,1.942179740748949460e-02,3.337906315333411672e-02,8.839890557175227037e-02,3.935679003577059742e-02,nan
9.688000000000000000e+03,1.280000000000000000e+03,1.013732377144444842e-02,1.263755944276724608e-02,2.071035066863476277e-03,7.213403845748134896e-03,1.206108601148520188e-15,6.814906599291036850e-03,8.033163762952436754e-03,2.071035066863475410e-03,9.533467489790123719e-03,1.627859321710160498e-02,4.252241849049904349e-02,1.931884420562571605e-02,nan
3.857600000000000000e+04,5.120000000000000000e+03,5.072365786456402292e-03,6.283333012597106838e-03,1.024925517106729063e-03,3.611050667887854160e-03,8.555868346844497441e-15,3.411559090788794223e-03,3.990509161800420007e-03,1.024925517106729063e-03,4.744010756855747239e-03,8.086462232184291749e-03,2.104373434198178808e-02,9.613375701932208930e-03,nan
1.539520000000000000e+05,2.048000000000000000e+04,2.536646469537692621e-03,3.137206113684618165e-03,5.111261365429741778e-04,1.806069727726358765e-03,2.657103927709445774e-14,1.706293504222274396e-03,1.991972001078187293e-03,5.111261365429746115e-04,2.369147487853670922e-03,4.036579218632683347e-03,1.049442369531132513e-02,4.800896553852007469e-03,nan
6.151040000000000000e+05,8.192000000000000000e+04,1.268381199960041429e-03,1.568044582276415010e-03,2.553955922981460797e-04,9.031029379560264608e-04,3.664836313382066605e-13,8.532110192103400112e-04,9.955749416699069196e-04,2.553955922981473808e-04,1.184215959595232449e-03,2.017456629892954139e-03,5.243773236914821279e-03,2.399723254286747782e-03,nan
2.459008000000000000e+06,3.276800000000000000e+05,6.341979304066042719e-04,7.839525440657059852e-04,1.276768292509199002e-04,4.515600249987416766e-04,1.985424065830113388e-12,4.266136261625820200e-04,4.977361597186193316e-04,1.276768292509187076e-04,5.920632854079752342e-04,1.008624336885656590e-03,2.621456126848881966e-03,1.199771057374161278e-03,nan
        }\tableUniformTwo

        \pgfplotstableread[col sep=comma]{
ndof,nelem,res,err,res1,res2,res3,res4,err1,err2,err3,relerr1,relerr2,relerr3,cond
6.202000000000000000e+03,8.000000000000000000e+02,1.238656959839968229e-03,7.762877477272346269e-03,1.708612881202346368e-04,8.917730325790300784e-04,5.713425732562170139e-16,8.425071766736988259e-04,5.468202834259527911e-03,1.708612881202346639e-04,5.507434149471079898e-03,1.108089567127128056e-01,3.508117903709686125e+00,1.116039493713288228e-01,nan
2.440400000000000000e+04,3.200000000000000000e+03,6.263637285749151405e-04,2.088273626172271885e-03,4.490359234860994802e-05,4.541325799465852704e-04,1.616210039877291779e-15,4.290440393473823476e-04,1.458248300455522537e-03,4.490359234860990059e-05,1.494115892218357642e-03,2.955028876931632686e-02,9.219589644448213939e-01,3.027711813967845153e-02,nan
9.680800000000000000e+04,1.280000000000000000e+04,3.140792420671675118e-04,6.241948856904872131e-04,1.258515630026769681e-05,2.281209874672330141e-04,3.824230283217133500e-15,2.155184428659384123e-04,4.260897162517936535e-04,1.258515630026769851e-05,4.559697609135108815e-04,8.634382877691949842e-03,2.583979825019546706e-01,9.239879176174173661e-03,nan
3.856160000000000000e+05,5.120000000000000000e+04,1.571523692609416505e-04,2.304408936740267038e-04,4.211514315623506650e-06,1.141931015858111428e-04,1.021022712889236827e-14,1.078845024167478527e-04,1.524856627539077992e-04,4.211514315623526131e-06,1.727234531815377669e-04,3.090005567742464855e-03,8.647066245907751580e-02,3.500108943829074396e-03,nan
1.539232000000000000e+06,2.048000000000000000e+05,7.859031370888425975e-05,1.023778759932019393e-04,1.747949966963880733e-06,5.711320368422398094e-05,3.746340886978002031e-14,5.395798406814487496e-05,6.645778851542934662e-05,1.747949966963872051e-06,7.785587010258695107e-05,1.346716358927211926e-03,3.588884668585731047e-02,1.577689782459781170e-03,nan
        }\tableUniformThree

        \pgfplotstableread[col sep=comma]{
ndof,nelem,res,err,res1,res2,res3,res4,err1,err2,err3,relerr1,relerr2,relerr3,cond
6.200200000000000000e+04,8.000000000000000000e+03,3.916188097125352112e-05,2.402783230105645824e-03,5.337132420228304580e-06,2.820118924965690004e-05,2.047969591660005624e-15,2.664321421376580749e-05,1.698956772871641633e-03,5.337132420228307968e-06,1.699083473466806229e-03,3.442806223690441092e-01,1.095818134339288008e+02,3.443062972775846609e-01,nan
2.440040000000000000e+05,3.200000000000000000e+04,1.980349788773681010e-05,6.071643213631083141e-04,1.348277401467073180e-06,1.436177024591935668e-05,5.348405470136325120e-15,1.356835406531370384e-05,4.292673059413062465e-04,1.348277401467068733e-06,4.293905859870073786e-04,8.698774307386120119e-02,2.768278375570570304e+01,8.701272483416103154e-02,nan
9.680080000000000000e+05,1.280000000000000000e+05,9.931206261539703603e-06,1.172016869255915276e-04,2.597244425229698599e-07,7.216529257568802105e-06,1.243821187779956271e-14,6.817852044291882978e-06,8.279451811319131407e-05,2.597244425229688541e-07,8.295321411387725295e-05,1.677767714865114126e-02,5.332653036097023680e+00,1.680983568190875035e-02,nan
        }\tableUniformFour

        %
        %
        \addlegendimage{line1}
        \addlegendentry{\(\ell = \hphantom{0}1\hphantom{{}^2}\)}
        \addlegendimage{line2}
        \addlegendentry{\(\ell = 10\hphantom{{}^2}\)}
        \addlegendimage{line3}
        \addlegendentry{\(\ell = 10^2\)}
        \addlegendimage{line4}
        \addlegendentry{\(\ell = 10^3\)}

        %
        %
        \addplot+ [line1, forget plot] table [x=ndof, y=relerr1] {\tableUniformOne};
        \addplot+ [line2, forget plot] table [x=ndof, y=relerr1] {\tableUniformTwo};
        \addplot+ [line3, forget plot] table [x=ndof, y=relerr1] {\tableUniformThree};
        \addplot+ [line4, forget plot] table [x=ndof, y=relerr1] {\tableUniformFour};

        %
        %
        \drawslopetriangle[ST1]{0.5}{3e4}{2e-3}
    \end{loglogaxis}
\end{tikzpicture}

%% file: figures/plot_dls_rectangle_error_scalingDiam.tex
\begin{tikzpicture}[>=stealth]
    %
    %
    \colorlet{col1}{TUblue}
    \colorlet{col2}{TUgreen}
    \colorlet{col3}{TUmagenta}
    \colorlet{col4}{TUyellow}
    \colorlet{col5}{purple}
    \colorlet{col6}{green}
    \pgfplotsset{%
        linedefault/.style = {%
            mark = *,%
            mark size = 2pt,%
            every mark/.append style = {solid},%
            gray,%
            every mark/.append style = {fill = gray!60!white}%
        },%
        line1/.style = {%
            linedefault,%
            col1,%
            every mark/.append style = {fill = col1!60!white}%
        },%
        line2/.style = {%
            linedefault,%
            mark = triangle*,%
            mark size = 2.75pt,%
            col2,%
            every mark/.append style = {fill = col2!60!white}%
        },%
        line3/.style = {%
            linedefault,%
            mark = square*,%
            mark size = 1.66pt,%
            col3,%
            every mark/.append style = {fill = col3!60!white}%
        },%
        line4/.style = {%
            linedefault,%
            mark = pentagon*,%
            mark size = 2.2pt,%
            col4,%
            every mark/.append style = {fill = col4!60!white}%
        },%
        line5/.style = {%
            linedefault,%
            mark = diamond*,%
            mark size = 2.75pt,%
            col5,%
            every mark/.append style = {fill = col5!60!white}%
        },%
        line6/.style = {%
            linedefault,%
            mark = halfsquare*,%
            mark size = 1.66pt,%
            col6,%
            every mark/.append style = {fill = col6!60!white}%
        },%
        minorline/.style = {%
            dashed,%
            every mark/.append style = {fill = black!20!white}%
        },%
        majorline/.style = {%
            solid%
        }%
    }

    %
    %
    \begin{loglogaxis}[%
            width            = 0.36\textwidth,%
            xlabel           = ndof,%
            ylabel           = {relative energy error},%
            ymajorgrids      = true,%
            font             = \footnotesize,%
            grid style       = {densely dotted, semithick},%
            legend style     = {legend pos  = south west}%
        ]

        %
        %
        \pgfplotstableread[col sep=comma]{
ndof,nelem,res,err,res1,res2,res3,res4,err1,err2,err3,relerr1,relerr2,relerr3,cond
6.400000000000000000e+01,8.000000000000000000e+00,3.722067886794899749e+00,3.858209861988078693e+00,3.546453038765342569e+00,8.249163461604881453e-01,5.364020030174311572e-16,7.719930178305682844e-01,9.776814571253283903e-01,3.546453038765343013e+00,1.163010383172664808e+00,1.981196849222083556e-01,3.640782396303650359e-02,2.356751772227673181e-01,nan
2.480000000000000000e+02,3.200000000000000000e+01,1.913889957607891601e+00,1.972780774879858701e+00,1.813323729694704234e+00,4.452051819860775028e-01,5.365506825390410093e-16,4.202667808561212914e-01,4.997848788626008121e-01,1.813323729694704456e+00,5.949252994468784417e-01,1.012775909858026485e-01,1.861554923104392911e-02,1.205570710374604254e-01,nan
9.760000000000000000e+02,1.280000000000000000e+02,9.637909646808745512e-01,9.919625182277083120e-01,9.118051445812148659e-01,2.269819218745720846e-01,8.705685723247048759e-16,2.144378909098067776e-01,2.513390795755259255e-01,9.118051445812144218e-01,2.990479609728425858e-01,5.093194607633622084e-02,9.360575434000671657e-03,6.059978674319232983e-02,nan
3.872000000000000000e+03,5.120000000000000000e+02,4.827577551156869706e-01,4.966800181526019875e-01,4.565495660113308984e-01,1.140497728107247172e-01,3.406977180645327062e-14,1.077506065923764661e-01,1.258527533760981920e-01,4.565495660113308429e-01,1.497151250695450964e-01,2.550309987342706808e-02,4.686929743056178488e-03,3.033862736241399807e-02,nan
1.542400000000000000e+04,2.048000000000000000e+03,2.414868222561601940e-01,2.484273748324086484e-01,2.283557848911531574e-01,5.709518900769539052e-02,5.545166613499776628e-15,5.394126641356537105e-02,6.294935241799405612e-02,2.283557848911533794e-01,7.488107306270747054e-02,1.275620579301969815e-02,2.344296435447094235e-03,1.517407791024387095e-02,nan
6.156800000000000000e+04,8.192000000000000000e+03,1.207569068588527361e-01,1.242246061328546591e-01,1.141880216732089753e-01,2.855640164390871663e-02,1.274265991111335540e-12,2.697884700331281208e-02,3.147755048595732991e-02,1.141880216732090170e-01,3.744345457654284920e-02,6.378685346797381868e-03,1.172252204194698623e-03,7.587630274706979108e-03,nan
2.460160000000000000e+05,3.276800000000000000e+04,6.038014033528765756e-02,6.211366768362285712e-02,5.709527711492968366e-02,1.427930263484323767e-02,1.101666416109153519e-11,1.349044826315926558e-02,1.573913460067629527e-02,5.709527711492989877e-02,1.872209086354116592e-02,3.189415494493354816e-03,5.861391016881970948e-04,3.793888813106971979e-03,nan
9.835520000000000000e+05,1.310720000000000000e+05,3.019028167135499527e-02,3.105700593981719840e-02,2.854779684736348522e-02,7.139790046750342338e-03,4.883028430778086825e-11,6.745354004479048489e-03,7.869615688106635851e-03,2.854779684736345052e-02,9.361092993289406999e-03,1.594717552658489988e-03,2.930711758453616008e-04,1.896954044532083924e-03,nan
        }\tableUniformOne

        \pgfplotstableread[col sep=comma]{
ndof,nelem,res,err,res1,res2,res3,res4,err1,err2,err3,relerr1,relerr2,relerr3,cond
6.220000000000000000e+02,8.000000000000000000e+01,1.221553379393372912e-01,1.256724286451672001e-01,1.155185885973321336e-01,2.886962401360561314e-02,3.174411493714387635e-16,2.727502676198633597e-02,3.184408612285672541e-02,1.155185885973321336e-01,3.788098717107377111e-02,6.452960995952800038e-02,1.185911780626396766e-02,7.676292915426845276e-02,nan
2.444000000000000000e+03,3.200000000000000000e+02,6.109514827029376149e-02,6.285035888301825080e-02,5.777241168351655853e-02,1.444620555306419182e-02,5.569919060396432942e-16,1.364817700493643066e-02,1.592576446682585425e-02,5.777241168351657241e-02,1.894428901124706535e-02,3.227234612375771278e-02,5.930905531533020060e-03,3.838915571764184298e-02,nan
9.688000000000000000e+03,1.280000000000000000e+03,3.054975920204242371e-02,3.142694717029238111e-02,2.888784592857444183e-02,7.224529235218787523e-03,1.111032003122803974e-15,6.825417367398843409e-03,7.963347648499393902e-03,2.888784592857442796e-02,9.472616337656315777e-03,1.613711619002777820e-02,2.965621136788011356e-03,1.919553399042176964e-02,nan
3.857600000000000000e+04,5.120000000000000000e+03,1.527515272073266990e-02,1.571369451990896338e-02,1.444412798807396918e-02,3.612443036619554081e-03,1.321423627335117854e-12,3.412874538131135775e-03,3.981732008617083589e-03,1.444412798807401255e-02,4.736366998101835223e-03,8.068676001193377315e-03,1.482831616099583019e-03,9.597886208243750286e-03,nan
1.539520000000000000e+05,2.048000000000000000e+04,7.637610498832804462e-03,7.856874874713140264e-03,7.222089622409996051e-03,1.806243826848023554e-03,2.286978389651529381e-12,1.706457985340488221e-03,1.990873278005836442e-03,7.222089622409994317e-03,2.368190844019709210e-03,4.034352740188968875e-03,7.414184390540112281e-04,4.798957988140384150e-03,nan
6.151040000000000000e+05,8.192000000000000000e+04,3.819551518321039757e-03,3.930187738518249076e-03,3.611048031986086905e-03,9.042195873039576880e-04,6.880292887426301691e-11,8.553904662397312109e-04,9.993744004763839648e-04,3.611048031986086038e-03,1.186574303767392452e-03,2.025155943162405218e-03,3.707095501718437290e-04,2.404502258745933301e-03,nan
2.459008000000000000e+06,3.276800000000000000e+05,4.052394898331838566e-03,5.927018465058996564e-03,1.805579386298585432e-03,2.141867722659281797e-03,4.263855155089880023e-10,2.928171810264843748e-03,4.251489309065810927e-03,1.805579386298573289e-03,3.714063760872056419e-03,8.615318580745963753e-03,1.853604593916541668e-04,7.526266727493481461e-03,nan
        }\tableUniformTwo

        \pgfplotstableread[col sep=comma]{
ndof,nelem,res,err,res1,res2,res3,res4,err1,err2,err3,relerr1,relerr2,relerr3,cond
6.202000000000000000e+03,8.000000000000000000e+02,3.864350227374720951e-03,3.975291917771464449e-03,3.654113684756289931e-03,9.138892507299173828e-04,4.819766076707508302e-16,8.634015844832458424e-04,1.007309140773591204e-03,3.654113684756292533e-03,1.198218387978538428e-03,2.041235088738539358e-02,3.751306624430749120e-03,2.428098106639999323e-02,nan
2.440400000000000000e+04,3.200000000000000000e+03,1.932180640927049853e-03,1.987650429846026488e-03,1.827060991839746544e-03,4.569482374745135675e-04,1.693309527540593780e-15,4.317041458450007209e-04,5.036557479990890884e-04,1.827060991839746110e-03,5.991103828779256821e-04,1.020619930710698654e-02,1.875657572045293271e-03,1.214051462512117686e-02,nan
9.680800000000000000e+04,1.280000000000000000e+04,9.660910113535957746e-04,9.938257737764125653e-04,9.135310146054992971e-04,2.284745702888590999e-04,3.942118017673574062e-15,2.158524920017878819e-04,2.518280212053452616e-04,9.135310146055004897e-04,2.995553399823315814e-04,5.103102636566683899e-03,9.378293185043888201e-04,6.070260322678572296e-03,nan
3.856160000000000000e+05,5.120000000000000000e+04,4.830470152837723608e-04,4.969161064832957296e-04,4.567655721462122792e-04,1.142396601900203700e-04,6.985593332682173948e-14,1.079302141023219904e-04,1.259212957550226073e-04,4.567655721462122250e-04,1.497820291618376007e-04,2.551698946335162390e-03,4.689147258203777216e-04,3.035218496605717169e-03,nan
1.539232000000000000e+06,2.048000000000000000e+05,8.557897171853764366e-04,1.302505557373286440e-03,2.283837527829341582e-04,4.880635158770351063e-04,5.398188827681036487e-11,6.648390075610444237e-04,9.689324734408774707e-04,2.283837527829311496e-04,8.399591956546603733e-04,1.963467701570728238e-02,2.344583553327976610e-04,1.702113198300074168e-02,nan
        }\tableUniformThree

        \pgfplotstableread[col sep=comma]{
ndof,nelem,res,err,res1,res2,res3,res4,err1,err2,err3,relerr1,relerr2,relerr3,cond
6.200200000000000000e+04,8.000000000000000000e+03,1.222019453857944469e-04,1.257101414922296497e-04,1.155535671432979922e-04,2.890001717800891133e-05,1.191211432384352422e-15,2.730343538170253558e-05,3.185401024125669213e-05,1.155535671432979109e-04,3.789109164441753242e-05,6.454972042798589119e-03,1.186270869758577996e-03,7.678340509825398556e-03,nan
2.440040000000000000e+05,3.200000000000000000e+04,6.110097444785825456e-05,6.285507217548894066e-05,5.777678488372499328e-05,1.445000973807171838e-05,5.052919703696984138e-15,1.365171877630168224e-05,1.592700552766553620e-05,5.777678488372506104e-05,1.894554622087720307e-05,3.227486103882250788e-03,5.931354483490361255e-04,3.839170335699787669e-03,nan
9.680080000000000000e+05,1.280000000000000000e+05,3.055049976868640669e-05,3.142756348159447610e-05,2.888839260800900197e-05,7.225025466607229571e-06,1.223953359129590802e-14,6.825893029407646638e-06,7.963565775993432764e-06,2.888839260800881223e-05,9.472810514658185242e-06,1.613755820874543606e-03,2.965677258801684949e-04,1.919592747529512596e-03,nan
        }\tableUniformFour

        %
        %
        \addlegendimage{line1}
        \addlegendentry{\(\ell = \hphantom{0}1\hphantom{{}^2}\)}
        \addlegendimage{line2}
        \addlegendentry{\(\ell = 10\hphantom{{}^2}\)}
        \addlegendimage{line3}
        \addlegendentry{\(\ell = 10^2\)}
        \addlegendimage{line4}
        \addlegendentry{\(\ell = 10^3\)}

        %
        %
        \addplot+ [line1, forget plot] table [x=ndof, y=relerr1] {\tableUniformOne};
        \addplot+ [line2, forget plot] table [x=ndof, y=relerr1] {\tableUniformTwo};
        \addplot+ [line3, forget plot] table [x=ndof, y=relerr1] {\tableUniformThree};
        \addplot+ [line4, forget plot] table [x=ndof, y=relerr1] {\tableUniformFour};

        %
        %
        \drawslopetriangle[ST1]{0.5}{3e4}{2e-3}
    \end{loglogaxis}
\end{tikzpicture}

%% file: figures/plot_dls_rectangle_error_scalingFriedrichs.tex
\begin{tikzpicture}[>=stealth]
    %
    %
    \colorlet{col1}{TUblue}
    \colorlet{col2}{TUgreen}
    \colorlet{col3}{TUmagenta}
    \colorlet{col4}{TUyellow}
    \colorlet{col5}{purple}
    \colorlet{col6}{green}
    \pgfplotsset{%
        linedefault/.style = {%
            mark = *,%
            mark size = 2pt,%
            every mark/.append style = {solid},%
            gray,%
            every mark/.append style = {fill = gray!60!white}%
        },%
        line1/.style = {%
            linedefault,%
            col1,%
            every mark/.append style = {fill = col1!60!white}%
        },%
        line2/.style = {%
            linedefault,%
            mark = triangle*,%
            mark size = 2.75pt,%
            col2,%
            every mark/.append style = {fill = col2!60!white}%
        },%
        line3/.style = {%
            linedefault,%
            mark = square*,%
            mark size = 1.66pt,%
            col3,%
            every mark/.append style = {fill = col3!60!white}%
        },%
        line4/.style = {%
            linedefault,%
            mark = pentagon*,%
            mark size = 2.2pt,%
            col4,%
            every mark/.append style = {fill = col4!60!white}%
        },%
        line5/.style = {%
            linedefault,%
            mark = diamond*,%
            mark size = 2.75pt,%
            col5,%
            every mark/.append style = {fill = col5!60!white}%
        },%
        line6/.style = {%
            linedefault,%
            mark = halfsquare*,%
            mark size = 1.66pt,%
            col6,%
            every mark/.append style = {fill = col6!60!white}%
        },%
        minorline/.style = {%
            dashed,%
            every mark/.append style = {fill = black!20!white}%
        },%
        majorline/.style = {%
            solid%
        }%
    }

    %
    %
    \begin{loglogaxis}[%
            width            = 0.36\textwidth,%
            xlabel           = ndof,%
            ylabel           = {relative energy error},%
            ymajorgrids      = true,%
            font             = \footnotesize,%
            grid style       = {densely dotted, semithick},%
            legend style     = {legend pos  = south west}%
        ]

        %
        %
        \pgfplotstableread[col sep=comma]{
ndof,nelem,res,err,res1,res2,res3,res4,err1,err2,err3,relerr1,relerr2,relerr3,cond
6.400000000000000000e+01,8.000000000000000000e+00,1.280722770480134942e+00,1.792391592290530022e+00,8.894857821520940622e-01,6.737049103391067328e-01,2.105017920753999974e-16,6.286394451198286104e-01,9.989798873695437953e-01,8.894857821520940622e-01,1.193114348289464166e+00,2.024356492463437185e-01,1.802475045613622928e-01,2.417755159786908747e-01,nan
2.480000000000000000e+02,3.200000000000000000e+01,7.191371363979215392e-01,8.963633276407829698e-01,4.267363878139665001e-01,4.209340599987722586e-01,4.945613973072224502e-16,3.973270609901127237e-01,5.062917834591607047e-01,4.267363878139665001e-01,6.041786991085400516e-01,1.025961655369450887e-01,8.647487183313389003e-02,1.224322018503307546e-01,nan
9.760000000000000000e+02,1.280000000000000000e+02,3.714282895700019016e-01,4.440483087671714957e-01,2.078979355871830859e-01,2.237382049793653349e-01,7.879334521542772076e-16,2.113732204043269058e-01,2.523372609288793700e-01,2.078979355871831136e-01,3.004717218715987070e-01,5.113421990875918755e-02,4.212892982098394995e-02,6.088830102216363088e-02,nan
3.872000000000000000e+03,5.120000000000000000e+02,1.872714828489223104e-01,2.212996605603755884e-01,1.031057089006786331e-01,1.136373426576444273e-01,2.066862209954338975e-15,1.073609586126717808e-01,1.259840371362226097e-01,1.031057089006785915e-01,1.499025514895755973e-01,2.552970352536471435e-02,2.089358493219823157e-02,3.037660789586067028e-02,nan
1.542400000000000000e+04,2.048000000000000000e+03,9.383277824791008692e-02,1.105517858784104690e-01,5.144174944904259666e-02,5.704341417553795285e-02,2.546447620260902952e-15,5.389235183717101113e-02,6.296597269996583035e-02,5.144174944904259666e-02,7.490480663131705852e-02,1.275957376630339722e-02,1.042427788561918131e-02,1.517888733677103148e-02,nan
6.156800000000000000e+04,8.192000000000000000e+03,4.694108785041319065e-02,5.526338988213455594e-02,2.570680218296977904e-02,2.854992285465225363e-02,2.620850447153686411e-14,2.697272613399086250e-02,3.147963462776941224e-02,2.570680218296969230e-02,3.744643089124272128e-02,6.379107682227939252e-03,5.209287249675903161e-03,7.588233402163545305e-03,nan
2.460160000000000000e+05,3.276800000000000000e+04,2.347363306728738469e-02,2.763012390184953204e-02,1.285163620918834161e-02,1.427849256918922380e-02,3.482614415314440901e-13,1.348968294942230750e-02,1.573939532426155488e-02,1.285163620918828783e-02,1.872246320345832132e-02,3.189468328137878839e-03,2.604285984911431946e-03,3.793964264949219178e-03,nan
9.835520000000000000e+05,1.310720000000000000e+05,1.173720271002646511e-02,1.381486529302895436e-02,6.425597311925169627e-03,7.139687827635480578e-03,2.496707236166071693e-12,6.745256434350469475e-03,7.869644818656510929e-03,6.425597311925235547e-03,9.361137360262795568e-03,1.594723455742077037e-03,1.302098250506576651e-03,1.896963035160600551e-03,nan
        }\tableUniformOne

        \pgfplotstableread[col sep=comma]{
ndof,nelem,res,err,res1,res2,res3,res4,err1,err2,err3,relerr1,relerr2,relerr3,cond
6.220000000000000000e+02,8.000000000000000000e+01,4.742414920409909751e-02,5.596126349428220265e-02,2.605707538316522076e-02,2.880271079916747745e-02,2.754251134131591179e-16,2.721180995134684158e-02,3.186547212552481162e-02,2.605707538316522076e-02,3.791152226658708879e-02,6.457294706160217335e-02,5.280267440159865100e-02,7.682480619466929006e-02,nan
2.444000000000000000e+03,3.200000000000000000e+02,2.374402069461517897e-02,2.796460535005248702e-02,1.301043460142452951e-02,1.443781845160109316e-02,5.188564963563268675e-16,1.364025324497226271e-02,1.592845953801948211e-02,1.301043460142452951e-02,1.894813765550365514e-02,3.227780747982448434e-02,2.636465267035130455e-02,3.839695470147175965e-02,nan
9.688000000000000000e+03,1.280000000000000000e+03,1.187601047907508971e-02,1.398027332924454968e-02,6.502935223302781266e-03,7.223480128569470179e-03,1.248652637484792936e-15,6.824426218576279228e-03,7.963685217593651292e-03,6.502935223302780399e-03,9.473098415882464859e-03,1.613780024803160296e-02,1.317770187948980043e-02,1.919651088515140352e-02,nan
3.857600000000000000e+04,5.120000000000000000e+03,5.938505423827758652e-03,6.989882131830129491e-03,3.251181748815174920e-03,3.612311875821205484e-03,1.020931086287926428e-14,3.412750623295509808e-03,3.981774226183360660e-03,3.251181748815180992e-03,4.736427289080987965e-03,8.068761551869091039e-03,6.588271660526395443e-03,9.598008383311165681e-03,nan
1.539520000000000000e+05,2.048000000000000000e+04,2.969315237250907871e-03,3.494909219350837010e-03,1.625555122663377601e-03,1.806227430931709680e-03,1.215400159275948900e-13,1.706442495105857425e-03,1.990878555442226246e-03,1.625555122663380637e-03,2.368198381083517351e-03,4.034363434511000976e-03,3.294063382082338277e-03,4.798973261424909278e-03,nan
6.151040000000000000e+05,8.192000000000000000e+04,1.484665505125353901e-03,1.747450767606611739e-03,8.127730899603384735e-04,9.031226960927382641e-04,3.111419246471325919e-12,8.532297241350826796e-04,9.954383505072991484e-04,8.127730899603377145e-04,1.184097369216133636e-03,2.017179838327544123e-03,1.647022630148455707e-03,2.399482939935170152e-03,nan
2.459008000000000000e+06,3.276800000000000000e+05,7.425308801805824248e-04,8.741197565200067117e-04,4.063860510323177031e-04,4.516807716941025389e-04,2.994332340131104374e-11,4.268336514731800330e-04,4.981191954661490099e-04,4.063860510323235036e-04,5.922946845727268929e-04,1.009400529591997350e-03,8.235103141266836031e-04,1.200239970119490896e-03,nan
        }\tableUniformTwo

        \pgfplotstableread[col sep=comma]{
ndof,nelem,res,err,res1,res2,res3,res4,err1,err2,err3,relerr1,relerr2,relerr3,cond
6.202000000000000000e+03,8.000000000000000000e+02,1.502351428033907007e-03,1.768308115515131868e-03,8.224824707854619691e-04,9.138680142280965375e-04,4.477656897895831506e-16,8.633815211926632255e-04,1.007315976580122055e-03,8.224824707854619691e-04,1.198228150217542799e-03,2.041248940978716495e-02,1.666697949301156750e-02,2.428117889072217706e-02,nan
2.440400000000000000e+04,3.200000000000000000e+03,7.511858373764349278e-04,8.841489008655330694e-04,4.112354465085131354e-04,4.569455828390455777e-04,1.541902407393451544e-15,4.317016378631942223e-04,5.036566025409754512e-04,4.112354465085127559e-04,5.991116032539816434e-04,1.020621662374606754e-02,8.333372439185591021e-03,1.214053935511047282e-02,nan
9.680800000000000000e+04,1.280000000000000000e+04,3.755941840905516497e-04,4.420738056808043740e-04,2.056169995745993795e-04,2.284742384649128500e-04,3.848833264716785593e-15,2.158521785075715371e-04,2.518281280153401381e-04,2.056169995745991627e-04,2.995554925140554866e-04,5.103104800989704919e-03,4.166671554776906969e-03,6.070263413617534541e-03,nan
3.856160000000000000e+05,5.120000000000000000e+04,1.877972720982251545e-04,2.210368675591508779e-04,1.028084093639424791e-04,1.142373125180784298e-04,6.388220965344427362e-14,1.079262841027973279e-04,1.259140877454125265e-04,1.028084093639425198e-04,1.497777363016886186e-04,2.551552881522082660e-03,2.083333945028144427e-03,3.035131505071416054e-03,nan
1.539232000000000000e+06,2.048000000000000000e+05,9.593310344977072700e-05,1.146242852612322440e-04,5.140422725562012233e-05,5.832161910723230377e-05,2.445299208800629430e-12,5.620813555409823442e-05,6.704642751658717605e-05,5.140422725561987161e-05,7.746682353862168965e-05,1.358644678994283756e-03,1.041667429951916328e-03,1.569806050789191625e-03,nan
        }\tableUniformThree

        \pgfplotstableread[col sep=comma]{
ndof,nelem,res,err,res1,res2,res3,res4,err1,err2,err3,relerr1,relerr2,relerr3,cond
6.200200000000000000e+04,8.000000000000000000e+03,4.750936872456139055e-05,5.591838194641916967e-05,2.600869613452052405e-05,2.890001046263256199e-05,1.266395518909964077e-15,2.730342903827718122e-05,3.185401239680519793e-05,2.600869613452052066e-05,3.789109472042178001e-05,6.454972479604041170e-03,5.270463754687045033e-03,7.678341133154181243e-03,nan
2.440040000000000000e+05,3.200000000000000000e+04,2.375468756678578375e-05,2.795918932550366023e-05,1.300434622687315960e-05,1.445000889584315894e-05,5.034579697893436006e-15,1.365171796944510856e-05,1.592700575454558795e-05,1.300434622687315451e-05,1.894554656497509168e-05,3.227486149857761947e-03,2.635231504403115228e-03,3.839170405428598068e-03,nan
9.680080000000000000e+05,1.280000000000000000e+05,1.187734444775556584e-05,1.397959485453541220e-05,6.502172829826221971e-06,7.225005067328022113e-06,1.220444269130145116e-14,6.825859755282409218e-06,7.963503245999408971e-06,6.502172829826166067e-06,9.472773450208773506e-06,1.613743149648528327e-03,1.317615694730079733e-03,1.919585236701741669e-03,nan
        }\tableUniformFour

        %
        %
        \addlegendimage{line1}
        \addlegendentry{\(\ell = \hphantom{0}1\hphantom{{}^2}\)}
        \addlegendimage{line2}
        \addlegendentry{\(\ell = 10\hphantom{{}^2}\)}
        \addlegendimage{line3}
        \addlegendentry{\(\ell = 10^2\)}
        \addlegendimage{line4}
        \addlegendentry{\(\ell = 10^3\)}

        %
        %
        \addplot+ [line1, forget plot] table [x=ndof, y=relerr1] {\tableUniformOne};
        \addplot+ [line2, forget plot] table [x=ndof, y=relerr1] {\tableUniformTwo};
        \addplot+ [line3, forget plot] table [x=ndof, y=relerr1] {\tableUniformThree};
        \addplot+ [line4, forget plot] table [x=ndof, y=relerr1] {\tableUniformFour};

        %
        %
        \drawslopetriangle[ST1]{0.5}{3e4}{2e-3}
    \end{loglogaxis}
\end{tikzpicture}

%% file: figures/plot_dls_lshape_error_scalingOne.tex
\begin{tikzpicture}[>=stealth]
    %
    %
    \colorlet{col1}{TUblue}
    \colorlet{col2}{TUgreen}
    \colorlet{col3}{TUmagenta}
    \colorlet{col4}{TUyellow}
    \colorlet{col5}{purple}
    \colorlet{col6}{green}
    \pgfplotsset{%
        linedefault/.style = {%
            mark = *,%
            mark size = 2pt,%
            every mark/.append style = {solid},%
            gray,%
            every mark/.append style = {fill = gray!60!white}%
        },%
        line1/.style = {%
            linedefault,%
            col1,%
            every mark/.append style = {fill = col1!60!white}%
        },%
        line2/.style = {%
            linedefault,%
            mark = triangle*,%
            mark size = 2.75pt,%
            col2,%
            every mark/.append style = {fill = col2!60!white}%
        },%
        line3/.style = {%
            linedefault,%
            mark = square*,%
            mark size = 1.66pt,%
            col3,%
            every mark/.append style = {fill = col3!60!white}%
        },%
        line4/.style = {%
            linedefault,%
            mark = pentagon*,%
            mark size = 2.2pt,%
            col4,%
            every mark/.append style = {fill = col4!60!white}%
        },%
        line5/.style = {%
            linedefault,%
            mark = diamond*,%
            mark size = 2.75pt,%
            col5,%
            every mark/.append style = {fill = col5!60!white}%
        },%
        line6/.style = {%
            linedefault,%
            mark = halfsquare*,%
            mark size = 1.66pt,%
            col6,%
            every mark/.append style = {fill = col6!60!white}%
        },%
        minorline/.style = {%
            dashed,%
            every mark/.append style = {fill = black!20!white}%
        },%
        majorline/.style = {%
            solid%
        }%
    }

    %
    %
    \begin{loglogaxis}[%
            width            = 0.36\textwidth,%
            xlabel           = ndof,%
            ylabel           = {relative energy error},%
            ymajorgrids      = true,%
            font             = \footnotesize,%
            grid style       = {densely dotted, semithick},%
            legend style     = {legend pos  = south west}%
        ]

        %
        %
        \pgfplotstableread[col sep=comma]{
ndof,nelem,res,err,res1,res2,res3,res4,err1,err2,err3,relerr1,relerr2,relerr3,cond
1.880000000000000000e+02,2.400000000000000000e+01,1.116001969258323401e+00,1.358536384763648686e+00,9.053403074152326369e-01,4.510999260266945754e-01,2.592194975286843052e-16,4.715168924826115271e-01,6.337020075340840108e-01,9.053403074152325258e-01,7.901909909276669808e-01,3.717194239555495217e-01,4.301227823712110154e-02,4.635133492877546724e-01,nan
7.360000000000000000e+02,9.600000000000000000e+01,5.986383594040479705e-01,7.225585233697768839e-01,4.588610891554643079e-01,2.582450137035929294e-01,4.528125879746569000e-16,2.848225747834292609e-01,3.512316058520802420e-01,4.588610891554643634e-01,4.337899026071750930e-01,2.056026071782401654e-01,2.180008516201056498e-02,2.539302655501662898e-01,nan
2.912000000000000000e+03,3.840000000000000000e+02,3.125477040689194075e-01,3.815165110278415650e-01,2.296147084108457748e-01,1.386715341407504243e-01,5.401355383630869495e-16,1.604162043535683202e-01,1.932602221906152828e-01,2.296147084108459413e-01,2.355470661823096667e-01,1.130372111200031726e-01,1.090878307291716365e-02,1.377706345721051162e-01,nan
1.158400000000000000e+04,1.536000000000000000e+03,1.633208985589179596e-01,2.049448109276436059e-01,1.148670168563833693e-01,7.350503187128673566e-02,8.185266552561036151e-16,8.986820702678353923e-02,1.088693796725255375e-01,1.148670168563833693e-01,1.302129107823122300e-01,6.365658733920447065e-02,5.457225609692605006e-03,7.613627957501903354e-02,nan
4.620800000000000000e+04,6.144000000000000000e+03,8.638620716412509504e-02,1.126721668011009397e-01,5.745261158953472030e-02,3.936402024631573127e-02,2.210188962489658346e-15,5.111015671853142567e-02,6.293811415334193970e-02,5.745261158953463704e-02,7.370894369083090536e-02,3.679554436348801694e-02,2.729520356316131151e-03,4.309250037193553678e-02,nan
1.845760000000000000e+05,2.457600000000000000e+04,4.658148753552105586e-02,6.356240493517854862e-02,2.873105472623013398e-02,2.151376315444394546e-02,3.859178430675532782e-13,2.969039356887615694e-02,3.723714885487078297e-02,2.873105472623023460e-02,4.275629264356498732e-02,2.176885897990025326e-02,1.364985794399540056e-03,2.499535366385511811e-02,nan
        }\tableUniformOne

        \pgfplotstableread[col sep=comma]{
ndof,nelem,res,err,res1,res2,res3,res4,err1,err2,err3,relerr1,relerr2,relerr3,cond
1.880000000000000000e+02,2.400000000000000000e+01,1.116001969258323401e+00,1.358536384763648686e+00,9.053403074152326369e-01,4.510999260266945754e-01,2.592194975286843052e-16,4.715168924826115271e-01,6.337020075340840108e-01,9.053403074152325258e-01,7.901909909276669808e-01,3.717194239555495217e-01,4.301227823712110154e-02,4.635133492877546724e-01,nan
5.350000000000000000e+02,7.000000000000000000e+01,7.153139218985220538e-01,8.602305021674170327e-01,5.759416472402235643e-01,2.741226535963929289e-01,4.175399183816779752e-16,3.237622532373635176e-01,4.034589088261352763e-01,5.759416472402234533e-01,4.954882891877958606e-01,2.362380011185393758e-01,2.736249456833235449e-02,2.901241252942862170e-01,nan
1.043000000000000000e+03,1.370000000000000000e+02,5.229567905673402350e-01,6.262300368572821796e-01,3.929295502606200419e-01,2.324517624385508963e-01,4.915178169484336282e-16,2.550614660659099586e-01,3.011380695254584627e-01,3.929295502606200974e-01,3.835183055301793487e-01,1.761717927689535057e-01,1.866772297529781244e-02,2.243658782545623620e-01,nan
2.136000000000000000e+03,2.820000000000000000e+02,3.626859114427763431e-01,4.335543945651375575e-01,2.607304122080996711e-01,1.681611192813299294e-01,5.351886081988628012e-16,1.878365259232670781e-01,2.055526183098776161e-01,2.607304122080996711e-01,2.788138918750124873e-01,1.202029842874894616e-01,1.238706070323271173e-02,1.630446848099108836e-01,nan
4.466000000000000000e+03,5.910000000000000000e+02,2.605708802073160024e-01,3.073468762214730532e-01,1.923612284254204163e-01,1.152435404234555955e-01,5.995253979342197505e-16,1.327149871014800808e-01,1.439341064355047495e-01,1.923612284254203608e-01,1.916826364759560142e-01,8.415469974534103836e-02,9.138903812973397239e-03,1.120720801935679578e-01,nan
8.593000000000000000e+03,1.140000000000000000e+03,1.859154161662533722e-01,2.205159248950244177e-01,1.319051679346172901e-01,8.696247167902126929e-02,7.535526456506134306e-16,9.799539356463150830e-02,1.053853126848044719e-01,1.319051679346173456e-01,1.418528662903024984e-01,6.161018181173214675e-02,6.266692338110522141e-03,8.292978082059654366e-02,nan
1.653400000000000000e+04,2.196000000000000000e+03,1.358319961443323343e-01,1.587920828218915559e-01,9.953972818486001017e-02,6.098894084042932179e-02,9.833485784101849062e-16,6.944438627132294706e-02,7.473033780794224867e-02,9.953972818486001017e-02,9.860093656203013091e-02,4.368703571151071746e-02,4.729040275081377274e-03,5.764168560089136945e-02,nan
3.241900000000000000e+04,4.311000000000000000e+03,9.727692397768453181e-02,1.142379527723223487e-01,6.837054044861433366e-02,4.630438773935324770e-02,1.443645950273783485e-15,5.142152091691028909e-02,5.441158385870758368e-02,6.837054044861445856e-02,7.358776117237537939e-02,3.180829166932362673e-02,3.248221039892199274e-03,4.301843108889501222e-02,nan
6.022400000000000000e+04,8.013000000000000000e+03,7.160847367189916712e-02,8.337517455539034117e-02,5.149093800725884629e-02,3.295714097462159459e-02,2.336554340775063019e-15,3.728650779478578925e-02,3.947981927370357302e-02,5.149093800725897813e-02,5.235882834241160910e-02,2.307923939041970143e-02,2.446286764748885085e-03,3.060809182379648982e-02,nan
1.211100000000000000e+05,1.612500000000000000e+04,5.100114055799274959e-02,5.930487347705259288e-02,3.617990532758547789e-02,2.427306836067204959e-02,9.125310206851333887e-14,2.651318428746774619e-02,2.810396755232736493e-02,3.617990532758551259e-02,3.765965316401487240e-02,1.642906785212764134e-02,1.718873785912032170e-03,2.201514771774479989e-02,nan
        }\tableAdaptiveOne

        \pgfplotstableread[col sep=comma]{
ndof,nelem,res,err,res1,res2,res3,res4,err1,err2,err3,relerr1,relerr2,relerr3,cond
1.880000000000000000e+02,2.400000000000000000e+01,3.824219241449398621e-01,1.332769092826115997e+00,2.965011670901880891e-01,1.710590778104971732e-01,1.365120637865488100e-16,1.705062400016396640e-01,8.809536110624445193e-01,2.965011670901880891e-01,9.551341517623289645e-01,5.167532451253596948e-01,1.408662642330038350e+00,5.602663593807049036e-01,nan
7.360000000000000000e+02,9.600000000000000000e+01,2.982545500763314483e-01,7.943421750058304864e-01,1.665630577047977623e-01,1.685276677306475068e-01,3.041176754106017304e-16,1.811379299278082777e-01,5.222880684669048978e-01,1.665630577047977623e-01,5.748490343891440313e-01,3.057349816636577011e-01,7.913263792954213782e-01,3.365029178327341164e-01,nan
2.912000000000000000e+03,3.840000000000000000e+02,1.929351722057153773e-01,3.832570343264242729e-01,6.831175870905535630e-02,1.192155708892397165e-01,4.341921508169276576e-16,1.354442023627350222e-01,2.493472577893737219e-01,6.831175870905532854e-02,2.829229630406152052e-01,1.458423171692905251e-01,3.245428667196586514e-01,1.654806268015966420e-01,nan
1.158400000000000000e+04,1.536000000000000000e+03,1.124315768153995726e-01,1.870253732926547841e-01,2.525128721694005693e-02,6.979322610974154906e-02,7.419011890696472953e-16,8.445198368469022943e-02,1.209042939917907983e-01,2.525128721694005693e-02,1.404386501056282821e-01,7.069347481655163801e-02,1.199665274238650853e-01,8.211533144709352861e-02,nan
4.620800000000000000e+04,6.144000000000000000e+03,6.385172244855187584e-02,1.002278144470027665e-01,9.706122911441709619e-03,3.862715691625507458e-02,1.629824256601680846e-15,4.990767863006551064e-02,6.510405688594773921e-02,9.706122911441707884e-03,7.558351502703537783e-02,3.806182065693078376e-02,4.611289084117534132e-02,4.418843204532160035e-02,nan
1.845760000000000000e+05,2.457600000000000000e+04,3.656732860201902802e-02,5.733951492201205813e-02,4.024362311756695190e-03,2.135410420602475504e-02,8.223640331687765926e-15,2.941047878954007946e-02,3.761375745010356320e-02,4.024362311756703864e-03,4.309094719265445239e-02,2.198902458474268073e-02,1.911937253751327948e-02,2.519099290880652195e-02,nan
        }\tableUniformTwo

        \pgfplotstableread[col sep=comma]{
ndof,nelem,res,err,res1,res2,res3,res4,err1,err2,err3,relerr1,relerr2,relerr3,cond
1.880000000000000000e+02,2.400000000000000000e+01,3.824219241449398621e-01,1.332769092826115997e+00,2.965011670901880891e-01,1.710590778104971732e-01,1.365120637865488100e-16,1.705062400016396640e-01,8.809536110624445193e-01,2.965011670901880891e-01,9.551341517623289645e-01,5.167532451253596948e-01,1.408662642330038350e+00,5.602663593807049036e-01,nan
4.820000000000000000e+02,6.300000000000000000e+01,3.281484143535136067e-01,9.832591185376652332e-01,2.099987714729588673e-01,1.680514232326457691e-01,2.202968013052563629e-16,1.879910023699416510e-01,6.547971385831415914e-01,2.099987714729587840e-01,7.028084500563752668e-01,3.836084598771605791e-01,9.976876446295803724e-01,4.117355608763829111e-01,nan
1.290000000000000000e+03,1.700000000000000000e+02,2.554768819036772887e-01,6.434309481675525300e-01,1.258125714196088596e-01,1.466439289064705720e-01,3.562051441838978151e-16,1.671382427061063325e-01,4.253738754441341308e-01,1.258125714196088873e-01,4.660811603382181145e-01,2.488785233721308721e-01,5.977240428906207548e-01,2.726956159106660449e-01,nan
2.904000000000000000e+03,3.840000000000000000e+02,2.012524275574248955e-01,4.044179229108319507e-01,7.306133709739864868e-02,1.218388475701061013e-01,4.537702224117844805e-16,1.425478019545909636e-01,2.613457620707737017e-01,7.306133709739867643e-02,2.998571160414546521e-01,1.528296788825469765e-01,3.471076541885493860e-01,1.753503343316198393e-01,nan
5.833000000000000000e+03,7.730000000000000000e+02,1.502916688376476073e-01,2.677799078580558878e-01,4.336384631105677512e-02,9.370373552835424558e-02,5.457167873306973422e-16,1.092097637720566700e-01,1.718820745939341754e-01,4.336384631105676818e-02,2.007042807686602592e-01,1.004925288000911715e-01,2.060176172820250562e-01,1.173437123277424193e-01,nan
1.022500000000000000e+04,1.357000000000000000e+03,1.153902121752471838e-01,1.818942161823756676e-01,2.594735500681710424e-02,7.318876368940849253e-02,6.772867036732545908e-16,8.535244989814110006e-02,1.148834950066590366e-01,2.594735500681713200e-02,1.386146573195128351e-01,6.716228893166006908e-02,1.232734799707480428e-01,8.103581514922181261e-02,nan
2.064000000000000000e+04,2.742000000000000000e+03,8.334387849485096200e-02,1.218721172165450201e-01,1.519864552670263484e-02,5.419810410415875496e-02,9.578767501298473827e-16,6.146355642356609444e-02,7.560150344270399891e-02,1.519864552670263658e-02,9.437280755704870183e-02,4.419613153554312335e-02,7.220735693308666381e-02,5.516970994275036122e-02,nan
3.687300000000000000e+04,4.904000000000000000e+03,6.246798611267741913e-02,8.727046113663000648e-02,9.712560638248711117e-03,4.039852958455190934e-02,1.397815807014267248e-15,4.664626739866851307e-02,5.371584548725558128e-02,9.712560638248721526e-03,6.809117047340039475e-02,3.140152037172077443e-02,4.614347584398725832e-02,3.980513119284490725e-02,nan
7.028500000000000000e+04,9.354000000000000000e+03,4.584839862450781206e-02,6.202794761706455218e-02,6.342774144187867183e-03,3.013478507726721420e-02,2.729473017271563242e-15,3.396674257173484812e-02,3.800484433520572719e-02,6.342774144187870652e-03,4.860933355640567610e-02,2.221697925249707661e-02,3.013393237957108159e-02,2.841618151557491131e-02,nan
1.241860000000000000e+05,1.653400000000000000e+04,3.426915537880588503e-02,4.554958354907551188e-02,4.433412560038840246e-03,2.271845867107684033e-02,5.482980727929356551e-15,2.527036800404386385e-02,2.745968975435425508e-02,4.433412560038842848e-03,3.607041520638917809e-02,1.605243060534853466e-02,2.106273237197589521e-02,2.108610265397687158e-02,nan
        }\tableAdaptiveTwo

        \pgfplotstableread[col sep=comma]{
ndof,nelem,res,err,res1,res2,res3,res4,err1,err2,err3,relerr1,relerr2,relerr3,cond
1.880000000000000000e+02,2.400000000000000000e+01,4.577498863445853688e-02,1.837506524651467243e+00,4.561002024051869636e-02,2.807250827418132123e-03,2.245946089358149711e-18,2.682346162642668813e-03,1.298396358638674419e+00,4.561002024051869636e-02,1.299429432506514992e+00,7.616184590881104333e-01,2.166909906603930125e+01,7.622244436336723661e-01,nan
7.360000000000000000e+02,9.600000000000000000e+01,4.556332549979515839e-02,1.814944889636574521e+00,4.499748174388054633e-02,5.047709822023327061e-03,9.044635440995738051e-18,5.075814502330655911e-03,1.282493303504765647e+00,4.499748174388054633e-02,1.283437223053804255e+00,7.507409996589637924e-01,2.137790623951251945e+01,7.512935476558276759e-01,nan
2.912000000000000000e+03,3.840000000000000000e+02,4.479510350267749963e-02,1.731328066933060184e+00,4.284260582668095868e-02,9.148218438424421395e-03,3.153905211600853264e-17,9.350001240540252792e-03,1.223419133653169011e+00,4.284260582668095868e-02,1.224298578732931464e+00,7.155734653073373419e-01,2.035412698412676846e+01,7.160878495816734857e-01,nan
1.158400000000000000e+04,1.536000000000000000e+03,4.220768153198283662e-02,1.484049239257810404e+00,3.647719095537036577e-02,1.466614849583035979e-02,1.268272085582048315e-16,1.535600887629934297e-02,1.048669323234047157e+00,3.647719095537035189e-02,1.049459008055835030e+00,6.131633207168806221e-01,1.732997566222007180e+01,6.136250542270722530e-01,nan
4.620800000000000000e+04,6.144000000000000000e+03,3.561914217222418577e-02,1.017761544560576992e+00,2.474053739656770223e-02,1.726719005445012667e-02,5.005087091103923446e-16,1.893338971335173371e-02,7.191372694380995245e-01,2.474053739656766060e-02,7.197694457955093883e-01,4.204296181575715718e-01,1.175400013712042124e+01,4.207992077698047839e-01,nan
1.845760000000000000e+05,2.457600000000000000e+04,2.606942158321831626e-02,5.608985517312458091e-01,1.402633955014033658e-02,1.403828740841092312e-02,2.294138623600582933e-15,1.690570989774133406e-02,3.962580513349956424e-01,1.402633955014039729e-02,3.967240903767064730e-01,2.316526883618602872e-01,6.663783984194272314e+00,2.319251350579720206e-01,nan
        }\tableUniformThree

        \pgfplotstableread[col sep=comma]{
ndof,nelem,res,err,res1,res2,res3,res4,err1,err2,err3,relerr1,relerr2,relerr3,cond
1.880000000000000000e+02,2.400000000000000000e+01,4.577498863445853688e-02,1.837506524651467243e+00,4.561002024051869636e-02,2.807250827418132123e-03,2.245946089358149711e-18,2.682346162642668813e-03,1.298396358638674419e+00,4.561002024051869636e-02,1.299429432506514992e+00,7.616184590881104333e-01,2.166909906603930125e+01,7.622244436336723661e-01,nan
4.440000000000000000e+02,5.800000000000000000e+01,4.565148882314879614e-02,1.825391571048322836e+00,4.524423307562039365e-02,4.282574216457766571e-03,6.077019619313374399e-18,4.321731510580914017e-03,1.289926706690196934e+00,4.524423307562039365e-02,1.290773581382515856e+00,7.556948097214886140e-01,2.149517925064248303e+01,7.561909455144376624e-01,nan
1.190000000000000000e+03,1.570000000000000000e+02,4.539070442026094804e-02,1.805574224526966143e+00,4.452177601785630573e-02,5.991989894904961415e-03,1.165715476633400755e-17,6.497966226650362526e-03,1.275949718279049572e+00,4.452177601785631267e-02,1.276741323902373493e+00,7.471090801098471124e-01,2.115191456850405771e+01,7.475725903411507556e-01,nan
3.177000000000000000e+03,4.210000000000000000e+02,4.464315090453255946e-02,1.744729682680493665e+00,4.252052627731765094e-02,9.138444381326093005e-03,4.078955990343851800e-17,1.007494923432448140e-02,1.232954828866864139e+00,4.252052627731765094e-02,1.233732572499145030e+00,7.212270428195595162e-01,2.020111100945323557e+01,7.216819903381942458e-01,nan
7.335000000000000000e+03,9.740000000000000000e+02,4.331622126262166250e-02,1.636010280479123669e+00,3.914608948037617137e-02,1.210187641525496162e-02,7.817954949444760751e-17,1.405073984161241997e-02,1.156152462731736552e+00,3.914608948037621300e-02,1.156852931202197343e+00,6.763008701946635615e-01,1.859794711594728156e+01,6.767106149743419286e-01,nan
1.651200000000000000e+04,2.196000000000000000e+03,4.093072594791748176e-02,1.441549877558886639e+00,3.405907995280338385e-02,1.469199070886192428e-02,2.023042242912069307e-16,1.730458922953061196e-02,1.018709804795998375e+00,3.405907995280334222e-02,1.019380381494223275e+00,5.956707857566696518e-01,1.618115359431267919e+01,5.960628924654312044e-01,nan
3.731800000000000000e+04,4.966000000000000000e+03,3.711271250097432362e-02,1.174580245896212016e+00,2.755366115098983187e-02,1.590314626700776654e-02,4.579836353897164518e-16,1.911123034151821290e-02,8.300098130217743453e-01,2.755366115098982147e-02,8.306402711729917732e-01,4.852571692103366807e-01,1.309048918028886987e+01,4.856257604401819639e-01,nan
7.468600000000000000e+04,9.944000000000000000e+03,3.208627224805678940e-02,8.726039303828811988e-01,2.050551163534317448e-02,1.584572800565312717e-02,9.552414155213850516e-16,1.891998317441337049e-02,6.165704331467142252e-01,2.050551163534316407e-02,6.171369736788570037e-01,3.604498157556630944e-01,9.741978621110328973e+00,3.607810178689185299e-01,nan
1.432050000000000000e+05,1.907100000000000000e+04,2.731282645022846053e-02,5.994611925486680315e-01,1.425940660711556146e-02,1.481075444976859973e-02,1.842501733548478920e-15,1.798058298670841826e-02,4.234991016148065279e-01,1.425940660711560830e-02,4.240270057750283317e-01,2.475734349624148234e-01,6.774511983389375658e+00,2.478820425740454469e-01,nan
        }\tableAdaptiveThree

        \pgfplotstableread[col sep=comma]{
ndof,nelem,res,err,res1,res2,res3,res4,err1,err2,err3,relerr1,relerr2,relerr3,cond
1.880000000000000000e+02,2.400000000000000000e+01,4.587751266962628374e-03,1.846411939302358585e+00,4.587584853438215517e-03,2.826338085077509453e-05,1.902789623034838220e-20,2.698298212484678455e-05,1.305601186203874553e+00,4.587584853438216384e-03,1.305611560245444780e+00,7.658446952691458920e-01,2.179539279719547835e+02,7.658507805153826897e-01,nan
7.360000000000000000e+02,9.600000000000000000e+01,4.587555557380855117e-03,1.848064312406804133e+00,4.586973169971440834e-03,5.162869528989567317e-05,7.105871516398743209e-20,5.174560697088573664e-05,1.306770013020325694e+00,4.586973169971440834e-03,1.306779551240234927e+00,7.649520065471367980e-01,2.179230449138295000e+02,7.649575899936390533e-01,nan
2.912000000000000000e+03,3.840000000000000000e+02,4.586750187581442248e-03,1.847889057347600561e+00,4.584583559544845935e-03,9.925077916757581013e-05,3.028569406215499728e-19,1.001007076105767623e-04,1.306646314883771698e+00,4.584583559544846802e-03,1.306655408919753514e+00,7.642527452390731613e-01,2.178093375501495927e+02,7.642580643081167135e-01,nan
1.158400000000000000e+04,1.536000000000000000e+03,4.583654773497868561e-03,1.844662981784001454e+00,4.575404532949241092e-03,1.938213018909195435e-04,1.358570886707254163e-18,1.949301039107084667e-04,1.304365264805805946e+00,4.575404532949241959e-03,1.304374117347351536e+00,7.626702903157017710e-01,2.173732327630013970e+02,7.626754664504953318e-01,nan
4.620800000000000000e+04,6.144000000000000000e+03,4.571674185598987965e-03,1.831409285934217257e+00,4.540187649075093368e-03,3.786951611774312838e-04,9.447763926941196247e-18,3.788019870050088724e-04,1.294993607342450170e+00,4.540187649075095103e-03,1.295002284262224412e+00,7.570928263485685594e-01,2.157001094777540970e+02,7.570978991409580416e-01,nan
1.845760000000000000e+05,2.457600000000000000e+04,4.525711517644245147e-03,1.782099099032659506e+00,4.409315721621172995e-03,7.216267285779588746e-04,7.234257652693278874e-17,7.205931394372311338e-04,1.260126275532391649e+00,4.409315721621148709e-03,1.260134725503883057e+00,7.366705570247588097e-01,2.094825052677455233e+02,7.366754968830304673e-01,nan
        }\tableUniformFour

        \pgfplotstableread[col sep=comma]{
ndof,nelem,res,err,res1,res2,res3,res4,err1,err2,err3,relerr1,relerr2,relerr3,cond
1.880000000000000000e+02,2.400000000000000000e+01,4.587751266962628374e-03,1.846411939302358585e+00,4.587584853438215517e-03,2.826338085077509453e-05,1.902789623034838220e-20,2.698298212484678455e-05,1.305601186203874553e+00,4.587584853438216384e-03,1.305611560245444780e+00,7.658446952691458920e-01,2.179539279719547835e+02,7.658507805153826897e-01,nan
4.440000000000000000e+02,5.800000000000000000e+01,4.587640395916463276e-03,1.847442733864194997e+00,4.587222973652636988e-03,4.364481915323519017e-05,4.722151007694971035e-20,4.387393302520870047e-05,1.306330982099157767e+00,4.587222973652636988e-03,1.306339533773991812e+00,7.653051431764816837e-01,2.179353552452201086e+02,7.653101531171651040e-01,nan
1.182000000000000000e+03,1.560000000000000000e+02,4.587370343224993902e-03,1.847721006721452452e+00,4.586451349645368283e-03,6.301319898793535633e-05,1.145669312063715579e-19,6.678338108903627787e-05,1.306527971239690666e+00,4.586451349645367416e-03,1.306536085805536818e+00,7.650136182842727228e-01,2.178983764739970752e+02,7.650183696202681460e-01,nan
3.244000000000000000e+03,4.300000000000000000e+02,4.586469180280322057e-03,1.847882944251481918e+00,4.583742665344172143e-03,1.076804999725816453e-04,3.973934577805861268e-19,1.157913187178011234e-04,1.306642435289810722e+00,4.583742665344171276e-03,1.306650646216736877e+00,7.643312127604837691e-01,2.177693987545484049e+02,7.643360158095551249e-01,nan
7.136000000000000000e+03,9.480000000000000000e+02,4.584892377795935295e-03,1.846607237928051370e+00,4.578359306448622282e-03,1.623594182957200239e-04,8.612486292137674880e-19,1.830398763035895543e-04,1.305740670046821261e+00,4.578359306448626619e-03,1.305748303609852057e+00,7.638037205693070097e-01,2.175136402938013873e+02,7.638081858846470507e-01,nan
1.626600000000000000e+04,2.164000000000000000e+03,4.582246390204705241e-03,1.843574565172628299e+00,4.569033069591245215e-03,2.275526032117964367e-04,3.360416357613630975e-18,2.629422033397989418e-04,1.303596393429920619e+00,4.569033069591247817e-03,1.303603752804033178e+00,7.625494221426482033e-01,2.170705586797046180e+02,7.625537270690081471e-01,nan
3.694100000000000000e+04,4.918000000000000000e+03,4.574707372849125611e-03,1.838449784167014212e+00,4.545396657116004130e-03,3.432033212396142690e-04,7.257609503322005539e-18,3.866888373988873690e-04,1.299972446633564038e+00,4.545396657116013671e-03,1.299980225336722528e+00,7.601336535735209177e-01,2.159475874577181287e+02,7.601382020191954636e-01,nan
8.094300000000000000e+04,1.078300000000000000e+04,4.562647930356713080e-03,1.825789798462336133e+00,4.506602071079534796e-03,4.721010564974645286e-04,2.680471561650149116e-17,5.342419880391544766e-04,1.291020875496719622e+00,4.506602071079533929e-03,1.291027953895110159e+00,7.548994009068346056e-01,2.141044925496049984e+02,7.549035398629164728e-01,nan
1.795090000000000000e+05,2.391900000000000000e+04,4.525372672841480176e-03,1.801745173338613837e+00,4.402569224725668036e-03,7.051790942445263595e-04,7.056882941610768440e-17,7.740184071708765658e-04,1.274018845377357989e+00,4.402569224725649821e-03,1.274026007841502661e+00,7.448427333352630741e-01,2.091619857221740233e+02,7.448469208002995723e-01,nan
        }\tableAdaptiveFour

        %
        %
        \addplot+ [line1, minorline, forget plot] table [x=ndof, y=relerr1] {\tableUniformOne};
        \label{leg:scal:unif:len1}
        \addplot+ [line2, minorline, forget plot] table [x=ndof, y=relerr1] {\tableUniformTwo};
        \label{leg:scal:unif:len10}
        \addplot+ [line3, minorline, forget plot] table [x=ndof, y=relerr1] {\tableUniformThree};
        \label{leg:scal:unif:len100}
        \addplot+ [line4, minorline, forget plot] table [x=ndof, y=relerr1] {\tableUniformFour};
        \label{leg:scal:unif:len1000}

        \addplot+ [line1, majorline, forget plot] table [x=ndof, y=relerr1] {\tableAdaptiveOne};
        \label{leg:scal:adap:len1}
        \addplot+ [line2, majorline, forget plot] table [x=ndof, y=relerr1] {\tableAdaptiveTwo};
        \label{leg:scal:adap:len10}
        \addplot+ [line3, majorline, forget plot] table [x=ndof, y=relerr1] {\tableAdaptiveThree};
        \label{leg:scal:adap:len100}
        \addplot+ [line4, majorline, forget plot] table [x=ndof, y=relerr1] {\tableAdaptiveFour};
        \label{leg:scal:adap:len1000}

        %
        %
        \drawslopetriangle[ST1]{0.5}{2e3}{3e-2}
        \drawswappedslopetriangle[ST2]{0.33}{1e5}{8e-2}
    \end{loglogaxis}
\end{tikzpicture}

%% file: figures/plot_dls_lshape_error_scalingDiam.tex
\begin{tikzpicture}[>=stealth]
    %
    %
    \colorlet{col1}{TUblue}
    \colorlet{col2}{TUgreen}
    \colorlet{col3}{TUmagenta}
    \colorlet{col4}{TUyellow}
    \colorlet{col5}{purple}
    \colorlet{col6}{green}
    \pgfplotsset{%
        linedefault/.style = {%
            mark = *,%
            mark size = 2pt,%
            every mark/.append style = {solid},%
            gray,%
            every mark/.append style = {fill = gray!60!white}%
        },%
        line1/.style = {%
            linedefault,%
            col1,%
            every mark/.append style = {fill = col1!60!white}%
        },%
        line2/.style = {%
            linedefault,%
            mark = triangle*,%
            mark size = 2.75pt,%
            col2,%
            every mark/.append style = {fill = col2!60!white}%
        },%
        line3/.style = {%
            linedefault,%
            mark = square*,%
            mark size = 1.66pt,%
            col3,%
            every mark/.append style = {fill = col3!60!white}%
        },%
        line4/.style = {%
            linedefault,%
            mark = pentagon*,%
            mark size = 2.2pt,%
            col4,%
            every mark/.append style = {fill = col4!60!white}%
        },%
        line5/.style = {%
            linedefault,%
            mark = diamond*,%
            mark size = 2.75pt,%
            col5,%
            every mark/.append style = {fill = col5!60!white}%
        },%
        line6/.style = {%
            linedefault,%
            mark = halfsquare*,%
            mark size = 1.66pt,%
            col6,%
            every mark/.append style = {fill = col6!60!white}%
        },%
        minorline/.style = {%
            dashed,%
            every mark/.append style = {fill = black!20!white}%
        },%
        majorline/.style = {%
            solid%
        }%
    }

    %
    %
    \begin{loglogaxis}[%
            width            = 0.36\textwidth,%
            xlabel           = ndof,%
            ylabel           = {relative energy error},%
            ymajorgrids      = true,%
            font             = \footnotesize,%
            grid style       = {densely dotted, semithick},%
            legend style     = {legend pos  = south west}%
        ]

        %
        %
        \pgfplotstableread[col sep=comma]{
ndof,nelem,res,err,res1,res2,res3,res4,err1,err2,err3,relerr1,relerr2,relerr3,cond
1.880000000000000000e+02,2.400000000000000000e+01,2.637678459993161528e+00,2.750392879609750985e+00,2.553030859086822701e+00,4.579952767477925013e-01,4.005459716866321505e-16,4.791882900687461366e-01,6.418937525993497761e-01,2.553030859086822701e+00,7.966597988679622322e-01,3.765245700978197552e-01,1.516165699786879185e-02,4.673078481731259837e-01,nan
7.360000000000000000e+02,9.600000000000000000e+01,1.352495157812587445e+00,1.411952400623623038e+00,1.296119271831964959e+00,2.595052255562113053e-01,4.391219064039562401e-16,2.863128062125283035e-01,3.527774664616804556e-01,1.296119271831965181e+00,4.350085906846744011e-01,2.065075171190665082e-01,7.697185263100086473e-03,2.546436564919155043e-01,nan
2.912000000000000000e+03,3.840000000000000000e+02,6.829860147169556361e-01,7.172309939193187711e-01,6.491247111150968419e-01,1.388780734094493263e-01,4.468094589450129788e-16,1.606856759165892767e-01,1.935544514767829904e-01,6.491247111150967308e-01,2.357839739502193832e-01,1.132093047746685621e-01,3.854914559825585278e-03,1.379092010762394682e-01,nan
1.158400000000000000e+04,1.536000000000000000e+03,3.449751470007734189e-01,3.665358620670863266e-01,3.248292525244516615e-01,7.354035969115059468e-02,9.608551864213432281e-16,8.992009924987545699e-02,1.089308190752326932e-01,3.248292525244518281e-01,1.302634697114541107e-01,6.369251133102164786e-02,1.929042126459683607e-03,7.616584168787737807e-02,nan
4.620800000000000000e+04,6.144000000000000000e+03,1.748290156893090863e-01,1.892074091304264682e-01,1.624861831661684464e-01,3.937073664767414688e-02,2.300181377963522425e-15,5.112113166590244451e-02,6.295201802064741281e-02,1.624861831661685296e-01,7.372060347736586472e-02,3.680367298909342083e-02,9.649460057202964699e-04,4.309931703936333014e-02,nan
1.845760000000000000e+05,2.457600000000000000e+04,8.915023094665179559e-02,9.908776620677903213e-02,8.126002673082503847e-02,2.151519220761478196e-02,7.096806902358540885e-12,2.969289877759232615e-02,3.724045781732550353e-02,8.126002673082513561e-02,4.275911328723839416e-02,2.177079340128485960e-02,4.825735741213454922e-04,2.499700261379481786e-02,nan
        }\tableUniformOne

        \pgfplotstableread[col sep=comma]{
ndof,nelem,res,err,res1,res2,res3,res4,err1,err2,err3,relerr1,relerr2,relerr3,cond
1.880000000000000000e+02,2.400000000000000000e+01,2.637678459993161528e+00,2.750392879609750985e+00,2.553030859086822701e+00,4.579952767477925013e-01,4.005459716866321505e-16,4.791882900687461366e-01,6.418937525993497761e-01,2.553030859086822701e+00,7.966597988679622322e-01,3.765245700978197552e-01,1.516165699786879185e-02,4.673078481731259837e-01,nan
5.200000000000000000e+02,6.800000000000000000e+01,1.689060681479619408e+00,1.757803092581543769e+00,1.630947981683745107e+00,2.837324446315512283e-01,4.923442844522535891e-16,3.353066754860668519e-01,4.078447169276278039e-01,1.630947981683745551e+00,5.133644711169433883e-01,2.388060315636147679e-01,9.685617784371764252e-03,3.005911980832293717e-01,nan
9.000000000000000000e+02,1.180000000000000000e+02,1.263430120662887735e+00,1.325017645518503207e+00,1.200179888991718569e+00,2.640082617615709948e-01,5.850588788134736119e-16,2.934681270252340046e-01,3.353768321888945247e-01,1.200179888991718347e+00,4.502914338805108496e-01,1.962886822822894617e-01,7.127434533727455661e-03,2.635456707684143418e-01,nan
1.723000000000000000e+03,2.270000000000000000e+02,8.817818537659222189e-01,9.236947205056689203e-01,8.367568040791399930e-01,1.828630788291901732e-01,3.335624343948499099e-14,2.096148440813274949e-01,2.343186791376746791e-01,8.367568040791401041e-01,3.132806156219768434e-01,1.370593953144042532e-01,4.969193147531007597e-03,1.832463886314588952e-01,nan
3.410000000000000000e+03,4.510000000000000000e+02,6.222732026803593497e-01,6.534537770504526000e-01,5.853447576434641642e-01,1.397300466056134727e-01,5.292392838745792098e-16,1.583381430640732057e-01,1.759368619942915102e-01,5.853447576434642752e-01,2.311267488466410835e-01,1.028735084769498570e-01,3.476148437939181966e-03,1.351440470587469134e-01,nan
7.179000000000000000e+03,9.520000000000000000e+02,4.420308523935467804e-01,4.665915239805026316e-01,4.137015820784105125e-01,1.028247852007110191e-01,2.260660466682103045e-13,1.169159484502152041e-01,1.318998757006352718e-01,4.137015820784106235e-01,1.707661384069598942e-01,7.712062092400341629e-02,2.456822382783710110e-03,9.984536040526098555e-02,nan
1.302300000000000000e+04,1.729000000000000000e+03,3.192232477880173058e-01,3.371310205387265668e-01,2.953117921888241404e-01,8.267433862206402528e-02,1.906734864943826890e-13,8.865315028784324924e-02,9.612511436878681814e-02,2.953117921888243069e-01,1.311801539536205929e-01,5.619708949005023518e-02,1.753748710882951532e-03,7.669112176831861805e-02,nan
2.668700000000000000e+04,3.548000000000000000e+03,2.298334059295063792e-01,2.414539837427214231e-01,2.147945028229316933e-01,5.421567551532285173e-02,4.429211180971685560e-15,6.121581957814559161e-02,6.686418488577217101e-02,2.147945028229315823e-01,8.770706129976821308e-02,3.908867896396239772e-02,1.275585981166355064e-03,5.127338601189946471e-02,nan
4.675600000000000000e+04,6.220000000000000000e+03,1.682210136073811757e-01,1.767918070850687529e-01,1.559065505904853599e-01,4.322707042206907352e-02,1.351834311583489381e-14,4.607469245471563363e-02,4.973900021825260892e-02,1.559065505904853044e-01,6.689187090124840140e-02,2.907680912099311690e-02,9.258719738118348667e-04,3.910416681893739976e-02,nan
9.396200000000000000e+04,1.250800000000000000e+04,1.230120986083583084e-01,1.288813332074744689e-01,1.145281042121732012e-01,3.020699082420572057e-02,3.024979650444091077e-12,3.320884603451938755e-02,3.548885460090002442e-02,1.145281042121734233e-01,4.726788120931361969e-02,2.074620205526092262e-02,6.801405168773881973e-04,2.763202772589966327e-02,nan
1.621760000000000000e+05,2.159700000000000000e+04,9.003791519252228592e-02,9.443995829810095111e-02,8.337978974709539526e-02,2.324831684122019917e-02,3.136872724745661449e-12,2.478210237240253619e-02,2.654376681777643326e-02,8.337978974709515934e-02,3.552667769720047053e-02,1.551700693786968975e-02,4.951620712294489343e-04,2.076825448668930332e-02,nan
        }\tableAdaptiveOne

        \pgfplotstableread[col sep=comma]{
ndof,nelem,res,err,res1,res2,res3,res4,err1,err2,err3,relerr1,relerr2,relerr3,cond
1.880000000000000000e+02,2.400000000000000000e+01,2.637678459993160196e+00,2.750392879609750985e+00,2.553030859086821813e+00,4.579952767477920017e-01,5.020370520385580885e-16,4.791882900687464697e-01,6.418937525993527737e-01,2.553030859086821813e+00,7.966597988679640086e-01,3.765245700978214760e-01,1.516165699786878318e-02,4.673078481731269274e-01,nan
7.360000000000000000e+02,9.600000000000000000e+01,1.352495157812586779e+00,1.411952400623619708e+00,1.296119271831964070e+00,2.595052255562130816e-01,5.185660576815422199e-16,2.863128062125276374e-01,3.527774664616750155e-01,1.296119271831964070e+00,4.350085906846711814e-01,2.065075171190632330e-01,7.697185263100078667e-03,2.546436564919135059e-01,nan
2.912000000000000000e+03,3.840000000000000000e+02,6.829860147169559692e-01,7.172309939193495243e-01,6.491247111150963978e-01,1.388780734094345326e-01,6.032160903052297910e-16,1.606856759166056248e-01,1.935544514768442192e-01,6.491247111150961757e-01,2.357839739502642085e-01,1.132093047747042835e-01,3.854914559825582676e-03,1.379092010762655585e-01,nan
1.158400000000000000e+04,1.536000000000000000e+03,3.449751470004121523e-01,3.665358620618785479e-01,3.248292525244454443e-01,7.354035969338121315e-02,2.533418932519373093e-15,8.992009924668746546e-02,1.089308190664071696e-01,3.248292525244452222e-01,1.302634697041971656e-01,6.369251132586128961e-02,1.929042126459646744e-03,7.616584168363416119e-02,nan
4.620800000000000000e+04,6.144000000000000000e+03,1.748290157115283405e-01,1.892074090638827533e-01,1.624861831663068079e-01,3.937073664115061516e-02,2.780036758991093740e-15,5.112113174647426456e-02,6.295201788210842697e-02,1.624861831663072242e-01,7.372060342457474602e-02,3.680367290809929970e-02,9.649460057211210057e-04,4.309931700850003905e-02,nan
1.845760000000000000e+05,2.457600000000000000e+04,8.915023095980512124e-02,9.908776703203080183e-02,8.126002673089063877e-02,2.151519197001579287e-02,2.522700851801951401e-12,2.969289898906642133e-02,3.724045900595646802e-02,8.126002673089054162e-02,4.275911416428832684e-02,2.177079409615920169e-02,4.825735741217341787e-04,2.499700312651870551e-02,nan
        }\tableUniformTwo

        \pgfplotstableread[col sep=comma]{
ndof,nelem,res,err,res1,res2,res3,res4,err1,err2,err3,relerr1,relerr2,relerr3,cond
1.880000000000000000e+02,2.400000000000000000e+01,2.637678459993160196e+00,2.750392879609750985e+00,2.553030859086821813e+00,4.579952767477920017e-01,5.020370520385580885e-16,4.791882900687464697e-01,6.418937525993527737e-01,2.553030859086821813e+00,7.966597988679640086e-01,3.765245700978214760e-01,1.516165699786878318e-02,4.673078481731269274e-01,nan
5.200000000000000000e+02,6.800000000000000000e+01,1.689824826457941631e+00,1.759701883427533753e+00,1.630948953464978812e+00,2.806935971607616986e-01,5.058609746318669376e-16,3.416497619171644518e-01,4.083758600305483211e-01,1.630948953464979034e+00,5.194086894481163563e-01,2.391170327151225516e-01,9.685623555433564569e-03,3.041302798309221678e-01,nan
9.000000000000000000e+02,1.180000000000000000e+02,1.263891925840282537e+00,1.325758663945737492e+00,1.200179818558539946e+00,2.618361895077994661e-01,5.387326554093745299e-16,2.973768874561096287e-01,3.347483001861777385e-01,1.200179818558539946e+00,4.529326811728668911e-01,1.959208163275056869e-01,7.127434115450265867e-03,2.650915347066432060e-01,nan
1.723000000000000000e+03,2.270000000000000000e+02,8.906942016673927709e-01,9.297990173143562131e-01,8.446794854811457975e-01,1.911801237931389508e-01,3.098035108199653518e-14,2.080934596705551198e-01,2.320120116081652861e-01,8.446794854811459086e-01,3.117903235759276193e-01,1.356939601896667846e-01,5.016242846073174812e-03,1.823529026000923203e-01,nan
3.372000000000000000e+03,4.460000000000000000e+02,6.238004355393209632e-01,6.552380315584509507e-01,5.866868935644197203e-01,1.391928537687150136e-01,6.691871871236574840e-16,1.598462441175860227e-01,1.765815731282333045e-01,5.866868935644196092e-01,2.322806814005059506e-01,1.032504827635975370e-01,3.484118892378558659e-03,1.358187724029570054e-01,nan
7.201000000000000000e+03,9.550000000000000000e+02,4.415806363266800427e-01,4.664179319108768818e-01,4.129728456291434924e-01,1.031067878911342828e-01,9.940371101108820958e-16,1.175409607005828555e-01,1.326863938481806782e-01,4.129728456291441585e-01,1.714451424466494900e-01,7.758049071226119642e-02,2.452494683550029888e-03,1.002423677024738280e-01,nan
1.296300000000000000e+04,1.721000000000000000e+03,3.196068135654367426e-01,3.375640632483467019e-01,2.956555428319420087e-01,8.255466523858523420e-02,3.096899969564585994e-14,8.900023886230117032e-02,9.617629515576284205e-02,2.956555428319417866e-01,1.314816224405342682e-01,5.622701102808692203e-02,1.755790120210913617e-03,7.686736760842785565e-02,nan
2.668700000000000000e+04,3.548000000000000000e+03,2.298141948031776627e-01,2.414091800387630493e-01,2.148352820954512465e-01,5.398701923219731413e-02,6.278657352319402989e-14,6.120267522276404831e-02,6.671443779333820212e-02,2.148352820954514408e-01,8.759781711820521788e-02,3.900113709035792076e-02,1.275828154348928939e-03,5.120952206516789229e-02,nan
4.654600000000000000e+04,6.192000000000000000e+03,1.684157028019938884e-01,1.770292190340424232e-01,1.560207223595257708e-01,4.339477538105462884e-02,3.102263960297704741e-13,4.624150315819253176e-02,4.984245239935004151e-02,1.560207223595259374e-01,6.717595198456974281e-02,2.913728599647756601e-02,9.265499981844791191e-04,3.927023713394924009e-02,nan
9.390200000000000000e+04,1.250000000000000000e+04,1.231086998170064228e-01,1.289539888947086044e-01,1.146686865395184912e-01,3.011900222586957193e-02,4.310083634969340504e-13,3.316157419142688051e-02,3.538399744873522990e-02,1.146686865395180888e-01,4.720377425276788691e-02,2.068490428472943496e-02,6.809753838948912497e-04,2.759455185104814010e-02,nan
1.615080000000000000e+05,2.150800000000000000e+04,9.008420606167649725e-02,9.451191856274919478e-02,8.339523705931979614e-02,2.328674145298595052e-02,2.979708756081232458e-12,2.486214572566559575e-02,2.659882974568377953e-02,8.339523705931979614e-02,3.564041866184678581e-02,1.554919573157866747e-02,4.952538071661705685e-04,2.083474540147361304e-02,nan
        }\tableAdaptiveTwo

        \pgfplotstableread[col sep=comma]{
ndof,nelem,res,err,res1,res2,res3,res4,err1,err2,err3,relerr1,relerr2,relerr3,cond
1.880000000000000000e+02,2.400000000000000000e+01,2.637678459993160640e+00,2.750392879609751429e+00,2.553030859086822257e+00,4.579952767477917797e-01,5.004525578964334788e-16,4.791882900687466917e-01,6.418937525993521076e-01,2.553030859086822257e+00,7.966597988679637865e-01,3.765245700978211429e-01,1.516165699786878665e-02,4.673078481731269274e-01,nan
7.360000000000000000e+02,9.600000000000000000e+01,1.352495157812586779e+00,1.411952400623622150e+00,1.296119271831964292e+00,2.595052255562119714e-01,5.621200894023069475e-16,2.863128062125279150e-01,3.527774664616803446e-01,1.296119271831964292e+00,4.350085906846742900e-01,2.065075171190664527e-01,7.697185263100082137e-03,2.546436564919153933e-01,nan
2.912000000000000000e+03,3.840000000000000000e+02,6.829860147169548590e-01,7.172309939193241002e-01,6.491247111150965088e-01,1.388780734094461622e-01,6.504249206369788802e-16,1.606856759165900261e-01,1.935544514767936763e-01,6.491247111150962867e-01,2.357839739502278487e-01,1.132093047746747377e-01,3.854914559825581375e-03,1.379092010762443254e-01,nan
1.158400000000000000e+04,1.536000000000000000e+03,3.449751470004120413e-01,3.665358620618432983e-01,3.248292525244456108e-01,7.354035969339363377e-02,2.648444129092359662e-14,8.992009924667634935e-02,1.089308190663421799e-01,3.248292525244457218e-01,1.302634697041513689e-01,6.369251132582334773e-02,1.929042126459646093e-03,7.616584168360747420e-02,nan
4.620800000000000000e+04,6.144000000000000000e+03,1.748290156892542135e-01,1.892074091293416693e-01,1.624861831661733869e-01,3.937073664806076123e-02,4.399615804803324763e-14,5.112113166540122738e-02,6.295201801894122207e-02,1.624861831661734979e-01,7.372060347602767127e-02,3.680367298809589932e-02,9.649460057203239002e-04,4.309931703858094210e-02,nan
1.845760000000000000e+05,2.457600000000000000e+04,8.915023094092536238e-02,9.908776622904345854e-02,8.126002673077765970e-02,2.151519221210334773e-02,2.449554446858866788e-12,2.969289875727648687e-02,3.724045785912266648e-02,8.126002673077761806e-02,4.275911330252044062e-02,2.177079342571947712e-02,4.825735741210615938e-04,2.499700262272867846e-02,nan
        }\tableUniformThree

        \pgfplotstableread[col sep=comma]{
ndof,nelem,res,err,res1,res2,res3,res4,err1,err2,err3,relerr1,relerr2,relerr3,cond
1.880000000000000000e+02,2.400000000000000000e+01,2.637678459993160640e+00,2.750392879609751429e+00,2.553030859086822257e+00,4.579952767477917797e-01,5.004525578964334788e-16,4.791882900687466917e-01,6.418937525993521076e-01,2.553030859086822257e+00,7.966597988679637865e-01,3.765245700978211429e-01,1.516165699786878665e-02,4.673078481731269274e-01,nan
5.200000000000000000e+02,6.800000000000000000e+01,1.688978641714293571e+00,1.757702758526648923e+00,1.630948078058457318e+00,2.841251316104520530e-01,5.137195209832153088e-16,3.345596037058675698e-01,4.080274928491081088e-01,1.630948078058457318e+00,5.128751496113709285e-01,2.389130526691053968e-01,9.685618356706766560e-03,3.003046848048798734e-01,nan
9.000000000000000000e+02,1.180000000000000000e+02,1.263329595463807520e+00,1.324957713430328754e+00,1.200180001947346309e+00,2.644275273484393285e-01,8.334221645407579645e-16,2.926563043585515556e-01,3.357114340032527622e-01,1.200180001947346531e+00,4.498652446940882288e-01,1.964845173636842746e-01,7.127435204530103464e-03,2.632962316128833469e-01,nan
1.723000000000000000e+03,2.270000000000000000e+02,8.817263207959226134e-01,9.235619555830786309e-01,8.367568484937282225e-01,1.831099970947996647e-01,6.517862351909793540e-16,2.091650311287733821e-01,2.343056824086197087e-01,8.367568484937283335e-01,3.128985610230816672e-01,1.370517931726035754e-01,4.969193411293031104e-03,1.830229144616044579e-01,nan
3.410000000000000000e+03,4.510000000000000000e+02,6.222259194082589984e-01,6.533807357177185970e-01,5.853447775479007076e-01,1.399314618475997474e-01,6.505067464500867665e-16,1.579739604050163504e-01,1.759931651479205161e-01,5.853447775479007076e-01,2.308772033476729213e-01,1.029064299630243762e-01,3.476148556144349842e-03,1.349981332308402282e-01,nan
7.179000000000000000e+03,9.520000000000000000e+02,4.419989282280771259e-01,4.665459814626957691e-01,4.137015887080490173e-01,1.029395419231413644e-01,7.275470869615779686e-16,1.166940391102257585e-01,1.319384582346428414e-01,4.137015887080491838e-01,1.706118154117456931e-01,7.714317976997411841e-02,2.456822422154707457e-03,9.975512919654641597e-02,nan
1.302300000000000000e+04,1.729000000000000000e+03,3.192300289154831883e-01,3.371580470795030871e-01,2.953117924376065240e-01,8.265645498220142162e-02,4.502087848028790887e-14,8.869423351291533253e-02,9.615016135403124664e-02,2.953117924376066350e-01,1.312312479150683775e-01,5.621173256933791051e-02,1.753748712360379955e-03,7.672099254603098606e-02,nan
2.667200000000000000e+04,3.546000000000000000e+03,2.297695165385952032e-01,2.413777765719066770e-01,2.147981952267574157e-01,5.421997038051509094e-02,2.807833977053529018e-13,6.095868165933100713e-02,6.683600802852537881e-02,2.147981952267577765e-01,8.750951023137908125e-02,3.907220682526793554e-02,1.275607909002009297e-03,5.115789802226138050e-02,nan
4.669600000000000000e+04,6.212000000000000000e+03,1.682145124978070561e-01,1.767957191535490802e-01,1.558629110481047075e-01,4.336454938412211135e-02,2.478426490300356920e-12,4.606941535549930605e-02,4.979406839246039440e-02,1.558629110481051516e-01,6.696290035204591207e-02,2.910900129982781775e-02,9.256128145328520020e-04,3.914568976418661417e-02,nan
9.395400000000000000e+04,1.250700000000000000e+04,1.230103792380004141e-01,1.289036483592094773e-01,1.144872685712351734e-01,3.024112170034352204e-02,2.246080445151376548e-12,3.331205928308185460e-02,3.547318831265362099e-02,1.144872685712352289e-01,4.743910629794360900e-02,2.073704380022357158e-02,6.798980089434098839e-04,2.773212310304249209e-02,nan
1.620630000000000000e+05,2.158200000000000000e+04,9.004534625957121008e-02,9.443147570354921971e-02,8.342756733644610123e-02,2.314886709381501484e-02,2.773672047735868335e-12,2.474136907290530596e-02,2.645487915477457916e-02,8.342756733644629552e-02,3.545820047293711530e-02,1.546504481459757066e-02,4.954458048557061275e-04,2.072822393747312739e-02,nan
        }\tableAdaptiveThree

        \pgfplotstableread[col sep=comma]{
ndof,nelem,res,err,res1,res2,res3,res4,err1,err2,err3,relerr1,relerr2,relerr3,cond
1.880000000000000000e+02,2.400000000000000000e+01,2.637678459993160196e+00,2.750392879609749652e+00,2.553030859086821369e+00,4.579952767477925568e-01,3.423839138177623525e-16,4.791882900687461921e-01,6.418937525993504423e-01,2.553030859086820925e+00,7.966597988679627873e-01,3.765245700978202548e-01,1.516165699786877798e-02,4.673078481731264833e-01,nan
7.360000000000000000e+02,9.600000000000000000e+01,1.352495157812586779e+00,1.411952400623620818e+00,1.296119271831964292e+00,2.595052255562128041e-01,5.080381130965519131e-16,2.863128062125264162e-01,3.527774664616770139e-01,1.296119271831964515e+00,4.350085906846714034e-01,2.065075171190645376e-01,7.697185263100082137e-03,2.546436564919137835e-01,nan
2.912000000000000000e+03,3.840000000000000000e+02,6.829860147169550810e-01,7.172309939193147743e-01,6.491247111150965088e-01,1.388780734094502145e-01,6.170700030074672055e-16,1.606856759165875836e-01,1.935544514767764679e-01,6.491247111150961757e-01,2.357839739502141652e-01,1.132093047746646902e-01,3.854914559825580941e-03,1.379092010762363318e-01,nan
1.158400000000000000e+04,1.536000000000000000e+03,3.449751470003383780e-01,3.665358620610047469e-01,3.248292525244442785e-01,7.354035969377580029e-02,8.891830131497487818e-14,8.992009924608598825e-02,1.089308190649652675e-01,3.248292525244445006e-01,1.302634697029461941e-01,6.369251132501826951e-02,1.929042126459639154e-03,7.616584168290280177e-02,nan
4.620800000000000000e+04,6.144000000000000000e+03,1.748290156891226799e-01,1.892074091287157811e-01,1.624861831661739697e-01,3.937073664854022492e-02,2.566422253678933757e-14,5.112113166458035624e-02,6.295201801819204357e-02,1.624861831661737199e-01,7.372060347506048661e-02,3.680367298765792328e-02,9.649460057203253097e-04,4.309931703801551939e-02,nan
1.845760000000000000e+05,2.457600000000000000e+04,8.915023094082820398e-02,9.908776624225530683e-02,8.126002673084596617e-02,2.151519220597389315e-02,3.300342342724167996e-12,2.969289876123922672e-02,3.724045787768594240e-02,8.126002673084584127e-02,4.275911331683984046e-02,2.177079343657157065e-02,4.825735741214668144e-04,2.499700263109980169e-02,nan
        }\tableUniformFour

        \pgfplotstableread[col sep=comma]{
ndof,nelem,res,err,res1,res2,res3,res4,err1,err2,err3,relerr1,relerr2,relerr3,cond
1.880000000000000000e+02,2.400000000000000000e+01,2.637678459993160196e+00,2.750392879609749652e+00,2.553030859086821369e+00,4.579952767477925568e-01,3.423839138177623525e-16,4.791882900687461921e-01,6.418937525993504423e-01,2.553030859086820925e+00,7.966597988679627873e-01,3.765245700978202548e-01,1.516165699786877798e-02,4.673078481731264833e-01,nan
5.200000000000000000e+02,6.800000000000000000e+01,1.688978641714293794e+00,1.757702758526647813e+00,1.630948078058457318e+00,2.841251316104528857e-01,5.583715369529872414e-16,3.345596037058685690e-01,4.080274928491037789e-01,1.630948078058457540e+00,5.128751496113698183e-01,2.389130526691028988e-01,9.685618356706770030e-03,3.003046848048792627e-01,nan
9.000000000000000000e+02,1.180000000000000000e+02,1.263329595463806410e+00,1.324957713430321204e+00,1.200180001947346753e+00,2.644275273484429922e-01,5.485228224208189461e-16,2.926563043585435619e-01,3.357114340032371635e-01,1.200180001947347197e+00,4.498652446940760719e-01,1.964845173636751707e-01,7.127435204530107801e-03,2.632962316128762970e-01,nan
1.723000000000000000e+03,2.270000000000000000e+02,8.817263207959222804e-01,9.235619555830751892e-01,8.367568484937282225e-01,1.831099970948007472e-01,5.063773554897659096e-16,2.091650311287714115e-01,2.343056824086135748e-01,8.367568484937280004e-01,3.128985610230771153e-01,1.370517931725999672e-01,4.969193411293029369e-03,1.830229144616017933e-01,nan
3.410000000000000000e+03,4.510000000000000000e+02,6.222259194082613298e-01,6.533807357177491282e-01,5.853447775479007076e-01,1.399314618475899497e-01,5.566929368075919770e-16,1.579739604050340862e-01,1.759931651479767489e-01,5.853447775479009296e-01,2.308772033477157204e-01,1.029064299630572527e-01,3.476148556144352010e-03,1.349981332308652360e-01,nan
7.179000000000000000e+03,9.520000000000000000e+02,4.419989282280665233e-01,4.665459814625657065e-01,4.137015887080489618e-01,1.029395419232138065e-01,6.910757484777406886e-16,1.166940391101221747e-01,1.319384582344221291e-01,4.137015887080490173e-01,1.706118154115611463e-01,7.714317976984506886e-02,2.456822422154706156e-03,9.975512919643850229e-02,nan
1.302300000000000000e+04,1.729000000000000000e+03,3.192300289154760273e-01,3.371580470795760287e-01,2.953117924376067460e-01,8.265645498219099940e-02,1.004178339429988057e-13,8.869423351289847102e-02,9.615016135420399734e-02,2.953117924376069681e-01,1.312312479151283018e-01,5.621173256943891305e-02,1.753748712360381906e-03,7.672099254606602747e-02,nan
2.667200000000000000e+04,3.546000000000000000e+03,2.297695165400830131e-01,2.413777765782426643e-01,2.147981952267612182e-01,5.421997037749910764e-02,1.324807054110225325e-12,6.095868166760803897e-02,6.683600803840127058e-02,2.147981952267613570e-01,8.750951024130411426e-02,3.907220683104136588e-02,1.275607909002030331e-03,5.115789802806353787e-02,nan
4.669600000000000000e+04,6.212000000000000000e+03,1.682145124990488683e-01,1.767957191584127452e-01,1.558629110481120350e-01,4.336454938087539596e-02,1.962386506050320988e-12,4.606941536306492085e-02,4.979406839910980459e-02,1.558629110481122015e-01,6.696290035992612244e-02,2.910900130371497224e-02,9.256128145328938522e-04,3.914568976879327644e-02,nan
9.395400000000000000e+04,1.250700000000000000e+04,1.230103792036966737e-01,1.289036482757920654e-01,1.144872685712903376e-01,3.024112172696897735e-02,1.404079678567263673e-12,3.331205913204904334e-02,3.547318821076148043e-02,1.144872685712901433e-01,4.743910614733674175e-02,2.073704374065909323e-02,6.798980089437362288e-04,2.773212301500018601e-02,nan
1.620630000000000000e+05,2.158200000000000000e+04,9.004534629613955565e-02,9.443147578784025220e-02,8.342756733648137857e-02,2.314886705174556969e-02,7.124665124869677640e-13,2.474136924523717002e-02,2.645487924127759224e-02,8.342756733648115652e-02,3.545820063279834372e-02,1.546504486516567071e-02,4.954458048559133186e-04,2.072822403092511312e-02,nan
        }\tableAdaptiveFour

        %
        %
        \addplot+ [line1, minorline, forget plot] table [x=ndof, y=relerr1] {\tableUniformOne};
        \addplot+ [line2, minorline, forget plot] table [x=ndof, y=relerr1] {\tableUniformTwo};
        \addplot+ [line3, minorline, forget plot] table [x=ndof, y=relerr1] {\tableUniformThree};
        \addplot+ [line4, minorline, forget plot] table [x=ndof, y=relerr1] {\tableUniformFour};

        \addplot+ [line1, majorline, forget plot] table [x=ndof, y=relerr1] {\tableAdaptiveOne};
        \addplot+ [line2, majorline, forget plot] table [x=ndof, y=relerr1] {\tableAdaptiveTwo};
        \addplot+ [line3, majorline, forget plot] table [x=ndof, y=relerr1] {\tableAdaptiveThree};
        \addplot+ [line4, majorline, forget plot] table [x=ndof, y=relerr1] {\tableAdaptiveFour};

        %
        %
        \drawslopetriangle[ST1]{0.5}{2e3}{3e-2}
        \drawswappedslopetriangle[ST2]{0.33}{1e5}{8e-2}
    \end{loglogaxis}
\end{tikzpicture}

%% file: figures/plot_dls_lshape_error_scalingWidth.tex
\begin{tikzpicture}[>=stealth]
    %
    %
    \colorlet{col1}{TUblue}
    \colorlet{col2}{TUgreen}
    \colorlet{col3}{TUmagenta}
    \colorlet{col4}{TUyellow}
    \colorlet{col5}{purple}
    \colorlet{col6}{green}
    \pgfplotsset{%
        linedefault/.style = {%
            mark = *,%
            mark size = 2pt,%
            every mark/.append style = {solid},%
            gray,%
            every mark/.append style = {fill = gray!60!white}%
        },%
        line1/.style = {%
            linedefault,%
            col1,%
            every mark/.append style = {fill = col1!60!white}%
        },%
        line2/.style = {%
            linedefault,%
            mark = triangle*,%
            mark size = 2.75pt,%
            col2,%
            every mark/.append style = {fill = col2!60!white}%
        },%
        line3/.style = {%
            linedefault,%
            mark = square*,%
            mark size = 1.66pt,%
            col3,%
            every mark/.append style = {fill = col3!60!white}%
        },%
        line4/.style = {%
            linedefault,%
            mark = pentagon*,%
            mark size = 2.2pt,%
            col4,%
            every mark/.append style = {fill = col4!60!white}%
        },%
        line5/.style = {%
            linedefault,%
            mark = diamond*,%
            mark size = 2.75pt,%
            col5,%
            every mark/.append style = {fill = col5!60!white}%
        },%
        line6/.style = {%
            linedefault,%
            mark = halfsquare*,%
            mark size = 1.66pt,%
            col6,%
            every mark/.append style = {fill = col6!60!white}%
        },%
        minorline/.style = {%
            dashed,%
            every mark/.append style = {fill = black!20!white}%
        },%
        majorline/.style = {%
            solid%
        }%
    }

    %
    %
    \begin{loglogaxis}[%
            width            = 0.36\textwidth,%
            xlabel           = ndof,%
            ylabel           = {relative energy error},%
            ymajorgrids      = true,%
            font             = \footnotesize,%
            grid style       = {densely dotted, semithick},%
            legend style     = {legend pos  = south west}%
        ]

        %
        %
        \pgfplotstableread[col sep=comma]{
ndof,nelem,res,err,res1,res2,res3,res4,err1,err2,err3,relerr1,relerr2,relerr3,cond
1.880000000000000000e+02,2.400000000000000000e+01,1.922740820674140982e+00,2.074393516215635280e+00,1.805416594311698431e+00,4.569966636030289298e-01,3.419920311564913281e-16,4.780765983880790326e-01,6.406661758091988235e-01,1.805416594311698431e+00,7.956922974840411333e-01,3.758044932606454669e-01,2.144361635409711836e-02,4.667403273939085828e-01,nan
7.360000000000000000e+02,9.600000000000000000e+01,9.945582442483172914e-01,1.073969236858786758e+00,9.165393591960622199e-01,2.593243740668844111e-01,4.688092239875740022e-16,2.860988872078772394e-01,3.525482443367516372e-01,9.165393591960622199e-01,4.348278511992884621e-01,2.063733359527845923e-01,1.088599390808184762e-02,2.545378559067851820e-01,nan
2.912000000000000000e+03,3.840000000000000000e+02,5.057426818516665978e-01,5.511020234172754328e-01,4.590089898125948031e-01,1.388485219457296638e-01,5.025444303102423968e-16,1.606471147304659730e-01,1.935111322685440227e-01,4.590089898125948031e-01,2.357490809610559568e-01,1.131839675250692623e-01,5.451773465470997719e-03,1.378887922919693754e-01,nan
1.158400000000000000e+04,1.536000000000000000e+03,2.573898347617985527e-01,2.856368399154572968e-01,2.296906505850448976e-01,7.353530990727477623e-02,1.662762042172394203e-15,8.991268136011652812e-02,1.089218311749720197e-01,2.296906505850450086e-01,1.302560710454535586e-01,6.368725605116554167e-02,2.728097531789443284e-03,7.616151564294285081e-02,nan
4.620800000000000000e+04,6.144000000000000000e+03,1.317731341735342010e-01,1.503268060996296174e-01,1.148954733162138031e-01,3.936977693556280300e-02,2.718548984073484596e-15,5.111956343126843744e-02,6.294999423637180624e-02,1.148954733162136782e-01,7.371890577863747518e-02,3.680248982297704141e-02,1.364644376419905372e-03,4.309832451282012700e-02,nan
1.845760000000000000e+05,2.457600000000000000e+04,6.816267135034016766e-02,8.072627075469120117e-02,5.745962124801556187e-02,2.151498805031592382e-02,1.262841745268080489e-12,2.969254082836034550e-02,3.723997783411574225e-02,5.745962124801558962e-02,4.275870399005184813e-02,2.177051280281995516e-02,6.824633441277266552e-04,2.499676333842039411e-02,nan
        }\tableUniformOne

        \pgfplotstableread[col sep=comma]{
ndof,nelem,res,err,res1,res2,res3,res4,err1,err2,err3,relerr1,relerr2,relerr3,cond
1.880000000000000000e+02,2.400000000000000000e+01,1.922740820674140982e+00,2.074393516215635280e+00,1.805416594311698431e+00,4.569966636030289298e-01,3.419920311564913281e-16,4.780765983880790326e-01,6.406661758091988235e-01,1.805416594311698431e+00,7.956922974840411333e-01,3.758044932606454669e-01,2.144361635409711836e-02,4.667403273939085828e-01,nan
5.200000000000000000e+02,6.800000000000000000e+01,1.233971684640085220e+00,1.326435976805148576e+00,1.153308656744947092e+00,2.834708747609424595e-01,4.989095628090362944e-16,3.349769005960285817e-01,4.075001973272314171e-01,1.153308656744947314e+00,5.130839424679570770e-01,2.386043043984664158e-01,1.369817671941523166e-02,3.004269396519593238e-01,nan
8.770000000000000000e+02,1.150000000000000000e+02,9.359957544856275513e-01,1.013208873017421352e+00,8.497143683802269587e-01,2.632790063170487582e-01,6.121725242636753425e-16,2.911317732849139084e-01,3.271078009082773463e-01,8.497143683802267367e-01,4.444976930360565182e-01,1.914169014929841539e-01,1.009229115230205268e-02,2.601111036957495792e-01,nan
1.715000000000000000e+03,2.260000000000000000e+02,6.616005979300529916e-01,7.146421876731671174e-01,6.007375170420488075e-01,1.816358976035203077e-01,5.084966840632932341e-16,2.093757089571785523e-01,2.300311708889703211e-01,6.007375170420489185e-01,3.113094159191037713e-01,1.345425569513482023e-01,7.135121424008948957e-03,1.820812573310048865e-01,nan
3.417000000000000000e+03,4.520000000000000000e+02,4.617347310036772345e-01,5.040347363132365865e-01,4.104711140168099681e-01,1.398337792563900395e-01,5.327515549690772887e-16,1.586157008360533083e-01,1.781509461630758073e-01,4.104711140168100236e-01,2.320058627478058799e-01,1.041703079043878599e-01,4.875275661719261945e-03,1.356609250671005351e-01,nan
7.457000000000000000e+03,9.890000000000000000e+02,3.258381667997057307e-01,3.539714403851788038e-01,2.899327674542118705e-01,9.905942005602334788e-02,2.476045302132410668e-15,1.108906334994037152e-01,1.235364171850229786e-01,2.899327674542117594e-01,1.611630373062883370e-01,7.222556721636196053e-02,3.443609306036928691e-03,9.422397094716433719e-02,nan
1.340500000000000000e+04,1.780000000000000000e+03,2.344976320950872273e-01,2.546685133084671815e-01,2.056302319867577832e-01,7.676116173804349574e-02,9.791458719554679528e-16,8.254133025279507263e-02,8.979435491929205204e-02,2.056302319867576722e-01,1.204542784148631002e-01,5.249452562884123624e-02,2.442325458201559125e-03,7.041857153532588010e-02,nan
2.799500000000000000e+04,3.722000000000000000e+03,1.689162967236020330e-01,1.823428508174387841e-01,1.501393028138544861e-01,5.155587365768200120e-02,5.061221657799226072e-14,5.773124752487863998e-02,6.270078383963474666e-02,1.501393028138545971e-01,8.231474169731570201e-02,3.665436136691811730e-02,1.783244797841191786e-03,4.812051944541541337e-02,nan
4.921700000000000000e+04,6.548000000000000000e+03,1.227077290724349901e-01,1.331457550763252562e-01,1.069430629415106110e-01,4.079830489985179776e-02,2.658166839654471566e-15,4.422517812401593335e-02,4.731534887527338168e-02,1.069430629415105277e-01,6.365713741870644704e-02,2.765985219340425655e-02,1.270191462714505542e-03,3.721301974752628389e-02,nan
1.006230000000000000e+05,1.339600000000000000e+04,8.964105170874783757e-02,9.645182359452665610e-02,7.942468585447533258e-02,2.778457675172269284e-02,6.084898859411802570e-13,3.090719532885407780e-02,3.294408907682101761e-02,7.942468585447522156e-02,4.369623034686281188e-02,1.925854070131595777e-02,9.433483119529113056e-04,2.554405522237433029e-02,nan
        }\tableAdaptiveOne

        \pgfplotstableread[col sep=comma]{
ndof,nelem,res,err,res1,res2,res3,res4,err1,err2,err3,relerr1,relerr2,relerr3,cond
1.880000000000000000e+02,2.400000000000000000e+01,1.922740820674140316e+00,2.074393516215634392e+00,1.805416594311697764e+00,4.569966636030293738e-01,5.108373094066665780e-16,4.780765983880793102e-01,6.406661758091976022e-01,1.805416594311697764e+00,7.956922974840409113e-01,3.758044932606447452e-01,2.144361635409710795e-02,4.667403273939084163e-01,nan
7.360000000000000000e+02,9.600000000000000000e+01,9.945582442483166252e-01,1.073969236858786314e+00,9.165393591960616648e-01,2.593243740668843555e-01,5.157809322350766461e-16,2.860988872078767953e-01,3.525482443367511931e-01,9.165393591960619979e-01,4.348278511992877404e-01,2.063733359527842592e-01,1.088599390808184415e-02,2.545378559067846824e-01,nan
2.912000000000000000e+03,3.840000000000000000e+02,5.057426818516671529e-01,5.511020234172910870e-01,4.590089898125944146e-01,1.388485219457243347e-01,6.439428133145974467e-16,1.606471147304736058e-01,1.935111322685674762e-01,4.590089898125945256e-01,2.357490809610736926e-01,1.131839675250828903e-01,5.451773465470995117e-03,1.378887922919796449e-01,nan
1.158400000000000000e+04,1.536000000000000000e+03,2.573898347616879745e-01,2.856368399137348968e-01,2.296906505850423719e-01,7.353530990782412846e-02,8.866323155729837443e-16,8.991268135935708006e-02,1.089218311726728727e-01,2.296906505850423441e-01,1.302560710436038438e-01,6.368725604982120037e-02,2.728097531789414661e-03,7.616151564186127154e-02,nan
4.620800000000000000e+04,6.144000000000000000e+03,1.317731341735766115e-01,1.503268061000431755e-01,1.148954733162138170e-01,3.936977693545282847e-02,2.661214932982521863e-13,5.111956343146243503e-02,6.294999423687241968e-02,1.148954733162138725e-01,7.371890577905299002e-02,3.680248982326971702e-02,1.364644376419908841e-03,4.309832451306305767e-02,nan
1.845760000000000000e+05,2.457600000000000000e+04,6.816267135295160651e-02,8.072627079249453108e-02,5.745962124799783993e-02,2.151498803894785414e-02,5.722308153279232554e-13,2.969254084262671484e-02,3.723997787600204917e-02,5.745962124799772891e-02,4.275870402496643324e-02,2.177051282730670112e-02,6.824633441275148021e-04,2.499676335883146278e-02,nan
        }\tableUniformTwo

        \pgfplotstableread[col sep=comma]{
ndof,nelem,res,err,res1,res2,res3,res4,err1,err2,err3,relerr1,relerr2,relerr3,cond
1.880000000000000000e+02,2.400000000000000000e+01,1.922740820674140316e+00,2.074393516215634392e+00,1.805416594311697764e+00,4.569966636030293738e-01,5.108373094066665780e-16,4.780765983880793102e-01,6.406661758091976022e-01,1.805416594311697764e+00,7.956922974840409113e-01,3.758044932606447452e-01,2.144361635409710795e-02,4.667403273939084163e-01,nan
5.200000000000000000e+02,6.800000000000000000e+01,1.235016042953527515e+00,1.328948313283806026e+00,1.153307442407722316e+00,2.804413281366299748e-01,5.862970435103269203e-16,3.413198369748098870e-01,4.080399724601615485e-01,1.153307442407722316e+00,5.191232450406568777e-01,2.389203598776273030e-01,1.369816229638475338e-02,3.039631430669819112e-01,nan
8.770000000000000000e+02,1.150000000000000000e+02,9.366091646854119990e-01,1.014225942798006930e+00,8.497124465888514777e-01,2.612080181341585594e-01,6.635216765958875773e-16,2.949506005200166481e-01,3.266590197423728248e-01,8.497124465888513667e-01,4.471430486739635390e-01,1.911542837871779421e-01,1.009226832665755663e-02,2.616591125727899536e-01,nan
1.723000000000000000e+03,2.270000000000000000e+02,6.598872744058504880e-01,7.084846107375587509e-01,6.005792965644394243e-01,1.807112833141506503e-01,5.950128800845396133e-16,2.051807874629054951e-01,2.202391615222041577e-01,6.005792965644397574e-01,3.045482948927958344e-01,1.287999465051444803e-01,7.133241893833849184e-03,1.781054914090307584e-01,nan
3.523000000000000000e+03,4.660000000000000000e+02,4.563381727380848574e-01,4.973298408528868397e-01,4.076199079816866244e-01,1.347649954515781900e-01,7.110307253536626205e-16,1.546898009432799426e-01,1.736936404202813478e-01,4.076199079816870130e-01,2.258616844383056488e-01,1.015639847567716969e-01,4.841411121967370557e-03,1.320682358889193608e-01,nan
7.563000000000000000e+03,1.003000000000000000e+03,3.222425425890999007e-01,3.493817624815786504e-01,2.878211126171053569e-01,9.587731752706893373e-02,3.106314932098341611e-14,1.086591154480516241e-01,1.196456864173669754e-01,2.878211126171055789e-01,1.578338772516938160e-01,6.994968746851078223e-02,3.418528612211391496e-03,9.227604200611382035e-02,nan
1.382500000000000000e+04,1.836000000000000000e+03,2.322508705998394540e-01,2.523357978060247198e-01,2.038399560629323981e-01,7.528815376261098480e-02,8.895793678126615357e-16,8.198434673559697183e-02,8.850968287698544879e-02,2.038399560629324259e-01,1.195351965058000510e-01,5.174315394445572203e-02,2.421061870361466018e-03,6.988080708838076416e-02,nan
2.830400000000000000e+04,3.763000000000000000e+03,1.670461567461012653e-01,1.802195249714255709e-01,1.485825221977489508e-01,5.089874609013781553e-02,1.644569644663043629e-14,5.689437781934299915e-02,6.160349948441946194e-02,1.485825221977488675e-01,8.128542384434543633e-02,3.601280358633008466e-02,1.764754496602040628e-03,4.751866416417415556e-02,nan
5.017800000000000000e+04,6.676000000000000000e+03,1.217269549754323443e-01,1.319251162030213864e-01,1.064443588589211037e-01,3.986612007070925584e-02,6.797187621692474833e-14,4.356308649776723801e-02,4.642517923633771720e-02,1.064443588589214368e-01,6.259822269975214015e-02,2.713944397404877834e-02,1.264268220469573403e-03,3.659395581838079115e-02,nan
1.021470000000000000e+05,1.359900000000000000e+04,8.838656133855403441e-02,9.511529162848525454e-02,7.829377981732400793e-02,2.750854801379344627e-02,1.980816919564575270e-13,3.042282125877189894e-02,3.252408292297988973e-02,7.829377981732420222e-02,4.311828815603983844e-02,1.901300457115144907e-02,9.299162374076761243e-04,2.520618987943163206e-02,nan
        }\tableAdaptiveTwo

        \pgfplotstableread[col sep=comma]{
ndof,nelem,res,err,res1,res2,res3,res4,err1,err2,err3,relerr1,relerr2,relerr3,cond
1.880000000000000000e+02,2.400000000000000000e+01,1.922740820674140760e+00,2.074393516215635280e+00,1.805416594311698208e+00,4.569966636030288742e-01,5.567945878842655693e-16,4.780765983880796433e-01,6.406661758091991565e-01,1.805416594311698208e+00,7.956922974840416884e-01,3.758044932606457444e-01,2.144361635409711142e-02,4.667403273939089714e-01,nan
7.360000000000000000e+02,9.600000000000000000e+01,9.945582442483169583e-01,1.073969236858788090e+00,9.165393591960617758e-01,2.593243740668839670e-01,5.323076011906775651e-16,2.860988872078779610e-01,3.525482443367539132e-01,9.165393591960619979e-01,4.348278511992900719e-01,2.063733359527859246e-01,1.088599390808184762e-02,2.545378559067861257e-01,nan
2.912000000000000000e+03,3.840000000000000000e+02,5.057426818516661537e-01,5.511020234172834265e-01,4.590089898125945256e-01,1.388485219457269992e-01,6.921863866546610637e-16,1.606471147304678881e-01,1.935111322685568458e-01,4.590089898125946921e-01,2.357490809610642557e-01,1.131839675250766869e-01,5.451773465470994250e-03,1.378887922919741493e-01,nan
1.158400000000000000e+04,1.536000000000000000e+03,2.573898347616963012e-01,2.856368399138619618e-01,2.296906505850419000e-01,7.353530990778243959e-02,1.025599385522540530e-15,8.991268135941624107e-02,1.089218311728425564e-01,2.296906505850419278e-01,1.302560710437411229e-01,6.368725604992048206e-02,2.728097531789404687e-03,7.616151564194162393e-02,nan
4.620800000000000000e+04,6.144000000000000000e+03,1.317731341735725314e-01,1.503268060999927436e-01,1.148954733162138586e-01,3.936977693544589652e-02,1.305742550691600950e-14,5.111956343145716147e-02,6.294999423680064377e-02,1.148954733162139558e-01,7.371890577901127339e-02,3.680248982322771589e-02,1.364644376419905806e-03,4.309832451303862583e-02,nan
1.845760000000000000e+05,2.457600000000000000e+04,6.816267134849723908e-02,8.072627073698593336e-02,5.745962124797754367e-02,2.151498805755292301e-02,1.315240085124818734e-12,2.969254081895937938e-02,3.723997781559997694e-02,5.745962124797746040e-02,4.275870397280241708e-02,2.177051279199559985e-02,6.824633441272712903e-04,2.499676332833633205e-02,nan
        }\tableUniformThree

        \pgfplotstableread[col sep=comma]{
ndof,nelem,res,err,res1,res2,res3,res4,err1,err2,err3,relerr1,relerr2,relerr3,cond
1.880000000000000000e+02,2.400000000000000000e+01,1.922740820674140760e+00,2.074393516215635280e+00,1.805416594311698208e+00,4.569966636030288742e-01,5.567945878842655693e-16,4.780765983880796433e-01,6.406661758091991565e-01,1.805416594311698208e+00,7.956922974840416884e-01,3.758044932606457444e-01,2.144361635409711142e-02,4.667403273939089714e-01,nan
5.200000000000000000e+02,6.800000000000000000e+01,1.233859515912456128e+00,1.326303261044500248e+00,1.153309048347520127e+00,2.838625367224141272e-01,6.401569480446092463e-16,3.342298673861656688e-01,4.076823716116002028e-01,1.153309048347520349e+00,5.125950283939103569e-01,2.387109732263252015e-01,1.369818137059097135e-02,3.001406649377047486e-01,nan
8.770000000000000000e+02,1.150000000000000000e+02,9.357311440333944130e-01,1.013091018592872317e+00,8.497150080943488559e-01,2.640620546825845905e-01,5.618921994934989812e-16,2.895658995981357231e-01,3.277791561936707065e-01,8.497150080943488559e-01,4.437326241530245574e-01,1.918097650938341003e-01,1.009229875036254942e-02,2.596634007837064217e-01,nan
1.715000000000000000e+03,2.260000000000000000e+02,6.613859788090241665e-01,7.142029354812076125e-01,6.007378416650720476e-01,1.822383530609973934e-01,6.710898663824066043e-16,2.081697413899942783e-01,2.300599022588999121e-01,6.007378416650719366e-01,3.102778110245565379e-01,1.345593615955192612e-01,7.135125279644100385e-03,1.814778836241305593e-01,nan
3.417000000000000000e+03,4.520000000000000000e+02,4.616709289348528422e-01,5.039399250231259009e-01,4.104711776139444157e-01,1.400351043146789709e-01,7.228037472001230531e-16,1.582517883986422069e-01,1.782064407697249442e-01,4.104711776139445267e-01,2.317570384016890139e-01,1.042027572984898093e-01,4.875276417079508519e-03,1.355154298603258867e-01,nan
7.457000000000000000e+03,9.890000000000000000e+02,3.258040156981619750e-01,3.539447121925007878e-01,2.899327842041494740e-01,9.915412363500258897e-02,7.502181055109759115e-16,1.107054517846925434e-01,1.235999547413750149e-01,2.899327842041496406e-01,1.610555529020725185e-01,7.226271444915076392e-02,3.443609504980434458e-03,9.416113018944854973e-02,nan
1.340500000000000000e+04,1.780000000000000000e+03,2.344597135345746486e-01,2.546077764917436537e-01,2.056302451682783528e-01,7.687006687256488224e-02,1.831527572914297140e-13,8.233196439874339301e-02,8.981331083523547143e-02,2.056302451682784083e-01,1.203116424911436499e-01,5.250560741473086351e-02,2.442325614762019034e-03,7.033518538972674539e-02,nan
2.799500000000000000e+04,3.722000000000000000e+03,1.689662562822069058e-01,1.824081779361495237e-01,1.501658747375557990e-01,5.167495811786863152e-02,1.830784846703171260e-15,5.770185898084784909e-02,6.283780352122321500e-02,1.501658747375554104e-01,8.230652398758531163e-02,3.673446194327088105e-02,1.783560399711054409e-03,4.811571544004866752e-02,nan
4.917200000000000000e+04,6.542000000000000000e+03,1.228061642224922284e-01,1.332665317004846639e-01,1.070430651456867055e-01,4.081492137117150032e-02,2.188098626747792901e-13,4.424113887235602699e-02,4.734603614201917526e-02,1.070430651456868720e-01,6.371894174679244849e-02,2.767779150708252200e-02,1.271379215735010150e-03,3.724914963165162213e-02,nan
1.002780000000000000e+05,1.335000000000000000e+04,8.980141505746799813e-02,9.656579002963894265e-02,7.965208719240317314e-02,2.786085160132932559e-02,1.314066177459198557e-12,3.071826981237952420e-02,3.299040147516606875e-02,7.965208719240307600e-02,4.349862320176919872e-02,1.928561412229099978e-02,9.460492186793526426e-04,2.542853752633915387e-02,nan
        }\tableAdaptiveThree

        \pgfplotstableread[col sep=comma]{
ndof,nelem,res,err,res1,res2,res3,res4,err1,err2,err3,relerr1,relerr2,relerr3,cond
1.880000000000000000e+02,2.400000000000000000e+01,1.922740820674140094e+00,2.074393516215634392e+00,1.805416594311697542e+00,4.569966636030289853e-01,4.048319820644552521e-16,4.780765983880792547e-01,6.406661758091989345e-01,1.805416594311697542e+00,7.956922974840413554e-01,3.758044932606456889e-01,2.144361635409710448e-02,4.667403273939089159e-01,nan
7.360000000000000000e+02,9.600000000000000000e+01,9.945582442483168473e-01,1.073969236858786314e+00,9.165393591960618869e-01,2.593243740668847996e-01,5.132991490655417332e-16,2.860988872078767398e-01,3.525482443367510821e-01,9.165393591960621089e-01,4.348278511992879070e-01,2.063733359527843148e-01,1.088599390808184589e-02,2.545378559067849045e-01,nan
2.912000000000000000e+03,3.840000000000000000e+02,5.057426818516661537e-01,5.511020234172716581e-01,4.590089898125945256e-01,1.388485219457311626e-01,5.590723396417386262e-16,1.606471147304641689e-01,1.935111322685389990e-01,4.590089898125945811e-01,2.357490809610517102e-01,1.131839675250662647e-01,5.451773465470993382e-03,1.378887922919668219e-01,nan
1.158400000000000000e+04,1.536000000000000000e+03,2.573898347616826454e-01,2.856368399138178860e-01,2.296906505850421776e-01,7.353530990778030241e-02,6.312731524087112855e-15,8.991268135937809103e-02,1.089218311727937760e-01,2.296906505850422053e-01,1.302560710436848901e-01,6.368725604989194933e-02,2.728097531789408590e-03,7.616151564190874745e-02,nan
4.620800000000000000e+04,6.144000000000000000e+03,1.317731341735329798e-01,1.503268060994621957e-01,1.148954733162134978e-01,3.936977693563536995e-02,6.892211638702929232e-14,5.111956343121003971e-02,6.294999423616391698e-02,1.148954733162134423e-01,7.371890577847396708e-02,3.680248982285548587e-02,1.364644376419899951e-03,4.309832451272452292e-02,nan
1.845760000000000000e+05,2.457600000000000000e+04,6.816267135232487173e-02,8.072627077357832914e-02,5.745962124798765364e-02,2.151498804330214537e-02,3.218966083324820124e-12,2.969254083805260638e-02,3.723997785406450017e-02,5.745962124798746629e-02,4.275870400837353952e-02,2.177051281448197720e-02,6.824633441273903357e-04,2.499676334913122563e-02,nan
        }\tableUniformFour

        \pgfplotstableread[col sep=comma]{
ndof,nelem,res,err,res1,res2,res3,res4,err1,err2,err3,relerr1,relerr2,relerr3,cond
1.880000000000000000e+02,2.400000000000000000e+01,1.922740820674140094e+00,2.074393516215634392e+00,1.805416594311697542e+00,4.569966636030289853e-01,4.048319820644552521e-16,4.780765983880792547e-01,6.406661758091989345e-01,1.805416594311697542e+00,7.956922974840413554e-01,3.758044932606456889e-01,2.144361635409710448e-02,4.667403273939089159e-01,nan
5.200000000000000000e+02,6.800000000000000000e+01,1.233859515912456128e+00,1.326303261044498694e+00,1.153309048347520127e+00,2.838625367224147933e-01,4.371881803543146497e-16,3.342298673861651137e-01,4.076823716115979268e-01,1.153309048347520127e+00,5.125950283939090246e-01,2.387109732263238970e-01,1.369818137059097135e-02,3.001406649377040270e-01,nan
8.770000000000000000e+02,1.150000000000000000e+02,9.357311440333947461e-01,1.013091018592870984e+00,8.497150080943488559e-01,2.640620546825850345e-01,4.531204528829544712e-16,2.895658995981358896e-01,3.277791561936683196e-01,8.497150080943488559e-01,4.437326241530234472e-01,1.918097650938327681e-01,1.009229875036254942e-02,2.596634007837058111e-01,nan
1.715000000000000000e+03,2.260000000000000000e+02,6.613859788090239444e-01,7.142029354812099440e-01,6.007378416650719366e-01,1.822383530609964220e-01,5.263561771952491639e-16,2.081697413899946947e-01,2.300599022589041309e-01,6.007378416650719366e-01,3.102778110245588694e-01,1.345593615955217315e-01,7.135125279644102120e-03,1.814778836241319471e-01,nan
3.417000000000000000e+03,4.520000000000000000e+02,4.616709289348539524e-01,5.039399250231446636e-01,4.104711776139444157e-01,1.400351043146745855e-01,5.630311665358704590e-16,1.582517883986493401e-01,1.782064407697523389e-01,4.104711776139445267e-01,2.317570384017086649e-01,1.042027572985058104e-01,4.875276417079508519e-03,1.355154298603373775e-01,nan
7.457000000000000000e+03,9.890000000000000000e+02,3.258040156981599211e-01,3.539447121924711448e-01,2.899327842041496406e-01,9.915412363501184545e-02,2.177784693130258687e-14,1.107054517846779024e-01,1.235999547413314525e-01,2.899327842041497516e-01,1.610555529020405163e-01,7.226271444912529818e-02,3.443609504980436192e-03,9.416113018942984247e-02,nan
1.340500000000000000e+04,1.780000000000000000e+03,2.344597135345476979e-01,2.546077764915832820e-01,2.056302451682791854e-01,7.687006687261040139e-02,8.582748165329948389e-16,8.233196439862210114e-02,8.981331083502748502e-02,2.056302451682792687e-01,1.203116424909581317e-01,5.250560741460927328e-02,2.442325614762029442e-03,7.033518538961829047e-02,nan
2.799500000000000000e+04,3.722000000000000000e+03,1.689662562858043060e-01,1.824081779611311238e-01,1.501658747375813618e-01,5.167495811125098065e-02,4.551806782786996245e-14,5.770185899724195044e-02,6.283780355536792983e-02,1.501658747375810843e-01,8.230652401683460906e-02,3.673446196323160162e-02,1.783560399711358203e-03,4.811571545714756620e-02,nan
4.917200000000000000e+04,6.542000000000000000e+03,1.228061642226106753e-01,1.332665317004405048e-01,1.070430651456817511e-01,4.081492137114797747e-02,6.722935983735019369e-14,4.424113887271848011e-02,4.734603614175320746e-02,1.070430651456818899e-01,6.371894174690610757e-02,2.767779150692704221e-02,1.271379215734950953e-03,3.724914963171806898e-02,nan
1.002780000000000000e+05,1.335000000000000000e+04,8.980141505785363409e-02,9.656579003105450476e-02,7.965208719235870871e-02,2.786085160031089025e-02,1.072988118339445967e-12,3.071826981454588157e-02,3.299040147655359773e-02,7.965208719235872259e-02,4.349862320394057985e-02,1.928561412310211831e-02,9.460492186788257204e-04,2.542853752760849961e-02,nan
        }\tableAdaptiveFour

        %
        %
        \addplot+ [line1, minorline, forget plot] table [x=ndof, y=relerr1] {\tableUniformOne};
        \addplot+ [line2, minorline, forget plot] table [x=ndof, y=relerr1] {\tableUniformTwo};
        \addplot+ [line3, minorline, forget plot] table [x=ndof, y=relerr1] {\tableUniformThree};
        \addplot+ [line4, minorline, forget plot] table [x=ndof, y=relerr1] {\tableUniformFour};

        \addplot+ [line1, majorline, forget plot] table [x=ndof, y=relerr1] {\tableAdaptiveOne};
        \addplot+ [line2, majorline, forget plot] table [x=ndof, y=relerr1] {\tableAdaptiveTwo};
        \addplot+ [line3, majorline, forget plot] table [x=ndof, y=relerr1] {\tableAdaptiveThree};
        \addplot+ [line4, majorline, forget plot] table [x=ndof, y=relerr1] {\tableAdaptiveFour};

        %
        %
        \drawslopetriangle[ST1]{0.5}{2e3}{3e-2}
        \drawswappedslopetriangle[ST2]{0.33}{1e5}{8e-2}
    \end{loglogaxis}
\end{tikzpicture}

%% file: figures/plot_dls_lshape_error_scalingFriedrichs.tex
\begin{tikzpicture}[>=stealth]
    %
    %
    \colorlet{col1}{TUblue}
    \colorlet{col2}{TUgreen}
    \colorlet{col3}{TUmagenta}
    \colorlet{col4}{TUyellow}
    \colorlet{col5}{purple}
    \colorlet{col6}{green}
    \pgfplotsset{%
        linedefault/.style = {%
            mark = *,%
            mark size = 2pt,%
            every mark/.append style = {solid},%
            gray,%
            every mark/.append style = {fill = gray!60!white}%
        },%
        line1/.style = {%
            linedefault,%
            col1,%
            every mark/.append style = {fill = col1!60!white}%
        },%
        line2/.style = {%
            linedefault,%
            mark = triangle*,%
            mark size = 2.75pt,%
            col2,%
            every mark/.append style = {fill = col2!60!white}%
        },%
        line3/.style = {%
            linedefault,%
            mark = square*,%
            mark size = 1.66pt,%
            col3,%
            every mark/.append style = {fill = col3!60!white}%
        },%
        line4/.style = {%
            linedefault,%
            mark = pentagon*,%
            mark size = 2.2pt,%
            col4,%
            every mark/.append style = {fill = col4!60!white}%
        },%
        line5/.style = {%
            linedefault,%
            mark = diamond*,%
            mark size = 2.75pt,%
            col5,%
            every mark/.append style = {fill = col5!60!white}%
        },%
        line6/.style = {%
            linedefault,%
            mark = halfsquare*,%
            mark size = 1.66pt,%
            col6,%
            every mark/.append style = {fill = col6!60!white}%
        },%
        minorline/.style = {%
            dashed,%
            every mark/.append style = {fill = black!20!white}%
        },%
        majorline/.style = {%
            solid%
        }%
    }

    %
    %
    \begin{loglogaxis}[%
            width            = 0.36\textwidth,%
            xlabel           = ndof,%
            ylabel           = {relative energy error},%
            ymajorgrids      = true,%
            font             = \footnotesize,%
            grid style       = {densely dotted, semithick},%
            legend style     = {legend pos  = south west}%
        ]

        %
        %
        \pgfplotstableread[col sep=comma]{
ndof,nelem,res,err,res1,res2,res3,res4,err1,err2,err3,relerr1,relerr2,relerr3,cond
1.880000000000000000e+02,2.400000000000000000e+01,6.691468253750407769e-01,1.025606713479946253e+00,3.566616625857724765e-01,3.931095612333450617e-01,2.477924273941170541e-16,4.074491443333447815e-01,5.936417091598806817e-01,3.566616625857724765e-01,7.564728087131913359e-01,3.482206960706684118e-01,1.633434078715055859e-01,4.437348049237161418e-01,nan
7.360000000000000000e+02,9.600000000000000000e+01,4.018328538784070036e-01,5.720407191928682522e-01,1.652282559352703240e-01,2.464800014124796934e-01,3.822651718777364391e-16,2.709554849054693348e-01,3.427995284376484930e-01,1.652282559352703240e-01,4.271050118626665837e-01,2.006666701171943290e-01,7.567037209895097205e-02,2.500170899051719076e-01,nan
2.912000000000000000e+03,3.840000000000000000e+02,2.226953929342509741e-01,3.122001557070244737e-01,7.752238586897751416e-02,1.366714899126683369e-01,5.032912734890589204e-16,1.578113480342056563e-01,1.914562321967724912e-01,7.752238586897754191e-02,2.341019608326406487e-01,1.119820638399216389e-01,3.550326745488059654e-02,1.369253976337967971e-01,nan
1.158400000000000000e+04,1.536000000000000000e+03,1.214925828018452447e-01,1.733462263358528099e-01,3.772328042619803312e-02,7.315877376198799686e-02,9.027751265676533212e-16,8.935992172888192198e-02,1.084461173911593773e-01,3.772328042619802618e-02,1.298665003823157305e-01,6.340910331328074301e-02,1.727629497202878839e-02,7.593373131076956717e-02,nan
4.620800000000000000e+04,6.144000000000000000e+03,6.703768583937516623e-02,9.857353874140875605e-02,1.866840647801303738e-02,3.929790042375387593e-02,2.034371188725006860e-15,5.100212686305510140e-02,6.283363054746113341e-02,1.866840647801300962e-02,7.362179033275986773e-02,3.673446005539963977e-02,8.549651322052077840e-03,4.304154785617510942e-02,nan
1.845760000000000000e+05,2.457600000000000000e+04,3.779800493525656280e-02,5.742160326159143763e-02,9.295180846622706078e-03,2.149966973848542462e-02,2.553853175909822800e-14,2.966568710227882044e-02,3.721084230050389352e-02,9.295180846622702608e-03,4.273398354292771983e-02,2.175348015284526634e-02,4.256954402829724732e-03,2.498231175058635706e-02,nan
        }\tableUniformOne

        \pgfplotstableread[col sep=comma]{
ndof,nelem,res,err,res1,res2,res3,res4,err1,err2,err3,relerr1,relerr2,relerr3,cond
1.880000000000000000e+02,2.400000000000000000e+01,6.691468253750407769e-01,1.025606713479946253e+00,3.566616625857724765e-01,3.931095612333450617e-01,2.477924273941170541e-16,4.074491443333447815e-01,5.936417091598806817e-01,3.566616625857724765e-01,7.564728087131913359e-01,3.482206960706684118e-01,1.633434078715055859e-01,4.437348049237161418e-01,nan
5.280000000000000000e+02,6.900000000000000000e+01,4.645504417268896402e-01,6.830605420970882191e-01,2.257866553728564707e-01,2.594642995420956000e-01,3.353784779502473084e-16,3.122591494693420788e-01,4.085020935114554996e-01,2.257866553728564429e-01,4.987164826053627209e-01,2.393181784222497832e-01,1.034048007908551087e-01,2.921696658646869005e-01,nan
1.134000000000000000e+03,1.490000000000000000e+02,3.512009102352414702e-01,4.899211964668722108e-01,1.303821863579617546e-01,2.161225937455042156e-01,4.358878004473377961e-16,2.441998961983719296e-01,2.886888337220095968e-01,1.303821863579617268e-01,3.737405804997301351e-01,1.689065895492492808e-01,5.971171407087753902e-02,2.186688207315767296e-01,nan
2.340000000000000000e+03,3.090000000000000000e+02,2.510111415316980121e-01,3.427687390905339027e-01,9.750850128206714484e-02,1.474123774849668278e-01,5.605968120035818230e-16,1.782365740105011642e-01,2.035671964278465973e-01,9.750850128206713097e-02,2.579590998859046524e-01,1.190419505586394061e-01,4.465639372516524985e-02,1.508492279386147161e-01,nan
4.956000000000000000e+03,6.560000000000000000e+02,1.819728532728017989e-01,2.425924512650964038e-01,6.148037562983391102e-02,1.130779319583521264e-01,5.774334875450577239e-16,1.286377240341318451e-01,1.442584950513282060e-01,6.148037562983391796e-02,1.850965894573657344e-01,8.434213409878976719e-02,2.815643512046980129e-02,1.082185237249779497e-01,nan
8.565000000000000000e+03,1.136000000000000000e+03,1.349300842170371439e-01,1.783259907176502901e-01,4.924296255606274569e-02,8.125570377033265734e-02,7.422863016000237500e-16,9.580589162628197486e-02,1.055700104768983877e-01,4.924296255606273875e-02,1.350194893044750766e-01,6.171751256163563598e-02,2.255201375430425717e-02,7.893403618670732924e-02,nan
1.652100000000000000e+04,2.194000000000000000e+03,9.966027545601267856e-02,1.289258107307255374e-01,3.282813032108918005e-02,6.346799554720226888e-02,3.376475797373012287e-15,6.947156184231227660e-02,7.654966312051884880e-02,3.282813032108913842e-02,9.840898133977166173e-02,4.475041932600895556e-02,1.503444163212453930e-02,5.752923005639901805e-02,nan
2.984500000000000000e+04,3.968000000000000000e+03,7.319005132368566757e-02,9.479901759857156662e-02,2.610373410740913711e-02,4.473276420333967274e-02,1.540268356267351765e-15,5.171420003385761299e-02,5.587635724236576784e-02,2.610373410740908853e-02,7.199501027590377022e-02,3.266452486906211428e-02,1.195484064946541157e-02,4.208726050993499834e-02,nan
5.687800000000000000e+04,7.567000000000000000e+03,5.413624159055402396e-02,6.907310236099371836e-02,1.789149090802056100e-02,3.457666917427426601e-02,2.095138213468278924e-15,3.761756445401324978e-02,4.056823527703840720e-02,1.789149090802054365e-02,5.296419837368814693e-02,2.371549357086034926e-02,8.193843911607373923e-03,3.096196069262863171e-02,nan
1.041770000000000000e+05,1.386800000000000000e+04,3.981748606707781468e-02,5.093995968901984805e-02,1.387211700638491842e-02,2.474252451444111309e-02,5.822000710554586147e-14,2.794287113228314384e-02,2.983779583704416691e-02,1.387211700638494617e-02,3.888637039461259604e-02,1.744262777094747821e-02,6.353073763287618754e-03,2.273225836991357024e-02,nan
        }\tableAdaptiveOne

        \pgfplotstableread[col sep=comma]{
ndof,nelem,res,err,res1,res2,res3,res4,err1,err2,err3,relerr1,relerr2,relerr3,cond
1.880000000000000000e+02,2.400000000000000000e+01,6.691468253750408879e-01,1.025606713479946475e+00,3.566616625857720324e-01,3.931095612333453393e-01,3.041234422105361890e-16,4.074491443333450591e-01,5.936417091598809037e-01,3.566616625857720324e-01,7.564728087131916690e-01,3.482206960706685228e-01,1.633434078715053916e-01,4.437348049237163083e-01,nan
7.360000000000000000e+02,9.600000000000000000e+01,4.018328538784071147e-01,5.720407191928684743e-01,1.652282559352701019e-01,2.464800014124798322e-01,5.069207530403312868e-16,2.709554849054694459e-01,3.427995284376487151e-01,1.652282559352700742e-01,4.271050118626667502e-01,2.006666701171943845e-01,7.567037209895086103e-02,2.500170899051719076e-01,nan
2.912000000000000000e+03,3.840000000000000000e+02,2.226953929342508076e-01,3.122001557070245847e-01,7.752238586897820805e-02,1.366714899126684757e-01,6.569518267649839023e-16,1.578113480342049624e-01,1.914562321967724634e-01,7.752238586897820805e-02,2.341019608326406209e-01,1.119820638399215418e-01,3.550326745488090879e-02,1.369253976337966583e-01,nan
1.158400000000000000e+04,1.536000000000000000e+03,1.214925828018439541e-01,1.733462263358526434e-01,3.772328042620228666e-02,7.315877376198844095e-02,8.112737495350342221e-16,8.935992172887800844e-02,1.084461173911589332e-01,3.772328042620225891e-02,1.298665003823146480e-01,6.340910331328046545e-02,1.727629497203073475e-02,7.593373131076891491e-02,nan
4.620800000000000000e+04,6.144000000000000000e+03,6.703768583937343151e-02,9.857353874140298289e-02,1.866840647801719030e-02,3.929790042375734538e-02,3.575439697243149286e-15,5.100212686304862741e-02,6.283363054745608189e-02,1.866840647801715214e-02,7.362179033275538520e-02,3.673446005539669074e-02,8.549651322053980831e-03,4.304154785617249346e-02,nan
1.845760000000000000e+05,2.457600000000000000e+04,3.779800493522718352e-02,5.742160326143515986e-02,9.295180846712273320e-03,2.149966973851244814e-02,6.181989757149400140e-14,2.966568710219373226e-02,3.721084230037273455e-02,9.295180846712295872e-03,4.273398354281245093e-02,2.175348015276857422e-02,4.256954402870758748e-03,2.498231175051895611e-02,nan
        }\tableUniformTwo

        \pgfplotstableread[col sep=comma]{
ndof,nelem,res,err,res1,res2,res3,res4,err1,err2,err3,relerr1,relerr2,relerr3,cond
1.880000000000000000e+02,2.400000000000000000e+01,6.691468253750408879e-01,1.025606713479946475e+00,3.566616625857720324e-01,3.931095612333453393e-01,3.041234422105361890e-16,4.074491443333450591e-01,5.936417091598809037e-01,3.566616625857720324e-01,7.564728087131916690e-01,3.482206960706685228e-01,1.633434078715053916e-01,4.437348049237163083e-01,nan
5.280000000000000000e+02,6.900000000000000000e+01,4.645624378874185223e-01,6.832668232616508908e-01,2.257670682367768356e-01,2.593329842899136373e-01,5.546873237131485948e-16,3.124002126359064491e-01,4.085803248249906816e-01,2.257670682367767800e-01,4.989437852461374723e-01,2.393640097062214867e-01,1.033958303585634786e-01,2.923028295737712390e-01,nan
1.149000000000000000e+03,1.510000000000000000e+02,3.478243066800156269e-01,4.845122293408223868e-01,1.333485847715212003e-01,2.134248369352593488e-01,5.116994883147742028e-16,2.401036072949007949e-01,2.861271747751423233e-01,1.333485847715212003e-01,3.675615529070488385e-01,1.674078093646074461e-01,6.107024884929990183e-02,2.150535838729011329e-01,nan
2.309000000000000000e+03,3.050000000000000000e+02,2.501052988851075187e-01,3.432873846834183507e-01,9.575799494931515332e-02,1.476593980156844510e-01,5.837575396969208164e-16,1.777069754155463088e-01,2.055318894307546496e-01,9.575799494931513944e-02,2.577465369565462017e-01,1.201908630152851709e-01,4.385470670041589220e-02,1.507249254692682727e-01,nan
4.768000000000000000e+03,6.310000000000000000e+02,1.848449590170492818e-01,2.446796373789327794e-01,6.218422338522493209e-02,1.172054516975935129e-01,7.689342782673883796e-16,1.286998963805105201e-01,1.459982741680137808e-01,6.218422338522488352e-02,1.862411105256408517e-01,8.535931290295911689e-02,2.847877933937318848e-02,1.088876791136458161e-01,nan
8.113000000000000000e+03,1.076000000000000000e+03,1.376219804433256033e-01,1.812277480934352536e-01,5.018360990588965204e-02,8.434036183668822129e-02,8.314390473836068581e-16,9.647858914398103891e-02,1.071448330200349963e-01,5.018360990588964510e-02,1.372774079480667542e-01,6.263817300273620547e-02,2.298280611263963172e-02,8.025404292533930284e-02,nan
1.571800000000000000e+04,2.087000000000000000e+03,1.017585707064303929e-01,1.314450757595068964e-01,3.419969509520793399e-02,6.426160897857313703e-02,1.829068745353436298e-15,7.110297588072787833e-02,7.790749118939764573e-02,3.419969509520792011e-02,1.001928720687816493e-01,4.554419650897505167e-02,1.566258311747242837e-02,5.857208061296409474e-02,nan
2.910100000000000000e+04,3.869000000000000000e+03,7.424270185928677490e-02,9.580819269278005590e-02,2.610512140004428785e-02,4.610559439404571941e-02,1.472204409577290003e-15,5.200745698097181124e-02,5.607192816284454345e-02,2.610512140004428439e-02,7.316878634941695325e-02,3.277885284872231880e-02,1.195547599390749983e-02,4.277343332837994089e-02,nan
5.466400000000000000e+04,7.272000000000000000e+03,5.510682960853015983e-02,7.036346662244638561e-02,1.859894103721107528e-02,3.478193849727023029e-02,2.563988152673700875e-15,3.848452697084937812e-02,4.130798675615700299e-02,1.859894103721106487e-02,5.384001353617614422e-02,2.414793958008080998e-02,8.517837924383335826e-03,3.147394719420765208e-02,nan
1.022990000000000000e+05,1.361800000000000000e+04,4.022992597714670354e-02,5.136979560401189887e-02,1.388065729318348727e-02,2.525979458233612585e-02,1.251713687766807458e-14,2.806629784802544894e-02,3.003943215831059860e-02,1.388065729318349768e-02,3.929142106253485228e-02,1.756050066336165699e-02,6.356984995579354619e-03,2.296904355571112097e-02,nan
        }\tableAdaptiveTwo

        \pgfplotstableread[col sep=comma]{
ndof,nelem,res,err,res1,res2,res3,res4,err1,err2,err3,relerr1,relerr2,relerr3,cond
1.880000000000000000e+02,2.400000000000000000e+01,6.691468253750407769e-01,1.025606713479946253e+00,3.566616625857722545e-01,3.931095612333451728e-01,5.363595415509453929e-16,4.074491443333450036e-01,5.936417091598807927e-01,3.566616625857722545e-01,7.564728087131916690e-01,3.482206960706685228e-01,1.633434078715054472e-01,4.437348049237164194e-01,nan
7.360000000000000000e+02,9.600000000000000000e+01,4.018328538784069481e-01,5.720407191928681412e-01,1.652282559352704905e-01,2.464800014124797767e-01,4.561343960297221054e-16,2.709554849054690573e-01,3.427995284376483265e-01,1.652282559352704905e-01,4.271050118626665837e-01,2.006666701171942457e-01,7.567037209895105532e-02,2.500170899051719076e-01,nan
2.912000000000000000e+03,3.840000000000000000e+02,2.226953929342505856e-01,3.122001557070244737e-01,7.752238586897911010e-02,1.366714899126683647e-01,1.923886787109607308e-15,1.578113480342043240e-01,1.914562321967721303e-01,7.752238586897912398e-02,2.341019608326404267e-01,1.119820638399213614e-01,3.550326745488131819e-02,1.369253976337965750e-01,nan
1.158400000000000000e+04,1.536000000000000000e+03,1.214925828018438014e-01,1.733462263358582778e-01,3.772328042620470140e-02,7.315877376198681725e-02,1.137825005451366394e-14,8.935992172887811946e-02,1.084461173911630411e-01,3.772328042620473609e-02,1.298665003823180064e-01,6.340910331328292182e-02,1.727629497203184844e-02,7.593373131077095495e-02,nan
4.620800000000000000e+04,6.144000000000000000e+03,6.703768583937258496e-02,9.857353874137066152e-02,1.866840647801323860e-02,3.929790042376421488e-02,1.715917990038685395e-15,5.100212686304367304e-02,6.283363054742975573e-02,1.866840647801323166e-02,7.362179033273559547e-02,3.673446005538125864e-02,8.549651322052180188e-03,4.304154785616087775e-02,nan
1.845760000000000000e+05,2.457600000000000000e+04,3.779800493522055688e-02,5.742160326144190446e-02,9.295180846769196537e-03,2.149966973850532190e-02,4.536360263089410609e-14,2.966568710217262761e-02,3.721084230037267904e-02,9.295180846769198271e-03,4.273398354280918271e-02,2.175348015276852912e-02,4.256954402896819499e-03,2.498231175051702710e-02,nan
        }\tableUniformThree

        \pgfplotstableread[col sep=comma]{
ndof,nelem,res,err,res1,res2,res3,res4,err1,err2,err3,relerr1,relerr2,relerr3,cond
1.880000000000000000e+02,2.400000000000000000e+01,6.691468253750407769e-01,1.025606713479946253e+00,3.566616625857722545e-01,3.931095612333451728e-01,5.363595415509453929e-16,4.074491443333450036e-01,5.936417091598807927e-01,3.566616625857722545e-01,7.564728087131916690e-01,3.482206960706685228e-01,1.633434078715054472e-01,4.437348049237164194e-01,nan
5.280000000000000000e+02,6.900000000000000000e+01,4.645504417268895847e-01,6.830605420970884412e-01,2.257866553728564984e-01,2.594642995420956555e-01,5.194272732835963824e-16,3.122591494693419678e-01,4.085020935114557772e-01,2.257866553728564984e-01,4.987164826053628874e-01,2.393181784222499497e-01,1.034048007908551364e-01,2.921696658646869560e-01,nan
1.134000000000000000e+03,1.490000000000000000e+02,3.518103722154240498e-01,4.904108312726485619e-01,1.301980667762059063e-01,2.153182715706491168e-01,5.709838553646765492e-16,2.458801401779371898e-01,2.883968379328984422e-01,1.301980667762058508e-01,3.746712034672192160e-01,1.687357481202074183e-01,5.962739200106293574e-02,2.192133113153246216e-01,nan
2.340000000000000000e+03,3.090000000000000000e+02,2.508498778676586105e-01,3.427315548724161332e-01,9.752950316736931635e-02,1.476779233623533027e-01,6.476971636419217874e-16,1.777776369218812824e-01,2.037485691937860610e-01,9.752950316736931635e-02,2.577584824378725825e-01,1.191480136582670352e-01,4.466601205020044596e-02,1.507319109408424007e-01,nan
4.978000000000000000e+03,6.590000000000000000e+02,1.808319377800240990e-01,2.413511010653371647e-01,6.134251940479320403e-02,1.117696545863764296e-01,6.547309484603931063e-16,1.282373944993939197e-01,1.434405665489926207e-01,6.134251940479314158e-02,1.841527984219278880e-01,8.386392423390900963e-02,2.809330050529166775e-02,1.076667271043160873e-01,nan
8.699000000000000000e+03,1.154000000000000000e+03,1.344346977620795303e-01,1.779605361967838417e-01,4.884125103040942023e-02,8.066437834016768060e-02,1.279946473836052365e-15,9.581481214906434341e-02,1.054508137401535123e-01,4.884125103040943411e-02,1.347761496807606540e-01,6.164782869910297169e-02,2.236804017957296289e-02,7.879177688204436203e-02,nan
1.653500000000000000e+04,2.196000000000000000e+03,9.952843663252358697e-02,1.285986834228695364e-01,3.272287639142416882e-02,6.360375860508357682e-02,9.660385394344731918e-16,6.920754980847172722e-02,7.644982682337005120e-02,3.272287639142414800e-02,9.809311248728795618e-02,4.469205569670890721e-02,1.498623803214425149e-02,5.734457524507118736e-02,nan
3.021900000000000000e+04,4.018000000000000000e+03,7.299357091303500245e-02,9.459466400747636172e-02,2.602412150952513117e-02,4.461038869502792492e-02,8.154122566641933895e-15,5.158216469645555968e-02,5.589007203276139674e-02,2.602412150952516240e-02,7.174395728225169133e-02,3.267254234074398744e-02,1.191838012173095963e-02,4.194049849538335206e-02,nan
5.741700000000000000e+04,7.639000000000000000e+03,5.384308662052388061e-02,6.863735113876702487e-02,1.763472364020083824e-02,3.450045980945933011e-02,1.883350083294694088e-14,3.738733437803316850e-02,4.041851604337635961e-02,1.763472364020082783e-02,5.259701564092276321e-02,2.362797027833591842e-02,8.076251089134359185e-03,3.074731197353178583e-02,nan
1.072600000000000000e+05,1.427900000000000000e+04,3.943988739115177750e-02,5.042059190781043865e-02,1.382792212250474884e-02,2.449631931583004171e-02,9.127650768473686462e-15,2.764459489986367813e-02,2.961721065122177934e-02,1.382792212250475057e-02,3.839069537464850912e-02,1.731367772017761153e-02,6.332833640087852899e-03,2.244249585140313877e-02,nan
        }\tableAdaptiveThree

        \pgfplotstableread[col sep=comma]{
ndof,nelem,res,err,res1,res2,res3,res4,err1,err2,err3,relerr1,relerr2,relerr3,cond
1.880000000000000000e+02,2.400000000000000000e+01,6.691468253750406658e-01,1.025606713479946253e+00,3.566616625857723655e-01,3.931095612333451728e-01,3.355125713439767860e-16,4.074491443333447260e-01,5.936417091598806817e-01,3.566616625857723655e-01,7.564728087131914469e-01,3.482206960706685783e-01,1.633434078715055027e-01,4.437348049237164194e-01,nan
7.360000000000000000e+02,9.600000000000000000e+01,4.018328538784069481e-01,5.720407191928683632e-01,1.652282559352701019e-01,2.464800014124798044e-01,3.834987345788232190e-16,2.709554849054692238e-01,3.427995284376486596e-01,1.652282559352701019e-01,4.271050118626666947e-01,2.006666701171944955e-01,7.567037209895088878e-02,2.500170899051720186e-01,nan
2.912000000000000000e+03,3.840000000000000000e+02,2.226953929342508076e-01,3.122001557070244182e-01,7.752238586897775008e-02,1.366714899126685867e-01,5.315586666598378399e-16,1.578113480342051012e-01,1.914562321967724356e-01,7.752238586897775008e-02,2.341019608326405654e-01,1.119820638399215557e-01,3.550326745488069369e-02,1.369253976337966583e-01,nan
1.158400000000000000e+04,1.536000000000000000e+03,1.214925828018438292e-01,1.733462263358544753e-01,3.772328042620518712e-02,7.315877376198780258e-02,9.109714234920839285e-16,8.935992172887714802e-02,1.084461173911596688e-01,3.772328042620522182e-02,1.298665003823156194e-01,6.340910331328095118e-02,1.727629497203206702e-02,7.593373131076955329e-02,nan
4.620800000000000000e+04,6.144000000000000000e+03,6.703768583937237679e-02,9.857353874141237815e-02,1.866840647803051298e-02,3.929790042375266162e-02,1.114832512314736760e-14,5.100212686304598370e-02,6.283363054746174403e-02,1.866840647803048175e-02,7.362179033275977058e-02,3.673446005539997977e-02,8.549651322060078384e-03,4.304154785617503309e-02,nan
1.845760000000000000e+05,2.457600000000000000e+04,3.779800493524770877e-02,5.742160326147267152e-02,9.295180846695616506e-03,2.149966973849944465e-02,7.922690462875794558e-14,2.966568710223453295e-02,3.721084230039533453e-02,9.295180846695616506e-03,4.273398354284679845e-02,2.175348015278176853e-02,4.256954402863119893e-03,2.498231175053901298e-02,nan
        }\tableUniformFour

        \pgfplotstableread[col sep=comma]{
ndof,nelem,res,err,res1,res2,res3,res4,err1,err2,err3,relerr1,relerr2,relerr3,cond
1.880000000000000000e+02,2.400000000000000000e+01,6.691468253750406658e-01,1.025606713479946253e+00,3.566616625857723655e-01,3.931095612333451728e-01,3.355125713439767860e-16,4.074491443333447260e-01,5.936417091598806817e-01,3.566616625857723655e-01,7.564728087131914469e-01,3.482206960706685783e-01,1.633434078715055027e-01,4.437348049237164194e-01,nan
5.280000000000000000e+02,6.900000000000000000e+01,4.645504417268894737e-01,6.830605420970886632e-01,2.257866553728564429e-01,2.594642995420956555e-01,3.551904851646315441e-16,3.122591494693418568e-01,4.085020935114558327e-01,2.257866553728564429e-01,4.987164826053629429e-01,2.393181784222500053e-01,1.034048007908551087e-01,2.921696658646870670e-01,nan
1.134000000000000000e+03,1.490000000000000000e+02,3.518103722154241608e-01,4.904108312726487284e-01,1.301980667762053234e-01,2.153182715706493111e-01,4.167527449078184146e-16,2.458801401779374673e-01,2.883968379328986087e-01,1.301980667762053234e-01,3.746712034672193825e-01,1.687357481202075016e-01,5.962739200106269982e-02,2.192133113153246771e-01,nan
2.340000000000000000e+03,3.090000000000000000e+02,2.508498778676588881e-01,3.427315548724158001e-01,9.752950316736823388e-02,1.476779233623536081e-01,5.454254263550114924e-16,1.777776369218819486e-01,2.037485691937860055e-01,9.752950316736828940e-02,2.577584824378725270e-01,1.191480136582670213e-01,4.466601205019997411e-02,1.507319109408424007e-01,nan
4.978000000000000000e+03,6.590000000000000000e+02,1.808319377800244598e-01,2.413511010653376920e-01,6.134251940479252402e-02,1.117696545863762908e-01,6.097016127272388116e-16,1.282373944993948633e-01,1.434405665489932313e-01,6.134251940479248932e-02,1.841527984219282488e-01,8.386392423390937045e-02,2.809330050529137285e-02,1.076667271043163093e-01,nan
8.699000000000000000e+03,1.154000000000000000e+03,1.344346977620793637e-01,1.779605361967694366e-01,4.884125103040672100e-02,8.066437834016962349e-02,7.681004107247820519e-16,9.581481214906387156e-02,1.054508137401407863e-01,4.884125103040665855e-02,1.347761496807525494e-01,6.164782869909550544e-02,2.236804017957169308e-02,7.879177688203960195e-02,nan
1.653500000000000000e+04,2.196000000000000000e+03,9.952843663252258777e-02,1.285986834228668441e-01,3.272287639142192062e-02,6.360375860508493684e-02,8.899438143247716728e-16,6.920754980847008964e-02,7.644982682336945445e-02,3.272287639142189980e-02,9.809311248728563859e-02,4.469205569670856026e-02,1.498623803214321759e-02,5.734457524506983428e-02,nan
3.021900000000000000e+04,4.018000000000000000e+03,7.299357091303204648e-02,9.459466400746124881e-02,2.602412150952882267e-02,4.461038869502938903e-02,1.951253959521169695e-15,5.158216469644823915e-02,5.589007203274805324e-02,2.602412150952884348e-02,7.174395728224083890e-02,3.267254234073619507e-02,1.191838012173264405e-02,4.194049849537701685e-02,nan
5.741700000000000000e+04,7.639000000000000000e+03,5.384308662052199324e-02,6.863735113878710603e-02,1.763472364020871042e-02,3.450045980945515983e-02,2.709471205164165117e-15,3.738733437803059417e-02,4.041851604339428278e-02,1.763472364020870695e-02,5.259701564093254011e-02,2.362797027834641003e-02,8.076251089137965675e-03,3.074731197353751735e-02,nan
1.072600000000000000e+05,1.427900000000000000e+04,3.943988739104545282e-02,5.042059190746867037e-02,1.382792212254327704e-02,2.449631931588659717e-02,8.700876694294158122e-15,2.764459489964260150e-02,2.961721065096148062e-02,1.382792212254333949e-02,3.839069537438655894e-02,1.731367772002544506e-02,6.332833640105525394e-03,2.244249585125000432e-02,nan
        }\tableAdaptiveFour

        %
        %
        \addplot+ [line1, minorline, forget plot] table [x=ndof, y=relerr1] {\tableUniformOne};
        \addplot+ [line2, minorline, forget plot] table [x=ndof, y=relerr1] {\tableUniformTwo};
        \addplot+ [line3, minorline, forget plot] table [x=ndof, y=relerr1] {\tableUniformThree};
        \addplot+ [line4, minorline, forget plot] table [x=ndof, y=relerr1] {\tableUniformFour};

        \addplot+ [line1, majorline, forget plot] table [x=ndof, y=relerr1] {\tableAdaptiveOne};
        \addplot+ [line2, majorline, forget plot] table [x=ndof, y=relerr1] {\tableAdaptiveTwo};
        \addplot+ [line3, majorline, forget plot] table [x=ndof, y=relerr1] {\tableAdaptiveThree};
        \addplot+ [line4, majorline, forget plot] table [x=ndof, y=relerr1] {\tableAdaptiveFour};

        %
        %
        \drawslopetriangle[ST1]{0.5}{2e3}{3e-2}
        \drawswappedslopetriangle[ST2]{0.33}{1e5}{8e-2}
    \end{loglogaxis}
\end{tikzpicture}

%% file: figures/legend_scal.tex
\begin{tikzpicture}[>=stealth]
    %
    %
    \colorlet{col1}{TUblue}
    \colorlet{col2}{TUgreen}
    \colorlet{col3}{TUmagenta}
    \colorlet{col4}{TUyellow}
    \colorlet{col5}{purple}
    \colorlet{col6}{green}
    \pgfplotsset{%
        linedefault/.style = {%
            mark = *,%
            mark size = 2pt,%
            every mark/.append style = {solid},%
            gray,%
            every mark/.append style = {fill = gray!60!white}%
        },%
        line1/.style = {%
            linedefault,%
            col1,%
            every mark/.append style = {fill = col1!60!white}%
        },%
        line2/.style = {%
            linedefault,%
            mark = triangle*,%
            mark size = 2.75pt,%
            col2,%
            every mark/.append style = {fill = col2!60!white}%
        },%
        line3/.style = {%
            linedefault,%
            mark = square*,%
            mark size = 1.66pt,%
            col3,%
            every mark/.append style = {fill = col3!60!white}%
        },%
        line4/.style = {%
            linedefault,%
            mark = pentagon*,%
            mark size = 2.2pt,%
            col4,%
            every mark/.append style = {fill = col4!60!white}%
        },%
        line5/.style = {%
            linedefault,%
            mark = diamond*,%
            mark size = 2.75pt,%
            col5,%
            every mark/.append style = {fill = col5!60!white}%
        },%
        line6/.style = {%
            linedefault,%
            mark = halfsquare*,%
            mark size = 1.66pt,%
            col6,%
            every mark/.append style = {fill = col6!60!white}%
        },%
        minorline/.style = {%
            dashed,%
            every mark/.append style = {fill = black!20!white}%
        },%
        majorline/.style = {%
            solid%
        }%
    }

    \matrix [
        matrix of nodes,
        anchor = south,
        font = \scriptsize,
        column 1/.style={anchor=base east},
    ] at (0,0) {
        & uniform
        & adaptive
        \\
        & \(\theta = 1\)
        & \(\theta = 0.5\)
        \\
        \hline \\
        \(\ell=\phantom{0}1\phantom{{}^2}\)
        & \ref*{leg:scal:unif:len1}
        & \ref*{leg:scal:adap:len1}
        \\
        \(\ell =10\phantom{{}^2}\)
        & \ref*{leg:scal:unif:len10}
        & \ref*{leg:scal:adap:len10}
        \\
        \(\ell = 10^2\)
        & \ref*{leg:scal:unif:len100}
        & \ref*{leg:scal:adap:len100}
        \\
        \(\ell = 10^3\)
        & \ref*{leg:scal:unif:len1000}
        & \ref*{leg:scal:adap:len1000}
        \\
    };
\end{tikzpicture}

%% file: figures/plot_dls_lshape_mesh.tex
\begin{tikzpicture}
    \begin{axis}[%
        axis equal image,%
        width=4.5cm,%
        xmin=-1.12, xmax=1.12,%
        ymin=-1.12, ymax=1.12,%
        font=\footnotesize%
    ]

        \addplot graphics [xmin=-1, xmax=1, ymin=-1, ymax=1]
        {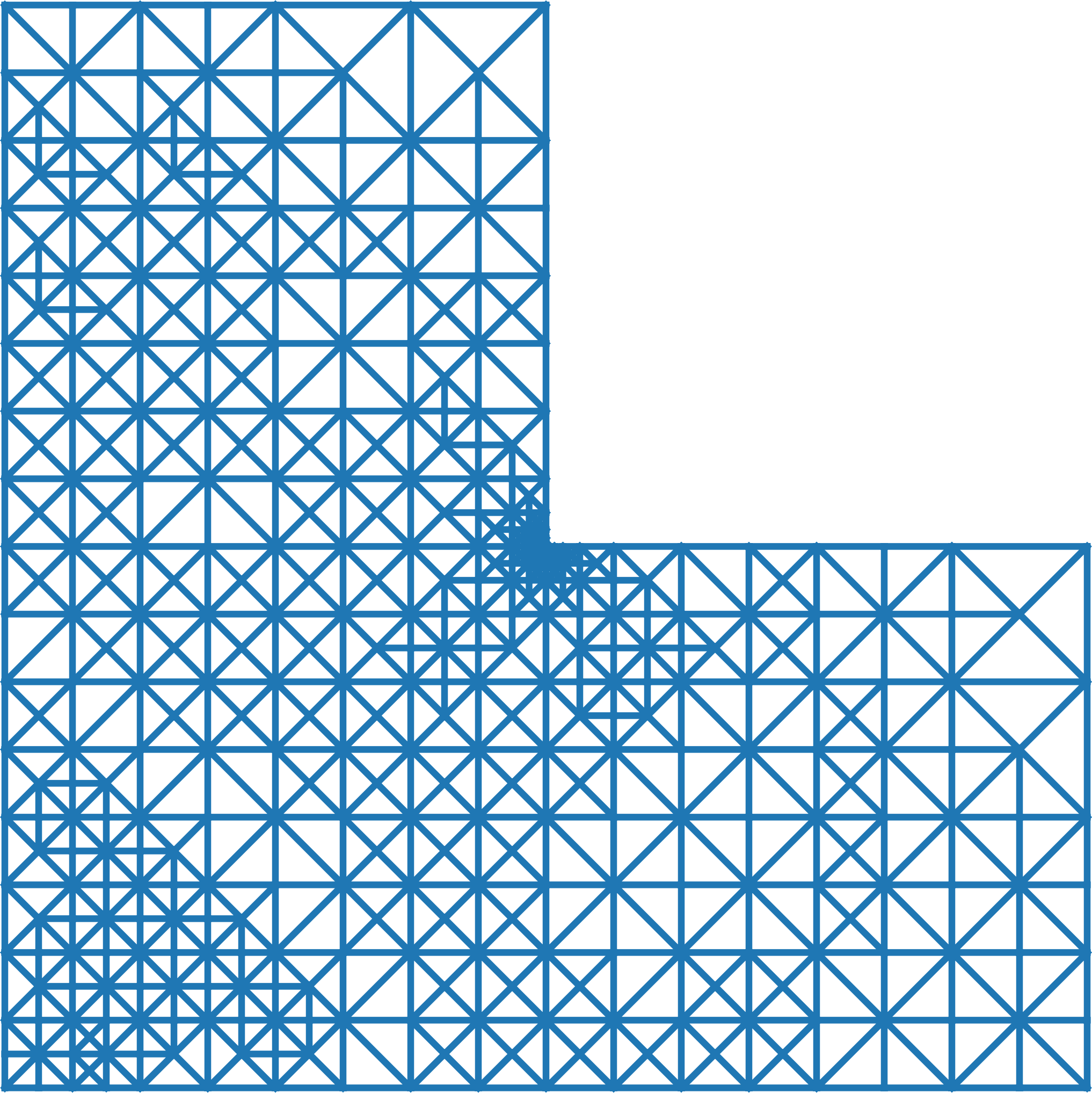};
    \end{axis}
\end{tikzpicture}

%% file: figures/legend.tex
\begin{tikzpicture}[>=stealth]
    %
    %
    \colorlet{col1}{TUblue}
    \colorlet{col2}{TUgreen}
    \colorlet{col3}{TUmagenta}
    \colorlet{col4}{TUyellow}
    \colorlet{col5}{purple}
    \colorlet{col6}{green}
    \pgfplotsset{%
        linedefault/.style = {%
            mark = *,%
            mark size = 2pt,%
            every mark/.append style = {solid},%
            gray,%
            every mark/.append style = {fill = gray!60!white}%
        },%
        line1/.style = {%
            linedefault,%
            col1,%
            every mark/.append style = {fill = col1!60!white}%
        },%
        line2/.style = {%
            linedefault,%
            mark = triangle*,%
            mark size = 2.75pt,%
            col2,%
            every mark/.append style = {fill = col2!60!white}%
        },%
        line3/.style = {%
            linedefault,%
            mark = square*,%
            mark size = 1.66pt,%
            col3,%
            every mark/.append style = {fill = col3!60!white}%
        },%
        line4/.style = {%
            linedefault,%
            mark = pentagon*,%
            mark size = 2.2pt,%
            col4,%
            every mark/.append style = {fill = col4!60!white}%
        },%
        line5/.style = {%
            linedefault,%
            mark = diamond*,%
            mark size = 2.75pt,%
            col5,%
            every mark/.append style = {fill = col5!60!white}%
        },%
        line6/.style = {%
            linedefault,%
            mark = halfsquare*,%
            mark size = 1.66pt,%
            col6,%
            every mark/.append style = {fill = col6!60!white}%
        },%
        minorline/.style = {%
            dashed,%
            every mark/.append style = {fill = black!20!white}%
        },%
        majorline/.style = {%
            solid%
        }%
    }

    \matrix [
        matrix of nodes,
        anchor = south,
        font = \scriptsize,
        column 1/.style={anchor=base east},
    ] at (0,0) {
        & uniform
        & adaptive
        \\
        & \(\theta = 1\)
        & \(\theta = 0.5\)
        \\
        \hline \\
        \(k=0\)
        & \ref*{leg:est:unif:0}
        & \ref*{leg:est:adap:0}
        \\
        \(k=1\)
        & \ref*{leg:est:unif:1}
        & \ref*{leg:est:adap:1}
        \\
        \(k=2\)
        & \ref*{leg:est:unif:2}
        & \ref*{leg:est:adap:2}
        \\
        \(k=3\)
        & \ref*{leg:est:unif:3}
        & \ref*{leg:est:adap:3}
        \\
    };
\end{tikzpicture}

%% file: figures/plot_dls_lshape_estimator_degree.tex
\begin{tikzpicture}[>=stealth]
    %
    %
    \colorlet{col1}{TUblue}
    \colorlet{col2}{TUgreen}
    \colorlet{col3}{TUmagenta}
    \colorlet{col4}{TUyellow}
    \colorlet{col5}{purple}
    \colorlet{col6}{green}
    \pgfplotsset{%
        linedefault/.style = {%
            mark = *,%
            mark size = 2pt,%
            every mark/.append style = {solid},%
            gray,%
            every mark/.append style = {fill = gray!60!white}%
        },%
        line1/.style = {%
            linedefault,%
            col1,%
            every mark/.append style = {fill = col1!60!white}%
        },%
        line2/.style = {%
            linedefault,%
            mark = triangle*,%
            mark size = 2.75pt,%
            col2,%
            every mark/.append style = {fill = col2!60!white}%
        },%
        line3/.style = {%
            linedefault,%
            mark = square*,%
            mark size = 1.66pt,%
            col3,%
            every mark/.append style = {fill = col3!60!white}%
        },%
        line4/.style = {%
            linedefault,%
            mark = pentagon*,%
            mark size = 2.2pt,%
            col4,%
            every mark/.append style = {fill = col4!60!white}%
        },%
        line5/.style = {%
            linedefault,%
            mark = diamond*,%
            mark size = 2.75pt,%
            col5,%
            every mark/.append style = {fill = col5!60!white}%
        },%
        line6/.style = {%
            linedefault,%
            mark = halfsquare*,%
            mark size = 1.66pt,%
            col6,%
            every mark/.append style = {fill = col6!60!white}%
        },%
        minorline/.style = {%
            dashed,%
            every mark/.append style = {fill = black!20!white}%
        },%
        majorline/.style = {%
            solid%
        }%
    }

    %
    %
    \begin{loglogaxis}[%
            width            = 0.36\textwidth,%
            xlabel           = ndof,%
            ylabel           = {built-in error estimator},%
            ymax             = 1e0,%
            ymin             = 3e-6,%
            ymajorgrids      = true,%
            font             = \footnotesize,%
            grid style       = {densely dotted, semithick},%
            legend style     = {legend pos  = south west}%
        ]

        %
        %
        \pgfplotstableread[col sep=comma]{
ndof,nelem,res,err,res1,res2,res3,res4,err1,err2,cond
1.880000000000000000e+02,2.400000000000000000e+01,6.691468253750407769e-01,1.025606713479946253e+00,3.566616625857725320e-01,3.931095612333450617e-01,2.393892856402309924e-16,4.074491443333447815e-01,6.925445981398684614e-01,7.564728087131914469e-01,1.861844080681640321e+02
7.360000000000000000e+02,9.600000000000000000e+01,4.018328538784070036e-01,5.720407191928681412e-01,1.652282559352703517e-01,2.464800014124796934e-01,3.594997419828183740e-16,2.709554849054693904e-01,3.805415788800027288e-01,4.271050118626665282e-01,2.158354393424532702e+02
2.912000000000000000e+03,3.840000000000000000e+02,2.226953929342510019e-01,3.122001557070245847e-01,7.752238586897748640e-02,1.366714899126683647e-01,4.720870189556533655e-16,1.578113480342057118e-01,2.065555837003762440e-01,2.341019608326406487e-01,2.354149255634362362e+03
1.158400000000000000e+04,1.536000000000000000e+03,1.214925828018451198e-01,1.733462263358546418e-01,3.772328042619833843e-02,7.315877376198784421e-02,7.058621229607447708e-16,8.935992172888182483e-02,1.148198861841061846e-01,1.298665003823167019e-01,3.452532552691214369e+04
4.620800000000000000e+04,6.144000000000000000e+03,6.703768583937476377e-02,9.857353874140176164e-02,1.866840647801305472e-02,3.929790042375498615e-02,1.332844287066715975e-15,5.100212686305377607e-02,6.554826106162134347e-02,7.362179033275560724e-02,5.364942042558327084e+05
1.845760000000000000e+05,2.457600000000000000e+04,3.779800493524370503e-02,5.742160326142143473e-02,9.295180846666025593e-03,2.149966973851240304e-02,5.681306352149327939e-15,2.966568710222893326e-02,3.835423277390073021e-02,4.273398354281396361e-02,8.490368912294194102e+06
        }\tableUniformZero

        \pgfplotstableread[col sep=comma]{
ndof,nelem,res,err,res1,res2,res3,res4,err1,err2,cond
1.880000000000000000e+02,2.400000000000000000e+01,6.691468253750407769e-01,1.025606713479946253e+00,3.566616625857725320e-01,3.931095612333450617e-01,2.393892856402309924e-16,4.074491443333447815e-01,6.925445981398684614e-01,7.564728087131914469e-01,1.861844080681634352e+02
5.280000000000000000e+02,6.900000000000000000e+01,4.645485304210772437e-01,6.830599490857015565e-01,2.257879949253541185e-01,2.594660624333418597e-01,4.221377968243194286e-16,3.122538725289212835e-01,4.667482632440184132e-01,4.987153023556160636e-01,2.128427996102370798e+02
1.134000000000000000e+03,1.490000000000000000e+02,3.514277916105422728e-01,4.898292054561634856e-01,1.302675895177524967e-01,2.157932520629477768e-01,3.997275055233180919e-16,2.448777658375585753e-01,3.163755221592782219e-01,3.739507714877192690e-01,1.106585141813846803e+03
2.302000000000000000e+03,3.040000000000000000e+02,2.501417772339131207e-01,3.416849949294442634e-01,9.756774828116387965e-02,1.481552895554908156e-01,5.184071061221481349e-16,1.763560415461084552e-01,2.258489202343929025e-01,2.563998849237131838e-01,4.320766697536862011e+03
4.850000000000000000e+03,6.420000000000000000e+02,1.827162881006695605e-01,2.434834612677447552e-01,6.173245639320502215e-02,1.138115980365454472e-01,5.550654254689614395e-16,1.289234886175351369e-01,1.572760924728729692e-01,1.858720706491045760e-01,3.002650222874309475e+04
8.519000000000000000e+03,1.130000000000000000e+03,1.351938722597004761e-01,1.787860591587285553e-01,4.926933115701741944e-02,8.141275564107232354e-02,6.863274354147775619e-16,9.603061659095710190e-02,1.167796488402617627e-01,1.353771345769053391e-01,2.195661886394199682e+05
1.637800000000000000e+04,2.175000000000000000e+03,1.001953397017019687e-01,1.293968314633148209e-01,3.316026559924474537e-02,6.387072948000693806e-02,1.377493817389158267e-15,6.971393547217655251e-02,8.359794490892365204e-02,9.876701676040477995e-02,1.526050124626174103e+06
2.965100000000000000e+04,3.942000000000000000e+03,7.366358106374841508e-02,9.533914920550538852e-02,2.635424251199278492e-02,4.496798373161311341e-02,1.768221127331582023e-14,5.205437076622775744e-02,6.185394011254979896e-02,7.255097148751847580e-02,1.134321577537557483e+07
5.533200000000000000e+04,7.361000000000000000e+03,5.476430034674915903e-02,6.993875022571281741e-02,1.815836851153471154e-02,3.483691541012369497e-02,7.663785009417625751e-15,3.815483678612911073e-02,4.490324281781864907e-02,5.362021603442729961e-02,7.710566014080482721e+07
1.020060000000000000e+05,1.357800000000000000e+04,4.041410942974788673e-02,5.163485719064491580e-02,1.401685756070034962e-02,2.523231564524061263e-02,2.868700345761843395e-15,2.828706758046184433e-02,3.331128810879161600e-02,3.945271298188963299e-02,5.550838271046789885e+08
        }\tableAdaptiveZero

        \pgfplotstableread[col sep=comma]{
ndof,nelem,res,err,res1,res2,res3,res4,err1,err2,cond
3.800000000000000000e+02,2.400000000000000000e+01,1.429544905823191270e-01,2.649557076419757906e-01,5.816283848669231982e-02,7.748609674522376201e-02,3.101857475521297217e-02,1.004331789209442632e-01,1.843098314644586733e-01,1.903455096334110419e-01,7.220823296696265061e+03
1.504000000000000000e+03,9.600000000000000000e+01,6.043023319836864776e-02,1.290671872306016765e-01,1.468726269570548582e-02,2.915852172614483498e-02,1.384234833024492080e-02,4.893125179644047112e-02,8.961392183059026240e-02,9.288532625641167384e-02,6.952862152031728328e+03
5.984000000000000000e+03,3.840000000000000000e+02,3.299879050228264610e-02,7.757485919926983275e-02,3.983931087607786060e-03,1.348187200302555962e-02,7.148377534395108056e-03,2.898600168395982710e-02,5.374083191374660895e-02,5.594445249535458242e-02,1.313286401135960114e+04
2.387200000000000000e+04,1.536000000000000000e+03,2.016212199610784711e-02,4.857279130894310432e-02,1.137728940393586529e-03,7.747887232526820024e-03,4.262654840198620149e-03,1.808360454637036180e-02,3.361484377301553195e-02,3.506220662847531200e-02,5.043913348106794001e+04
9.536000000000000000e+04,6.144000000000000000e+03,1.263390460609820944e-02,3.058304462746423355e-02,3.414287908117189093e-04,4.792162396347206775e-03,2.655905513844827003e-03,1.137894184207646243e-02,2.115795523764861172e-02,2.208310550731320404e-02,5.673349497205432272e+05
3.811840000000000000e+05,2.457600000000000000e+04,7.952201830990024076e-03,1.926739112649915667e-02,1.078196987798035983e-04,3.009227955488317108e-03,1.669733750245612882e-03,7.168153541083874779e-03,1.332828378032429109e-02,1.391363404336342334e-02,8.972909084775183350e+06
        }\tableUniformOne

        \pgfplotstableread[col sep=comma]{
ndof,nelem,res,err,res1,res2,res3,res4,err1,err2,cond
3.800000000000000000e+02,2.400000000000000000e+01,1.429544905823191270e-01,2.649557076419757906e-01,5.816283848669231982e-02,7.748609674522376201e-02,3.101857475521297217e-02,1.004331789209442632e-01,1.843098314644586733e-01,1.903455096334110419e-01,7.220823296696225043e+03
9.410000000000000000e+02,6.000000000000000000e+01,9.048552861572219108e-02,1.723899622309089663e-01,3.181784498837661457e-02,4.657744344487177446e-02,1.919519511623793764e-02,6.809803136472918073e-02,1.177978602354333021e-01,1.258648608704095451e-01,6.575537898597228377e+03
1.657000000000000000e+03,1.060000000000000000e+02,5.461672803270744420e-02,9.889981584593049435e-02,1.658444428552899308e-02,3.009022628867067187e-02,1.358916670968292502e-02,4.022258095612405576e-02,6.755965415739066671e-02,7.222788038211234996e-02,1.251520602726039579e+04
2.682000000000000000e+03,1.720000000000000000e+02,3.906428975490600780e-02,6.727301167620299005e-02,1.348658806205335625e-02,2.185195330908825881e-02,9.470620031590409882e-03,2.787346712008392097e-02,4.586435459046426383e-02,4.921502898492190886e-02,4.766783496311864292e+04
4.142000000000000000e+03,2.660000000000000000e+02,2.736237289062587452e-02,4.587644150401347376e-02,8.801945939295292842e-03,1.605758390699676616e-02,6.791470739562765167e-03,1.916390153780805905e-02,3.102208294304604000e-02,3.379760723699156560e-02,1.891961357416887186e+05
5.528000000000000000e+03,3.550000000000000000e+02,1.789686763959938420e-02,2.870725154843154631e-02,5.745194310199929454e-03,1.075726856481438488e-02,4.467533532744634092e-03,1.231312027089981158e-02,1.955759948688892344e-02,2.101443774587767252e-02,7.753271398940539220e+05
8.291000000000000000e+03,5.330000000000000000e+02,1.249295063885906901e-02,1.957246137199181862e-02,4.041650490638224423e-03,7.515471054698831532e-03,3.204933307625098659e-03,8.543124387257416269e-03,1.328974480978349117e-02,1.436885266988788788e-02,3.163988872234159615e+06
1.327300000000000000e+04,8.540000000000000000e+02,8.248363375542640263e-03,1.274733939861605672e-02,2.359825262671173731e-03,5.086567941530156095e-03,2.232341650879941623e-03,5.622294943757080200e-03,8.595629524193595244e-03,9.413278643826525696e-03,1.264977298992297426e+07
1.897100000000000000e+04,1.221000000000000000e+03,5.447211056598747804e-03,8.169662857762239119e-03,1.640798505876144529e-03,3.400699778372295954e-03,1.461327451629655010e-03,3.644125636519470173e-03,5.506712451709151399e-03,6.034857842873450964e-03,5.094934276248360425e+07
2.705600000000000000e+04,1.742000000000000000e+03,3.795429211931716871e-03,5.552369508253003627e-03,1.104586523556202798e-03,2.397035874636590201e-03,1.060227549986266339e-03,2.513027670503841828e-03,3.756613471384596080e-03,4.088601519198158292e-03,2.037979702250770628e+08
4.257200000000000000e+04,2.742000000000000000e+03,2.521971172035490042e-03,3.633188819655926768e-03,7.325091624243566281e-04,1.608443598106142937e-03,7.118381972641818880e-04,1.652260418992561146e-03,2.439987182271654599e-03,2.691936765531991772e-03,1.021559285691922307e+09
6.434700000000000000e+04,4.146000000000000000e+03,1.654286607367381470e-03,2.369957368777833584e-03,4.713215955699466271e-04,1.054314049770528616e-03,4.576109251289303898e-04,1.092489935204713674e-03,1.582971914287373098e-03,1.763773752044664120e-03,5.675444554343406677e+09
9.747400000000000000e+04,6.282000000000000000e+03,1.065468831961948956e-03,1.513301870864532713e-03,3.138129662588048040e-04,6.784876866757804433e-04,3.005773986150441900e-04,6.971749712671648433e-04,1.023324527077856773e-03,1.114849525560725153e-03,3.158204548207907104e+10
1.539850000000000000e+05,9.926000000000000000e+03,7.039262121589830992e-04,9.880828360568061873e-04,2.096265624133581828e-04,4.513101998226642887e-04,1.994464939610840460e-04,4.561896711653064566e-04,6.640247352533644150e-04,7.316958670662034297e-04,2.459662173605694580e+11
        }\tableAdaptiveOne

        \pgfplotstableread[col sep=comma]{
ndof,nelem,res,err,res1,res2,res3,res4,err1,err2,cond
6.440000000000000000e+02,2.400000000000000000e+01,3.847702620941702983e-02,1.131196847568739644e-01,5.743021451740769855e-03,2.150090107081665050e-02,8.633463227888469455e-03,3.017737258875701281e-02,8.061772967664361522e-02,7.935266058096705766e-02,7.006228862074787321e+04
2.560000000000000000e+03,9.600000000000000000e+01,2.273908253830363504e-02,7.156810052386690379e-02,1.763513034254757640e-03,1.086925323657852820e-02,4.748736181433878267e-03,1.932006043520741873e-02,5.084724758081463070e-02,5.036417800430836866e-02,6.799243010727243382e+04
1.020800000000000000e+04,3.840000000000000000e+02,1.425692585739921775e-02,4.518126249798776195e-02,5.772156993588011854e-04,6.730218126621295531e-03,2.949944689420091245e-03,1.220363663984028700e-02,3.207450208855686102e-02,3.182094902235373807e-02,7.929658784635017219e+04
4.076800000000000000e+04,1.536000000000000000e+03,8.976453842815396156e-03,2.848202059630345478e-02,1.929824724181147897e-04,4.233353394718327931e-03,1.853763013109814258e-03,7.692968418249724466e-03,2.021604058194495512e-02,2.006332974452172332e-02,3.042834676585858106e+05
1.629440000000000000e+05,6.144000000000000000e+03,5.654525003601752232e-03,1.794713443837424663e-02,6.590397029447043066e-05,2.666520798019271239e-03,1.167026285081042801e-03,4.847373120228893120e-03,1.273806090224119739e-02,1.264284141322175904e-02,1.205355947972324211e+06
        }\tableUniformTwo

        \pgfplotstableread[col sep=comma]{
ndof,nelem,res,err,res1,res2,res3,res4,err1,err2,cond
6.440000000000000000e+02,2.400000000000000000e+01,3.847702620941702983e-02,1.131196847568739644e-01,5.743021451740769855e-03,2.150090107081665050e-02,8.633463227888469455e-03,3.017737258875701281e-02,8.061772967664361522e-02,7.935266058096705766e-02,7.006228862074839708e+04
1.228000000000000000e+03,4.600000000000000000e+01,2.766348712976239577e-02,7.671488518218108621e-02,3.618531229729282823e-03,1.529989583881437819e-02,6.874136424179911123e-03,2.169871396330199168e-02,5.410487767489381705e-02,5.438598919115111147e-02,7.266397462063950661e+04
1.812000000000000000e+03,6.800000000000000000e+01,2.191528436873128871e-02,5.319009979447197889e-02,2.469583004351433731e-03,1.128592104956859139e-02,4.423009160546887732e-03,1.808993711364413268e-02,3.632871834469516420e-02,3.885113820183003219e-02,8.583580719178859727e+04
2.608000000000000000e+03,9.800000000000000000e+01,1.577085774849998692e-02,3.565247261401977008e-02,2.173384840600311139e-03,9.417338822730016851e-03,3.355471988001987142e-03,1.200212021659951764e-02,2.412096496923229680e-02,2.625410163091661947e-02,3.320778297554539749e+05
3.537000000000000000e+03,1.330000000000000000e+02,1.017091125036426527e-02,2.321000291408364946e-02,1.539937484950214437e-03,5.852907755644591245e-03,2.028254811531791588e-03,7.918691901358779553e-03,1.576822328544433618e-02,1.703136429333078003e-02,1.311948035012585111e+06
4.651000000000000000e+03,1.750000000000000000e+02,8.092902219110827480e-03,1.695784923961482851e-02,1.510838194245316471e-03,5.193220960524437026e-03,1.645009350504783195e-03,5.791099599456047374e-03,1.140067022256757191e-02,1.255361976920468109e-02,5.234429868500921875e+06
5.873000000000000000e+03,2.210000000000000000e+02,4.191811761451745555e-03,9.126518237599393416e-03,5.782971443977294197e-04,2.680996191051443923e-03,9.722850580446663286e-04,3.017246997794071663e-03,6.301598558113834063e-03,6.601756641501717489e-03,2.190076274810091034e+07
7.464000000000000000e+03,2.810000000000000000e+02,3.322800901646483923e-03,6.465390771769898536e-03,5.532126236288504662e-04,2.412354060919971559e-03,8.221866556008502244e-04,2.059009133820241446e-03,4.479860580459394531e-03,4.661772947209397513e-03,8.759603005152300000e+07
8.924000000000000000e+03,3.360000000000000000e+02,2.228831305062019610e-03,4.396070553795494024e-03,4.331369736184238197e-04,1.610807235518265821e-03,5.242480816484671206e-04,1.382224781791673756e-03,3.011409317499701057e-03,3.202631767222673901e-03,3.501953099590717554e+08
1.096800000000000000e+04,4.130000000000000000e+02,1.578989731468655004e-03,2.960255868633515988e-03,3.251040863768998781e-04,1.139572940980586140e-03,3.645787740928859333e-04,9.777380706356283372e-04,2.002006257492097274e-03,2.180615911329101490e-03,1.407959411232272863e+09
1.375300000000000000e+04,5.180000000000000000e+02,1.047917380671772305e-03,1.868951829161593107e-03,2.366059575763547734e-04,7.740218403094632424e-04,2.656985115482563128e-04,6.102810411919737214e-04,1.282711104065344413e-03,1.359276705911616027e-03,5.358758530384290695e+09
1.709700000000000000e+04,6.440000000000000000e+02,7.276033096020212724e-04,1.268979328044669711e-03,9.938593963253733031e-05,5.573056455429883319e-04,1.982301981505454013e-04,4.118788864876185172e-04,8.885892362114583380e-04,9.059347130416393543e-04,2.140637755832029343e+10
2.290300000000000000e+04,8.630000000000000000e+02,4.581817818729878996e-04,8.247713294497135712e-04,7.161598234707830025e-05,3.473513552351447688e-04,1.181162347278577753e-04,2.649477069589236625e-04,5.742169281090586880e-04,5.920495463685818702e-04,8.598030341170127869e+10
3.033000000000000000e+04,1.143000000000000000e+03,2.840047592173665738e-04,5.204948069465009739e-04,5.710321851807570307e-05,2.129595130612663286e-04,7.199030880598332786e-05,1.639010888545852987e-04,3.593620015964348502e-04,3.765286122818289677e-04,3.379660005989685059e+11
4.061900000000000000e+04,1.531000000000000000e+03,1.793346164369370461e-04,3.139209865195773759e-04,4.072699651197968045e-05,1.342704307694953241e-04,4.689435312368020790e-05,1.013636398013714675e-04,2.173399509075229575e-04,2.265165148878560224e-04,1.327079637776505615e+12
5.459100000000000000e+04,2.058000000000000000e+03,1.165553772241248052e-04,2.086956456297543640e-04,1.931101311277747416e-05,8.821300265701417971e-05,3.105715290295322904e-05,6.682933598512244815e-05,1.462634026751136499e-04,1.488653403002777897e-04,5.488579746794000977e+12
7.214200000000000000e+04,2.720000000000000000e+03,7.463836027364426377e-05,1.317633647283836089e-04,1.080184212952104163e-05,5.832684088722381500e-05,2.003876343289723555e-05,4.062797832952164535e-05,9.277248989718397552e-05,9.356735222726956158e-05,2.322314511983621094e+13
9.388500000000000000e+04,3.540000000000000000e+03,4.938128655286498717e-05,9.086149518680140412e-05,8.854870207830906019e-06,3.761878423627395472e-05,1.265094726745719018e-05,2.801576950646391534e-05,6.345258671675551708e-05,6.503522542848397497e-05,9.153998881275926562e+13
        }\tableAdaptiveTwo

        \pgfplotstableread[col sep=comma]{
ndof,nelem,res,err,res1,res2,res3,res4,err1,err2,cond
9.800000000000000000e+02,2.400000000000000000e+01,1.579640099765125191e-02,7.745260268776205104e-02,1.370050608843645466e-03,8.724146622972304282e-03,2.376964189468375958e-03,1.287977300850719382e-02,5.590315601602467710e-02,5.360730183992037240e-02,3.403034996212002006e+05
3.904000000000000000e+03,9.600000000000000000e+01,9.878890492758007555e-03,4.874987673361046159e-02,4.597026614927667422e-04,5.485531215102422073e-03,1.487285267298336976e-03,8.067098651481540567e-03,3.519224418065177579e-02,3.373509198255723179e-02,3.247330292022909853e+05
1.558400000000000000e+04,3.840000000000000000e+02,6.221388043142065012e-03,3.073272744655416305e-02,1.566223248604904021e-04,3.455887714721620364e-03,9.331394800182685537e-04,5.085983626081959363e-03,2.218452971733632373e-02,2.126845498678276922e-02,3.194800256842822419e+05
6.227200000000000000e+04,1.536000000000000000e+03,3.919371882127006476e-03,1.936717337819044912e-02,5.457295902936510284e-05,2.177670252939919483e-03,5.869537918074706944e-04,3.204954798167057497e-03,1.398004532577559744e-02,1.340319877305964046e-02,9.179296387696918100e+05
2.489600000000000000e+05,6.144000000000000000e+03,2.469152167699198777e-03,1.220237508991658774e-02,1.951126756177166505e-05,1.372042380208051858e-03,3.695455996560514716e-04,2.019224478451975540e-03,8.808159391221696832e-03,8.444778621956361420e-03,3.638998359430677257e+06
        }\tableUniformThree

        \pgfplotstableread[col sep=comma]{
ndof,nelem,res,err,res1,res2,res3,res4,err1,err2,cond
9.800000000000000000e+02,2.400000000000000000e+01,1.579640099765125191e-02,7.745260268776205104e-02,1.370050608843645466e-03,8.724146622972304282e-03,2.376964189468375958e-03,1.287977300850719382e-02,5.590315601602467710e-02,5.360730183992037240e-02,3.403034996212003171e+05
1.872000000000000000e+03,4.600000000000000000e+01,1.039356507247008593e-02,5.031147165047349101e-02,5.929456298192313773e-04,6.092947895662580160e-03,1.403739324888181447e-03,8.281311022069621705e-03,3.626497685769639739e-02,3.487256275565568664e-02,3.339195624512357172e+05
2.885000000000000000e+03,7.100000000000000000e+01,6.632988268975648430e-03,3.409387589654797740e-02,2.335733114688661556e-04,3.886403303307548428e-03,9.616437310822971126e-04,5.283283787858489115e-03,2.462999996533182107e-02,2.357446659750691426e-02,3.451049727763903211e+05
3.939000000000000000e+03,9.700000000000000000e+01,4.226750430316881559e-03,2.150737668553739240e-02,1.203941624656627895e-04,2.493166838520854967e-03,6.034425663264664490e-04,3.357216202491790803e-03,1.553265980299863534e-02,1.487628083016409065e-02,1.003872316667611594e+06
4.912000000000000000e+03,1.210000000000000000e+02,2.730742758590313374e-03,1.358305083504622070e-02,9.909063793046497847e-05,1.644824769234428568e-03,3.821538681097262332e-04,2.143746010764720038e-03,9.806494877595418880e-03,9.398506700702226324e-03,3.981125962346578017e+06
5.885000000000000000e+03,1.450000000000000000e+02,1.824054917509096560e-03,8.607757977415784328e-03,9.615546400788242698e-05,1.145561788503714597e-03,2.441427577275650813e-04,1.394995688667135256e-03,6.209528915401064594e-03,5.961144860390762794e-03,1.589067077477816120e+07
7.506000000000000000e+03,1.850000000000000000e+02,1.286101238448659687e-03,5.495209525334464329e-03,9.555393546061036636e-05,8.604566953766562163e-04,1.588935429746196777e-04,9.377062217665108187e-04,3.957253206433091093e-03,3.812804058367233288e-03,6.354630405612101406e+07
9.775000000000000000e+03,2.410000000000000000e+02,8.714848727395625369e-04,3.510754538010496008e-03,7.141875551132739492e-05,5.950943062094559083e-04,1.025179910807744090e-04,6.242900551802448147e-04,2.524334443650425898e-03,2.439904310164926413e-03,2.541440411169292629e+08
1.115300000000000000e+04,2.750000000000000000e+02,5.610684253399463852e-04,2.239490469807899845e-03,3.614356074781348098e-05,3.858755426603923040e-04,6.532449535521635893e-05,4.004050409436056015e-04,1.610502489033887588e-03,1.556148867292606692e-03,1.020233728323160887e+09
1.261300000000000000e+04,3.110000000000000000e+02,4.138554201614976901e-04,1.451800461600342420e-03,3.051776134266695630e-05,2.937800656097847452e-04,4.361380643614078362e-05,2.865939356967334996e-04,1.041361405257881919e-03,1.011578471470351648e-03,4.082341180966885090e+09
1.451800000000000000e+04,3.580000000000000000e+02,2.598944594397424785e-04,9.173026361383381259e-04,2.215181140932124192e-05,1.864393711282921284e-04,2.760272554380993208e-05,1.775749918614720013e-04,6.586792533326811267e-04,6.384244414928414834e-04,1.580540114787137222e+10
1.682900000000000000e+04,4.150000000000000000e+02,1.885224102100716080e-04,5.838493699383948046e-04,1.390655347084193207e-05,1.392944025682014693e-04,1.861903918888980664e-05,1.248907831460253296e-04,4.178830836031713391e-04,4.077423392484108306e-04,6.474824001145267487e+10
2.120500000000000000e+04,5.230000000000000000e+02,1.268342996471843022e-04,3.755186949483942734e-04,1.122595519158902115e-05,9.256984126512383471e-05,1.206985784804454359e-05,8.512379524378567680e-05,2.683836638821142893e-04,2.626490038377398036e-04,2.589929274680123596e+11
2.562300000000000000e+04,6.320000000000000000e+02,8.227510996148546787e-05,2.390920068307111207e-04,4.763044222881068083e-06,6.012641199763711045e-05,7.742045455098618334e-06,5.542005450526912808e-05,1.705147784840336936e-04,1.675998151815139719e-04,1.001990135899522339e+12
3.166000000000000000e+04,7.810000000000000000e+02,5.399157872754615616e-05,1.525528944145721627e-04,2.641670925701570216e-06,4.053399144016093251e-05,4.997030440247689441e-06,3.521558423370705428e-05,1.085717843393566614e-04,1.071660078552512634e-04,4.011303716025214355e+12
3.899300000000000000e+04,9.620000000000000000e+02,3.439634968241429830e-05,9.765917907600624257e-05,1.390516986193124636e-06,2.507788333204449734e-05,3.231105258688908502e-06,2.327735080239308273e-05,6.967625526231181862e-05,6.842904873239610743e-05,1.634752065090042578e+13
4.912300000000000000e+04,1.212000000000000000e+03,2.454634069079307011e-05,6.715404851809319002e-05,1.202436218302408555e-06,1.653312396436414989e-05,2.274664911216767049e-06,1.795991926401804454e-05,4.787890315557513248e-05,4.708796942945679936e-05,6.539011166224096094e+13
5.678000000000000000e+04,1.401000000000000000e+03,1.629489535089447263e-05,3.973776052138151669e-05,1.154489986382009330e-06,1.269897126152423009e-05,1.392120848858520415e-06,1.004932315178996441e-05,2.821126958365616750e-05,2.798595861736635142e-05,2.459732738189519062e+14
7.087800000000000000e+04,1.749000000000000000e+03,9.575786073179552835e-06,2.439309361035165650e-05,6.551446549960352244e-07,7.426700503813102715e-06,8.305738689725898131e-07,5.951531825939585947e-06,1.734842094755755993e-05,1.714804089421662474e-05,9.843618264263920000e+14
        }\tableAdaptiveThree

        %
        %
        \addplot+ [line1, minorline, forget plot] table [x=ndof, y=res] {\tableUniformZero};
        \label{leg:est:unif:0}
        \addplot+ [line2, minorline, forget plot] table [x=ndof, y=res] {\tableUniformOne};
        \label{leg:est:unif:1}
        \addplot+ [line3, minorline, forget plot] table [x=ndof, y=res] {\tableUniformTwo};
        \label{leg:est:unif:2}
        \addplot+ [line4, minorline, forget plot] table [x=ndof, y=res] {\tableUniformThree};
        \label{leg:est:unif:3}

        \addplot+ [line1, majorline, forget plot] table [x=ndof, y=res] {\tableAdaptiveZero};
        \label{leg:est:adap:0}
        \addplot+ [line2, majorline, forget plot] table [x=ndof, y=res] {\tableAdaptiveOne};
        \label{leg:est:adap:1}
        \addplot+ [line3, majorline, forget plot] table [x=ndof, y=res] {\tableAdaptiveTwo};
        \label{leg:est:adap:2}
        \addplot+ [line4, majorline, forget plot] table [x=ndof, y=res] {\tableAdaptiveThree};
        \label{leg:est:adap:3}

        %
        %
        \drawslopetriangle[ST1]{2}{3e3}{2e-5}
        \drawswappedslopetriangle[ST3]{0.33}{2e5}{1.5e-1}
    \end{loglogaxis}
\end{tikzpicture}

%% file: figures/plot_opdls_lshape_estimator_degree.tex
\begin{tikzpicture}[>=stealth]
    %
    %
    \colorlet{col1}{TUblue}
    \colorlet{col2}{TUgreen}
    \colorlet{col3}{TUmagenta}
    \colorlet{col4}{TUyellow}
    \colorlet{col5}{purple}
    \colorlet{col6}{green}
    \pgfplotsset{%
        linedefault/.style = {%
            mark = *,%
            mark size = 2pt,%
            every mark/.append style = {solid},%
            gray,%
            every mark/.append style = {fill = gray!60!white}%
        },%
        line1/.style = {%
            linedefault,%
            col1,%
            every mark/.append style = {fill = col1!60!white}%
        },%
        line2/.style = {%
            linedefault,%
            mark = triangle*,%
            mark size = 2.75pt,%
            col2,%
            every mark/.append style = {fill = col2!60!white}%
        },%
        line3/.style = {%
            linedefault,%
            mark = square*,%
            mark size = 1.66pt,%
            col3,%
            every mark/.append style = {fill = col3!60!white}%
        },%
        line4/.style = {%
            linedefault,%
            mark = pentagon*,%
            mark size = 2.2pt,%
            col4,%
            every mark/.append style = {fill = col4!60!white}%
        },%
        line5/.style = {%
            linedefault,%
            mark = diamond*,%
            mark size = 2.75pt,%
            col5,%
            every mark/.append style = {fill = col5!60!white}%
        },%
        line6/.style = {%
            linedefault,%
            mark = halfsquare*,%
            mark size = 1.66pt,%
            col6,%
            every mark/.append style = {fill = col6!60!white}%
        },%
        minorline/.style = {%
            dashed,%
            every mark/.append style = {fill = black!20!white}%
        },%
        majorline/.style = {%
            solid%
        }%
    }

    %
    %
    \begin{loglogaxis}[%
            width            = 0.36\textwidth,%
            xlabel           = ndof,%
            ymax             = 1e0,%
            ymin             = 3e-6,%
            ymajorgrids      = true,%
            font             = \footnotesize,%
            grid style       = {densely dotted, semithick},%
            legend style     = {legend pos  = south west}%
        ]

        %
        %
        \pgfplotstableread[col sep=comma]{
ndof,nelem,res,err,res1,res2,res3,res4,err1,err2,cond
1.440000000000000000e+02,2.400000000000000000e+01,6.306198149120191276e-01,1.008644855299557497e+00,3.366965743021015922e-01,3.651053388672639288e-01,2.652373984422513042e-01,2.840140521366953386e-01,6.762132392680514359e-01,7.483983559314346667e-01,2.296586155100305859e+02
5.760000000000000000e+02,9.600000000000000000e+01,3.891634116065782467e-01,5.745268224882443819e-01,1.620156561466426781e-01,2.310425831275144948e-01,1.314317966611835697e-01,2.335467741263464181e-01,3.810417180298759643e-01,4.299863705738576902e-01,4.384474670547972437e+02
2.304000000000000000e+03,3.840000000000000000e+02,2.199346201832386161e-01,3.132415464677545947e-01,7.722439043103720979e-02,1.329370244949013291e-01,4.796426422386291383e-02,1.497825342074489141e-01,2.070719725645633613e-01,2.350350284781635712e-01,5.219306858866519178e+03
9.216000000000000000e+03,1.536000000000000000e+03,1.209213915613219781e-01,1.733866068509177238e-01,3.769733884062748119e-02,7.238555346200944984e-02,1.632586479936865356e-02,8.771938993183199529e-02,1.148033304170861940e-01,1.299350251488029240e-01,7.644441501019104908e+04
3.686400000000000000e+04,6.144000000000000000e+03,6.691209920672000488e-02,9.851845448837180663e-02,1.866593364274058559e-02,3.912977304488995250e-02,5.704439904109297362e-03,5.064713386096535847e-02,6.549793757790635385e-02,7.359283965045819442e-02,1.190822091908842791e+06
1.474560000000000000e+05,2.457600000000000000e+04,3.776848089809536901e-02,5.739465368140826695e-02,9.294922804628987303e-03,2.146013391463778147e-02,2.079038829906257931e-03,2.958382681999287625e-02,3.833217236807209827e-02,4.271757053898535966e-02,1.888587027434254438e+07
        }\tableUniformZero

        \pgfplotstableread[col sep=comma]{
ndof,nelem,res,err,res1,res2,res3,res4,err1,err2,cond
1.440000000000000000e+02,2.400000000000000000e+01,6.306198149120191276e-01,1.008644855299557497e+00,3.366965743021015922e-01,3.651053388672639288e-01,2.652373984422513042e-01,2.840140521366953386e-01,6.762132392680514359e-01,7.483983559314346667e-01,2.296586155100324049e+02
4.200000000000000000e+02,7.000000000000000000e+01,4.394954169616491213e-01,6.715707909097310413e-01,2.134046978046295895e-01,2.398934244640127944e-01,1.626743030107819787e-01,2.521960993435343745e-01,4.599266805425582527e-01,4.893616001774407365e-01,4.248331973432046880e+02
1.044000000000000000e+03,1.740000000000000000e+02,3.080995764558332461e-01,4.448715612364610994e-01,1.206147436548664076e-01,1.831248710847279693e-01,8.866217124494756019e-02,1.974379234232068736e-01,2.920082058547046611e-01,3.356216824200766724e-01,2.851195280339358760e+03
2.352000000000000000e+03,3.920000000000000000e+02,2.283813563786274703e-01,3.195496515583942654e-01,8.852265376466934843e-02,1.299399484381094927e-01,4.840039087579297433e-02,1.584133696212980114e-01,2.069620642343262262e-01,2.434721416896762913e-01,2.390779644215623557e+04
4.500000000000000000e+03,7.500000000000000000e+02,1.636713346895881938e-01,2.212306526684754138e-01,5.883268021431856815e-02,9.802389768926991764e-02,2.600950890058896514e-02,1.142008864520090011e-01,1.419485766226470380e-01,1.696867799061674009e-01,1.575015820655770076e+05
9.180000000000000000e+03,1.530000000000000000e+03,1.187388127792072257e-01,1.588100720996714343e-01,4.122529757414088170e-02,6.988946933720374677e-02,1.314419176578006734e-02,8.568589614824374812e-02,1.006913074464998586e-01,1.228083857275931723e-01,1.226046921408378752e+06
1.662600000000000000e+04,2.771000000000000000e+03,8.648401513452842748e-02,1.137322684093443892e-01,3.016879753258562905e-02,5.267870751491721126e-02,7.324607844276403469e-03,6.116070980386480171e-02,7.275750577325909574e-02,8.741495427665293971e-02,8.826925579184673727e+06
3.108000000000000000e+04,5.180000000000000000e+03,6.407539385855630376e-02,8.380844066854448160e-02,2.230113949688580441e-02,3.842147594042780168e-02,3.825850255794115921e-03,4.601595735882178101e-02,5.311393808999341493e-02,6.482873057422372520e-02,6.463194480318483710e+07
5.716800000000000000e+04,9.528000000000000000e+03,4.660186986402903991e-02,6.027330625225215538e-02,1.622121547324421181e-02,2.855032938168600154e-02,2.134798690145196885e-03,3.299890559009913127e-02,3.851564025972416261e-02,4.636180434324470301e-02,4.729308135390356779e+08
1.078380000000000000e+05,1.797300000000000000e+04,3.443687955697830733e-02,4.437535905440612161e-02,1.192594613064655411e-02,2.113007189858816595e-02,1.108489517543145104e-03,2.441232875128364344e-02,2.799514469590566254e-02,3.443028266893495559e-02,3.438310353276770592e+09
        }\tableAdaptiveZero

        \pgfplotstableread[col sep=comma]{
ndof,nelem,res,err,res1,res2,res3,res4,err1,err2,cond
3.360000000000000000e+02,2.400000000000000000e+01,1.410774541754422629e-01,2.532432142128129837e-01,5.793087820617605355e-02,7.321032798409224385e-02,4.723246674586529387e-02,9.463721988967961629e-02,1.753891392472060329e-01,1.826766963215638306e-01,7.230740466181046941e+03
1.344000000000000000e+03,9.600000000000000000e+01,6.044335677174061255e-02,1.291863839669573910e-01,1.469509362887662200e-02,2.959645169782602198e-02,1.298982141402372885e-02,4.891592979392801049e-02,8.974391242259822932e-02,9.292551848412342053e-02,6.953307961995453297e+03
5.376000000000000000e+03,3.840000000000000000e+02,3.319076797554559172e-02,7.916865575463447047e-02,3.994376469488723075e-03,1.409646253926915768e-02,3.964753671979430981e-03,2.951681736952689247e-02,5.494290076540755269e-02,5.699959394134658519e-02,1.673272813861177201e+04
2.150400000000000000e+04,1.536000000000000000e+03,2.031558977735474503e-02,4.987345534912807771e-02,1.143409419566413041e-03,8.165126979841825713e-03,1.352857110533361605e-03,1.851799886441691731e-02,3.458845579109927348e-02,3.593049226560430381e-02,7.729624800803052494e+04
8.601600000000000000e+04,6.144000000000000000e+03,1.273734677141928444e-02,3.146054049942184627e-02,3.442315075950557412e-04,5.061888876694767628e-03,4.836574720980640050e-04,1.167325421184009504e-02,2.181382519503181558e-02,2.266986190686552660e-02,1.194891384845554363e+06
3.440640000000000000e+05,2.457600000000000000e+04,8.018827900752740698e-03,1.983188318091227628e-02,1.091357260475482187e-04,3.180962498518550077e-03,1.777413066239565630e-04,7.357960037002686700e-03,1.375004510734577703e-02,1.429125082164982644e-02,1.890740388982083276e+07
        }\tableUniformOne

        \pgfplotstableread[col sep=comma]{
ndof,nelem,res,err,res1,res2,res3,res4,err1,err2,cond
3.360000000000000000e+02,2.400000000000000000e+01,1.410774541754422629e-01,2.532432142128129837e-01,5.793087820617605355e-02,7.321032798409224385e-02,4.723246674586529387e-02,9.463721988967961629e-02,1.753891392472060329e-01,1.826766963215638306e-01,7.230740466181054217e+03
8.540000000000000000e+02,6.100000000000000000e+01,8.748715718044233758e-02,1.684296540347229310e-01,3.311526801709450962e-02,4.457825042334847415e-02,2.055699479589141038e-02,6.440164011751835227e-02,1.152801018908801434e-01,1.227967689570240856e-01,6.443302249040961215e+03
1.358000000000000000e+03,9.700000000000000000e+01,5.759116983626711567e-02,1.050673055722659566e-01,1.654478567884715492e-02,3.150369769104819501e-02,1.230241574477199118e-02,4.357958821036068964e-02,7.151586571125451153e-02,7.697155092490802486e-02,1.550524917010031277e+04
2.170000000000000000e+03,1.550000000000000000e+02,4.141981427373113422e-02,7.172515584232455610e-02,1.294810477267752807e-02,2.391145890572530727e-02,7.650618415257546753e-03,3.029286662120440204e-02,4.865311310762630848e-02,5.270078334846894919e-02,6.252838230085997930e+04
3.262000000000000000e+03,2.330000000000000000e+02,3.033955551758293553e-02,5.015740412317001157e-02,9.967159336669190531e-03,1.830299842149288661e-02,5.276672396546587714e-03,2.140797377063736670e-02,3.396670332767071782e-02,3.690566695542835118e-02,2.499584164417928259e+05
4.634000000000000000e+03,3.310000000000000000e+02,2.063765791858325135e-02,3.282160479639936734e-02,7.336568817536011020e-03,1.261348304823302173e-02,3.043348420066842829e-03,1.427325324290278780e-02,2.259776723547589003e-02,2.380333290491641607e-02,1.018860286316490616e+06
6.524000000000000000e+03,4.660000000000000000e+02,1.400410033005189750e-02,2.208387186187773729e-02,4.174915921989151221e-03,8.904496261487114539e-03,1.627868997453578047e-03,9.835898126494229191e-03,1.505494962739059583e-02,1.615691456091065958e-02,4.076691744258814491e+06
9.702000000000000000e+03,6.930000000000000000e+02,9.742231738923525300e-03,1.471505667422966061e-02,2.913031071305545623e-03,6.418956689802151130e-03,9.715999900890765536e-04,6.654195496542055638e-03,9.977065324930664211e-03,1.081624058662977381e-02,1.630663577492151223e+07
1.527400000000000000e+04,1.091000000000000000e+03,6.534829400284570926e-03,9.775121048782925079e-03,1.886694431174944231e-03,4.297115553398902013e-03,5.309073156346052081e-04,4.516338645090722109e-03,6.609974557783129812e-03,7.201474006328045102e-03,6.521203520578209311e+07
2.233000000000000000e+04,1.595000000000000000e+03,4.280859669157644799e-03,6.393246211616394545e-03,1.192913015934106988e-03,2.833009672566197228e-03,2.609563914990130720e-04,2.967941374137594585e-03,4.345083261514961180e-03,4.689759969641294282e-03,2.928693184765223861e+08
3.235400000000000000e+04,2.311000000000000000e+03,2.939255400262257793e-03,4.278937844654475721e-03,8.798031802792400692e-04,1.972485081737179041e-03,1.611462988583335393e-04,1.987084080916036097e-03,2.897366575165785836e-03,3.148741972205469815e-03,1.572087604532146931e+09
5.154800000000000000e+04,3.682000000000000000e+03,1.911687575861123782e-03,2.750815083306280481e-03,5.489416556833490342e-04,1.295184838722848340e-03,8.661533450177437985e-05,1.291590671896096498e-03,1.856267786951105829e-03,2.030087073421481432e-03,1.000104831224785614e+10
7.842800000000000000e+04,5.602000000000000000e+03,1.228069321714765324e-03,1.768564637063715613e-03,3.573866704758233630e-04,8.268487638269265494e-04,4.109147874345383264e-05,8.337035677240038590e-04,1.206569594144262653e-03,1.293062523607757617e-03,6.375636762051554108e+10
1.176980000000000000e+05,8.407000000000000000e+03,8.248451186374084311e-04,1.167742506468436893e-03,2.522381977433809210e-04,5.591888267446263857e-04,2.435729739990704620e-05,5.508719809960150844e-04,7.943734678354903104e-04,8.559166752740630891e-04,3.713117678216314087e+11
        }\tableAdaptiveOne

        \pgfplotstableread[col sep=comma]{
ndof,nelem,res,err,res1,res2,res3,res4,err1,err2,cond
6.000000000000000000e+02,2.400000000000000000e+01,3.819561679950280592e-02,1.091662736805051048e-01,5.670447324352385562e-03,2.007928134485041649e-02,1.297152897598293463e-02,2.924573586928157914e-02,7.767238907163920292e-02,7.670903001075458916e-02,7.007366623479264672e+04
2.400000000000000000e+03,9.600000000000000000e+01,2.277932785179561520e-02,7.183010594498503987e-02,1.768816172240626370e-03,1.105710855070412418e-02,4.088785592132335514e-03,1.941111149851879394e-02,5.104347467606239891e-02,5.053837960463263695e-02,6.798193358015047852e+04
9.600000000000000000e+03,3.840000000000000000e+02,1.431082609520344801e-02,4.584306043342896431e-02,5.823146197496271688e-04,7.021356604078357946e-03,1.366474247951473273e-03,1.238119351188859789e-02,3.256631628527141553e-02,3.226486066157242816e-02,1.010560841474677291e+05
3.840000000000000000e+04,1.536000000000000000e+03,9.014753329646749483e-03,2.898781487936142806e-02,1.952364080309571383e-04,4.446171942711974086e-03,4.613407076057559155e-04,7.826006653888250772e-03,2.059176009742792746e-02,2.040276519421086940e-02,4.009493434865098679e+05
1.536000000000000000e+05,6.144000000000000000e+03,5.679343816742773998e-03,1.828203560329835775e-02,6.681547802816383604e-05,2.805554438313187862e-03,1.588234668496425425e-04,4.934989491448724466e-03,1.298679488682430092e-02,1.286763320769688211e-02,1.601304117351882625e+06
        }\tableUniformTwo

        \pgfplotstableread[col sep=comma]{
ndof,nelem,res,err,res1,res2,res3,res4,err1,err2,cond
6.000000000000000000e+02,2.400000000000000000e+01,3.819561679950280592e-02,1.091662736805051048e-01,5.670447324352385562e-03,2.007928134485041649e-02,1.297152897598293463e-02,2.924573586928157914e-02,7.767238907163920292e-02,7.670903001075458916e-02,7.007366623479271948e+04
1.150000000000000000e+03,4.600000000000000000e+01,2.756964904917026324e-02,7.703019285530524063e-02,3.615022502163543674e-03,1.521574048550653965e-02,6.269939773021809305e-03,2.182169236616771438e-02,5.433206919464296514e-02,5.460473302154351255e-02,7.276368591044633649e+04
1.700000000000000000e+03,6.800000000000000000e+01,2.198641813043256815e-02,5.357698815880844456e-02,2.465229996057532232e-03,1.128756341673838731e-02,3.997719197554329938e-03,1.827387141804841242e-02,3.660315366322466901e-02,3.912419714288668110e-02,1.063015404990217794e+05
2.350000000000000000e+03,9.400000000000000000e+01,1.585428130065361443e-02,3.607674668813766472e-02,2.167285357792936035e-03,9.343712303961203974e-03,3.318469608709834814e-03,1.217965145539988414e-02,2.441775181880920936e-02,2.655757985425426237e-02,4.247983726905161748e+05
3.325000000000000000e+03,1.330000000000000000e+02,1.022165500908132733e-02,2.346448107127773949e-02,1.538460393819226319e-03,5.897545309268036225e-03,1.801554614198701797e-03,8.005543769071345184e-03,1.595606973065698611e-02,1.720423525457506417e-02,1.689795423082960071e+06
4.375000000000000000e+03,1.750000000000000000e+02,8.128557037679385450e-03,1.706780558628640324e-02,1.509294112177854029e-03,5.187770089703944328e-03,1.748342259932253624e-03,5.815996186059663696e-03,1.148245800000523634e-02,1.262787178464306985e-02,6.751156253198359162e+06
5.525000000000000000e+03,2.210000000000000000e+02,4.207805289181119381e-03,9.264352108815279707e-03,5.784715725420058827e-04,2.734282212934087235e-03,6.544591268083478153e-04,3.076748287200448129e-03,6.403523852472148802e-03,6.695005770492692618e-03,2.825060354364511371e+07
7.025000000000000000e+03,2.810000000000000000e+02,3.331159766031210456e-03,6.545093333746404736e-03,5.533849120563855747e-04,2.447544432006348511e-03,6.419410708908593029e-04,2.094714405543551556e-03,4.540244928450749146e-03,4.714278601984501008e-03,1.130028801320911646e+08
8.350000000000000000e+03,3.340000000000000000e+02,2.198225575137922796e-03,4.409293946533826905e-03,4.301334292020647003e-04,1.610680439374087645e-03,3.441330669457995373e-04,1.390849332863017956e-03,3.034411508149149374e-03,3.199096732853174644e-03,4.517641988687487841e+08
1.032500000000000000e+04,4.130000000000000000e+02,1.538033157362616492e-03,2.946448251106259866e-03,3.212353954825196271e-04,1.131930137895854633e-03,2.075089912616295753e-04,9.685184537760034530e-04,2.010119562325146398e-03,2.154292608167448185e-03,1.816508657687121153e+09
1.315000000000000000e+04,5.260000000000000000e+02,9.890414840953529619e-04,1.842006646037166887e-03,1.984798353354147937e-04,7.491959651551121074e-04,1.105492321722191008e-04,6.043948091175163281e-04,1.278220901324935745e-03,1.326325680766665205e-03,6.908659052490003586e+09
1.712500000000000000e+04,6.850000000000000000e+02,6.485146494104617040e-04,1.174625050322526849e-03,1.219014788500662313e-04,5.082158015828334334e-04,5.862216188487281185e-05,3.794620153738755735e-04,8.263951438649462739e-04,8.347544998630616929e-04,2.726510556326704788e+10
2.295000000000000000e+04,9.180000000000000000e+02,4.077883613519849753e-04,7.576825104850408226e-04,6.942790158930151870e-05,3.183220337087589889e-04,3.430147752301317224e-05,2.428283460189559467e-04,5.287449788117095728e-04,5.426891689321968523e-04,1.122120581666528778e+11
2.960000000000000000e+04,1.184000000000000000e+03,2.741965634972319912e-04,5.184784302531738600e-04,5.084198902757097190e-05,2.095272353955302046e-04,1.707900678980735733e-05,1.685392913833202164e-04,3.588055606703724519e-04,3.742705602499554296e-04,4.489444048630886230e+11
4.047500000000000000e+04,1.619000000000000000e+03,1.692209133215686657e-04,3.129745288966800216e-04,3.588779727677518207e-05,1.297913811140073280e-04,9.379447807776055638e-06,1.020490414425664178e-04,2.183110099257196200e-04,2.242618083475457600e-04,1.750882096834728027e+12
5.402500000000000000e+04,2.161000000000000000e+03,1.090935958472342309e-04,1.981103946960378827e-04,1.738611838880867725e-05,8.629806808459527308e-05,5.246929398166668794e-06,6.422031391240002609e-05,1.395930987571483229e-04,1.405755855972050470e-04,6.855385731618833984e+12
7.090000000000000000e+04,2.836000000000000000e+03,7.219652495714911895e-05,1.283887793841192320e-04,1.078768863498282686e-05,5.900769926371286420e-05,3.575501259532920860e-06,4.001588694692670645e-05,9.080658885280398546e-05,9.076255887017426935e-05,2.752245437487217188e+13
9.117500000000000000e+04,3.647000000000000000e+03,4.798932552939210887e-05,8.880366444441850887e-05,8.549489185352279410e-06,3.790147359310555705e-05,1.704783625335345602e-06,2.811500687706529060e-05,6.212129314302242643e-05,6.345892968681016264e-05,1.100902615881208438e+14
1.274500000000000000e+05,5.098000000000000000e+03,2.991038790260966415e-05,5.386633579543171908e-05,6.027534678034814011e-06,2.378441082139401839e-05,1.024344002718999720e-06,1.707491292607125270e-05,3.729533393569878284e-05,3.886695509878440022e-05,4.807198707087767500e+14
        }\tableAdaptiveTwo

        \pgfplotstableread[col sep=comma]{
ndof,nelem,res,err,res1,res2,res3,res4,err1,err2,cond
9.360000000000000000e+02,2.400000000000000000e+01,1.587287669783075603e-02,7.595987805906435908e-02,1.358141089144095960e-03,8.544501383718315288e-03,3.874310374247793946e-03,1.273125613104020108e-02,5.481240670708409368e-02,5.258805135888882515e-02,3.403802350252678152e+05
3.744000000000000000e+03,9.600000000000000000e+01,9.880771211355181619e-03,4.888935784724834227e-02,4.608958037078263225e-04,5.514618161488502902e-03,1.299875972795064586e-03,8.081863880170288977e-03,3.529654494296664835e-02,3.382784689875251161e-02,3.246558266713191988e+05
1.497600000000000000e+04,3.840000000000000000e+02,6.213848900663823151e-03,3.101060628718815479e-02,1.575809511539809273e-04,3.499670015181694741e-03,4.363194349487324749e-04,5.113611398966347484e-03,2.238955745375374382e-02,2.145612778028817627e-02,3.192977679735824349e+05
5.990400000000000000e+04,1.536000000000000000e+03,3.912984282204073905e-03,1.957581702869206144e-02,5.497430347534909090e-05,2.209542810322811992e-03,1.491579503757459234e-04,3.225538140206056417e-03,1.413379017229841363e-02,1.354431938881652843e-02,1.209546226606750395e+06
2.396160000000000000e+05,6.144000000000000000e+03,2.464803052673721507e-03,1.233995624749361329e-02,1.966621898909293986e-05,1.392847920339726530e-03,5.224274491968399392e-05,2.032759871383960082e-03,8.909497725342371832e-03,8.537878569772276716e-03,4.834927912092171609e+06
        }\tableUniformThree

        \pgfplotstableread[col sep=comma]{
ndof,nelem,res,err,res1,res2,res3,res4,err1,err2,cond
9.360000000000000000e+02,2.400000000000000000e+01,1.587287669783075603e-02,7.595987805906435908e-02,1.358141089144095960e-03,8.544501383718315288e-03,3.874310374247793946e-03,1.273125613104020108e-02,5.481240670708409368e-02,5.258805135888882515e-02,3.403802350252678152e+05
1.794000000000000000e+03,4.600000000000000000e+01,1.044915960577807480e-02,5.018233837711233436e-02,5.928998277704849995e-04,6.095577088923767248e-03,1.808404555794396401e-03,8.270853594208021214e-03,3.617391115849837024e-02,3.478096083336495309e-02,3.339776074129303452e+05
2.769000000000000000e+03,7.100000000000000000e+01,6.635266756571738325e-03,3.431627056434442824e-02,2.343898911081658217e-04,3.923572513949794352e-03,6.649669619238980818e-04,5.304264697122472103e-03,2.479432431245664342e-02,2.372441584810916493e-02,3.452127746248840704e+05
3.783000000000000000e+03,9.700000000000000000e+01,4.221917769610821630e-03,2.170500740314325674e-02,1.204489754034804478e-04,2.521598716824297298e-03,2.384367398600577596e-04,3.375613948613621646e-03,1.567808216511660385e-02,1.501016608816691834e-02,1.299883519846176263e+06
4.719000000000000000e+03,1.210000000000000000e+02,2.725873354535585970e-03,1.371783137792322721e-02,9.884956520065375069e-05,1.661408201999718453e-03,1.083837079515970263e-04,2.156058920021773368e-03,9.905557097018662799e-03,9.489933419726310007e-03,5.196154601200701669e+06
5.655000000000000000e+03,1.450000000000000000e+02,1.819834332671921252e-03,8.693959554343065610e-03,9.586001748769681756e-05,1.152889520420732807e-03,7.842489721117074867e-05,1.402605840790699995e-03,6.272830438739948958e-03,6.019678647516872247e-03,2.078143583830968291e+07
7.215000000000000000e+03,1.850000000000000000e+02,1.281645461708655342e-03,5.548619558195228181e-03,9.524621882136063782e-05,8.607615494533052631e-04,7.319812976987359805e-05,9.419526719000080942e-04,3.996384121383923273e-03,3.849167826418727414e-03,8.314545740445955098e+07
8.307000000000000000e+03,2.130000000000000000e+02,9.966580641392046736e-04,3.610273636436359557e-03,9.519906774514425973e-05,7.174788920416048802e-04,7.257209391171580187e-05,6.813382162877577516e-04,2.589364890799469160e-03,2.515803090912021017e-03,3.325787006826478839e+08
1.002300000000000000e+04,2.570000000000000000e+02,5.877306847250159030e-04,2.271458219838518987e-03,3.663291604626303205e-05,4.182903696924609787e-04,2.336163683546792045e-05,4.105761656211856102e-04,1.633552922775719906e-03,1.578298860470625856e-03,1.335089161239186764e+09
1.181700000000000000e+04,3.030000000000000000e+02,3.831165054876371225e-04,1.443372331996678470e-03,2.708890123921533854e-05,2.705521397650076398e-04,1.387692044826759351e-05,2.695429816451521028e-04,1.037688083674370753e-03,1.003258256768285814e-03,5.342453017069465637e+09
1.353300000000000000e+04,3.470000000000000000e+02,2.741880337021723171e-04,9.129440701841704930e-04,1.290109051990256363e-05,2.031636466690750790e-04,6.034013376168046660e-06,1.835776757766593840e-04,6.552165774825424943e-04,6.357343091862382393e-04,2.117872461866301727e+10
1.548300000000000000e+04,3.970000000000000000e+02,1.775337249757544603e-04,5.794029171933716199e-04,1.201643916437640030e-05,1.309859489858040493e-04,5.184404817014923740e-06,1.191202409691875542e-04,4.154488993188299199e-04,4.038687317767514348e-04,8.471489527912284851e+10
2.063100000000000000e+04,5.290000000000000000e+02,1.206242561955739540e-04,3.740274611655121113e-04,1.036028941375277699e-05,8.801327873889524557e-05,4.081937120081584325e-06,8.173050924211808384e-05,2.674512911938693488e-04,2.614695939199258401e-04,3.450490087660043335e+11
2.655900000000000000e+04,6.810000000000000000e+02,7.536952986073883816e-05,2.334359605725026333e-04,4.883640297637097618e-06,5.643651737489691186e-05,2.480672192711689101e-06,4.965361881008645192e-05,1.668595716373047076e-04,1.632489848097747509e-04,1.265135852051228271e+12
3.174600000000000000e+04,8.140000000000000000e+02,5.132304931960395271e-05,1.509797757820412331e-04,2.491197114071798622e-06,3.942300938049305776e-05,1.301756094422553341e-06,3.274112232128794270e-05,1.078271337470901000e-04,1.056797138673387930e-04,5.218576744270708984e+12
3.868800000000000000e+04,9.920000000000000000e+02,3.149745786349974575e-05,9.430840625361631090e-05,1.310130367562511299e-06,2.347964704294816849e-05,4.579761585551246562e-07,2.094922057881156937e-05,6.744224356288509327e-05,6.592131122256004679e-05,2.087430866446758594e+13
4.578600000000000000e+04,1.174000000000000000e+03,2.079312900766588045e-05,5.949765744365865313e-05,1.170176942825198382e-06,1.586918368930536052e-05,3.534898175050502724e-07,1.338017017215858154e-05,4.253500105738561608e-05,4.160222261287443600e-05,8.048023902612678125e+13
5.951400000000000000e+04,1.526000000000000000e+03,1.292739499076695692e-05,3.738336084600879885e-05,7.916873589213317810e-07,9.901393078428031436e-06,2.203662857167556843e-07,8.270708928478118392e-06,2.674158372548604874e-05,2.612285145223855400e-05,3.228427580781133125e+14
7.215000000000000000e+04,1.850000000000000000e+03,8.609161170931742225e-06,2.379273493406855371e-05,5.774810469590360504e-07,6.749398004320083672e-06,1.550863546335971917e-07,5.310908256538430609e-06,1.700316015683268081e-05,1.664291982567792680e-05,1.393205347293677750e+15
9.126000000000000000e+04,2.340000000000000000e+03,5.802633111431548545e-06,1.582601220917260784e-05,3.486582543694925725e-07,4.218547487638094829e-06,7.386616570535245681e-08,3.968298039669914082e-06,1.126140265581457409e-05,1.111950865229631463e-05,5.579504527340372000e+15
        }\tableAdaptiveThree

        %
        %
        \addplot+ [line1, minorline, forget plot] table [x=ndof, y=res] {\tableUniformZero};
        \addplot+ [line2, minorline, forget plot] table [x=ndof, y=res] {\tableUniformOne};
        \addplot+ [line3, minorline, forget plot] table [x=ndof, y=res] {\tableUniformTwo};
        \addplot+ [line4, minorline, forget plot] table [x=ndof, y=res] {\tableUniformThree};

        \addplot+ [line1, majorline, forget plot] table [x=ndof, y=res] {\tableAdaptiveZero};
        \addplot+ [line2, majorline, forget plot] table [x=ndof, y=res] {\tableAdaptiveOne};
        \addplot+ [line3, majorline, forget plot] table [x=ndof, y=res] {\tableAdaptiveTwo};
        \addplot+ [line4, majorline, forget plot] table [x=ndof, y=res] {\tableAdaptiveThree};

        %
        %
        \drawslopetriangle[ST1]{2}{4e3}{1e-5}
        \drawswappedslopetriangle[ST3]{0.33}{1.5e5}{1.2e-1}

    \end{loglogaxis}
\end{tikzpicture}

%% file: figures/plot_dls_lshape_efficiency_degree.tex
\begin{tikzpicture}[>=stealth]
    %
    %
    \colorlet{col1}{TUblue}
    \colorlet{col2}{TUgreen}
    \colorlet{col3}{TUmagenta}
    \colorlet{col4}{TUyellow}
    \colorlet{col5}{purple}
    \colorlet{col6}{green}
    \pgfplotsset{%
        linedefault/.style = {%
            mark = *,%
            mark size = 2pt,%
            every mark/.append style = {solid},%
            gray,%
            every mark/.append style = {fill = gray!60!white}%
        },%
        line1/.style = {%
            linedefault,%
            col1,%
            every mark/.append style = {fill = col1!60!white}%
        },%
        line2/.style = {%
            linedefault,%
            mark = triangle*,%
            mark size = 2.75pt,%
            col2,%
            every mark/.append style = {fill = col2!60!white}%
        },%
        line3/.style = {%
            linedefault,%
            mark = square*,%
            mark size = 1.66pt,%
            col3,%
            every mark/.append style = {fill = col3!60!white}%
        },%
        line4/.style = {%
            linedefault,%
            mark = pentagon*,%
            mark size = 2.2pt,%
            col4,%
            every mark/.append style = {fill = col4!60!white}%
        },%
        line5/.style = {%
            linedefault,%
            mark = diamond*,%
            mark size = 2.75pt,%
            col5,%
            every mark/.append style = {fill = col5!60!white}%
        },%
        line6/.style = {%
            linedefault,%
            mark = halfsquare*,%
            mark size = 1.66pt,%
            col6,%
            every mark/.append style = {fill = col6!60!white}%
        },%
        minorline/.style = {%
            dashed,%
            every mark/.append style = {fill = black!20!white}%
        },%
        majorline/.style = {%
            solid%
        }%
    }

    %
    %
    \pgfplotstableset{%
        create on use/effRes/.style={create col/expr={\thisrow{res}/\thisrow{err}}}
    }

    %
    %
    \begin{semilogxaxis}[%
            width            = 0.36\textwidth,%
            xlabel           = ndof,%
            ylabel           = {efficiency index},%
            ymax             = 1,%
            ymin             = 0,%
            ymajorgrids      = true,%
            font             = \footnotesize,%
            grid style       = {densely dotted, semithick},%
            legend style     = {legend pos  = south west}%
        ]

        %
        %
        \pgfplotstableread[col sep=comma]{
ndof,nelem,res,err,res1,res2,res3,res4,err1,err2,cond
1.880000000000000000e+02,2.400000000000000000e+01,6.691468253750407769e-01,1.025606713479946253e+00,3.566616625857725320e-01,3.931095612333450617e-01,2.393892856402309924e-16,4.074491443333447815e-01,6.925445981398684614e-01,7.564728087131914469e-01,1.861844080681640321e+02
7.360000000000000000e+02,9.600000000000000000e+01,4.018328538784070036e-01,5.720407191928681412e-01,1.652282559352703517e-01,2.464800014124796934e-01,3.594997419828183740e-16,2.709554849054693904e-01,3.805415788800027288e-01,4.271050118626665282e-01,2.158354393424532702e+02
2.912000000000000000e+03,3.840000000000000000e+02,2.226953929342510019e-01,3.122001557070245847e-01,7.752238586897748640e-02,1.366714899126683647e-01,4.720870189556533655e-16,1.578113480342057118e-01,2.065555837003762440e-01,2.341019608326406487e-01,2.354149255634362362e+03
1.158400000000000000e+04,1.536000000000000000e+03,1.214925828018451198e-01,1.733462263358546418e-01,3.772328042619833843e-02,7.315877376198784421e-02,7.058621229607447708e-16,8.935992172888182483e-02,1.148198861841061846e-01,1.298665003823167019e-01,3.452532552691214369e+04
4.620800000000000000e+04,6.144000000000000000e+03,6.703768583937476377e-02,9.857353874140176164e-02,1.866840647801305472e-02,3.929790042375498615e-02,1.332844287066715975e-15,5.100212686305377607e-02,6.554826106162134347e-02,7.362179033275560724e-02,5.364942042558327084e+05
1.845760000000000000e+05,2.457600000000000000e+04,3.779800493524370503e-02,5.742160326142143473e-02,9.295180846666025593e-03,2.149966973851240304e-02,5.681306352149327939e-15,2.966568710222893326e-02,3.835423277390073021e-02,4.273398354281396361e-02,8.490368912294194102e+06
        }\tableUniformZero

        \pgfplotstableread[col sep=comma]{
ndof,nelem,res,err,res1,res2,res3,res4,err1,err2,cond
1.880000000000000000e+02,2.400000000000000000e+01,6.691468253750407769e-01,1.025606713479946253e+00,3.566616625857725320e-01,3.931095612333450617e-01,2.393892856402309924e-16,4.074491443333447815e-01,6.925445981398684614e-01,7.564728087131914469e-01,1.861844080681634352e+02
5.280000000000000000e+02,6.900000000000000000e+01,4.645485304210772437e-01,6.830599490857015565e-01,2.257879949253541185e-01,2.594660624333418597e-01,4.221377968243194286e-16,3.122538725289212835e-01,4.667482632440184132e-01,4.987153023556160636e-01,2.128427996102370798e+02
1.134000000000000000e+03,1.490000000000000000e+02,3.514277916105422728e-01,4.898292054561634856e-01,1.302675895177524967e-01,2.157932520629477768e-01,3.997275055233180919e-16,2.448777658375585753e-01,3.163755221592782219e-01,3.739507714877192690e-01,1.106585141813846803e+03
2.302000000000000000e+03,3.040000000000000000e+02,2.501417772339131207e-01,3.416849949294442634e-01,9.756774828116387965e-02,1.481552895554908156e-01,5.184071061221481349e-16,1.763560415461084552e-01,2.258489202343929025e-01,2.563998849237131838e-01,4.320766697536862011e+03
4.850000000000000000e+03,6.420000000000000000e+02,1.827162881006695605e-01,2.434834612677447552e-01,6.173245639320502215e-02,1.138115980365454472e-01,5.550654254689614395e-16,1.289234886175351369e-01,1.572760924728729692e-01,1.858720706491045760e-01,3.002650222874309475e+04
8.519000000000000000e+03,1.130000000000000000e+03,1.351938722597004761e-01,1.787860591587285553e-01,4.926933115701741944e-02,8.141275564107232354e-02,6.863274354147775619e-16,9.603061659095710190e-02,1.167796488402617627e-01,1.353771345769053391e-01,2.195661886394199682e+05
1.637800000000000000e+04,2.175000000000000000e+03,1.001953397017019687e-01,1.293968314633148209e-01,3.316026559924474537e-02,6.387072948000693806e-02,1.377493817389158267e-15,6.971393547217655251e-02,8.359794490892365204e-02,9.876701676040477995e-02,1.526050124626174103e+06
2.965100000000000000e+04,3.942000000000000000e+03,7.366358106374841508e-02,9.533914920550538852e-02,2.635424251199278492e-02,4.496798373161311341e-02,1.768221127331582023e-14,5.205437076622775744e-02,6.185394011254979896e-02,7.255097148751847580e-02,1.134321577537557483e+07
5.533200000000000000e+04,7.361000000000000000e+03,5.476430034674915903e-02,6.993875022571281741e-02,1.815836851153471154e-02,3.483691541012369497e-02,7.663785009417625751e-15,3.815483678612911073e-02,4.490324281781864907e-02,5.362021603442729961e-02,7.710566014080482721e+07
1.020060000000000000e+05,1.357800000000000000e+04,4.041410942974788673e-02,5.163485719064491580e-02,1.401685756070034962e-02,2.523231564524061263e-02,2.868700345761843395e-15,2.828706758046184433e-02,3.331128810879161600e-02,3.945271298188963299e-02,5.550838271046789885e+08
        }\tableAdaptiveZero

        \pgfplotstableread[col sep=comma]{
ndof,nelem,res,err,res1,res2,res3,res4,err1,err2,cond
3.800000000000000000e+02,2.400000000000000000e+01,1.429544905823191270e-01,2.649557076419757906e-01,5.816283848669231982e-02,7.748609674522376201e-02,3.101857475521297217e-02,1.004331789209442632e-01,1.843098314644586733e-01,1.903455096334110419e-01,7.220823296696265061e+03
1.504000000000000000e+03,9.600000000000000000e+01,6.043023319836864776e-02,1.290671872306016765e-01,1.468726269570548582e-02,2.915852172614483498e-02,1.384234833024492080e-02,4.893125179644047112e-02,8.961392183059026240e-02,9.288532625641167384e-02,6.952862152031728328e+03
5.984000000000000000e+03,3.840000000000000000e+02,3.299879050228264610e-02,7.757485919926983275e-02,3.983931087607786060e-03,1.348187200302555962e-02,7.148377534395108056e-03,2.898600168395982710e-02,5.374083191374660895e-02,5.594445249535458242e-02,1.313286401135960114e+04
2.387200000000000000e+04,1.536000000000000000e+03,2.016212199610784711e-02,4.857279130894310432e-02,1.137728940393586529e-03,7.747887232526820024e-03,4.262654840198620149e-03,1.808360454637036180e-02,3.361484377301553195e-02,3.506220662847531200e-02,5.043913348106794001e+04
9.536000000000000000e+04,6.144000000000000000e+03,1.263390460609820944e-02,3.058304462746423355e-02,3.414287908117189093e-04,4.792162396347206775e-03,2.655905513844827003e-03,1.137894184207646243e-02,2.115795523764861172e-02,2.208310550731320404e-02,5.673349497205432272e+05
3.811840000000000000e+05,2.457600000000000000e+04,7.952201830990024076e-03,1.926739112649915667e-02,1.078196987798035983e-04,3.009227955488317108e-03,1.669733750245612882e-03,7.168153541083874779e-03,1.332828378032429109e-02,1.391363404336342334e-02,8.972909084775183350e+06
        }\tableUniformOne

        \pgfplotstableread[col sep=comma]{
ndof,nelem,res,err,res1,res2,res3,res4,err1,err2,cond
3.800000000000000000e+02,2.400000000000000000e+01,1.429544905823191270e-01,2.649557076419757906e-01,5.816283848669231982e-02,7.748609674522376201e-02,3.101857475521297217e-02,1.004331789209442632e-01,1.843098314644586733e-01,1.903455096334110419e-01,7.220823296696225043e+03
9.410000000000000000e+02,6.000000000000000000e+01,9.048552861572219108e-02,1.723899622309089663e-01,3.181784498837661457e-02,4.657744344487177446e-02,1.919519511623793764e-02,6.809803136472918073e-02,1.177978602354333021e-01,1.258648608704095451e-01,6.575537898597228377e+03
1.657000000000000000e+03,1.060000000000000000e+02,5.461672803270744420e-02,9.889981584593049435e-02,1.658444428552899308e-02,3.009022628867067187e-02,1.358916670968292502e-02,4.022258095612405576e-02,6.755965415739066671e-02,7.222788038211234996e-02,1.251520602726039579e+04
2.682000000000000000e+03,1.720000000000000000e+02,3.906428975490600780e-02,6.727301167620299005e-02,1.348658806205335625e-02,2.185195330908825881e-02,9.470620031590409882e-03,2.787346712008392097e-02,4.586435459046426383e-02,4.921502898492190886e-02,4.766783496311864292e+04
4.142000000000000000e+03,2.660000000000000000e+02,2.736237289062587452e-02,4.587644150401347376e-02,8.801945939295292842e-03,1.605758390699676616e-02,6.791470739562765167e-03,1.916390153780805905e-02,3.102208294304604000e-02,3.379760723699156560e-02,1.891961357416887186e+05
5.528000000000000000e+03,3.550000000000000000e+02,1.789686763959938420e-02,2.870725154843154631e-02,5.745194310199929454e-03,1.075726856481438488e-02,4.467533532744634092e-03,1.231312027089981158e-02,1.955759948688892344e-02,2.101443774587767252e-02,7.753271398940539220e+05
8.291000000000000000e+03,5.330000000000000000e+02,1.249295063885906901e-02,1.957246137199181862e-02,4.041650490638224423e-03,7.515471054698831532e-03,3.204933307625098659e-03,8.543124387257416269e-03,1.328974480978349117e-02,1.436885266988788788e-02,3.163988872234159615e+06
1.327300000000000000e+04,8.540000000000000000e+02,8.248363375542640263e-03,1.274733939861605672e-02,2.359825262671173731e-03,5.086567941530156095e-03,2.232341650879941623e-03,5.622294943757080200e-03,8.595629524193595244e-03,9.413278643826525696e-03,1.264977298992297426e+07
1.897100000000000000e+04,1.221000000000000000e+03,5.447211056598747804e-03,8.169662857762239119e-03,1.640798505876144529e-03,3.400699778372295954e-03,1.461327451629655010e-03,3.644125636519470173e-03,5.506712451709151399e-03,6.034857842873450964e-03,5.094934276248360425e+07
2.705600000000000000e+04,1.742000000000000000e+03,3.795429211931716871e-03,5.552369508253003627e-03,1.104586523556202798e-03,2.397035874636590201e-03,1.060227549986266339e-03,2.513027670503841828e-03,3.756613471384596080e-03,4.088601519198158292e-03,2.037979702250770628e+08
4.257200000000000000e+04,2.742000000000000000e+03,2.521971172035490042e-03,3.633188819655926768e-03,7.325091624243566281e-04,1.608443598106142937e-03,7.118381972641818880e-04,1.652260418992561146e-03,2.439987182271654599e-03,2.691936765531991772e-03,1.021559285691922307e+09
6.434700000000000000e+04,4.146000000000000000e+03,1.654286607367381470e-03,2.369957368777833584e-03,4.713215955699466271e-04,1.054314049770528616e-03,4.576109251289303898e-04,1.092489935204713674e-03,1.582971914287373098e-03,1.763773752044664120e-03,5.675444554343406677e+09
9.747400000000000000e+04,6.282000000000000000e+03,1.065468831961948956e-03,1.513301870864532713e-03,3.138129662588048040e-04,6.784876866757804433e-04,3.005773986150441900e-04,6.971749712671648433e-04,1.023324527077856773e-03,1.114849525560725153e-03,3.158204548207907104e+10
1.539850000000000000e+05,9.926000000000000000e+03,7.039262121589830992e-04,9.880828360568061873e-04,2.096265624133581828e-04,4.513101998226642887e-04,1.994464939610840460e-04,4.561896711653064566e-04,6.640247352533644150e-04,7.316958670662034297e-04,2.459662173605694580e+11
        }\tableAdaptiveOne

        \pgfplotstableread[col sep=comma]{
ndof,nelem,res,err,res1,res2,res3,res4,err1,err2,cond
6.440000000000000000e+02,2.400000000000000000e+01,3.847702620941702983e-02,1.131196847568739644e-01,5.743021451740769855e-03,2.150090107081665050e-02,8.633463227888469455e-03,3.017737258875701281e-02,8.061772967664361522e-02,7.935266058096705766e-02,7.006228862074787321e+04
2.560000000000000000e+03,9.600000000000000000e+01,2.273908253830363504e-02,7.156810052386690379e-02,1.763513034254757640e-03,1.086925323657852820e-02,4.748736181433878267e-03,1.932006043520741873e-02,5.084724758081463070e-02,5.036417800430836866e-02,6.799243010727243382e+04
1.020800000000000000e+04,3.840000000000000000e+02,1.425692585739921775e-02,4.518126249798776195e-02,5.772156993588011854e-04,6.730218126621295531e-03,2.949944689420091245e-03,1.220363663984028700e-02,3.207450208855686102e-02,3.182094902235373807e-02,7.929658784635017219e+04
4.076800000000000000e+04,1.536000000000000000e+03,8.976453842815396156e-03,2.848202059630345478e-02,1.929824724181147897e-04,4.233353394718327931e-03,1.853763013109814258e-03,7.692968418249724466e-03,2.021604058194495512e-02,2.006332974452172332e-02,3.042834676585858106e+05
1.629440000000000000e+05,6.144000000000000000e+03,5.654525003601752232e-03,1.794713443837424663e-02,6.590397029447043066e-05,2.666520798019271239e-03,1.167026285081042801e-03,4.847373120228893120e-03,1.273806090224119739e-02,1.264284141322175904e-02,1.205355947972324211e+06
        }\tableUniformTwo

        \pgfplotstableread[col sep=comma]{
ndof,nelem,res,err,res1,res2,res3,res4,err1,err2,cond
6.440000000000000000e+02,2.400000000000000000e+01,3.847702620941702983e-02,1.131196847568739644e-01,5.743021451740769855e-03,2.150090107081665050e-02,8.633463227888469455e-03,3.017737258875701281e-02,8.061772967664361522e-02,7.935266058096705766e-02,7.006228862074839708e+04
1.228000000000000000e+03,4.600000000000000000e+01,2.766348712976239577e-02,7.671488518218108621e-02,3.618531229729282823e-03,1.529989583881437819e-02,6.874136424179911123e-03,2.169871396330199168e-02,5.410487767489381705e-02,5.438598919115111147e-02,7.266397462063950661e+04
1.812000000000000000e+03,6.800000000000000000e+01,2.191528436873128871e-02,5.319009979447197889e-02,2.469583004351433731e-03,1.128592104956859139e-02,4.423009160546887732e-03,1.808993711364413268e-02,3.632871834469516420e-02,3.885113820183003219e-02,8.583580719178859727e+04
2.608000000000000000e+03,9.800000000000000000e+01,1.577085774849998692e-02,3.565247261401977008e-02,2.173384840600311139e-03,9.417338822730016851e-03,3.355471988001987142e-03,1.200212021659951764e-02,2.412096496923229680e-02,2.625410163091661947e-02,3.320778297554539749e+05
3.537000000000000000e+03,1.330000000000000000e+02,1.017091125036426527e-02,2.321000291408364946e-02,1.539937484950214437e-03,5.852907755644591245e-03,2.028254811531791588e-03,7.918691901358779553e-03,1.576822328544433618e-02,1.703136429333078003e-02,1.311948035012585111e+06
4.651000000000000000e+03,1.750000000000000000e+02,8.092902219110827480e-03,1.695784923961482851e-02,1.510838194245316471e-03,5.193220960524437026e-03,1.645009350504783195e-03,5.791099599456047374e-03,1.140067022256757191e-02,1.255361976920468109e-02,5.234429868500921875e+06
5.873000000000000000e+03,2.210000000000000000e+02,4.191811761451745555e-03,9.126518237599393416e-03,5.782971443977294197e-04,2.680996191051443923e-03,9.722850580446663286e-04,3.017246997794071663e-03,6.301598558113834063e-03,6.601756641501717489e-03,2.190076274810091034e+07
7.464000000000000000e+03,2.810000000000000000e+02,3.322800901646483923e-03,6.465390771769898536e-03,5.532126236288504662e-04,2.412354060919971559e-03,8.221866556008502244e-04,2.059009133820241446e-03,4.479860580459394531e-03,4.661772947209397513e-03,8.759603005152300000e+07
8.924000000000000000e+03,3.360000000000000000e+02,2.228831305062019610e-03,4.396070553795494024e-03,4.331369736184238197e-04,1.610807235518265821e-03,5.242480816484671206e-04,1.382224781791673756e-03,3.011409317499701057e-03,3.202631767222673901e-03,3.501953099590717554e+08
1.096800000000000000e+04,4.130000000000000000e+02,1.578989731468655004e-03,2.960255868633515988e-03,3.251040863768998781e-04,1.139572940980586140e-03,3.645787740928859333e-04,9.777380706356283372e-04,2.002006257492097274e-03,2.180615911329101490e-03,1.407959411232272863e+09
1.375300000000000000e+04,5.180000000000000000e+02,1.047917380671772305e-03,1.868951829161593107e-03,2.366059575763547734e-04,7.740218403094632424e-04,2.656985115482563128e-04,6.102810411919737214e-04,1.282711104065344413e-03,1.359276705911616027e-03,5.358758530384290695e+09
1.709700000000000000e+04,6.440000000000000000e+02,7.276033096020212724e-04,1.268979328044669711e-03,9.938593963253733031e-05,5.573056455429883319e-04,1.982301981505454013e-04,4.118788864876185172e-04,8.885892362114583380e-04,9.059347130416393543e-04,2.140637755832029343e+10
2.290300000000000000e+04,8.630000000000000000e+02,4.581817818729878996e-04,8.247713294497135712e-04,7.161598234707830025e-05,3.473513552351447688e-04,1.181162347278577753e-04,2.649477069589236625e-04,5.742169281090586880e-04,5.920495463685818702e-04,8.598030341170127869e+10
3.033000000000000000e+04,1.143000000000000000e+03,2.840047592173665738e-04,5.204948069465009739e-04,5.710321851807570307e-05,2.129595130612663286e-04,7.199030880598332786e-05,1.639010888545852987e-04,3.593620015964348502e-04,3.765286122818289677e-04,3.379660005989685059e+11
4.061900000000000000e+04,1.531000000000000000e+03,1.793346164369370461e-04,3.139209865195773759e-04,4.072699651197968045e-05,1.342704307694953241e-04,4.689435312368020790e-05,1.013636398013714675e-04,2.173399509075229575e-04,2.265165148878560224e-04,1.327079637776505615e+12
5.459100000000000000e+04,2.058000000000000000e+03,1.165553772241248052e-04,2.086956456297543640e-04,1.931101311277747416e-05,8.821300265701417971e-05,3.105715290295322904e-05,6.682933598512244815e-05,1.462634026751136499e-04,1.488653403002777897e-04,5.488579746794000977e+12
7.214200000000000000e+04,2.720000000000000000e+03,7.463836027364426377e-05,1.317633647283836089e-04,1.080184212952104163e-05,5.832684088722381500e-05,2.003876343289723555e-05,4.062797832952164535e-05,9.277248989718397552e-05,9.356735222726956158e-05,2.322314511983621094e+13
9.388500000000000000e+04,3.540000000000000000e+03,4.938128655286498717e-05,9.086149518680140412e-05,8.854870207830906019e-06,3.761878423627395472e-05,1.265094726745719018e-05,2.801576950646391534e-05,6.345258671675551708e-05,6.503522542848397497e-05,9.153998881275926562e+13
1.297270000000000000e+05,4.892000000000000000e+03,2.529581014309899224e-04,5.156777610713853596e-04,6.684943304108737542e-06,1.110013994403902455e-04,3.174833342094263107e-05,2.249752103229582208e-04,3.703016402337507441e-04,3.588875151127302564e-04,3.712418396069915000e+14
        }\tableAdaptiveTwo

        \pgfplotstableread[col sep=comma]{
ndof,nelem,res,err,res1,res2,res3,res4,err1,err2,cond
9.800000000000000000e+02,2.400000000000000000e+01,1.579640099765125191e-02,7.745260268776205104e-02,1.370050608843645466e-03,8.724146622972304282e-03,2.376964189468375958e-03,1.287977300850719382e-02,5.590315601602467710e-02,5.360730183992037240e-02,3.403034996212002006e+05
3.904000000000000000e+03,9.600000000000000000e+01,9.878890492758007555e-03,4.874987673361046159e-02,4.597026614927667422e-04,5.485531215102422073e-03,1.487285267298336976e-03,8.067098651481540567e-03,3.519224418065177579e-02,3.373509198255723179e-02,3.247330292022909853e+05
1.558400000000000000e+04,3.840000000000000000e+02,6.221388043142065012e-03,3.073272744655416305e-02,1.566223248604904021e-04,3.455887714721620364e-03,9.331394800182685537e-04,5.085983626081959363e-03,2.218452971733632373e-02,2.126845498678276922e-02,3.194800256842822419e+05
6.227200000000000000e+04,1.536000000000000000e+03,3.919371882127006476e-03,1.936717337819044912e-02,5.457295902936510284e-05,2.177670252939919483e-03,5.869537918074706944e-04,3.204954798167057497e-03,1.398004532577559744e-02,1.340319877305964046e-02,9.179296387696918100e+05
2.489600000000000000e+05,6.144000000000000000e+03,2.469152167699198777e-03,1.220237508991658774e-02,1.951126756177166505e-05,1.372042380208051858e-03,3.695455996560514716e-04,2.019224478451975540e-03,8.808159391221696832e-03,8.444778621956361420e-03,3.638998359430677257e+06
        }\tableUniformThree

        \pgfplotstableread[col sep=comma]{
ndof,nelem,res,err,res1,res2,res3,res4,err1,err2,cond
9.800000000000000000e+02,2.400000000000000000e+01,1.579640099765125191e-02,7.745260268776205104e-02,1.370050608843645466e-03,8.724146622972304282e-03,2.376964189468375958e-03,1.287977300850719382e-02,5.590315601602467710e-02,5.360730183992037240e-02,3.403034996212003171e+05
1.872000000000000000e+03,4.600000000000000000e+01,1.039356507247008593e-02,5.031147165047349101e-02,5.929456298192313773e-04,6.092947895662580160e-03,1.403739324888181447e-03,8.281311022069621705e-03,3.626497685769639739e-02,3.487256275565568664e-02,3.339195624512357172e+05
2.885000000000000000e+03,7.100000000000000000e+01,6.632988268975648430e-03,3.409387589654797740e-02,2.335733114688661556e-04,3.886403303307548428e-03,9.616437310822971126e-04,5.283283787858489115e-03,2.462999996533182107e-02,2.357446659750691426e-02,3.451049727763903211e+05
3.939000000000000000e+03,9.700000000000000000e+01,4.226750430316881559e-03,2.150737668553739240e-02,1.203941624656627895e-04,2.493166838520854967e-03,6.034425663264664490e-04,3.357216202491790803e-03,1.553265980299863534e-02,1.487628083016409065e-02,1.003872316667611594e+06
4.912000000000000000e+03,1.210000000000000000e+02,2.730742758590313374e-03,1.358305083504622070e-02,9.909063793046497847e-05,1.644824769234428568e-03,3.821538681097262332e-04,2.143746010764720038e-03,9.806494877595418880e-03,9.398506700702226324e-03,3.981125962346578017e+06
5.885000000000000000e+03,1.450000000000000000e+02,1.824054917509096560e-03,8.607757977415784328e-03,9.615546400788242698e-05,1.145561788503714597e-03,2.441427577275650813e-04,1.394995688667135256e-03,6.209528915401064594e-03,5.961144860390762794e-03,1.589067077477816120e+07
7.506000000000000000e+03,1.850000000000000000e+02,1.286101238448659687e-03,5.495209525334464329e-03,9.555393546061036636e-05,8.604566953766562163e-04,1.588935429746196777e-04,9.377062217665108187e-04,3.957253206433091093e-03,3.812804058367233288e-03,6.354630405612101406e+07
9.775000000000000000e+03,2.410000000000000000e+02,8.714848727395625369e-04,3.510754538010496008e-03,7.141875551132739492e-05,5.950943062094559083e-04,1.025179910807744090e-04,6.242900551802448147e-04,2.524334443650425898e-03,2.439904310164926413e-03,2.541440411169292629e+08
1.115300000000000000e+04,2.750000000000000000e+02,5.610684253399463852e-04,2.239490469807899845e-03,3.614356074781348098e-05,3.858755426603923040e-04,6.532449535521635893e-05,4.004050409436056015e-04,1.610502489033887588e-03,1.556148867292606692e-03,1.020233728323160887e+09
1.261300000000000000e+04,3.110000000000000000e+02,4.138554201614976901e-04,1.451800461600342420e-03,3.051776134266695630e-05,2.937800656097847452e-04,4.361380643614078362e-05,2.865939356967334996e-04,1.041361405257881919e-03,1.011578471470351648e-03,4.082341180966885090e+09
1.451800000000000000e+04,3.580000000000000000e+02,2.598944594397424785e-04,9.173026361383381259e-04,2.215181140932124192e-05,1.864393711282921284e-04,2.760272554380993208e-05,1.775749918614720013e-04,6.586792533326811267e-04,6.384244414928414834e-04,1.580540114787137222e+10
1.682900000000000000e+04,4.150000000000000000e+02,1.885224102100716080e-04,5.838493699383948046e-04,1.390655347084193207e-05,1.392944025682014693e-04,1.861903918888980664e-05,1.248907831460253296e-04,4.178830836031713391e-04,4.077423392484108306e-04,6.474824001145267487e+10
2.120500000000000000e+04,5.230000000000000000e+02,1.268342996471843022e-04,3.755186949483942734e-04,1.122595519158902115e-05,9.256984126512383471e-05,1.206985784804454359e-05,8.512379524378567680e-05,2.683836638821142893e-04,2.626490038377398036e-04,2.589929274680123596e+11
2.562300000000000000e+04,6.320000000000000000e+02,8.227510996148546787e-05,2.390920068307111207e-04,4.763044222881068083e-06,6.012641199763711045e-05,7.742045455098618334e-06,5.542005450526912808e-05,1.705147784840336936e-04,1.675998151815139719e-04,1.001990135899522339e+12
3.166000000000000000e+04,7.810000000000000000e+02,5.399157872754615616e-05,1.525528944145721627e-04,2.641670925701570216e-06,4.053399144016093251e-05,4.997030440247689441e-06,3.521558423370705428e-05,1.085717843393566614e-04,1.071660078552512634e-04,4.011303716025214355e+12
3.899300000000000000e+04,9.620000000000000000e+02,3.439634968241429830e-05,9.765917907600624257e-05,1.390516986193124636e-06,2.507788333204449734e-05,3.231105258688908502e-06,2.327735080239308273e-05,6.967625526231181862e-05,6.842904873239610743e-05,1.634752065090042578e+13
4.912300000000000000e+04,1.212000000000000000e+03,2.454634069079307011e-05,6.715404851809319002e-05,1.202436218302408555e-06,1.653312396436414989e-05,2.274664911216767049e-06,1.795991926401804454e-05,4.787890315557513248e-05,4.708796942945679936e-05,6.539011166224096094e+13
5.678000000000000000e+04,1.401000000000000000e+03,1.629489535089447263e-05,3.973776052138151669e-05,1.154489986382009330e-06,1.269897126152423009e-05,1.392120848858520415e-06,1.004932315178996441e-05,2.821126958365616750e-05,2.798595861736635142e-05,2.459732738189519062e+14
7.087800000000000000e+04,1.749000000000000000e+03,9.575786073179552835e-06,2.439309361035165650e-05,6.551446549960352244e-07,7.426700503813102715e-06,8.305738689725898131e-07,5.951531825939585947e-06,1.734842094755755993e-05,1.714804089421662474e-05,9.843618264263920000e+14
        }\tableAdaptiveThree

        %
        %
        \addplot+ [line1, minorline, forget plot] table [x=ndof, y=effRes] {\tableUniformZero};
        \addplot+ [line2, minorline, forget plot] table [x=ndof, y=effRes] {\tableUniformOne};
        \addplot+ [line3, minorline, forget plot] table [x=ndof, y=effRes] {\tableUniformTwo};
        \addplot+ [line4, minorline, forget plot] table [x=ndof, y=effRes] {\tableUniformThree};

        \addplot+ [line1, majorline, forget plot] table [x=ndof, y=effRes] {\tableAdaptiveZero};
        \addplot+ [line2, majorline, forget plot] table [x=ndof, y=effRes] {\tableAdaptiveOne};
        \addplot+ [line3, majorline, forget plot] table [x=ndof, y=effRes] {\tableAdaptiveTwo};
        \addplot+ [line4, majorline, forget plot] table [x=ndof, y=effRes] {\tableAdaptiveThree};
    \end{semilogxaxis}
\end{tikzpicture}

%% file: figures/plot_opdls_lshape_efficiency_degree.tex
\begin{tikzpicture}[>=stealth]
    \colorlet{col0}{TUblue}
    \colorlet{col1}{TUgreen}
    \colorlet{col2}{TUmagenta}
    \colorlet{col3}{TUyellow}
    \colorlet{col4}{purple}
    \colorlet{col5}{green}
    \pgfplotsset{%
        degdefault/.style = {%
            mark = *,%
            mark size = 2pt,%
            every mark/.append style = {solid},%
            gray,%
            every mark/.append style = {fill = gray!60!white}%
        },%
        deg0/.style = {%
            degdefault,%
            col0,%
            every mark/.append style = {fill = col0!60!white}%
        },%
        deg1/.style = {%
            degdefault,%
            mark = triangle*,%
            mark size = 2.75pt,%
            col1,%
            every mark/.append style = {fill = col1!60!white}%
        },%
        deg2/.style = {%
            degdefault,%
            mark = square*,%
            mark size = 1.66pt,%
            col2,%
            every mark/.append style = {fill = col2!60!white}%
        },%
        deg3/.style = {%
            degdefault,%
            mark = pentagon*,%
            mark size = 2.2pt,%
            col3,%
            every mark/.append style = {fill = col3!60!white}%
        },%
        deg4/.style = {%
            degdefault,%
            mark = diamond*,%
            mark size = 2.75pt,%
            col4,%
            every mark/.append style = {fill = col4!60!white}%
        },%
        deg5/.style = {%
            degdefault,%
            mark = halfsquare*,%
            mark size = 1.66pt,%
            col5,%
            every mark/.append style = {fill = col5!60!white}%
        },%
        uniform/.style = {%
            dashed,%
            every mark/.append style = {%
                fill = black!20!white
            }%
        },%
        adaptive/.style = {%
            solid%
        }%
    }

    %
    %
    \pgfplotstableset{%
        create on use/effRes/.style={%
            create col/expr={%
                \thisrow{res} / \thisrow{err}
            }
        }
    }

    \begin{semilogxaxis}[%
            width            = 0.36\textwidth,%
            xlabel           = ndof,%
            ymax             = 1,%
            ymin             = 0,%
            ymajorgrids      = true,%
            font             = \footnotesize,%
            grid style       = {%
                densely dotted,%
                semithick%
            },%
            legend style     = {%
                legend pos  = south west %
            }%
        ]

        \pgfplotstableread[col sep=comma]{
ndof,nelem,res,err,res1,res2,res3,res4,err1,err2,cond
1.440000000000000000e+02,2.400000000000000000e+01,6.306198149120191276e-01,1.008644855299557497e+00,3.366965743021015922e-01,3.651053388672639288e-01,2.652373984422513042e-01,2.840140521366953386e-01,6.762132392680514359e-01,7.483983559314346667e-01,2.296586155100305859e+02
5.760000000000000000e+02,9.600000000000000000e+01,3.891634116065782467e-01,5.745268224882443819e-01,1.620156561466426781e-01,2.310425831275144948e-01,1.314317966611835697e-01,2.335467741263464181e-01,3.810417180298759643e-01,4.299863705738576902e-01,4.384474670547972437e+02
2.304000000000000000e+03,3.840000000000000000e+02,2.199346201832386161e-01,3.132415464677545947e-01,7.722439043103720979e-02,1.329370244949013291e-01,4.796426422386291383e-02,1.497825342074489141e-01,2.070719725645633613e-01,2.350350284781635712e-01,5.219306858866519178e+03
9.216000000000000000e+03,1.536000000000000000e+03,1.209213915613219781e-01,1.733866068509177238e-01,3.769733884062748119e-02,7.238555346200944984e-02,1.632586479936865356e-02,8.771938993183199529e-02,1.148033304170861940e-01,1.299350251488029240e-01,7.644441501019104908e+04
3.686400000000000000e+04,6.144000000000000000e+03,6.691209920672000488e-02,9.851845448837180663e-02,1.866593364274058559e-02,3.912977304488995250e-02,5.704439904109297362e-03,5.064713386096535847e-02,6.549793757790635385e-02,7.359283965045819442e-02,1.190822091908842791e+06
1.474560000000000000e+05,2.457600000000000000e+04,3.776848089809536901e-02,5.739465368140826695e-02,9.294922804628987303e-03,2.146013391463778147e-02,2.079038829906257931e-03,2.958382681999287625e-02,3.833217236807209827e-02,4.271757053898535966e-02,1.888587027434254438e+07
        }\tableUniformZero

        \pgfplotstableread[col sep=comma]{
ndof,nelem,res,err,res1,res2,res3,res4,err1,err2,cond
1.440000000000000000e+02,2.400000000000000000e+01,6.306198149120191276e-01,1.008644855299557497e+00,3.366965743021015922e-01,3.651053388672639288e-01,2.652373984422513042e-01,2.840140521366953386e-01,6.762132392680514359e-01,7.483983559314346667e-01,2.296586155100324049e+02
4.200000000000000000e+02,7.000000000000000000e+01,4.394954169616491213e-01,6.715707909097310413e-01,2.134046978046295895e-01,2.398934244640127944e-01,1.626743030107819787e-01,2.521960993435343745e-01,4.599266805425582527e-01,4.893616001774407365e-01,4.248331973432046880e+02
1.044000000000000000e+03,1.740000000000000000e+02,3.080995764558332461e-01,4.448715612364610994e-01,1.206147436548664076e-01,1.831248710847279693e-01,8.866217124494756019e-02,1.974379234232068736e-01,2.920082058547046611e-01,3.356216824200766724e-01,2.851195280339358760e+03
2.352000000000000000e+03,3.920000000000000000e+02,2.283813563786274703e-01,3.195496515583942654e-01,8.852265376466934843e-02,1.299399484381094927e-01,4.840039087579297433e-02,1.584133696212980114e-01,2.069620642343262262e-01,2.434721416896762913e-01,2.390779644215623557e+04
4.500000000000000000e+03,7.500000000000000000e+02,1.636713346895881938e-01,2.212306526684754138e-01,5.883268021431856815e-02,9.802389768926991764e-02,2.600950890058896514e-02,1.142008864520090011e-01,1.419485766226470380e-01,1.696867799061674009e-01,1.575015820655770076e+05
9.180000000000000000e+03,1.530000000000000000e+03,1.187388127792072257e-01,1.588100720996714343e-01,4.122529757414088170e-02,6.988946933720374677e-02,1.314419176578006734e-02,8.568589614824374812e-02,1.006913074464998586e-01,1.228083857275931723e-01,1.226046921408378752e+06
1.662600000000000000e+04,2.771000000000000000e+03,8.648401513452842748e-02,1.137322684093443892e-01,3.016879753258562905e-02,5.267870751491721126e-02,7.324607844276403469e-03,6.116070980386480171e-02,7.275750577325909574e-02,8.741495427665293971e-02,8.826925579184673727e+06
3.108000000000000000e+04,5.180000000000000000e+03,6.407539385855630376e-02,8.380844066854448160e-02,2.230113949688580441e-02,3.842147594042780168e-02,3.825850255794115921e-03,4.601595735882178101e-02,5.311393808999341493e-02,6.482873057422372520e-02,6.463194480318483710e+07
5.716800000000000000e+04,9.528000000000000000e+03,4.660186986402903991e-02,6.027330625225215538e-02,1.622121547324421181e-02,2.855032938168600154e-02,2.134798690145196885e-03,3.299890559009913127e-02,3.851564025972416261e-02,4.636180434324470301e-02,4.729308135390356779e+08
1.078380000000000000e+05,1.797300000000000000e+04,3.443687955697830733e-02,4.437535905440612161e-02,1.192594613064655411e-02,2.113007189858816595e-02,1.108489517543145104e-03,2.441232875128364344e-02,2.799514469590566254e-02,3.443028266893495559e-02,3.438310353276770592e+09
        }\tableAdaptiveZero

        \pgfplotstableread[col sep=comma]{
ndof,nelem,res,err,res1,res2,res3,res4,err1,err2,cond
3.360000000000000000e+02,2.400000000000000000e+01,1.410774541754422629e-01,2.532432142128129837e-01,5.793087820617605355e-02,7.321032798409224385e-02,4.723246674586529387e-02,9.463721988967961629e-02,1.753891392472060329e-01,1.826766963215638306e-01,7.230740466181046941e+03
1.344000000000000000e+03,9.600000000000000000e+01,6.044335677174061255e-02,1.291863839669573910e-01,1.469509362887662200e-02,2.959645169782602198e-02,1.298982141402372885e-02,4.891592979392801049e-02,8.974391242259822932e-02,9.292551848412342053e-02,6.953307961995453297e+03
5.376000000000000000e+03,3.840000000000000000e+02,3.319076797554559172e-02,7.916865575463447047e-02,3.994376469488723075e-03,1.409646253926915768e-02,3.964753671979430981e-03,2.951681736952689247e-02,5.494290076540755269e-02,5.699959394134658519e-02,1.673272813861177201e+04
2.150400000000000000e+04,1.536000000000000000e+03,2.031558977735474503e-02,4.987345534912807771e-02,1.143409419566413041e-03,8.165126979841825713e-03,1.352857110533361605e-03,1.851799886441691731e-02,3.458845579109927348e-02,3.593049226560430381e-02,7.729624800803052494e+04
8.601600000000000000e+04,6.144000000000000000e+03,1.273734677141928444e-02,3.146054049942184627e-02,3.442315075950557412e-04,5.061888876694767628e-03,4.836574720980640050e-04,1.167325421184009504e-02,2.181382519503181558e-02,2.266986190686552660e-02,1.194891384845554363e+06
3.440640000000000000e+05,2.457600000000000000e+04,8.018827900752740698e-03,1.983188318091227628e-02,1.091357260475482187e-04,3.180962498518550077e-03,1.777413066239565630e-04,7.357960037002686700e-03,1.375004510734577703e-02,1.429125082164982644e-02,1.890740388982083276e+07
        }\tableUniformOne

        \pgfplotstableread[col sep=comma]{
ndof,nelem,res,err,res1,res2,res3,res4,err1,err2,cond
3.360000000000000000e+02,2.400000000000000000e+01,1.410774541754422629e-01,2.532432142128129837e-01,5.793087820617605355e-02,7.321032798409224385e-02,4.723246674586529387e-02,9.463721988967961629e-02,1.753891392472060329e-01,1.826766963215638306e-01,7.230740466181054217e+03
8.540000000000000000e+02,6.100000000000000000e+01,8.748715718044233758e-02,1.684296540347229310e-01,3.311526801709450962e-02,4.457825042334847415e-02,2.055699479589141038e-02,6.440164011751835227e-02,1.152801018908801434e-01,1.227967689570240856e-01,6.443302249040961215e+03
1.358000000000000000e+03,9.700000000000000000e+01,5.759116983626711567e-02,1.050673055722659566e-01,1.654478567884715492e-02,3.150369769104819501e-02,1.230241574477199118e-02,4.357958821036068964e-02,7.151586571125451153e-02,7.697155092490802486e-02,1.550524917010031277e+04
2.170000000000000000e+03,1.550000000000000000e+02,4.141981427373113422e-02,7.172515584232455610e-02,1.294810477267752807e-02,2.391145890572530727e-02,7.650618415257546753e-03,3.029286662120440204e-02,4.865311310762630848e-02,5.270078334846894919e-02,6.252838230085997930e+04
3.262000000000000000e+03,2.330000000000000000e+02,3.033955551758293553e-02,5.015740412317001157e-02,9.967159336669190531e-03,1.830299842149288661e-02,5.276672396546587714e-03,2.140797377063736670e-02,3.396670332767071782e-02,3.690566695542835118e-02,2.499584164417928259e+05
4.634000000000000000e+03,3.310000000000000000e+02,2.063765791858325135e-02,3.282160479639936734e-02,7.336568817536011020e-03,1.261348304823302173e-02,3.043348420066842829e-03,1.427325324290278780e-02,2.259776723547589003e-02,2.380333290491641607e-02,1.018860286316490616e+06
6.524000000000000000e+03,4.660000000000000000e+02,1.400410033005189750e-02,2.208387186187773729e-02,4.174915921989151221e-03,8.904496261487114539e-03,1.627868997453578047e-03,9.835898126494229191e-03,1.505494962739059583e-02,1.615691456091065958e-02,4.076691744258814491e+06
9.702000000000000000e+03,6.930000000000000000e+02,9.742231738923525300e-03,1.471505667422966061e-02,2.913031071305545623e-03,6.418956689802151130e-03,9.715999900890765536e-04,6.654195496542055638e-03,9.977065324930664211e-03,1.081624058662977381e-02,1.630663577492151223e+07
1.527400000000000000e+04,1.091000000000000000e+03,6.534829400284570926e-03,9.775121048782925079e-03,1.886694431174944231e-03,4.297115553398902013e-03,5.309073156346052081e-04,4.516338645090722109e-03,6.609974557783129812e-03,7.201474006328045102e-03,6.521203520578209311e+07
2.233000000000000000e+04,1.595000000000000000e+03,4.280859669157644799e-03,6.393246211616394545e-03,1.192913015934106988e-03,2.833009672566197228e-03,2.609563914990130720e-04,2.967941374137594585e-03,4.345083261514961180e-03,4.689759969641294282e-03,2.928693184765223861e+08
3.235400000000000000e+04,2.311000000000000000e+03,2.939255400262257793e-03,4.278937844654475721e-03,8.798031802792400692e-04,1.972485081737179041e-03,1.611462988583335393e-04,1.987084080916036097e-03,2.897366575165785836e-03,3.148741972205469815e-03,1.572087604532146931e+09
5.154800000000000000e+04,3.682000000000000000e+03,1.911687575861123782e-03,2.750815083306280481e-03,5.489416556833490342e-04,1.295184838722848340e-03,8.661533450177437985e-05,1.291590671896096498e-03,1.856267786951105829e-03,2.030087073421481432e-03,1.000104831224785614e+10
7.842800000000000000e+04,5.602000000000000000e+03,1.228069321714765324e-03,1.768564637063715613e-03,3.573866704758233630e-04,8.268487638269265494e-04,4.109147874345383264e-05,8.337035677240038590e-04,1.206569594144262653e-03,1.293062523607757617e-03,6.375636762051554108e+10
1.176980000000000000e+05,8.407000000000000000e+03,8.248451186374084311e-04,1.167742506468436893e-03,2.522381977433809210e-04,5.591888267446263857e-04,2.435729739990704620e-05,5.508719809960150844e-04,7.943734678354903104e-04,8.559166752740630891e-04,3.713117678216314087e+11
        }\tableAdaptiveOne

        \pgfplotstableread[col sep=comma]{
ndof,nelem,res,err,res1,res2,res3,res4,err1,err2,cond
6.000000000000000000e+02,2.400000000000000000e+01,3.819561679950280592e-02,1.091662736805051048e-01,5.670447324352385562e-03,2.007928134485041649e-02,1.297152897598293463e-02,2.924573586928157914e-02,7.767238907163920292e-02,7.670903001075458916e-02,7.007366623479264672e+04
2.400000000000000000e+03,9.600000000000000000e+01,2.277932785179561520e-02,7.183010594498503987e-02,1.768816172240626370e-03,1.105710855070412418e-02,4.088785592132335514e-03,1.941111149851879394e-02,5.104347467606239891e-02,5.053837960463263695e-02,6.798193358015047852e+04
9.600000000000000000e+03,3.840000000000000000e+02,1.431082609520344801e-02,4.584306043342896431e-02,5.823146197496271688e-04,7.021356604078357946e-03,1.366474247951473273e-03,1.238119351188859789e-02,3.256631628527141553e-02,3.226486066157242816e-02,1.010560841474677291e+05
3.840000000000000000e+04,1.536000000000000000e+03,9.014753329646749483e-03,2.898781487936142806e-02,1.952364080309571383e-04,4.446171942711974086e-03,4.613407076057559155e-04,7.826006653888250772e-03,2.059176009742792746e-02,2.040276519421086940e-02,4.009493434865098679e+05
1.536000000000000000e+05,6.144000000000000000e+03,5.679343816742773998e-03,1.828203560329835775e-02,6.681547802816383604e-05,2.805554438313187862e-03,1.588234668496425425e-04,4.934989491448724466e-03,1.298679488682430092e-02,1.286763320769688211e-02,1.601304117351882625e+06
        }\tableUniformTwo

        \pgfplotstableread[col sep=comma]{
ndof,nelem,res,err,res1,res2,res3,res4,err1,err2,cond
6.000000000000000000e+02,2.400000000000000000e+01,3.819561679950280592e-02,1.091662736805051048e-01,5.670447324352385562e-03,2.007928134485041649e-02,1.297152897598293463e-02,2.924573586928157914e-02,7.767238907163920292e-02,7.670903001075458916e-02,7.007366623479271948e+04
1.150000000000000000e+03,4.600000000000000000e+01,2.756964904917026324e-02,7.703019285530524063e-02,3.615022502163543674e-03,1.521574048550653965e-02,6.269939773021809305e-03,2.182169236616771438e-02,5.433206919464296514e-02,5.460473302154351255e-02,7.276368591044633649e+04
1.700000000000000000e+03,6.800000000000000000e+01,2.198641813043256815e-02,5.357698815880844456e-02,2.465229996057532232e-03,1.128756341673838731e-02,3.997719197554329938e-03,1.827387141804841242e-02,3.660315366322466901e-02,3.912419714288668110e-02,1.063015404990217794e+05
2.350000000000000000e+03,9.400000000000000000e+01,1.585428130065361443e-02,3.607674668813766472e-02,2.167285357792936035e-03,9.343712303961203974e-03,3.318469608709834814e-03,1.217965145539988414e-02,2.441775181880920936e-02,2.655757985425426237e-02,4.247983726905161748e+05
3.325000000000000000e+03,1.330000000000000000e+02,1.022165500908132733e-02,2.346448107127773949e-02,1.538460393819226319e-03,5.897545309268036225e-03,1.801554614198701797e-03,8.005543769071345184e-03,1.595606973065698611e-02,1.720423525457506417e-02,1.689795423082960071e+06
4.375000000000000000e+03,1.750000000000000000e+02,8.128557037679385450e-03,1.706780558628640324e-02,1.509294112177854029e-03,5.187770089703944328e-03,1.748342259932253624e-03,5.815996186059663696e-03,1.148245800000523634e-02,1.262787178464306985e-02,6.751156253198359162e+06
5.525000000000000000e+03,2.210000000000000000e+02,4.207805289181119381e-03,9.264352108815279707e-03,5.784715725420058827e-04,2.734282212934087235e-03,6.544591268083478153e-04,3.076748287200448129e-03,6.403523852472148802e-03,6.695005770492692618e-03,2.825060354364511371e+07
7.025000000000000000e+03,2.810000000000000000e+02,3.331159766031210456e-03,6.545093333746404736e-03,5.533849120563855747e-04,2.447544432006348511e-03,6.419410708908593029e-04,2.094714405543551556e-03,4.540244928450749146e-03,4.714278601984501008e-03,1.130028801320911646e+08
8.350000000000000000e+03,3.340000000000000000e+02,2.198225575137922796e-03,4.409293946533826905e-03,4.301334292020647003e-04,1.610680439374087645e-03,3.441330669457995373e-04,1.390849332863017956e-03,3.034411508149149374e-03,3.199096732853174644e-03,4.517641988687487841e+08
1.032500000000000000e+04,4.130000000000000000e+02,1.538033157362616492e-03,2.946448251106259866e-03,3.212353954825196271e-04,1.131930137895854633e-03,2.075089912616295753e-04,9.685184537760034530e-04,2.010119562325146398e-03,2.154292608167448185e-03,1.816508657687121153e+09
1.315000000000000000e+04,5.260000000000000000e+02,9.890414840953529619e-04,1.842006646037166887e-03,1.984798353354147937e-04,7.491959651551121074e-04,1.105492321722191008e-04,6.043948091175163281e-04,1.278220901324935745e-03,1.326325680766665205e-03,6.908659052490003586e+09
1.712500000000000000e+04,6.850000000000000000e+02,6.485146494104617040e-04,1.174625050322526849e-03,1.219014788500662313e-04,5.082158015828334334e-04,5.862216188487281185e-05,3.794620153738755735e-04,8.263951438649462739e-04,8.347544998630616929e-04,2.726510556326704788e+10
2.295000000000000000e+04,9.180000000000000000e+02,4.077883613519849753e-04,7.576825104850408226e-04,6.942790158930151870e-05,3.183220337087589889e-04,3.430147752301317224e-05,2.428283460189559467e-04,5.287449788117095728e-04,5.426891689321968523e-04,1.122120581666528778e+11
2.960000000000000000e+04,1.184000000000000000e+03,2.741965634972319912e-04,5.184784302531738600e-04,5.084198902757097190e-05,2.095272353955302046e-04,1.707900678980735733e-05,1.685392913833202164e-04,3.588055606703724519e-04,3.742705602499554296e-04,4.489444048630886230e+11
4.047500000000000000e+04,1.619000000000000000e+03,1.692209133215686657e-04,3.129745288966800216e-04,3.588779727677518207e-05,1.297913811140073280e-04,9.379447807776055638e-06,1.020490414425664178e-04,2.183110099257196200e-04,2.242618083475457600e-04,1.750882096834728027e+12
5.402500000000000000e+04,2.161000000000000000e+03,1.090935958472342309e-04,1.981103946960378827e-04,1.738611838880867725e-05,8.629806808459527308e-05,5.246929398166668794e-06,6.422031391240002609e-05,1.395930987571483229e-04,1.405755855972050470e-04,6.855385731618833984e+12
7.090000000000000000e+04,2.836000000000000000e+03,7.219652495714911895e-05,1.283887793841192320e-04,1.078768863498282686e-05,5.900769926371286420e-05,3.575501259532920860e-06,4.001588694692670645e-05,9.080658885280398546e-05,9.076255887017426935e-05,2.752245437487217188e+13
9.117500000000000000e+04,3.647000000000000000e+03,4.798932552939210887e-05,8.880366444441850887e-05,8.549489185352279410e-06,3.790147359310555705e-05,1.704783625335345602e-06,2.811500687706529060e-05,6.212129314302242643e-05,6.345892968681016264e-05,1.100902615881208438e+14
1.274500000000000000e+05,5.098000000000000000e+03,2.991038790260966415e-05,5.386633579543171908e-05,6.027534678034814011e-06,2.378441082139401839e-05,1.024344002718999720e-06,1.707491292607125270e-05,3.729533393569878284e-05,3.886695509878440022e-05,4.807198707087767500e+14
        }\tableAdaptiveTwo

        \pgfplotstableread[col sep=comma]{
ndof,nelem,res,err,res1,res2,res3,res4,err1,err2,cond
9.360000000000000000e+02,2.400000000000000000e+01,1.587287669783075603e-02,7.595987805906435908e-02,1.358141089144095960e-03,8.544501383718315288e-03,3.874310374247793946e-03,1.273125613104020108e-02,5.481240670708409368e-02,5.258805135888882515e-02,3.403802350252678152e+05
3.744000000000000000e+03,9.600000000000000000e+01,9.880771211355181619e-03,4.888935784724834227e-02,4.608958037078263225e-04,5.514618161488502902e-03,1.299875972795064586e-03,8.081863880170288977e-03,3.529654494296664835e-02,3.382784689875251161e-02,3.246558266713191988e+05
1.497600000000000000e+04,3.840000000000000000e+02,6.213848900663823151e-03,3.101060628718815479e-02,1.575809511539809273e-04,3.499670015181694741e-03,4.363194349487324749e-04,5.113611398966347484e-03,2.238955745375374382e-02,2.145612778028817627e-02,3.192977679735824349e+05
5.990400000000000000e+04,1.536000000000000000e+03,3.912984282204073905e-03,1.957581702869206144e-02,5.497430347534909090e-05,2.209542810322811992e-03,1.491579503757459234e-04,3.225538140206056417e-03,1.413379017229841363e-02,1.354431938881652843e-02,1.209546226606750395e+06
2.396160000000000000e+05,6.144000000000000000e+03,2.464803052673721507e-03,1.233995624749361329e-02,1.966621898909293986e-05,1.392847920339726530e-03,5.224274491968399392e-05,2.032759871383960082e-03,8.909497725342371832e-03,8.537878569772276716e-03,4.834927912092171609e+06
        }\tableUniformThree

        \pgfplotstableread[col sep=comma]{
ndof,nelem,res,err,res1,res2,res3,res4,err1,err2,cond
9.360000000000000000e+02,2.400000000000000000e+01,1.587287669783075603e-02,7.595987805906435908e-02,1.358141089144095960e-03,8.544501383718315288e-03,3.874310374247793946e-03,1.273125613104020108e-02,5.481240670708409368e-02,5.258805135888882515e-02,3.403802350252678152e+05
1.794000000000000000e+03,4.600000000000000000e+01,1.044915960577807480e-02,5.018233837711233436e-02,5.928998277704849995e-04,6.095577088923767248e-03,1.808404555794396401e-03,8.270853594208021214e-03,3.617391115849837024e-02,3.478096083336495309e-02,3.339776074129303452e+05
2.769000000000000000e+03,7.100000000000000000e+01,6.635266756571738325e-03,3.431627056434442824e-02,2.343898911081658217e-04,3.923572513949794352e-03,6.649669619238980818e-04,5.304264697122472103e-03,2.479432431245664342e-02,2.372441584810916493e-02,3.452127746248840704e+05
3.783000000000000000e+03,9.700000000000000000e+01,4.221917769610821630e-03,2.170500740314325674e-02,1.204489754034804478e-04,2.521598716824297298e-03,2.384367398600577596e-04,3.375613948613621646e-03,1.567808216511660385e-02,1.501016608816691834e-02,1.299883519846176263e+06
4.719000000000000000e+03,1.210000000000000000e+02,2.725873354535585970e-03,1.371783137792322721e-02,9.884956520065375069e-05,1.661408201999718453e-03,1.083837079515970263e-04,2.156058920021773368e-03,9.905557097018662799e-03,9.489933419726310007e-03,5.196154601200701669e+06
5.655000000000000000e+03,1.450000000000000000e+02,1.819834332671921252e-03,8.693959554343065610e-03,9.586001748769681756e-05,1.152889520420732807e-03,7.842489721117074867e-05,1.402605840790699995e-03,6.272830438739948958e-03,6.019678647516872247e-03,2.078143583830968291e+07
7.215000000000000000e+03,1.850000000000000000e+02,1.281645461708655342e-03,5.548619558195228181e-03,9.524621882136063782e-05,8.607615494533052631e-04,7.319812976987359805e-05,9.419526719000080942e-04,3.996384121383923273e-03,3.849167826418727414e-03,8.314545740445955098e+07
8.307000000000000000e+03,2.130000000000000000e+02,9.966580641392046736e-04,3.610273636436359557e-03,9.519906774514425973e-05,7.174788920416048802e-04,7.257209391171580187e-05,6.813382162877577516e-04,2.589364890799469160e-03,2.515803090912021017e-03,3.325787006826478839e+08
1.002300000000000000e+04,2.570000000000000000e+02,5.877306847250159030e-04,2.271458219838518987e-03,3.663291604626303205e-05,4.182903696924609787e-04,2.336163683546792045e-05,4.105761656211856102e-04,1.633552922775719906e-03,1.578298860470625856e-03,1.335089161239186764e+09
1.181700000000000000e+04,3.030000000000000000e+02,3.831165054876371225e-04,1.443372331996678470e-03,2.708890123921533854e-05,2.705521397650076398e-04,1.387692044826759351e-05,2.695429816451521028e-04,1.037688083674370753e-03,1.003258256768285814e-03,5.342453017069465637e+09
1.353300000000000000e+04,3.470000000000000000e+02,2.741880337021723171e-04,9.129440701841704930e-04,1.290109051990256363e-05,2.031636466690750790e-04,6.034013376168046660e-06,1.835776757766593840e-04,6.552165774825424943e-04,6.357343091862382393e-04,2.117872461866301727e+10
1.548300000000000000e+04,3.970000000000000000e+02,1.775337249757544603e-04,5.794029171933716199e-04,1.201643916437640030e-05,1.309859489858040493e-04,5.184404817014923740e-06,1.191202409691875542e-04,4.154488993188299199e-04,4.038687317767514348e-04,8.471489527912284851e+10
2.063100000000000000e+04,5.290000000000000000e+02,1.206242561955739540e-04,3.740274611655121113e-04,1.036028941375277699e-05,8.801327873889524557e-05,4.081937120081584325e-06,8.173050924211808384e-05,2.674512911938693488e-04,2.614695939199258401e-04,3.450490087660043335e+11
2.655900000000000000e+04,6.810000000000000000e+02,7.536952986073883816e-05,2.334359605725026333e-04,4.883640297637097618e-06,5.643651737489691186e-05,2.480672192711689101e-06,4.965361881008645192e-05,1.668595716373047076e-04,1.632489848097747509e-04,1.265135852051228271e+12
3.174600000000000000e+04,8.140000000000000000e+02,5.132304931960395271e-05,1.509797757820412331e-04,2.491197114071798622e-06,3.942300938049305776e-05,1.301756094422553341e-06,3.274112232128794270e-05,1.078271337470901000e-04,1.056797138673387930e-04,5.218576744270708984e+12
3.868800000000000000e+04,9.920000000000000000e+02,3.149745786349974575e-05,9.430840625361631090e-05,1.310130367562511299e-06,2.347964704294816849e-05,4.579761585551246562e-07,2.094922057881156937e-05,6.744224356288509327e-05,6.592131122256004679e-05,2.087430866446758594e+13
4.578600000000000000e+04,1.174000000000000000e+03,2.079312900766588045e-05,5.949765744365865313e-05,1.170176942825198382e-06,1.586918368930536052e-05,3.534898175050502724e-07,1.338017017215858154e-05,4.253500105738561608e-05,4.160222261287443600e-05,8.048023902612678125e+13
5.951400000000000000e+04,1.526000000000000000e+03,1.292739499076695692e-05,3.738336084600879885e-05,7.916873589213317810e-07,9.901393078428031436e-06,2.203662857167556843e-07,8.270708928478118392e-06,2.674158372548604874e-05,2.612285145223855400e-05,3.228427580781133125e+14
7.215000000000000000e+04,1.850000000000000000e+03,8.609161170931742225e-06,2.379273493406855371e-05,5.774810469590360504e-07,6.749398004320083672e-06,1.550863546335971917e-07,5.310908256538430609e-06,1.700316015683268081e-05,1.664291982567792680e-05,1.393205347293677750e+15
9.126000000000000000e+04,2.340000000000000000e+03,5.802633111431548545e-06,1.582601220917260784e-05,3.486582543694925725e-07,4.218547487638094829e-06,7.386616570535245681e-08,3.968298039669914082e-06,1.126140265581457409e-05,1.111950865229631463e-05,5.579504527340372000e+15
        }\tableAdaptiveThree






        %

        \addplot+ [deg0, uniform, forget plot]
        table [x=ndof, y=effRes] {\tableUniformZero};
        \addplot+ [deg1, uniform, forget plot]
        table [x=ndof, y=effRes] {\tableUniformOne};
        \addplot+ [deg2, uniform, forget plot]
        table [x=ndof, y=effRes] {\tableUniformTwo};
        \addplot+ [deg3, uniform, forget plot]
        table [x=ndof, y=effRes] {\tableUniformThree};

        \addplot+ [deg0, adaptive, forget plot]
        table [x=ndof, y=effRes] {\tableAdaptiveZero};
        \addplot+ [deg1, adaptive, forget plot]
        table [x=ndof, y=effRes] {\tableAdaptiveOne};
        \addplot+ [deg2, adaptive, forget plot]
        table [x=ndof, y=effRes] {\tableAdaptiveTwo};
        \addplot+ [deg3, adaptive, forget plot]
        table [x=ndof, y=effRes] {\tableAdaptiveThree};


    \end{semilogxaxis}
\end{tikzpicture}

%% file: Paper.bbl
\begin{thebibliography}{10}

\bibitem{MR3082295}
{\sc M.~Aurada, M.~Feischl, J.~Kemetm\"{u}ller, M.~Page, and D.~Praetorius}, {\em Each {$H^{1/2}$}-stable projection yields convergence and quasi-optimality of adaptive {FEM} with inhomogeneous {D}irichlet data in {$R^d$}}, ESAIM Math. Model. Numer. Anal., 47 (2013), pp.~1207--1235.

\bibitem{MR2194203}
{\sc C.~Bahriawati and C.~Carstensen}, {\em Three {MATLAB} implementations of the lowest-order {R}aviart-{T}homas {MFEM} with a posteriori error control}, Comput. Methods Appl. Math., 5 (2005), pp.~333--361.

\bibitem{MR3696081}
{\sc R.~E. Bank and J.~S. Ovall}, {\em Some remarks on interpolation and best approximation}, Numer. Math., 137 (2017), pp.~289--302.

\bibitem{MR4586821}
{\sc L.~Beir\~{a}o~da Veiga, F.~Brezzi, L.~D. Marini, and A.~Russo}, {\em The virtual element method}, Acta Numer., 32 (2023), pp.~123--202.

\bibitem{MR2195400}
{\sc R.~Bensow and M.~G. Larson}, {\em Discontinuous least-squares finite element method for the div-curl problem}, Numer. Math., 101 (2005), pp.~601--617.

\bibitem{MR2149925}
{\sc R.~E. Bensow and M.~G. Larson}, {\em Discontinuous/continuous least-squares finite element methods for elliptic problems}, Math. Models Methods Appl. Sci., 15 (2005), pp.~825--842.

\bibitem{MR2878612}
{\sc P.~Bochev, J.~Lai, and L.~Olson}, {\em A locally conservative, discontinuous least-squares finite element method for the {S}tokes equations}, Internat. J. Numer. Methods Fluids, 68 (2012), pp.~782--804.

\bibitem{MR3049438}
\leavevmode\vrule height 2pt depth -1.6pt width 23pt, {\em A non-conforming least-squares finite element method for incompressible fluid flow problems}, Internat. J. Numer. Methods Fluids, 72 (2013), pp.~375--402.

\bibitem{MR3097958}
{\sc D.~Boffi, F.~Brezzi, and M.~Fortin}, {\em Mixed finite element methods and applications}, vol.~44 of Springer Series in Computational Mathematics, Springer, Heidelberg, 2013.

\bibitem{MR2189548}
{\sc J.~H. Bramble, T.~V. Kolev, and J.~E. Pasciak}, {\em A least-squares approximation method for the time-harmonic {M}axwell equations}, J. Numer. Math., 13 (2005), pp.~237--263.

\bibitem{MR2373954}
{\sc S.~C. Brenner and L.~R. Scott}, {\em The mathematical theory of finite element methods}, vol.~15 of Texts in Applied Mathematics, Springer, New York, third~ed., 2008.

\bibitem{BringmannDissertation}
{\sc P.~Bringmann}, {\em Adaptive least-squares finite element method with optimal convergence rates}, PhD thesis,  (2021).
\newblock Humboldt-{U}niversit{\"a}t zu {B}erlin, Germany.

\bibitem{DLSFEM}
\leavevmode\vrule height 2pt depth -1.6pt width 23pt, {\em Discontinuous least-squares finite element method for the {Poisson} model problem using the {NGSolve} software library}.
\newblock \url{https://www.codeocean.com/}, 2024.
\newblock Python software package, available under DOI: \href{https://doi.org/10.24433/CO.1482259.v1}{10.24433/CO.1482259.v1}.

\bibitem{MR4782072}
\leavevmode\vrule height 2pt depth -1.6pt width 23pt, {\em Review and computational comparison of adaptive least-squares finite element schemes}, Comput. Math. Appl., 172 (2024), pp.~1--15.

\bibitem{MR3715170}
{\sc P.~Bringmann and C.~Carstensen}, {\em {$h$}-adaptive least-squares finite element methods for the 2{D} {S}tokes equations of any order with optimal convergence rates}, Comput. Math. Appl., 74 (2017), pp.~1923--1939.

\bibitem{MR3757107}
{\sc P.~Bringmann, C.~Carstensen, and G.~Starke}, {\em An adaptive least-squares {FEM} for linear elasticity with optimal convergence rates}, SIAM J. Numer. Anal., 56 (2018), pp.~428--447.

\bibitem{MR4796047}
{\sc P.~Bringmann, C.~Carstensen, and J.~Streitberger}, {\em Local parameter selection in the {${\rm C}^0$} interior penalty method for the biharmonic equation}, J. Numer. Math., 32 (2024), pp.~257--273.

\bibitem{MR1302685}
{\sc Z.~Cai, R.~Lazarov, T.~A. Manteuffel, and S.~F. McCormick}, {\em First-order system least squares for second-order partial differential equations. {I}}, SIAM J. Numer. Anal., 31 (1994), pp.~1785--1799.

\bibitem{MR2084237}
{\sc Z.~Cai and G.~Starke}, {\em Least-squares methods for linear elasticity}, SIAM J. Numer. Anal., 42 (2004), pp.~826--842.

\bibitem{MR1618488}
{\sc Y.~Cao and M.~D. Gunzburger}, {\em Least-squares finite element approximations to solutions of interface problems}, SIAM J. Numer. Anal., 35 (1998), pp.~393--405.

\bibitem{MR3215064}
{\sc C.~Carstensen, L.~Demkowicz, and J.~Gopalakrishnan}, {\em A posteriori error control for {DPG} methods}, SIAM J. Numer. Anal., 52 (2014), pp.~1335--1353.

\bibitem{MR3521055}
{\sc C.~Carstensen, L.~Demkowicz, and J.~Gopalakrishnan}, {\em Breaking spaces and forms for the {DPG} method and applications including {M}axwell equations}, Comput. Math. Appl., 72 (2016), pp.~494--522.

\bibitem{MR4332791}
{\sc C.~Carstensen, A.~Ern, and S.~Puttkammer}, {\em Guaranteed lower bounds on eigenvalues of elliptic operators with a hybrid high-order method}, Numer. Math., 149 (2021), pp.~273--304.

\bibitem{MR3246802}
{\sc C.~Carstensen and J.~Gedicke}, {\em Guaranteed lower bounds for eigenvalues}, Math. Comp., 83 (2014), pp.~2605--2629.

\bibitem{MR3671598}
{\sc C.~Carstensen, E.-J. Park, and P.~Bringmann}, {\em Convergence of natural adaptive least squares finite element methods}, Numer. Math., 136 (2017), pp.~1097--1115.

\bibitem{MR4700405}
{\sc C.~Carstensen and S.~Puttkammer}, {\em Adaptive guaranteed lower eigenvalue bounds with optimal convergence rates}, Numer. Math., 156 (2024), pp.~1--38.

\bibitem{MR3820383}
{\sc C.~Carstensen and J.~Storn}, {\em Asymptotic exactness of the least-squares finite element residual}, SIAM J. Numer. Anal., 56 (2018), pp.~2008--2028.

\bibitem{MR2431595}
{\sc D.~Day and P.~Bochev}, {\em Analysis and computation of a least-squares method for consistent mesh tying}, J. Comput. Appl. Math., 218 (2008), pp.~21--33.

\bibitem{DemkowiczGopalakrishnan2017}
{\sc L.~Demkowicz and J.~Gopalakrishnan}, {\em Discontinuous petrov–galerkin (dpg) method}, in Encyclopedia of Computational Mechanics Second Edition, E.~Stein, R.~de~Borst, and T.~Hughes, eds., John Wiley \& Sons, Ltd, Chichester, UK, 2017, pp.~1--15.

\bibitem{MR2882148}
{\sc D.~A. Di~Pietro and A.~Ern}, {\em Mathematical aspects of discontinuous {G}alerkin methods}, vol.~69 of Math\'{e}matiques \& Applications (Berlin) [Mathematics \& Applications], Springer, Heidelberg, 2012.

\bibitem{MR1393904}
{\sc W.~D\"{o}rfler}, {\em A convergent adaptive algorithm for {P}oisson's equation}, SIAM J. Numer. Anal., 33 (1996), pp.~1106--1124.

\bibitem{MR4242224}
{\sc A.~Ern and J.-L. Guermond}, {\em Finite elements {I}---{A}pproximation and interpolation}, vol.~72 of Texts in Applied Mathematics, Springer, Cham, 2021.

\bibitem{MR2970849}
{\sc V.~J. Ervin}, {\em Computational bases for {$RT_k$} and {$BDM_k$} on triangles}, Comput. Math. Appl., 64 (2012), pp.~2765--2774.

\bibitem{MR4253887}
{\sc T.~F\"{u}hrer and N.~Heuer}, {\em A robust {DPG} method for large domains}, Comput. Math. Appl., 94 (2021), pp.~15--27.

\bibitem{MR4138307}
{\sc T.~F\"{u}hrer and D.~Praetorius}, {\em A short note on plain convergence of adaptive least-squares finite element methods}, Comput. Math. Appl., 80 (2020), pp.~1619--1632.

\bibitem{MR4216839}
{\sc G.~Gantner and R.~Stevenson}, {\em Further results on a space-time {FOSLS} formulation of parabolic {PDE}s}, ESAIM Math. Model. Numer. Anal., 55 (2021), pp.~283--299.

\bibitem{MR2285779}
{\sc M.~I. Gerritsma}, {\em Direct minimization of the discontinuous least-squares spectral element method for viscoelastic fluids}, J. Sci. Comput., 27 (2006), pp.~245--256.

\bibitem{MR1910570}
{\sc M.~I. Gerritsma and M.~M.~J. Proot}, {\em Analysis of a discontinuous least squares spectral element method}, J. Sci. Comput., 17 (2002), pp.~297--306.

\bibitem{MR851383}
{\sc V.~Girault and P.-A. Raviart}, {\em Finite element methods for {N}avier-{S}tokes equations}, vol.~5 of Springer Series in Computational Mathematics, Springer-Verlag, Berlin, 1986.
\newblock Theory and algorithms.

\bibitem{MR2684360}
{\sc T.~Gudi}, {\em A new error analysis for discontinuous finite element methods for linear elliptic problems}, Math. Comp., 79 (2010), pp.~2169--2189.

\bibitem{IgelbuscherDissertation}
{\sc M.~Igelb\"uscher}, {\em Mixed and {Hybrid} {Least}-{Squares} {FEM} in {Nonlinear} {Solid} {Mechanics}}, PhD thesis,  (2021).
\newblock Universit\"at Duisburg-Essen, Germany.

\bibitem{MR4433564}
{\sc M.~Igelb\"{u}scher and J.~Schr\"{o}der}, {\em Hybrid mixed finite element formulations based on a least-squares approach}, in Non-standard discretisation methods in solid mechanics, J.~Schr{\"o}der and P.~Wriggers, eds., vol.~98 of Lect. Notes Appl. Comput. Mech., Springer, Cham, 2022, pp.~169--189.

\bibitem{MR0461948}
{\sc D.~C. Jespersen}, {\em A least squares decomposition method for solving elliptic equations}, Math. Comp., 31 (1977), pp.~873--880.

\bibitem{MR4414910}
{\sc N.~K. Kumar and S.~Mohapatra}, {\em Performance of nonconforming spectral element method for {S}tokes problems}, Comput. Appl. Math., 41 (2022), pp.~1--20.

\bibitem{MR4367673}
{\sc R.~Li, Q.~Liu, and F.~Yang}, {\em A discontinuous least squares finite element method for time-harmonic {M}axwell equations}, IMA J. Numer. Anal., 42 (2022), pp.~817--839.

\bibitem{MR4563176}
\leavevmode\vrule height 2pt depth -1.6pt width 23pt, {\em A discontinuous least squares finite element method for the {H}elmholtz equation}, Numer. Methods Partial Differential Equations, 39 (2023), pp.~1425--1448.

\bibitem{MR4065377}
{\sc R.~Li and F.~Yang}, {\em A least squares method for linear elasticity using a patch reconstructed space}, Comput. Methods Appl. Mech. Engrg., 363 (2020), pp.~112902, 19.

\bibitem{MR4593282}
{\sc Y.~Liang and S.~Zhang}, {\em Least-squares methods with nonconforming finite elements for general second-order elliptic equations}, J. Sci. Comput., 96 (2023), pp.~1--35, Paper No. 15.

\bibitem{MR4058304}
{\sc S.~Mohapatra, P.~Dutt, B.~V.~R. Kumar, and M.~I. Gerritsma}, {\em Non-conforming least-squares spectral element method for {S}tokes equations on non-smooth domains}, J. Comput. Appl. Math., 372 (2020), pp.~112696, 20.

\bibitem{MohapatraKumarJoshi2023}
{\sc S.~Mohapatra, N.~K. Kumar, and S.~Joshi}, {\em Least-squares formulations for {Stokes} equations with non-standard boundary conditions --- {A} unified approach}, Math. Meth. Appl. Sci., 46 (2023), pp.~16463--16482.

\bibitem{MR4708036}
{\sc H.~Monsuur, R.~Stevenson, and J.~Storn}, {\em Minimal residual methods in negative or fractional {S}obolev norms}, Math. Comp., 93 (2024), pp.~1027--1052.

\bibitem{ngsolve}
{\sc J.~Sch\"oberl}, {\em C++11 {Implementation} of {Finite} {Elements} in {NGSolve}}, ASC Report, 30 (2014).
\newblock TU Wien, Institute for Analysis and Scientic Computing, Austria.

\bibitem{MR1011446}
{\sc L.~R. Scott and S.~Zhang}, {\em Finite element interpolation of nonsmooth functions satisfying boundary conditions}, Math. Comp., 54 (1990), pp.~483--493.

\bibitem{MR4803195}
{\sc R.~Stevenson}, {\em A convenient inclusion of inhomogeneous boundary conditions in minimal residual methods}, Comput. Methods Appl. Math., 24 (2024), pp.~983--994.

\bibitem{MR3816182}
{\sc A.~Veeser and P.~Zanotti}, {\em Quasi-optimal nonconforming methods for symmetric elliptic problems. {I}---{A}bstract theory}, SIAM J. Numer. Anal., 56 (2018), pp.~1621--1642.

\bibitem{MR3059294}
{\sc R.~Verf\"{u}rth}, {\em A posteriori error estimation techniques for finite element methods}, Numerical Mathematics and Scientific Computation, Oxford University Press, Oxford, 2013.

\bibitem{MR3895838}
{\sc X.~Ye and S.~Zhang}, {\em A discontinuous least-squares finite-element method for second-order elliptic equations}, Int. J. Comput. Math., 96 (2019), pp.~557--567.

\bibitem{MR4018037}
{\sc X.~Ye, S.~Zhang, and P.~Zhu}, {\em A discontinuous {G}alerkin least-squares method for div-curl systems}, J. Comput. Appl. Math., 367 (2020), pp.~112474, 10.

\end{thebibliography}
